\documentclass[10pt]{amsart}
\usepackage{amsmath,amssymb}
\usepackage{amsthm}
\usepackage{graphicx}
\usepackage[pagebackref]{hyperref}

\renewcommand{\d}{\mathrm{d}}
\renewcommand{\L}{\Lambda}

\renewcommand{\ss}{\scriptstyle }

\newcommand{\cF}{{\mathcal F}}
\newcommand{\cI}{{\mathcal I}}
\newcommand{\cH}{{\mathcal H}}
\newcommand{\cZ}{{\mathcal Z}}

\newcommand{\bbR}{{\mathbb R}}
\newcommand{\bbC}{{\mathbb C}}
\newcommand{\C}[1]{\ensuremath{\mathbb C^{\,#1}}} 
\newcommand{\bbS}{{\mathbb S}}
\newcommand{\bbZ}{{\mathbb Z}}
\newcommand{\bbP}{{\mathbb P}}

\newcommand{\sfP}{{\mathsf P}}

\newcommand{\iC}{\mathrm{i}}
\newcommand{\ov}{\overline}
\newcommand{\ot}{\otimes}

\newcommand{\euso}{\operatorname{\mathfrak{so}}}
\newcommand{\eusu}{\operatorname{\mathfrak{su}}}
\newcommand{\eusl}{\operatorname{\mathfrak{sl}}}
\newcommand{\eusp}{\operatorname{\mathfrak{sp}}}
\newcommand{\euu}{\operatorname{\mathfrak u}}
\newcommand{\eug}{\operatorname{\mathfrak g}}

\newcommand{\euk}{\operatorname{\mathfrak k}}
\newcommand{\eum}{\operatorname{\mathfrak m}}
\newcommand{\eun}{\operatorname{\mathfrak n}}
\newcommand{\eue}{\operatorname{\mathfrak e}}
\newcommand{\eut}{\operatorname{\mathfrak t}}

\newcommand{\SO}{\operatorname{SO}}
\newcommand{\SL}{\operatorname{SL}}
\newcommand{\SU}{\operatorname{SU}}
\newcommand{\Un}{\operatorname{U}}
\newcommand{\GL}{\operatorname{GL}}
\newcommand{\Symp}{\operatorname{Sp}}
\newcommand{\Spin}{\operatorname{Spin}}
\newcommand{\Or}{\operatorname{O}}
\newcommand{\I}{\operatorname{I}}

\newcommand{\E}{\operatorname{E}}

\DeclareMathOperator{\tr}{tr}
\DeclareMathOperator{\sgn}{sgn}
\DeclareMathOperator{\Hom}{Hom}

\DeclareMathOperator{\Gr}{Gr}
\DeclareMathOperator{\Aut}{Aut}

\DeclareMathOperator{\vol}{vol}
\DeclareMathOperator{\ev}{ev}
\DeclareMathOperator{\rank}{rank}
\DeclareMathOperator{\Ad}{Ad}

\newcommand{\ab}{\mathsf{a}}
\newcommand{\bb}{\mathsf{b}}
\newcommand{\cb}{\mathsf{c}}
\newcommand{\eb}{\mathbf{e}}
\newcommand{\ub}{\mathbf{u}}
\newcommand{\vb}{\mathbf{v}}
\newcommand{\zb}{\mathbf{z}}

\newcommand{\us}{\mathsf{u}}
\newcommand{\vs}{\mathsf{v}}

\newcommand{\la}{\langle}
\newcommand{\ra}{\rangle}
\newcommand{\dl}{[\![}
\newcommand{\dr}{]\!]}

\newcommand{\w}{{\mathchoice{\,{\scriptstyle\wedge}\,}{{\scriptstyle\wedge}}
      {{\scriptscriptstyle\wedge}}{{\scriptscriptstyle\wedge}}}}

\numberwithin{equation}{section}

\newtheorem{theorem}{Theorem}
\newtheorem{lemma}{Lemma}
\newtheorem{proposition}{Proposition}
\newtheorem{corollary}{Corollary}

\theoremstyle{remark}
\newtheorem{definition}{Definition}

\newtheorem{remark}{Remark}
\newtheorem{example}{Example}

\begin{document}

\author[R. Bryant]{Robert L. Bryant}
\address{Duke University Mathematics Department\\
         P.O. Box 90320\\
         Durham, NC 27708-0320}
\email{\href{mailto:bryant@math.duke.edu}{bryant@math.duke.edu}}
\urladdr{\href{http://www.math.duke.edu/~bryant}%
         {http://www.math.duke.edu/\lower3pt\hbox{\symbol{'176}}bryant}}

\title[Rigidity of extremal cycles]
      {Rigidity and quasi-rigidity \\
                    of \\
              extremal cycles \\
                    in \\
        Hermitian symmetric spaces}

\date{March 5, 2001}

\begin{abstract}
Let~$M$ be a compact Hermitian symmetric space
and let~$W\not=\emptyset$ be a compact complex subvariety 
of~$M$ of codimension~$p$.  There exists a nontrivial 
holomorphic exterior differential system~$\cI$ on~$M$ 
with the property that any compact complex
subvariety~$V\subset M$ of dimension~$p$ that 
satisfies~$[V]\cap[W]=0$ is necessarily an integral 
variety of~$\cI$.  

The system~$\cI$ is almost never involutive.  
However, its $p$-dimensional integral varieties (when they exist) 
can sometimes be described explicitly by taking advantage
of this non-involutivity.  In this article, 
several of these ideals~$\cI$ will be analyzed,
particularly in the case where~$M$ is a Grassmannian, 
and the results applied to prove various results
about the rigidity of algebraic subvarieties 
with certain specified homology classes.

These rigidity results have implications for the classification 
of certain holomorphic bundles over compact K\"ahler manifolds 
that are generated by their global sections. For example,
if~$F\to M$ is generated by its global sections and~$M$ is
compact and K\"ahler, then, as is well-known, $c_2(F)\ge0$.
If equality holds, then either~$F$ is the pullback to~$M$ of
a holomorphic bundle~$F'\to C$ over a curve~$C$ via a holomorphic
map~$\kappa:M\to C$ or else~$F = L\oplus T$ where~$L$ is a line
bundle and~$T$ is trivial.  There is a similar
(though more complicated) characterization when~$c_3(F)=0$.

\end{abstract}



\thanks{
This article was begun at l'Institut
des Hautes \'Etudes Scientifiques during June of 1999.  
I thank l'IHES for its hospitality. Thanks are also due to
Duke University for its support via a research grant and the
National Science Foundation for its support via DMS-9870164.\hfill\break
\hspace*{\parindent} 
This is Version~$2$. The most recent version
can be found at arXiv:math.DG/0006186.
}

\maketitle

\setcounter{tocdepth}{2}
\tableofcontents

\listoffigures

\section[Introduction]{Introduction}\label{sec:intro}

\subsection{An overview}
\label{ssec:overiew}
This article is an account of some basic local geometric properties 
of the Hermitian symmetric spaces and how these properties can be used
to derive interesting topological and algebro-geometric consequences.

\subsubsection{A seed problem}
\label{sssec: seed prob}
The study that lead to this article was inspired by a simple
problem:  To understand, from a geometric point of view, why a certain
subvariety in~$\Gr(3,6)$, the Grassmannian of $3$-planes in~$\C6$, 
cannot be smoothed,%
\footnote{My use of `smoothed' and `smoothable'
is not always in agreement with the usage common in algebraic
geometry.  For more discussion on this, see~\S\ref{sssec:smoothability}.} 
i.e., is not homologous to a smooth subvariety of~$\Gr(3,6)$. 

The subvariety in question can be described as follows  
(where I will generalize the setting for the sake of exposition):  
The Grassmannian~$\Gr(m,n)$ of $m$-dimensional subspaces of~$\C{n}$
is a smooth compact algebraic variety of dimension~$m(n{-}m)$ that is
naturally embedded into the projective space~$\bbP\bigl(\L^m(\C{n})\bigr)$.
For any~$k\le n{-}m$ and any subspace~$W\subset\C{n}$ 
of codimension~$m{+}k{-}1$, the subvariety
$$
\sigma(W)= \bigl\{ \ E\in\Gr(m,n)\mid E\cap W\not=\{0\}\ \bigr\}
$$
has codimension~$k$ in~$\Gr(m,n)$.  
(The subvariety~$\sigma(W)$ is one of an important family of 
subvarieties of the Grassmannians known as \emph{Schubert cycles}
(see~\S\ref{ssec:Schubert-cycles}) whose associated
homology classes form a natural basis for the integral homology 
of~$\Gr(m,n)$. In particular,~$\sigma(W)$ is denoted~$\sigma_{(k)}$ 
in the standard notation for Schubert cycles.)

When~$m\ge 2$ and~$0< k < n{-}m$, the variety~$\sigma(W)$ 
is singular.  For example, when~$m=2$, 
the singular locus of~$\sigma(W)$ is $\Gr(2,W)\subset\Gr(2,n)$, 
i.e., the set of $2$-planes that lie completely in~$W$.  

Now, in some cases, $\sigma(W)$, though singular, is homologous 
to a smooth subvariety of~$\Gr(m,n)$.  
For example, when~$k=1$ the 
hypersurface~$\sigma(W)\subset\Gr(m,n)$
is a hyperplane section~$\Gr(m,n)\cap\bbP(H_W)$ 
where~$H_W\subset\L^m(\C{n})$ is the hyperplane of~$m$-vectors 
that are annihilated by the decomposable $m$-form~$\alpha_W$ 
(unique up to multiples) that has~$W\subset\C{n}$ as its kernel.
Meanwhile, for the generic hyperplane~$H\subset \L^m(\C{n})$,
the intersection~$\Gr(m,n)\cap \bbP(H)$ is a smooth hypersurface 
in~$\Gr(m,n)$ that is homologous to~$\sigma(W)$.

However, when~$k=2$, Hartshorne, Rees, and Thomas~\cite[Theorem 2]{MR50:9870} 
show that $\sigma(W)\subset\Gr(3,6)$ is not homologous to
a smooth subvariety.  They do this by using results of 
Thom~\cite{MR14:1005,MR15:890a} to prove the stronger result 
that the integral homology class of~$\Gr(3,6)$ that~$\sigma(W)$ 
represents is not representable as an integral linear combination
of homology classes of smooth, oriented submanifolds 
of~$\Gr(3,6)$ of (real) codimension~$4$.

A slightly different situation presents itself for the case $k=2$
in~$\Gr(2,5)$.  In this case, $\sigma(W)$ itself is singular, but 
its homology class can be written as the difference of the homology
classes of two nonsingular subvarieties.  However, this use of 
differences is essential, because it turns out that~$\sigma(W)$ is
not homologous to any nonsingular subvariety.

In fact, when~$k=2$ and~$n\ge 5$, it turns out 
(see Theorem~\ref{thm:codim2-sigma2}) that any 
codimension~$2$ subvariety of~$\Gr(m,n)$ that is 
homologous to~$\sigma(W)$ must actually be equal 
to~$\sigma(W')$ for some subspace~$W'\subset\C{n}$ 
of codimension~$m{+}1$ (and hence must be singular). 
Moreover, it turns out that no integral multiple of the
homology class of $\sigma(W)$ can be represented by a smooth
subvariety of~$\Gr(m,n)$ (see the
discussion in Example~\ref{ex:sigma2-must-be-singular}).  This
is in spite of the fact that results of Thom~\cite{MR15:890a} show
that there is an integral multiple of the homology class of~$\sigma(W)$
that can be represented by a smooth submanifold of~$\Gr(m,n)$, 
even one with a complex normal bundle.  Of course, such a submanifold
cannot be holomorphic.

My goal in this article is to explain these sorts of nonsmoothability 
and rigidity results
from a more geometric perspective, using techniques from differential 
and algebraic geometry,
along the lines of Griffiths and Harris~\cite{MR81k:53004}
rather than the topological techniques of Thom.%
\footnote{That this might be an interesting problem was suggested to
me by Robin Hartshorne.  I would like to thank him for a very stimulating
conversation.}

\subsubsection{The basic idea}
\label{sssec: basic idea}

I now want to explain why one might expect to be able to approach 
this problem by local, differential-geometric techniques.  

Recall that complex vector spaces are canonically oriented, so that
it makes sense to say whether a top-degree differential form is
positive or not.  More generally, one says that a (real-valued) 
$2p$-form on a complex manifold~$M$ is \emph{weakly positive} 
if it is non-negative on every complex $p$-plane~$E\subset T_mM$. 
The standard example of such a form is the $p$-th power of a K\"ahler
form (which is actually positive on each tangent complex $p$-plane).

It is an interesting feature of the Grassmannians~$\Gr(m,n)$ (which,
as will be seen, generalizes to other Hermitian symmetric spaces and
Schubert varieties) that there exists a closed, weakly positive
$\bigl(m(n{-}m){-}2,m(n{-}m){-}2\bigr)$-form~$\phi$ 
on~$\Gr(m,n)$ that is non-zero in cohomology
and yet vanishes identically when pulled back to the smooth locus 
of~$\sigma(W)$ where~$W\subset\C{n}$ is a subspace of codimension~$m{+}1$.   
It follows that~$\phi$ must vanish identically when pulled back 
to the smooth locus of any codimension~$2$ subvariety~$X\subset\Gr(m,n)$ 
that is homologous to~$W$.

The conditions on a complex $(m(n{-}m){-}2)$-plane~$E\subset T_V\Gr(m,n)$ 
that~$\phi$ vanish identically on~$E$ turn out to be
very restrictive.   An analysis of these conditions
shows that any codimension~$2$ complex submanifold~$X\subset\Gr(m,n)$ 
to which $\phi$ pulls back to be zero
must satisfy an overdetermined
system of holomorphic first order partial differential equations.
Fortunately, this system of equations is fairly simple,%
\footnote{In particular, one does not need to explicitly invoke
the machinery of exterior differential systems;  
elementary arguments using the moving frame suffice.}
and one can describe its local solutions explicitly 
in terms of local subvarieties of~$\bbP^{n-1}$:
One finds that there exists a rational 
map~$\lambda:X\dashrightarrow \bbP^{n-1}$ whose differential
generically has rank~$n{-}m{-}2$ with the property that, for
all~$V\in X$, the point~$\lambda(V)$ lies in~$\bbP(V)$.
Thus, letting~$Y\subset\bbP^{n-1}$ be the closure of the 
image of~$\lambda$, one finds that~$X$ is a subset of
the variety of~$m$-planes whose projectivizations meet~$Y$.
For dimension reasons, $X$ must be open in this variety.
(See Theorem~\ref{thm:codim2-sigma2} for details). 

Thus, when~$n\ge m{+}2$, 
for any codimension~$2$ algebraic variety~$X\subset\Gr(m,n)$
that satisfies the equation in integral 
homology~$[X] = r\,\bigl[\sigma(W)\bigr]$, 
there exists a codimension~$m{+}1$ subvariety~$Y\subset\bbP^{n-1}$ 
of degree~$r$ so that
$$
X = \bigl\{\ V\in\Gr(m,n)\mid \bbP(V)\cap Y\not=\emptyset\ \bigr\}.
$$
From this description,%
\footnote{When~$m=2$ and~$n=4$, this is a classical 
result~\cite[p.~143]{gFano}.
(I thank Igor Dolgachev for supplying me with this reference).  In a
private communication (24 July 2000), Chad Schoen has supplied a proof 
of a version of this result when~$m=2$ that is valid over any 
algebraically closed field.} 
it is easy to see that~$X$ is singular unless~$Y$ is a single point.  
In particular, $X$ is singular if~$n\ge m{+}3$ or if~$n=m{+}2$ and~$r>1$.  

Thus, this line of argument realizes the original goal 
of finding a geometric explanation for the fact that many 
of the varieties~$\sigma(W)$ are not homologous to smooth subvarieties.
It also provides sharper results, since it shows that~$\sigma(W)$ 
cannot even be deformed in any non-trivial way:  
Any codimension~$2$ subvariety~$X\subset\Gr(m,n)$ that is homologous
to~$\sigma(W)$ is of the form~$X = \sigma(W')$ 
for some subspace~$W'\subset\C{n}$ of codimension~$m{+}1$.  For this 
reason, the cycle~$\sigma(W)$ will be said to be~\emph{rigid}.  
Even when $[X]=r\,\bigl[\sigma(W)\bigr]$ for some~$r>1$, the variety~$X$
displays a form of rigidity: It is a union of an $(n{-}m{-}2)$-parameter
family of~$\Gr(m{-}1,n{-}1)$s linearly embedded into~$\Gr(m,n)$.  
Roughly speaking, it can only `deform' in $n{-}m{-}2$ 
of its~$m(n{-}m){-}2$ dimensions.
I refer to this (rather loosely defined) 
property as \emph{quasi-rigidity}.

\subsubsection{The general program}
\label{sssec: gen prog}

While the considerations above may seem very special, they actually
generalize to cover an enormous number of cases. 

Let~$M=U/K$ be an irreducible Hermitian symmetric space of compact
type, where~$U$ is compact, connected simple Lie group and~$K\subset U$ 
is a symmetric subgroup with a central subgroup of dimension~$1$.

By results of Kostant~\cite{MR26:265,MR26:266}, 
there is an essentially canonical basis~$\sfP(M)$ for the 
integral homology~$H_*(M,\bbZ)$ (which is all of even degree and
torsion free) that generalizes the well-known Schubert basis
of the Grassmannians~$\Gr(m,n)$ (see~\S\ref{sec:grass}).
Each of the elements of the basis~$\sfP(M)$
is representable by a (generalized) Schubert 
variety in~$M$ and, furthermore, the integral homology 
class of any compact subvariety of~$M$ is an integral combination
of elements of~$\sfP(M)$ with non-negative coefficients.

For each such Schubert variety~$\sigma\subset M$, 
there is a unique $U$-invariant harmonic form~$\phi_\sigma$
of type~$(p,p)$ where~$p = \dim M - \dim \sigma$ 
that represents intersection with~$\sigma$,
in the sense that
$$
\int_X\phi_\sigma = [\sigma]\cap[X]\in\bbZ
$$
when~$X\subset M$ is any subvariety of (complex) dimension~$p$
and the right hand side is interpreted as the homological 
intersection pairing.   
As Kostant shows, the form~$\phi_\sigma$ is weakly positive.%
\footnote{He actually proves the stronger result 
that it is \emph{positive} in the sense of~\S\ref{ssec:positivity}.}
In particular, $\phi_\sigma$ must vanish identically on any~$X$ 
that satisfies the homological condition~$[\sigma]\cap[X]=0$.

It turns out (and, in any case, follows easily from Kostant's
results) that any complex submanifold of~$M$ on which~$\phi_\sigma$
vanishes must satisfy a first order system of holomorphic
partial differential equations. This system depends only 
on the cohomology class~$[\sigma]$ and turns out to be invariant 
under~$G$, the identity component of the group of biholomorphisms 
of~$M$ (which contains~$U$ as a maximal compact subgroup).

Thus, one may expect to get global information about the 
complex subvarieties~$X$ that satisfy~$[X]\cap[\sigma]=0$
by studying the studying the local solutions of this system
of partial differential equations.  This expectation is amply 
borne out by the results in this article.  

The cases in which~$\sigma$ has low dimension or codimension 
turn out to be particularly accessible, and a complete description 
of the subvarieties~$X$ satisfying~$[X]\cap[\sigma]=0$ is 
available.  It often takes the form of saying that such
subvarieties are rigid or quasi-rigid in a sense analogous
to that of the examples discussed above in the Grassmannian
case.

These descriptions of rigid and quasi-rigid varieties in the
Grassmannians will be applied to the characterization of 
holomorphic bundles over compact complex manifolds that are 
generated by their sections and yet satisfy certain vanishing
conditions on polynomials in their Chern classes.

Of course, the idea of using (weak) positivity of a $(p,p)$ form 
representing a cohomology class on a complex manifold~$M$ to derive 
information about the subvarieties on which it vanishes is not new.
In fact, this already appears in the work of Kostant cited above.
Another place where this technique has been used to great effect
is in Griffiths and Harris~\cite[\S4]{MR81k:53004}, where they
combine these ideas with information coming from the geometry of 
Gauss maps to study the subvarieties of Abelian varieties 
that have degenerate Gauss maps.  

After the first version of this
article was posted to the arXiv, Dan Burns%
\footnote{private communication, 15 September 2000}
brought to my attention
the (unpublished) 1997 thesis%
\footnote{Also, see the preprints~\cite{mW98-1} and~\cite{mW98-2},
which contain expositions of some of Walters' thesis results.}
of Maria Walters~\cite{mW97},
in which she also investigated the consequences of positivity 
of certain of the forms on the Grassmannians to prove rigidity
results.  Some of her results anticipate mine.
I will discuss the relation between her results
and the results of the present article as the 
opportunity arises.

In concluding this overview, I would like to thank several people
for their very helpful comments and suggestions on the first 
version of this article or for references to the algebraic
geometry literature:  Dan Burns, Igor Dolgachev, Phillip Griffiths, 
Robin Hartshorne, Chad Schoen, and Maria Walters.  Any errors
or infelicities that remain are solely due to me.

\subsection{Background}\label{ssec:background}
In this section, I introduce some of the concepts that
will be important in this article.

\subsubsection{Differential ideals and systems}
\label{sssec: dif ideals}
The reader will probably be relieved to know that, in the cases
studied in this article, no essential use is made of the theory
of exterior differential systems as such.  For example, there
is no use of the concepts of polar spaces, regularity, 
involutivity, characteristic variety, and so on.  The Cartan-K\"ahler
theorem will not even be mentioned outside this sentence.

In fact, the reader only needs to know the following
exterior differential systems terminology to read this article:%
\footnote{In this article, I will only be concerned with 
differential systems in the holomorphic category, and so have
adopted definitions suitable for this purpose.  Of course, the
general theory is not restricted to this case.}

\begin{definition}[Differential ideals, integral elements, and
integral varieties]
\label{def:diff ideal}
A \emph{differential ideal}~$\cI$ on~$M$ is a sheaf of ideals 
of holomorphic differential forms on~$M$ that is closed 
under exterior derivative.
An \emph{integral element} of~$\cI$ 
is a (complex) subspace~$E\subset T_xM$
on which all of the forms in~$\cI_x$ vanish.  
An \emph{integral variety} of~$\cI$ is a subvariety~$X\subset M$
with the property that, at every smooth point~$x\in X$, 
the tangent space~$T_xX$ is an integral element of~$\cI$.
\end{definition}

The general procedure for extracting information 
about integral varieties of an ideal~$\cI$ is to, first, 
compute the space of integral elements, which is essentially an 
algebraic problem and, second, use the integral elements to describe 
the local submanifolds that are integral varieties of~$\cI$, 
which is a differential geometric problem.  
One then applies the local description from the second step
to deduce global results about algebraic integral varieties of~$\cI$.

As a guide to the reader, I generally call algebraic
results about integral elements `Lemmas', 
local differential geometric results `Propositions',
and global topological or algebro-geometric results `Theorems'.
Thus, the names do not always reflect the degree of difficulty
of the corresponding proofs.  In fact, the most difficult arguments
in the article tend to be the algebraic arguments that compute
the integral elements of a given ideal.  The differential geometric
arguments are usually straightforward (if somewhat involved)
applications of the method of the moving frame.%
\footnote{Undoubtedly, the reason that no deeper results are needed 
from the theory of exterior differential systems is that I only 
analyze ideals of low degree or codegree in each case.  It seems
unlikely to me that such elementary methods will suffice for all of the 
ideals in the midrange.}

I do not mean to suggest that some knowledge of exterior differential
systems (EDS) would not be helpful, and the interested reader might want 
to consult~\cite{MR92h:58007}.  Indeed, many of the results in this 
article were first found by doing an exterior differential systems 
analysis.  However, once the results were found, it was possible 
to prove them without invoking EDS theory, so I did.  
Meanwhile, for the reader familiar with EDS theory, I have included
comments from time to time that point out EDS features 
that may be of interest. 
Other readers can safely ignore these comments.

This EDS avoidance does not significantly lengthen any of the proofs, 
so I feel that the savings of not having to introduce and discuss
concepts from exterior differential systems justifies this strategy.
The main disadvantage to the reader is that it does not explain
why the rigidity results could have been anticipated, making them
seem somewhat miraculous.   

There is a more general notion that, while it will not play any
direct role in this article, will be needed in the discussions
of the work of Maria Walters:

\begin{definition}
\label{def:diff sys}
Let~$M$ be a complex manifold and let~$m$ be an integer satisfying
$0<m<\dim M$.  Let~$\Gr(m,TM)$ denote the complex manifold whose
elements are the complex $m$-planes tangent to~$M$, i.e., each 
$E\in \Gr(m,TM)$ is an $m$-dimensional subspace~$E\subset T_xM$
for some~$x\in M$.  A \emph{differential system} for $m$-dimensional
subvarieties of~$M$ is a subvariety~$\Sigma\subset \Gr(m,TM)$.  
A \emph{solution}%
\footnote{Synonyms:  \emph{integral} of~$\Sigma$ 
or \emph{$\Sigma$-variety}.}
 of~$\Sigma$ is a subvariety~$X\subset M$ with 
the property that~$T_xX$ lies in~$\Sigma$
for every smooth point~$x\in X$.
\end{definition}

Strictly speaking, such a~$\Sigma\subset \Gr(m,TM)$ should be 
called a \emph{first-order} differential system, but since no other
kind will appear in this article, I will leave this as understood.

For any differential ideal~$\cI$ on~$M$ and any integer~$m$
with~$0< m < \dim M$, the space of $m$-dimensional integral elements
of~$\cI$ defines a differential system~$\mathcal{V}_m(\cI)\subset\Gr(m,TM)$.
Not every differential system in the above sense is of the 
form~$\mathcal{V}_m(\cI)$ for some ideal~$\cI$, so this is a proper
generalization. 

\subsubsection{Effective cycles}
\label{sssec:effective cycles}
Let~$M$ be a compact complex manifold. 
For each integer~$p$ in the range~$0\le p\le \dim M$, let~$\cZ^+_p(M)$ 
denote the semigroup of effective $p$-cycles in~$M$.  Thus, an element~$X$ 
in~$\cZ^+_p(M)$ is a formal sum~$X_1+ X_2+\dots+X_k$, 
where each~$X_i\subset M$ is an irreducible, compact, complex, 
$p$-dimensional subvariety.  (Of course, the~$X_i$ need not be
distinct.)

Since a compact complex $p$-dimensional subvariety~$X\subset M$ is 
triangulable~\cite{MR30:3478}, its singular locus has codimension at 
least~$2$~\cite{MR80b:14001}, and its smooth locus~$X^\circ\subset X$ is 
canonically oriented, it follows that~$X$ defines a 
homology class~$[X]\in H_{2p}(M,\bbZ)$   
and this extends to a semigroup
homomorphism~$[\cdot]:\cZ^+_p(M)\to H_{2p}(M,\bbZ)$.   

It is a fundamental problem in complex geometry to describe the 
image semigroup~$\left[\cZ^+_p(M)\right]\subset H_{2p}(M,\bbZ)$.
Certain of these classes will play an important role in this article:

\begin{definition}[Atomic classes and extremal rays]
\label{def: atomic and extremal}
A class~$z\in \bigl[\cZ^+_p(M)\bigr]$ will be said to be \emph{atomic}
if it cannot be written as a sum~$z = z_1+z_2$ 
where~$z_1,z_2\in\bigl[\cZ^+_p(M)\bigr]$ are both nonzero.  When~$z$
is atomic, the ray~$R_z = \bbZ^+\cdot z$ will be said to be \emph{extremal}
if $z_1+z_2\in R_z$ for~$z_1,z_2\in\bigl[\cZ^+_p(M)\bigr]$ implies
$z_1,z_2\in R_z$.
\end{definition}

Obviously, these notions are not useful when $\bigl[\cZ^+_p(M)\bigr]$
contains torsion classes in~$H_{2p}(M,\bbZ)$.  However, in the cases
of interest in this article, there will be no torsion classes anyway.
For example, when~$M$ admits a K\"ahler form~$\omega$, there are no
torsion classes in~$\bigl[\cZ^+_p(M)\bigr]$ 
(see Example~\ref{ex:the-Kahler-form}).

Given~$z\in H_{2p}(M,\bbZ)$, 
one could ask for a description of the set
\begin{equation}
\cZ^+_p(M,z) = \bigl\{\ X\in \cZ^+_p(M)\ \vrule\ [X] = z\ \bigl\}.
\end{equation}
When $M$ is compact and K\"ahler, $\cZ^+_p(M,z)$
has the structure of a (possibly reducible) compact, 
complex analytic space~\cite{MR80h:32056}.  Furthermore, when $M$ 
is projective, $\cZ^+_p(M,z)$ is a finite union of irreducible 
projective varieties~\cite{MR56:662,MR80h:32056}.  The study of
these varieties is a large part of algebraic geometry,
with even relatively simple cases, 
such as~$\cZ^+_1\bigl(\bbP^3,r[\bbP^1]\bigr)$, 
not being fully understood~\cite[Chapter~IV, \S6]{MR57:3116}.

\begin{example}[Surfaces of degree~$2$ in~$\bbP^4$]
\label{ex:deg2surfsinP4}
As is well-known,~$H_4(\bbP^4,\bbZ)= \bbZ\cdot[\bbP^2]$.  
The variety~$\cZ_2^+\bigl(\bbP^4,[\bbP^2]\bigr)$ consists of the
linear $\bbP^2$s in~$\bbP^4$, and so can be identified with~$\Gr(3,5)$,
which has complex dimension~$6$ and is irreducible and smooth. 

On the other hand, the variety~$\cZ_2^+\bigl(\bbP^4,\,2[\bbP^2]\bigr)$
is neither irreducible nor smooth.  It has two irreducible components,
one of dimension~$12$ and the other of dimension~$13$.
The first component consists of pairs of $\bbP^2$s in~$\bbP^4$ and so is
identifiable with the symmetric product~$\Gr(3,5)^{(2)}$, a singular
variety of dimension~$12$.  The second component consists of the degree~$2$ 
surfaces that are degenerate, i.e., that lie in some hyperplane in~$\bbP^4$.
Since the quadric surfaces that lie in a given~$\bbP^3$ form a~$\bbP^9$, 
this second component has dimension~$13$.  These two irreducible components 
meet in a common subvariety of dimension~$10$ that consists of the pairs 
of $\bbP^2$s that meet in at least a line.
\end{example}

\subsubsection{Smoothability}
\label{sssec:smoothability}
In this article, I will adopt a rather na\"{\i}ve notion of 
what it means for a subvariety to be smoothable.

\begin{definition}[Smoothability]
\label{def: smoothability}
An class~$z\in H_{2p}(M,\bbZ)$ will be said to be 
\emph{smoothly representable} if there exists a smooth 
(i.e., nonsingular) subvariety~$X\subset M$ so that~$[X] = z$.  

An element~$X\in\cZ^+_p(M)$ will be said to be 
\emph{smoothable} if~$[X]$ is smoothly representable.
\end{definition}

As to why this definition may be viewed as na\"{\i}ve, 
the reader should compare the discussion in~\cite{MR50:9870},
where a class~$z\in H_{2p}(M,\bbZ)$ is regarded 
as smoothly representable if there exist nonsingular, 
irreducible subvarieties~$X_1,\dots,X_k\subset M$ 
of dimension~$p$ and (not necessarily positive) 
integers~$m_1,\dots,m_k$
so that~$z = m_1[X_1]+\dots+m_k[X_k]$.  
I might have called this latter notion 
\emph{virtual smooth representability} and the associated notion
of smoothability \emph{virtual smoothability}, 
but it turns out that algebraic geometers have found this version
of (homological) smooth representability to be the most useful, 
so, in most sources, `smoothability' means `virtual smoothability' 
and not the term as I have defined it.

Example~\ref{ex:deg2surfsinP4} shows why this article's notion
of `smoothability' is different from some sort of `smoothable
under small deformations' since the union of a pair of transversely 
intersecting $\bbP^2$s in~$\bbP^4$ is not `smoothable by a small 
deformation' in the strict sense, though it is smoothable in the
sense adopted in this article because it is homologous to a 
smooth quadric surface in any $\bbP^3$ in~$\bbP^4$.

\begin{example}[Curves of degree~$3$ in~$\bbP^3$]
\label{ex:deg1curvesinP3}
As is well-known,~$H_2(\bbP^3,\bbZ) = \bbZ\cdot[\bbP^1]$.  
The variety~$\cZ_1^+\bigl(\bbP^3,[\bbP^1]\bigr)$ consists of the
linear $\bbP^1$s in~$\bbP^3$, and so can be identified with~$\Gr(2,4)$,
which has complex dimension~$4$ and is irreducible and smooth. 

The variety~$\cZ_1^+\bigl(\bbP^3,\,2[\bbP^1]\bigr)$
is neither irreducible nor smooth.  It has two irreducible components,
both of dimension~$8$.  
The first component consists of pairs of $\bbP^1$s in~$\bbP^3$ and so is
identifiable with the symmetric product~$\Gr(2,4)^{(2)}$.  
The second component consists of the degree~$2$ 
curves that are degenerate, i.e., that lie in some hyperplane in~$\bbP^3$.
Since the conics that lie in a given~$\bbP^2$ form a~$\bbP^5$, 
this second component has dimension~$8$ as well.  
These two irreducible components meet in a common subvariety 
of dimension~$6$ that consists of the pairs 
of $\bbP^1$s that meet in at least a point.

The variety~$\cZ_1^+\bigl(\bbP^3,\,3[\bbP^1]\bigr)$ has four 
irreducible components, each of dimension~$12$.  In addition to the
two components that are the image of the natural map
\begin{equation}
\cZ_1^+\bigl(\bbP^3,[\bbP^1]\bigr) \times
\cZ_1^+\bigl(\bbP^3,2[\bbP^1]\bigr) 
\longrightarrow
\cZ_1^+\bigl(\bbP^3,3[\bbP^1]\bigr),
\end{equation}
there are two further components, one consisting of the degenerate
curves of degree~$3$ (the generic member of which is a nonsingular
plane cubic) and the other consisting of the closure of the space
of twisted cubic curves.  Note that the `generic' element of each
component is a smooth (though possibly reducible) curve, but that 
these represent four very different ways of smoothing~$3[\bbP^1]$.
\end{example}

\subsubsection{Positivity}\label{sssec:background-positivity}
An $\bbR$-valued $2p$-form~$\phi$ on~$M$ 
is said to be \emph{weakly positive} 
in the sense of Harvey and Knapp~\cite{MR50:7573}
if it evaluates on each complex $p$-plane~$E\subset T_xM$ to be
nonnegative (see \S\ref{ssec:positivity}).  If~$\phi$ is closed
and weakly positive, its deRham cohomology 
class~$[\phi]\in H^{2p}(M,\bbR)$ satisfies
\begin{equation}
\bigl\la [\phi],[X] \bigr\ra = \int_X\phi \ge 0
\end{equation}
for any compact $p$-dimensional subvariety~$X\subset M$.  
When $[\phi]\not=0$, this implies that the 
image~$\bigl[\cZ^+_p(M)\bigr]$ must lie in a 
closed `halfspace'~$H^+\bigl([\phi]\bigr)\subset H_{2p}(M,\bbZ)$.

The intersection of these halfspaces $H^+\bigl([\phi]\bigr)$ as~$\phi$ 
ranges over the closed, weakly positive $2p$-forms on~$M$ is a 
semigroup~$H^+_{2p}(M,\bbZ)$ that evidently satisfies
\begin{equation}
\bigl[\cZ^+_p(M)\bigr]\subseteq H^+_{2p}(M,\bbZ)\subset H_{2p}(M,\bbZ).
\end{equation}

\begin{example}[K\"ahler forms]
\label{ex:the-Kahler-form} 
When~$M$ admits a K\"ahler structure~$\omega$,
it defines an Hermitian metric on~$M$ and the 
Wirtinger theorem~\cite[p. 31]{MR80b:14001} implies 
\begin{equation}\label{eq: Wirtinger theorem}
\bigl\la[\omega^p], [X]\bigr\ra = \int_{X} \omega^p = p!\,\vol(X) > 0. 
\end{equation}
Thus, $\bigl[\cZ^+_p(M)\bigr]$ (if nonempty) lies 
strictly on one side of a hyperplane in~$H_{2p}(M,\bbZ)$. 

Note also that, in the K\"ahler case, 
$\bigl[\cZ^+_p(M)\bigr]$ cannot contain
any torsion classes.  Moreover, if a class~$z\in \bigl[\cZ^+_p(M)\bigr]$
is atomic, then any~$X\in \cZ^+_p(M,z)$ must be irreducible.  If, in 
addition, the ray~$R_z$ is extremal, then any $X\in\cZ^+_p(M,\,r\,z)$
is the sum~$X = X_1+\cdots+X_k$ 
of irreducible~$X_i\in \cZ^+_p(M,r_i z)$ where~$r = r_1+\cdots+r_k$.
\end{example}

Despite its fundamental importance, the K\"ahler form 
is not typical of the sort of positive form that will be
studied in this article.  Instead, I will be interested in closed 
weakly positive $2p$-forms~$\phi$ that vanish identically on certain 
$p$-dimensional compact subvarieties~$X$.  
In such a situation, any effective cycle~$X'$ whose homology class
is that of~$r[X]$ for any~$r\in\bbZ^+$ must necessarily be a union of
irreducible $p$-cycles on which~$\phi$ vanishes.

\begin{definition}[Zero planes of a weakly positive form]
\label{def: Z phi}
Let~$\phi\in\Omega^{2p}(M)$ be a weakly positive form.  
Then~$Z(\phi)\subset \Gr(p,TM)$ denotes the set of complex tangent
$p$-planes on which~$\phi$ vanishes.  A complex subvariety~$X\subset M$
to which~$\phi$ pulls back to become zero is known 
as a \emph{$\phi$-null subvariety}.
\end{definition} 

In many cases, the set~$Z(\phi)$ is rather small, 
and, consequently, this implies severe 
restrictions on the possible $\phi$-null subvarieties.
It is not generally true that~$Z(\phi)$ is a differential system
in the sense of Definition~\ref{def:diff sys}.

However, in many cases, when~$\phi$ satisfies a strengthened
condition known as \emph{positivity} (see~\S\ref{ssec:positivity}),
one can construct an exterior differential
system~$\cI_\phi$ whose complex integral varieties are exactly
the $\phi$-null subvarieties.

These exterior differential systems~$\cI_\phi$ are 
usually very far from being involutive, and their integral manifolds
display varying degrees of rigidity, as will be explored in
this article.  In some cases, this rigidity 
permits a complete description of the integral manifolds and, hence,
a complete description of the effective $p$-cycles whose homology 
classes lie on the boundary of the halfspace~
$H^+\bigl([\phi]\bigr)\subset H_{2p}(M,\bbZ)$.

\begin{example}[The $n$-quadric]
\label{ex:the-even-quadrics}
A simple example will illustrate these ideas.  The proofs
will be taken up in \S\ref{ssec:quadrics-as-symmetric}.
Let~$(,)$ be the standard complex inner product on~$\C{n+2}$
and let~$Q_n\subset\bbP^{n+1}$ be the space of null lines for
this inner product, i.e.,~$[v]$ lies in~$Q_n$ for~$v\not=0$
in~$\C{n+2}$ if and only if~$(v,v)=0$.  Then~$Q_n$ is a compact
complex manifold of dimension~$n$.  It can also be regarded as
an Hermitian symmetric space:
\begin{equation}
Q_n = \frac{\SO(n + 2)}{\SO(2)\times\SO(n)}
\end{equation}
and so carries an $\SO(n{+}2)$-invariant K\"ahler structure~$\omega$.

When~$n$ is odd, $H_{2p}(Q_n,\bbZ) \simeq \bbZ$ for~$0\le p\le n$ and
there is a unique generator~$a_p\in H_{2p}(Q_n,\bbZ)$ on which
$[\omega^p]$ is positive.  It is not hard to see that
\begin{equation*}
[\cZ^+_p(Q_n)] = \bbZ^+\cdot a_p = H^+_{2p}(Q_n,\bbZ).
\end{equation*}

When~$n=2m$, one still has~$H_{2p}(Q_{2m},\bbZ) \simeq \bbZ$ 
for~$0\le p < m$ and~$m < p\le 2m$, 
but~$H_{2m}(Q_{2m},\bbZ) \simeq \bbZ^2$. 
A pair of generators of~$H_{2m}(Q_{2m},\bbZ)$ 
can be described as follows:
The subvariety~$N_{m+1}\subset \Gr(m{+}1,\C{2m+2})$ 
of isotropic $(m{+}1)$-planes in~$\C{2m+2}$ has two components, 
say~$N^\pm_{m+1}$, each of dimension~$\frac12 m(m{+}1)$.  
(These components are exchanged by any 
orientation reversing element of~$\Or(2m{+}2)$.)  
If~$W^\pm\subset\C{2m+2}$
are two isotropic~$(m{+}1)$-dimensional subspaces with~$W^\pm\in
N^\pm_{m+1}$,
then it is not difficult to show that the two $m$-cycles~$P_\pm 
= \bbP(W^\pm)\subset Q_{2m}$ have the property that~$[P_+]$ and~$[P_-]$
are a basis for $H_{2m}(Q_{2m},\bbZ)$.

Moreover, there exist $\SO(2m{+}2)$-invariant $2m$-forms~$\phi_\pm$
with the properties
\begin{enumerate}
\item $\phi_+$ and $\phi_-$ are closed and weakly positive;
\item $\omega^m = \phi_+ + \phi_-$; and
\item $\phi_\pm$ vanishes on~$P_\mp$.
\end{enumerate}

It follows that~$a = r_+\,[P_+] + r_-\,[P_-]$ 
lies in~$H^+_{2m}(Q_{2m},\bbZ)$
if and only if~$r_+$ and $r_-$ are nonnegative integers.  
Moreover, $\phi_-$
(respectively, $\phi_+$) must vanish on any effective $m$-cycle~$X$ 
that is homologous to~$r[P_+]$ (respectively, $r[P_-]$). 
In particular, the two classes~$[P_\pm]\in H_{2m}(Q_{2m},\bbZ)$ are
atomic and generate extremal rays.

It will be shown in~\S\ref{ssec:quadrics-as-symmetric}
that any $m$-cycle~$X\subset Q_{2m}$ on which~$\phi_-$ vanishes 
has the property that, 
for any smooth point~$x\in X$, there is a unique~$W^+_x\in N^+_{m+1}$ 
so that~$X$ and~$\bbP(W^+_x)$ are tangent at~$x$.  Further,
it will be shown that this implies, when~$X$ is irreducible, that~$X$
must actually be equal to~$\bbP(W^+_x)$ for some (and hence any) $x\in X$.
Consequently, for every integer~$r\ge1$,
\begin{equation}
\label{eq: extremal Zm in Q2m}
\cZ^+_m\bigl(Q_{2m},r[P_+]\bigr) = \bigl(N^+_{m+1}\bigr)^{(r)}
\end{equation}
Of course, the analogous formula holds for the classes~$r[P_-]$.  

Thus, the extremal classes in~$H^+_{2m}(Q_{2m},\bbZ)$ satisfy a 
strong form of rigidity.

The rigidity of subvarieties representing~$2[P_+]$ or $2[P_-]$ 
should be contrasted with the `flexibility' of 
the subvarieties~$X\subset Q_{2m}$ that satisfy~$[X] = [P_+]+[P_-]$.  
An example of such a variety is~$X = Q_{2m}\cap P^{m+1}$, where~$P^{m+1}$
is any $(m{+}1)$-dimensional linear projective subspace 
of~$\bbP^{2m+1}$.  Thus,  $\cZ^+_m\bigl(Q_{2m},[P_+]+[P_-]\bigr)$ contains%
\footnote{This containment is proper since, when $P_-$ and $P_+$
are chosen generically, 
the union ~$P_-\cup P_+$ does not lie in a~$P^{m+1}$.  Moreover, the
rigidity results above imply that the locus of reducible elements 
of~$\cZ^+_m\bigl(Q_{2m},[P_+]+[P_-]\bigr)$
is $N^-_{m+1}\times N^+_{m+1}$.  When $m$ is even, this is 
an irreducible component of~$\cZ^+_m\bigl(Q_{2m},[P_+]+[P_-]\bigr)$.}
 $\Gr(m{+}2,2m{+}2)$, 
a space of dimension~$m(m{+}2)$.  This dimension is greater than that of
$\cZ^+_m\bigl(Q_{2m},2[P_+]\bigr)$, which, being the symmetric square 
of~$N^+_{m+1}$, has dimension~$m(m{+}1)$. 
\end{example}

\subsubsection{Grassmannians}\label{sssec:background-grassmannians}
Much of this article will deal with the case~$M=\Gr(m,n)$, 
the Grassmannian of $m$-planes in~$\C{n}$.  This is a complex manifold 
of dimension~$N=m(n{-}m)$. Its homology groups are 
described as follows~\cite{MR50:7573,MR58:1235}:   

Let~$\sfP(m,n)$ denote the set of $m$-tuples~$\ab 
= (a_1,\ldots, a_m)$ where~$a_1,\ldots,a_m$ are integers satisfying
\begin{equation}
 n{-}m\ge a_1\ge a_2\ge\cdots\ge a_m\ge  0.
\end{equation}
Define~$|\ab| = a_1{+}\cdots{+}a_m$ and~$d(\ab) = m(n{-}m)-|\ab|\ge0$.
Let~$\sigma_\ab\subset\Gr(m,n)$ denote the set 
of~$m$-planes~$E\subset\C{n}$ that satisfy~
\begin{equation}
 \dim\bigl(E\cap \C{n-m+i-a_i}\bigr) \ge i.
\end{equation}
Then $\sigma_\ab$ is an irreducible complex subvariety of~$\Gr(m,n)$
of dimension~$d(\ab)$ that is known as the
\emph{Schubert variety} (or \emph{Schubert cycle}) of type~$\ab$.  

It is known~\cite[Chapter~0, \S5]{MR80b:14001} that the set
\begin{equation}
\sfP_p(m,n) =  \{\ [\sigma_\ab]\ \mid\ \ab\in \sfP(m,n),\ d(\ab) = p\ \}
\end{equation}
is a basis for the free abelian group $H_{2p}\bigl(\Gr(m,n),\bbZ\bigr)$
and, that, furthermore, the semigroup generated by~$\sfP_p(m,n)$
is equal to $\left[\cZ^+_p\bigl(\Gr(m,n)\bigr)\right] 
= H^+_{2p}\bigl(\Gr(m,n),\bbZ\bigr)$.  Thus, each of the 
classes~$[\sigma_\ab]$ is atomic and each of the rays 
$R_\ab = \bigl\{r\,[\sigma_\ab]\,\mid\,r\in\bbZ^+\bigr\}$ is
extremal.

\begin{example}[The cycle~$\sigma_{(2)}$]
\label{ex:the-grassmannians-sigma2}
Consider the case~$\ab = (2)$.%
\footnote{As is standard practice, 
I will suppress trailing zeroes when~$m$
can be inferred from context.  Thus, $\ab = (2)$ is an abbreviated
way of writing~$a_1 = 2$ and $a_j=0$ for~$j>1$.  
Note, though, that $m$ is needed to compute the dual $\ab^*$
(defined in~\S\ref{ssec:Schubert-cycles}), 
so some care must be taken with this shorthand.}
The cycle~$\sigma_{(2)}$ consists of the~$m$-planes~$V\subset\C{n}$ 
that meet~$\C{n-m-1}$ in at least a line.  In other words~$V$ lies 
in~$\sigma_{(2)}$ if and only if~$\bbP V\cap \bbP^{n-m-2}\not=\emptyset$.
Note that~$\sigma_{(2)}$ has codimension~$2$ in~$\Gr(m,n)$.

More generally, for any subset~$Y\subset\bbP^{n-1}$, define
\begin{equation}
\Psi_m(Y) 
= \bigl\{\ V\in \Gr(m,n) \mid \bbP V \cap Y\not=\emptyset\ \bigr\}.
\end{equation}
If~$Y$ is an algebraic subvariety of~$\bbP^{n-1}$
of dimension~$n{-}m{-}2$ and degree~$r$, 
then it is easy to see that~$\Psi_m(Y)$ 
is an algebraic subvariety of~$\Gr(m,n)$ of codimension~$2$ that 
satisfies~$\bigl[\Psi_m(Y)\bigr] = r[\sigma_{(2)}]$. 
 
Conversely, by Theorem~\ref{thm:codim2-sigma2},
any codimension~$2$ subvariety~$X\subset\Gr(m,n)$ satisfying~$[X]
= r[\sigma_{(2)}]$ is of the form~$X=\Psi_m(Y)$ for some
algebraic variety~$Y\subset\bbP^{n-1}$ of dimension~$n{-}m{-}2$ 
and degree~$r$. 

Thus, when~$\Psi_m$ is extended additively to a
semigroup homomorphism~$\Psi_m:\cZ^+_{n-m-2}(\bbP^{n-1})
\to \cZ^+_{m(n-m)-2}\bigl(\Gr(m,n)\bigr)$ in the obvious way,
the map
$$
\Psi_m: \cZ^+_{n-m-2}\bigl(\bbP^{n-1},r\bigl[\bbP^{n-m-2}\bigr]\bigr)
    \longrightarrow
 \cZ^+_{m(n-m)-2}\bigl(\Gr(m,n),r[\sigma_{(2)}]\bigr)
$$
is a bijection for all~$r\ge1$.
\end{example}

\begin{example}[The cycle~$\sigma_{(1,1)}$]
\label{ex:the-grassmannians-sigma11} 
The cycle~$\sigma_{(1,1)}\subset\Gr(m,n)$ is the set of~$m$-planes~$V$
that meet~$\C{n-m+1}$ in a subspace of dimension at least~$2$.
Equivalently~$\sigma_{(1,1)}$ is the union of the~$\Gr(m,\xi)$
where~$\xi\in\Gr(n{-}1,n)$ is any hyperplane that contains~$\C{n-m+1}$.
The set of such hyperplanes is a~$\bbP^{m-2}$ 
in~$\Gr(n{-}1,n)\simeq\bbP^{n-1}$.

More generally, for any subset~$Y^*\subset\Gr(n{-}1,n)\simeq\bbP^{n-1}$, 
define
\begin{equation}
\Sigma_m(Y^*)
= \bigcup_{\xi\in Y^*}\Gr(m,\xi) \subset\Gr(m,n).
\end{equation}

If~$Y^*$ is an algebraic subvariety of~$\Gr(n{-}1,n)$ 
of dimension~$m{-}2$ and degree~$r$, 
then~$\Sigma_m(Y^*)$ is an algebraic subvariety of~$\Gr(m,n)$ 
of codimension~$2$ that satisfies~$[\Sigma_m(Y^*)]=r[\sigma_{(1,1)}]$.

Conversely, by Theorem~\ref{thm:codim2-sigma11}, any
codimension~$2$ algebraic variety~$X\subset\Gr(m,n)$ 
that satisfies~$[X] = r[\sigma_{(1,1)}]$ is of the form~$X = \Sigma_m(Y^*)$
for some algebraic subvariety~$Y^*\subset\Gr(n{-}1,n)$ 
of dimension~$m{-}2$ and degree~$r$. 

Thus, when~$\Sigma_m$ is extended additively to a
semigroup homomorphism~$\Sigma_m:\cZ^+_{m-2}\bigl(\Gr(n{-}1,n)\bigr)
\to \cZ^+_{m(n-m)-2}\bigl(\Gr(m,n)\bigr)$ in the obvious way,
the map
$$
\Sigma_m: \cZ^+_{m-2}\bigl(\bbP^{n-1},r\bigl[\bbP^{m-2}\bigr]\bigr)
    \longrightarrow
 \cZ^+_{m(n-m)-2}\bigl(\Gr(m,n),r[\sigma_{(1,1)}]\bigr)
$$
is a bijection for all~$r\ge1$. 

Thus, for example, 
any~$X\in\cZ^+_{2n{-}6}\bigl(\Gr(2,n),r[\sigma_{(1,1)}]\bigr)$ 
is of the form
\begin{equation}
X = \Gr(2,\xi_1)+\cdots+\Gr(2,\xi_r)
\end{equation}
for some unique~$\xi_1,\ldots,\xi_r\in\Gr(n{-}1,n)$.  
Of course, such a~$X$ will be singular when~$r>1$ and~$n\ge4$.
\end{example}

\subsection{Notation}\label{ssec:notation}
The notation used in this article is mostly standard.  The space 
of~$n$-by-$m$ matrices with complex entries will be denoted 
by~$\C{n,m}$ and it will be endowed with the Hermitian 
inner product~$\la \ub,\vb\ra = \tr\bigl(\ub^*\,\vb\bigr)$, 
where $\ub^*$ denotes the conjugate transpose of~$\ub$.  
As usual, $\C{n,1}$ will be abbreviated to~$\C{n}$.  When~$m<n$, 
I will regard~$\C{m}$ as the subspace of~$\C{n}$ defined by setting
the bottom~$n{-}m$ entries equal to zero.

Unless otherwise specified,
all projective spaces and Grassmannians
are meant to be taken in the complex category.  
For any nonzero vector~$v$
in a vector space~$V$, the symbol~$[v]$ denotes the line~$\bbC v$ 
spanned by~$v$. The Grassmannian of $m$-dimensional
subspaces of a vector space~$V$ will be denoted~$\Gr(m,V)$,
with the shorthand notation~$\Gr(m,n)$ for~$\Gr(m,\C{n})$.
The space~$\Gr(m,V)$ will frequently be identified with
the projectivization of the cone of simple $m$-vectors in~$\L^m(V)$.
As usual,~$\Gr(1,V)$ will be denoted~$\bbP V$ and $\Gr(\dim V{-}1,V)$
will be identified with~$\bbP V^*$.

If~$W\subset V$ is a pair of vector spaces, then $\Gr(m,W)$ 
is a submanifold of~$\Gr(m,V)$ in the obvious sense.
If $m$ satisfies $\dim W\le m\le\dim V$,
then $[W,V]_m$ will denote the set of~$E\in\Gr(m,V)$
satisfying~$W\subset E\subset V$.  This space is a Grassmannian in
its own right, naturally identified 
with~$\Gr\bigl(m{-}\dim W, V/W\bigr)$.

If~$V$ is a vector space and $E\subset V$ is a subspace, 
then~$\cdot_{|_E}:\L(V^*)\to\L(E^*)$  denotes the induced pullback 
homomorphism.  I.e., for~$\phi\in\L^m(V^*)$, its pullback to~$E$
will be denoted~$\phi_{|_E}\in\L^k(E^*)$.

\subsection{Orientation and positivity}\label{ssec:positivity}
A complex vector space~$E$ of dimension~$n$ 
carries a canonical orientation,
namely, the one for which~$(\eb_1,\iC\,\eb_1,\dots,
\eb_n,\iC\,\eb_n)$ is a positively oriented $\bbR$-basis of~$E$ 
whenever~$(\eb_1,\dots,\eb_n)$ is a $\bbC$-basis of~$E$.

\begin{definition}[Weak positivity]
\label{def: weak pos}
A real-valued $(p,p)$-form~$\phi$ on a complex vector space~$V$ 
of dimension~$n\ge p$ is \emph{weakly positive} \cite{MR50:7573} 
if $\phi(\eb_1,\iC\,\eb_1,\dots,\eb_p,\iC\,\eb_p)\ge0$
for all $\eb_1,\dots,\eb_p\in V$. 
\end{definition}
 
If~$\zeta$ is any $(p,0)$-form on~$V$,
the real-valued $(p,p)$-form $\phi = \iC^{p^2}\zeta\w\ov{\zeta}$
is weakly positive, as is any sum of the form
\begin{equation}\label{eq:positive}
\phi = \iC^{p^2} \sum_{k=1}^K \zeta_k\w\ov{\zeta_k}
\end{equation}
where~$\zeta_1,\dots,\zeta_K$ are $(p,0)$-forms on~$V$.

\begin{definition}[Positivity]
A $2p$-form $\phi$ that can be expressed
in the form~\eqref{eq:positive} is said to be \emph{positive}.
\end{definition}

\begin{remark}[Usage caveats]
This terminology can be misleading, since, for example, $\phi=0$ 
is positive according to this definition. 
In~\cite[p.~401]{MR81k:53004}, Griffiths and Harris adopted the
more suggestive terminology `non-negative' for what had, until then,
been called `positive'.  For various reasons, I have not followed suit.

Compare Griffiths~\cite{MR41:2717} and Harvey and Knapp~\cite{MR50:7573}.
To be accurate, Harvey and Knapp define positivity somewhat
differently, but prove that their definition is equivalent to the
one given above~\cite[Theorem 1.2]{MR50:7573}.  The reader should also 
keep in mind that the term `positive' in reference to~$(p,p)$-forms 
is used somewhat differently by some authors.  
See~\cite{MR50:7573} for a thorough discussion.

When~$2\le p\le n{-}2$, Harvey and Knapp show that there exist 
weakly positive forms that are not positive, so the two concepts 
really are different.
\end{remark}

Positivity is preserved under pullback and, if~$S:W\to V$
is a linear surjection, then~$S^*\phi$ is positive if 
\emph{and only if}~$\phi$ is positive.  Moreover, the wedge product of
positive forms is positive, as is the sum.

\subsubsection{The zero locus and ideal of a positive form}
\label{sssec: zero locus and ideal}
Given~$\phi$ of the form~\eqref{eq:positive},
the linear span of the forms~$\zeta_1,\dots,\zeta_K$ in~$\L^{(p,0)}(V)$ 
is well-defined, even though the representation~\eqref{eq:positive} 
is not unique~\cite[Theorem 1.2]{MR50:7573}.  Let~$\I_\phi\subset
\L^{(*,0)}(V)$ denote the ideal generated by~$\zeta_1,\dots,\zeta_K$.

It follows that a given complex $p$-plane~$E\in\Gr(p,V)$ 
satisfies~$\phi_{|_E}>0$ unless each of the~$\zeta_1,\dots,\zeta_K$ 
vanishes on~$E$, i.e., unless $E$ is an integral element of~$\I_\phi$.

Thus, when~$\phi$ is positive, the locus~$Z(\phi)\subset\Gr(p,V)$ 
consisting of the $p$-planes on which~$\phi$ vanishes 
is a complex subvariety of~$\Gr(p,V)$.  (This is not generally the case 
for forms~$\phi$ that are only weakly positive.)  

The subvariety~$Z(\phi)$ can be singular and/or reducible, as will be seen. 

If~$\phi$ is positive, with a representation as in~\eqref{eq:positive} 
and~$\psi$ is a positive $(q,q)$ form, with a representation of the form
\begin{equation}
\label{eq:positive psi}
\psi = \iC^{q^2} \sum_{j=1}^J \eta_j\w\ov{\eta_j},
\end{equation}
then
\begin{equation}
\label{eq:positive wedge}
\psi\w\phi = \iC^{(p+q)^2} \sum_{j=1}^J \sum_{k=1}^K 
       \bigl(\eta_j\w\zeta_k\bigr)\w\ov{\bigl(\eta_j\w\zeta_k\bigr)}\,,
\end{equation}
so that, not only is~$\psi\w\phi$ positive, 
but one also has the equality of ideals
\begin{equation}
\label{eq: ideals wedge}
\I_{\psi\w\phi} = \I_\psi\w\I_\phi\,.
\end{equation}  

\subsubsection{The generalized Wirtinger inequality}
\label{sssec: generalized Wirtinger}
The positive definite Hermitian inner products~$\la,\ra$ on~$V$ are 
in one-to-one correspondence 
with the positive $(1,1)$-forms~$\omega$ on~$V$
that are non-zero on each line.  As is shown in~\cite{MR50:7573}, a choice
of such an~$\omega$ defines several norms on the space of real-valued
$(p,p)$-forms.  For one of these norms, 
denoted~$|\!|\cdot|\!|_1$ in~\cite{MR50:7573}, Harvey and Knapp prove their
\emph{Generalized Wirtinger Inequality} \cite[Theorem 1.8\,(b)]{MR50:7573}:
\begin{equation}\label{eq:Generalized Wirtinger Inequality}
\phi\w\frac{\omega^{n-p}}{(n{-}p)!} 
\le  |\!|\phi|\!|_1\,\frac{\omega^n}{n!}
\end{equation}
and show that equality 
in~\eqref{eq:Generalized Wirtinger Inequality} holds
if and only if~$\phi$ is positive.  
In particular, if~$\phi$ is positive, 
then~$\phi\w\omega^{n-p}\ge0$, with equality only for~$\phi=0$.

As an application of this fact, consider a compact complex 
manifold~$M$ endowed with a K\"ahler structure~$\omega$.  
If~$\phi\in\Omega^{p,p}(M)$ is a nonzero, closed, positive $(p,p)$-form, 
then its cohomology class $[\phi]\in H^{p,p}(M,\bbR)$ is nonzero, since
$$
\int_M \phi\w\frac{\omega^{n-p}}{(n{-}p)!} =
\int_M  |\!|\phi|\!|_1\,\frac{\omega^n}{n!} > 0.
$$
This motivates the following definition:  A class~$a\in H^{p,p}(M,\bbR)$
will be said to be \emph{positive} if $a = [\phi]$ for 
some closed positive $(p,p)$-form~$\phi$.  Denote the set of positive 
classes by~$H^{p,p}_+(M,\bbR)\subset H^{p,p}(M,\bbR)$.  Since the positive
forms are closed under addition and scalar multiplication by non-negative
numbers,~$H^{p,p}_+(M,\bbR)$ is a convex cone 
in~$H^{p,p}(M,\bbR)$ that (except for~$0$) lies strictly 
on one side of the hyperplane 
$$
H^{p,p}_\omega(M,\bbR) 
=\{\>a\in H^{p,p}(M,\bbR) \mid a\cup[\omega^{n-p}]=0\>\}.
$$
Note that the cone~$H^{p,p}_+(M,\bbR)$ does not depend on the choice
of K\"ahler structure~$\omega$, even though~$H^{p,p}_\omega(M,\bbR)$ does.

\section[Grassmannians]{Geometry of Grassmannians}
\label{sec:grass}

To avoid trivialities, assume that~$0<m<n$ throughout this section.

\subsection{Partition posets and their operations}
\label{ssec: partition posets}
This is as good a place as any to collect the basic definitions 
and properties of partitions that will be needed in what follows.

\begin{definition}(The partition poset)
\label{def: P(m,n)}
$\sfP(m,n)$ is the set of \emph{partitions} 
$\ab = (a_1,\dots,a_m)$ where the~$a_i$ are integers satisfying%
\footnote{For notational convenience, adopt the 
convention, for~$\ab\in\sfP(m,n)$, that~$a_0 = n{-}m$ and~$a_{m+1}=0$.}
$$
n-m\ge a_1 \ge a_2 \ge \dots \ge a_m \ge 0.
$$
Set~$|\ab| = a_1 + \dots + a_m$ and~$d(\ab) = m(n{-}m)-|\ab|$.  

Also, define~$\ab\le\bb$ to mean that~$a_i\le b_i$ for~$1\le i\le m$.
When~$\ab \le \bb$, a \emph{chain} from~$\ab$ to~$\bb$ is a set of
elements~$\ab_p,\ab_{p+1},\dots,\ab_q\in \sfP(m,n)$ where~$|\ab_k|=k$ 
for~$p\le k\le q$, and
$$
\ab = \ab_p \le \ab_{p+1} \le \dots \le \ab_q = \bb.
$$
The number of distinct chains from~$\ab$ to~$\bb$ will be 
denoted~$\mu^\bb_\ab$.  By definition~$\mu^\bb_\ab = 0$ 
if $\ab \not\le \bb$.
\end{definition} 

The relation~$\le$ is a partial order on~$\sfP(m,n)$.  As
examples, the Hasse diagrams%
\footnote{See Remark~\ref{rem: Hasse diagrams}
in~\S\ref{ssec:hss-ideal-poset} for an
explanation of how these diagrams depict the poset structure.} 
of~$\sfP(3,6)$ and~$\sfP(2,5)$
are to be found in Figures~\ref{fig:Gr36poset} and~\ref{fig:Gr25poset}, 
respectively.  (The labeling above the nodes will be explained
later.)  This poset structure is sometimes referred to as the
\emph{Bruhat poset} of~$\Gr(m,n)$.

When~$m$ and~$n$ can be inferred from context, I will usually 
suppress trailing zeroes, e.g., writing~$(2,1)$ 
to denote~$(2,1,0)\in\sfP(3,6)$.  The partition~$(q,0,\dots,0)$
will often be denoted more simply by~$q$ when this will not cause 
confusion.

There are two operations on the partition posets 
that will be needed. 

\begin{definition}[Dual and conjugate partitions]
\label{def: dual and conjugate}
For each~$\ab=(a_1,\ldots,a_m)\in \sfP(m,n)$ define its 
\emph{dual partition}~$\ab^*\in \sfP(m,n)$ by
\begin{equation}
\ab^* = (n{-}m{-}a_m,\,n{-}m{-}a_{m-1},\ldots,n{-}m{-}a_1).
\end{equation}
For~$\ab = (a_1,\dots,a_m)\in\sfP(m,n)$, 
the \emph{conjugate partition} $\ab'\in\sfP(n{-}m,n)$ 
is defined as follows:
Set~$a_0 = n{-}m$ and then, when $1\le a\le n{-}m$,
set $\ab'_a = j$ where $j\in\{0,\dots,m\}$ 
is the largest integer for which~$a_j\ge a$.
\end{definition}

One can show that~$(\ab')' = \ab$
and~$|\ab'| = |\ab|$.  One also has~$(\ab^*)' = (\ab')^*$,
$(\ab^*)^* = \ab$, and $d(\ab) + d(\ab^*) = m(n{-}m)$. 
Moreover, $\ab \le \bb$ is equivalent to~$\bb^*\le \ab^*$
and~$\ab'\le\bb'$.

For more on the conjugate 
construction as well as its interpretation in terms of 
Young tableaux, see \cite[p.~45]{MR93a:20069}.

\subsection{Schubert cycles}\label{ssec:Schubert-cycles}

The group~$\SL(n,\bbC)$ acts transitively on~$\Gr(m,n)$ on the left
in the usual way:  $A\cdot E = A(E)\subset\C{n}$ for~$A\in\SL(n,\bbC)$
and~$E\in\Gr(m,n)$.  When~$0<m<n$, the stabilizer of~$\C{m}\in\Gr(m,n)$ 
is a maximal parabolic subgroup that will be 
denoted~$P_m\subset \SL(n,\bbC)$.  For notational convenience, set
$P_0 = P_n = \SL(n,\bbC)$.

The definition of Schubert cycles in a Grassmannian was already given
in~\S\ref{sssec:background-grassmannians},
but it is convenient to generalize this definition slightly 
and it will be necessary to discuss the geometry of these cycles 
in a bit more detail.  For proofs of the statements in this subsection, 
see~\cite[Chapter 1, Section 5]{MR80b:14001}.

Let~$V$ be a complex vector
space of dimension~$n$.  A \emph{flag}~$F$ in~$V$ is a nested sequence
of vector spaces~$V_i\subset V_{i+1}$ for~$0\le i\le n$ with~$\dim V_i=i$
and~$V_n = V$.  Since the connected group~$\SL(V)\simeq\SL(n,\bbC)$
acts transitively on the set of flags in~$V$, the choice of flag will
not materially affect the constructions to be made below.

For~$\ab = (a_1,\ldots,a_m)\in \sfP(m,n)$, 
the \emph{Schubert cell}~$W_\ab(F) \subset\Gr(m,V)$ is, 
by definition, the set of~$E\in\Gr(m,V)$ satisfying
\begin{equation}
i = \dim(E\cap V_{n-m+i-a_i}) > \dim(E\cap V_{n-m+i-a_i-1}) 
\end{equation}
for~$1\le i\le m$.  For~$\ab\not=\bb$ in~$\sfP(m,n)$, the sets~$W_\ab(F)$
and~$W_\bb(F)$ are disjoint and the union 
of the~$W_\ab(F)$ as~$\ab$ ranges
over~$\sfP(m,n)$ is the whole of~$\Gr(m,V)$.

When~$V=\C{n}$ and $F$ is the standard flag, i.e.,~$V_i = \C{i}\subset\C{n}$ 
for all~$i$, then $W_\ab(F)$ will be denoted~$W_{\ab}$.

The set~$W_\ab(F)$ is a subvariety of~$\Gr(m,V)$ 
that is biholomorphic with~$\C{d(\ab)}$. In fact, the description of
Schubert cells in~\cite[pp.~195--6]{MR80b:14001} shows that there exists
a closed, nilpotent subgroup of~$\SL(V)$ that acts simply transitively
on~$W_\ab(F)$.  I will need a description of this subgroup later, 
so I give it here:  

\begin{definition}[Subspace type]
\label{def:n-of-type-a}
For any partition~$\ab\in \sfP(m,n)$, let~$\eun_\ab\subset \C{n-m,m}$ 
denote the vector space of matrices~$\textsf{Z} = (z^a_i)$ 
that satisfy~$z^a_i=0$ when~$a>n{-}m-a_i$. 
(Note that~$\eun_\ab$ has dimension~$d(\ab)$.)

A subspace~$A\subset\C{n-m,m}$
is \emph{of type~$\ab$} if $A = q\,\eun_\ab\,s^{-1}$ for
some~$q\in \GL(n{-}m,\bbC)$ and $s\in\GL(m,\bbC)$.

More generally, if $Q$ and~$S$ are vector spaces of dimensions~$n{-}m$ 
and~$m$ respectively, a subspace~$A\subset Q\ot S^*$ will be said to be of 
\emph{type~$\ab$} if there exist isomorphisms~$q:\C{n-m}\to Q$ 
and~$s:\C{m}\to S$ so that $A = q{\ot}(s^{-1})^*(\eun_\ab)$.

Let~$P_\ab\subset\GL(n{-}m,\bbC){\times}\GL(m,\bbC)$ be the
pairs~$(q,s)\in\GL(n{-}m,\bbC){\times}\GL(m,\bbC)$ that satisfy~$\eun_\ab
= q\,\eun_\ab\,s^{-1}$.
\end{definition}  


Now, let~$N_\ab\subset\SL(n,\bbC)$
be the abelian nilpotent subgroup defined by
\begin{equation}
N_\ab = \left\{\ \begin{pmatrix}\I_m & 0\\ \textsf{Z} & \I_{n-m}
\end{pmatrix}
       \  \vrule\ \textsf{Z}\in \eun_\ab \right\}.
\end{equation}
Then~$N_\ab\cdot \C{m}\subset\Gr(m,n)$ is a 
Schubert cell~$W_\ab(F)$ for some flag~$F$ (that depends on~$\ab$).

\begin{remark}[Closure of type]
\label{rem: closure of type}
Since~$P_\ab$ contains the pairs~$(q,s)$ where~$q$ and~$s$ are
upper triangular matrices, it is a parabolic subgroup 
of~$\GL(n{-}m,\bbC){\times}\GL(m,\bbC)$.  It follows that, 
for any~$Q$ and~$S$, the subspaces of~$Q{\otimes}S^*$ of type~$\ab$ 
form a closed $\bigl(\GL(Q){\times}\GL(S)\bigr)$-orbit
in~$\Gr\bigl(d(\ab),Q{\otimes}S^*\bigr)$.  
This fact will be useful in~\S\ref{ssec: Walters diff sys}.  

In fact, the $\bigl(\GL(Q){\times}\GL(S)\bigr)$-orbit of a 
subspace~$A\in \Gr\bigl(d,\,Q{\otimes}S^*\bigr)$
is closed only when~$A$ has type~$\ab$ for some~$\ab$.   The reason for
this is simple:  If the orbit of~$A$ is closed, then its 
stabilizer~$P_A\subset \GL(Q){\times}\GL(S)$ must
be a parabolic subgroup of~$\GL(Q){\times}\GL(S)$.  
Every parabolic subgroup contains a Borel subgroup 
and all Borel subgroups are conjugate, so there is 
a~$\eun\in\Gr\bigl(d,\,\C{n-m,m}\bigr)$ such 
that~$A = q{\otimes}(s^{-1})^*(\eun)$ for some isomorphisms~$q:\C{n-m}\to Q$
and~$s:\C{m}\to S$ and so that~$\eun$ is stable under the action of the
Borel subgroup consisting of the pairs of upper triangular matrices
in~$\GL(n{-}m,\bbC){\times}\GL(m,\bbC)$.  It is easily proved that the
only subspaces of~$\C{n-m,m}$ that are stable under this Borel subgroup
are the subspaces of the form~$\eun_\ab$ for~$\ab\in\sfP(m,n)$.
\end{remark}

The closure~$\sigma_\ab(F) = \ov{W_\ab(F)}\subset\Gr(m,V)$ is an 
irreducible variety of dimension~$d(\ab)$, known as the \emph{Schubert cycle}
or \emph{Schubert variety} of type~$\ab$ associated to the flag~$F$.

Note that~$\sigma_\ab(F) = \sigma_\ab(F')$ if~$V_i = V'_i$ for
all~$i$ such that~$a_i > a_{i+1}$.  In other words, $\sigma_\ab(F)$
frequently depends only on partial flag information.  As usual, when~$F$
is the standard flag, one simply writes~$\sigma_\ab$.

Since the connected Lie group~$\SL(V)\simeq\SL(n,\bbC)$ acts transitively 
on the space of flags in~$V$, the homology class~$[\sigma_{\ab}(F)]$ is 
independent of the choice of~$F$ and will usually just be written 
as~$[\sigma_\ab]$.  

The classes~$\bigl\{[\sigma_\ab]\>\vrule\>d(\ab) 
= p\,\bigr\}$ form a basis for~$H_{2p}\bigl(\Gr(m,n),\bbZ\bigr)$ 
as a free abelian group and one has the homology intersection pairing
\begin{equation}
[\sigma_{\ab}]\cap[\sigma_{\bb^*}] = \delta^\bb_\ab\,.
\end{equation}

\begin{remark}[Singularity of Schubert cycles]
\label{rem: sing of Schubert}
For most~$n$, $m$, and~$\ab\in \sfP(m,n)$, 
the Schubert cycle~$\sigma_\ab$ is singular.  In fact, 
as shown in~\cite{MR85j:14095}, $\sigma_\ab$ is singular
unless~$\ab^* = (p,\dots,p)$ for some~$p$ with~$0\le p\le n{-}m$.
(As usual, I suppress trailing zeroes, so the length of~$\ab^*$
can be anywhere from~$0$ to~$m$.)  When~$\ab^*=(p,\dots,p)$ has
length~$q$, then~$\sigma_\ab=[W_{-},W_{+}]_m\subset\Gr(m,n)$,
where~$W_{-}$ has dimension~$m{-}q$, $W_{+}$ has dimension~$m{+}p$,
and~$W_{-}$ is a subspace of~$W_{+}$.  
Thus,~$\sigma_\ab\simeq\Gr(q,p{+}q)$ in this case.
\end{remark}

\begin{remark}[$A$-cycles]
\label{rem: A-cycles}
The Schubert cycles can be generalized in a way that 
will be used later on to produce examples of 
subvarieties of~$\Gr(m,n)$ that satisfy certain differential
systems or homological conditions, so I will describe it here.
If~$A\subset \C{n-m,m}$ is any complex subspace of dimension~$d$, 
define~$N_A\subset\SL(n,\bbC)$ to be the abelian nilpotent subgroup
\begin{equation}
N_A = \left\{\ \begin{pmatrix}\I_m & 0\\ \textsf{Z} & \I_{n-m}
\end{pmatrix}
       \  \vrule\ \textsf{Z}\in A \right\}.
\end{equation}
Then~$N_A{\cdot}\C{m}\subset\Gr(m,n)$ is biholomorphic to~$\C{d}$.
The closure~$\sigma_A = \overline{N_A{\cdot}\C{m}}$ is the image
of a rational map of~$\bbP^d$ into~$\Gr(m,n)$ and hence is an
irreducible, $d$-dimensional, algebraic subvariety of~$\Gr(m,n)$
\cite[pp.~492-3]{MR80b:14001}.  Note that, while~$\sigma_A$ will
generally be singular, it is `quasi-homogeneous', in the sense
that it contains a Zariski-open subset, namely~$N_A{\cdot}\C{m}$
that is homogeneous under a subgroup of~$\SL(n,\bbC)$.

What is not so obvious is how one expresses the homology 
class~$[\sigma_A]$ in terms of the homology classes of the Schubert
cycles.  For example, when~$\dim A = 1$, then one easily sees
that~$[\sigma_A] = r\,[\sigma_{(1)^*}]$ where~$r>0$ is the 
rank of a generator of~$A$.  

When~$\dim A > 1$, the homology class~$[\sigma_A]$ is more difficult 
to compute.
However, one can say that, for the generic~$A\in\Gr(d,\C{n-m,m})$, 
the class~$[\sigma_A]$ is a linear combination with strictly 
positive coefficients of all of the~$[\sigma_{\ab^*}]$ with~$|\ab|=d$.
This is because each of the corresponding forms~$\phi_\ab$ will be 
nonzero on the Zariski open subset~$N_A{\cdot}\C{m}\subset\sigma_A$.

I will usually refer to any subvariety of~$\Gr(m,n)$ that
is equivalent to~$\sigma_A$ under the action of~$\SL(n,\bbC)$
as an \emph{$A$-cycle}.  For any subspace~$B\subset\C{n-m,m}$, 
the cycle~$\sigma_B$ is an $A$-cycle when there are~$q\in\GL(n{-}m,\bbC)$
and~$s\in\GL(m,\bbC)$ so that~$B = q\,A\,s^{-1}$.  In this case,
one also says that~$B$ is a subspace of \emph{type~$A$}.   

More generally, if~$Q$ and~$S$ are vector spaces,
a subspace~$B\subset Q\otimes S^* = \Hom(S,Q)$ 
is said to be \emph{of type~$A$} 
if there are isomorphisms~$i:\C{n-m}\to Q$
and~$\pi:S\to\C{m}$ so that
$$
B = \left\{\,q\circ a\circ \pi \mid a\in A\,\right\}
$$
(where the elements of~$A$ are regarded as linear maps from~$\C{m}$ 
to~$\C{n-m}$).

\end{remark}

\subsection{The canonical bundles}\label{ssec: canonical bundles over
Gr(m,n)}
The trivial bundle~$\Gr(m,V)\times V$ over~$\Gr(m,V)$ contains
the subbundle~$S$ of rank~$m$ that consists of the 
pairs~$(E,v)$ with~$v\in E$.  The quotient construction then defines
a canonical bundle~$Q$ over~$\Gr(m,V)$ of rank~$n{-}m$
 whose fiber over~$E$
is canonically isomorphic to~$V/E$.  These fit into the exact sequence
\begin{equation}\label{eq:exact-seq}
0\longrightarrow S\longrightarrow \Gr(m,V)\times V
\longrightarrow Q\longrightarrow0.
\end{equation}
Moreover, there is a canonical bundle isomorphism
\begin{equation}
T\Gr(m,V) = Q\ot S^*
\end{equation}
corresponding to the canonical 
isomorphism~$T_E\Gr(m,V)\simeq V/E\otimes E^*$.

When~$0\le q\le m(n{-}m)$, there is a canonical 
decomposition of the $q$-th (complex) exterior power
of the cotangent bundle of the form~\cite[p.~80]{MR93a:20069}
\begin{equation}\label{eq: Gr ext alg decomp}
\L^{q,0}\bigl(T^*\Gr(m,V)\bigr) = \L^{q}(S\ot Q^*) 
= \bigoplus_{\substack{\ab\in \sfP(m,n)\\ |\ab|=q}}
\bbS_\ab(S)\ot\bbS_{\ab'}(Q^*)
\end{equation}
where~$\bbS_\bb$ denotes
the Schur functor associated to the partition~$\bb$
in the category of vector spaces and linear 
maps~\cite[Lecture~6]{MR93a:20069}.
The formula~\eqref{eq: Gr ext alg decomp} 
seems to be due to Ehresmann~\cite{cE34}.

For example, 
\begin{equation}
\begin{split}
\L^{2,0}\bigl(T^*\Gr(m,V)\bigr) 
&= \bigl(\bbS_{(2)}(S)\ot\bbS_{(1,1)}(Q^*)\bigr)\oplus
     \bigl(\bbS_{(1,1)}(S)\ot\bbS_{(2)}(Q^*)\bigr)\\
&=\bigl(S^2(S)\ot\L^2(Q^*)\bigr)\oplus\bigl(\L^2(S)\ot S^2(Q^*)\bigr),
\end{split}
\end{equation}
and, as long as~$2\le m\le n{-}2$, both of these summands will 
be nontrivial.  

\begin{definition}[The ideal~$\cI_\ab$]
\label{def:ideal-a}
For~$\ab\in \sfP(m,n)$, let~$\cI_\ab$ denote the exterior ideal on~$\Gr(m,n)$ 
generated by the sections of the subbundle 
$$
\I_\ab = \bbS_\ab(S)\ot\bbS_{\ab'}(Q^*)
\subset\L^{|\ab|,0}\bigl(T^*\Gr(m,n)\bigr).
$$ 
\end{definition}

The ideal~$\cI_\ab$ is invariant under the action of~$\SL(n,\bbC)$.
It is not hard to see that~$\cI_{\ab}$ is holomorphic and
differentially closed.  This will be proved below 
(see Proposition~\ref{prop: I-a is holo and closed}), when a
different description of~$\cI_{\ab}$ is given.

\subsection{Chern classes}\label{sec:chern-classes}

Let~$c(Q)$ and $c(S)$ denote, respectively, the total Chern classes
of the canonical quotient bundle and subbundle over~$\Gr(m,n)$.  In view
of~\eqref{eq:exact-seq}, these satisfy~$c(Q)c(S) = c(Q\oplus S) = 1$.
Writing~$c(Q) = 1 + q_1 + \dots+ q_{n-m}$ and~$c(S) =1+s_1+\dots + s_m$
with~$s_j,q_j\in H^{2j}\bigl(\Gr(m,n),\bbZ\bigr)$, this gives the relation
\begin{equation}\label{eq:whitney-sum}
(1 + s_1 + \dots + s_m)(1 + q_1 + \dots+ q_{n-m}) = 1,
\end{equation}
which allows one to compute the $s_i$ recursively in terms of the $q_j$
(or vice versa).  For example, $s_1 = -q_1$, $s_2 = {q_1}^2-q_2$, etc. 
In fact, comparing like degrees on both sides for degrees between $0$
and~$m$ gives a recursive formula for $c(S)$ in terms of $c(Q)$
and then the remaining degrees between $m{+}1$ and~$n$ yield graded 
polynomial relations~$R_{m+1}(q)= \dots = R_n(q) = 0$ on the~$q_j$.  

It is well-known~\cite{MR55:13428} that the classes~$q_1,\dots,q_{n-m}\in 
H^*\bigl(\Gr(m,n),\bbZ\bigr)$ generate the 
ring~$H^*\bigl(\Gr(m,n),\bbZ\bigr)$, i.e., that this ring 
is isomorphic to the polynomial ring on the classes~$q_1,\dots,q_{n-m}$
modulo the ideal generated 
by the relations~$R_{m+1}(q)= \dots = R_n(q) = 0$.

Certain polynomials in these classes, 
the so-called \emph{Schur classes}, 
will play an important role in this article.  These are defined 
for each~$\ab\in \sfP(m,n)$ by the Giambelli determinant formula:
\begin{equation}\label{eq:Schur-polys}
q_\ab = 
\begin{vmatrix}
q_{a_1} & q_{a_1+1} & \hdots & q_{a_1+m-1}\\
q_{a_2-1} & q_{a_2} & \hdots & q_{a_2+m-2}\\
\hdotsfor[2.0]4\\
q_{a_m-m+1} & q_{a_m-m+2} & \hdots & q_{a_m}\\
\end{vmatrix}.
\end{equation}
where, by convention,~$q_0=1$ and~$q_j=0$ unless~$0\le j\le n{-}m$.
These classes correspond naturally 
to the Schubert cycles~\cite[p.~205 and p.~411]{MR80b:14001}, 
i.e., using the natural pairing between cohomology and homology, 
they satisfy
\begin{equation}\label{eq:Schubert-cycle-pairing-relation}
\bigl\la q_\ab,[\sigma_{\bb^*}]\bigr\ra 
= [\sigma_\ab]\cap[\sigma_{\bb^*}] 
= \delta^\bb_\ab\,.
\end{equation}
Thus, $\{\,q_\ab \mid \ab\in \sfP(m,n),\>|\ab|=p\,\}$ is a basis
of the lattice~$H^{2p}\bigl(\Gr(m,n),\bbZ\bigr)$.

An explicit formula for the product~$q_\ab\,q_\bb$
in~$H^*\bigl(\Gr(m,n),\bbZ\bigr)$ is known, of course, as
this is the basis for the Schubert calculus.  However, I will not 
need to work with the full formula in what follows, only the simplest
\emph{Pieri formula}~\cite[p. 203]{MR80b:14001}:
\begin{equation}
\label{eq: Pieri 1}
q_1\,q_\ab 
= \sum_{\substack{\bb\in\sfP(m,n) \\ |\bb|= |\ab|+1 \\ \bb\ge \ab}} q_\bb
= \sum_{\substack{\bb\in\sfP(m,n) \\ |\bb|= |\ab|+1}} \mu^\bb_\ab\,q_\bb\,,
\end{equation}
which, by induction and the definition of~$\mu^\bb_\ab$, generalizes to
\begin{equation}
\label{eq: Pieri p}
(q_1)^p\,q_\ab 
= \sum_{\substack{\bb\in\sfP(m,n) \\ |\bb|= |\ab|+p}} \mu^\bb_\ab\,q_\bb\,.
\end{equation}

\subsection{$\Gr(m,n)$ as an Hermitian symmetric space}
\label{ssec:gr-mn-symmspace}

I will now briefly review the K\"ahler geometry of~$\Gr(m,n)$.
For more details, consult~\cite{MR80b:14001,MR81k:53004}. 

Let~$\SU(n)$ denote the group of special unitary 
$n$-by-$n$ matrices.  When a name is needed for the 
inclusion~$\SU(n)\subset\GL(n,\bbC)$, I will write it 
as~$\us:\SU(n)\to\GL(n,\bbC)$.  I will also write
\begin{equation}
\us = \begin{pmatrix}\us_1&\us_2&\cdots&\us_n\end{pmatrix}
\end{equation}
and regard each column as a function~$\us_k:\SU(n)\to\C{n}$. 

In what follows, the Hermitian summation convention will be assumed, 
i.e., when a subscript occurs both barred and unbarred in a single
term, a summation over that subscript is implied.  Adopt the index range 
conventions
$$
1\le i,j,k \le m < a,b,c\le n
$$
together with the comprehensive index range~$1\le A,B,C\le n$.  
The complex valued 1-forms~$\upsilon_{\bar{A}B}={\us_A}\!^*\,\d \us_B
=-\,\ov{\upsilon_{\bar{B}A}}$ satisfy the {\it structure equations}
\begin{equation}
\d \us_A = \us_B\,\upsilon_{\bar{B}A}\qquad\text{and}\qquad
\d \upsilon_{\bar{A}B} = -\upsilon_{\bar{A}C}\w\upsilon_{\bar{C}B}\,.
\end{equation}

 The map
\begin{equation}
\pi_m = [\us_1\w\cdots\w\us_{m}] : \SU(n)\to\Gr(m,n)
\end{equation}
makes~$\SU(n)$ into a principal right 
$S\bigl(\Un(m)\times\Un(n{-}m)\bigr)$-bundle over~$\Gr(m,n)$,
where~$S\bigl(\Un(m)\times\Un(n{-}m)\bigr)\subset\SU(n)$ is
the group of matrices of the form
\begin{equation}
\begin{pmatrix} A&0\\0&B\\ \end{pmatrix}
\end{equation}
with~$A\in\Un(m)$, $B\in\Un(n{-}m)$, and $\det A\,\det B = 1$.
In particular, $\Gr(m,n)$ is an Hermitian symmetric space
\begin{equation}
\Gr(m,n) = \SU(n)/S\bigl(\Un(m)\times\Un(n{-}m)\bigr).
\end{equation}

Write the left-invariant $\eusu(n)$-valued 
$1$-form~$\upsilon = \us^*\,\d \us 
= \us^{-1}\,\d \us = -\upsilon^*$
on~$\SU(n)$ in block form as
\begin{equation}
\upsilon = 
\begin{pmatrix} 
         \>\sigma&-\omega^*\\
         \>\omega&\theta\\ 
\end{pmatrix}
\end{equation}
where~$\sigma = -\sigma^*$ is $m$-by-$m$, $\omega$ is~$(n{-}m)$-by-$m$,
and $\theta = -\theta^*$ is $(n{-}m)$-by-$(n{-}m)$.  

By the structure equations, $d\theta + \theta\w\theta = \omega\w\omega^*$.  
For each~$1\le j\le n{-}m$, there is a unique 
form~$\phi_j\in\Omega^{j,j}\bigl(\Gr(m,n)\bigr)$ so that
\begin{equation}\label{eq:Chernforms-definition}
\det\left(I_{n-m} + \frac{\iC}{2\pi}\ \omega\w\omega^* \right)
 = \pi_m^*\bigl(1 + \phi_1 + \cdots + \phi_{n-m}\bigr).
\end{equation}
Each~$\phi_j$ is invariant under the action of~$\SU(n)$ and
satisfies~$ [\phi_j] = c_j(Q) = q_j $.  

It is well-known \cite{MR58:1235}
that the forms~$\phi_1,\ldots,\phi_{n-m}$
generate the ring of~$\SU(n)$-invariant forms on~$\Gr(m,n)$.
Moreover, the map
\[
[\cdot]:\bbZ[\phi_1,\ldots,\phi_{n-m}]
\to H^*\bigl(\Gr(m,n),\bbZ\bigr)
\] 
is a isomorphism of rings.

In particular, $\phi_1>0$
defines an $\SU(n)$-invariant K\"ahler form on~$\Gr(m,n)$.
The normalization is such that, if ~$E^-\subset\C{n}$ is 
an $(m{-}1)$-plane and~$E^+\subset\C{n}$ is an $(m{+}1)$-plane
containing~$E^-$, then the line
\begin{equation}
[E^-,E^+]_m = \{E\in\Gr(m,n)\ \vrule\ E^-\subset E\subset E^+\ \}
\simeq \bbP^1
\end{equation} 
has unit area.  When~$m=1$, this defines the usual Fubini-Study metric 
on~$\bbP^{n-1}$.

This is as good a place as any to prove the following result 
for future use.

\begin{proposition}\label{prop: I-a is holo and closed}
For each~$\ab\in \sfP(m,n)$, the ideal~$\cI_\ab$ on~$\Gr(m,n)$ is
holomorphic and differentially closed.
\end{proposition}

\begin{proof}
The ideal~$\cI_\ab$ is the sheaf of sections of the 
sub-bundle~$\I_\ab\subset\L^{|\ab|,0}(T^*)$.   By its construction,
this bundle is $\SU(n)$-invariant and its fiber over~$\C{m}\in\Gr(m,n)$
is the subspace~$\bbS_\ab(\C{m})\ot\bbS_{\ab'}\bigl((\C{m})^\perp\bigr)
\subset\L^{|\ab|}\bigl(\C{m}\ot(\C{m})^\perp\bigr)$.  This latter
subspace is a (minimal) $K$-invariant subspace of 
$\L^{|\ab|}\bigl(\C{m}\ot(\C{m})^\perp\bigr)$ where~$K 
= S\bigl(\Un(m)\times\Un(n{-}m)\bigr)$ is the stabilizer in~$\SU(n)$
of~$\C{m}$.  Since~$\Gr(m,n)$ is a symmetric space, it follows that
the sub-bundle~$\I_\ab$ is parallel with respect to the Levi-Civita
connection of the $\SU(n)$-invariant metric associated to
the K\"ahler form~$\phi_1$ on~$\Gr(m,n)$.  Equivalently,
~$\nabla \I_\ab\subset T^*\ot \L^{|\ab|}(T^*)$ 
is a subspace of~$T^*\ot\I_\ab$.  Since the Levi-Civita connection
is torsion-free, the exterior derivative~$d$ on $p$-forms 
is just~$d = W\circ\nabla$ where~$W: T^*\ot\L^{p}(T^*)\to \L^{p+1}(T^*)$
is the bundle map induced by wedge product.  The differential closure
and holomorphicity of~$\cI_\ab$ now follow immediately.
\end{proof}

\subsubsection{Schur forms}\label{sssec:schur-forms}
For~$\ab\in \sfP(m,n)$, define the \emph{Schur form}~$\phi_\ab$
on~$\Gr(m,n)$ to be the polynomial
\begin{equation}\label{eq:Schub-form-def}
\phi_\ab = 
\det \begin{vmatrix}
      \phi_{a_1} & \phi_{a_1+1} & \hdots & \phi_{a_1+m-1}\\
      \phi_{a_2-1} & \phi_{a_2} & \hdots & \phi_{a_2+m-2}\\
       \hdotsfor[2.0]4\\
      \phi_{a_m-m+1} & \phi_{a_m-m+2} & \hdots & \phi_{a_m}\\
\end{vmatrix},
\end{equation}
where, again, the convention is that~$\phi_0 = 1$ and~$\phi_j = 0$
unless~$0\le j\le n{-}m$.  Note that~$\phi_{(p)} = \phi_p$ 
for $0\le p\le n{-}m$, so that a potential notational confusion is avoided.

Then, by the above discussion, the 
set~$\{\phi_\ab\mid \ab\in \sfP(m,n), |\ab| = p\}$ is a basis for the
$\SU(n)$-invariant~$2p$-forms on~$\Gr(m,n)$.  Since~$[\phi_\ab] = q_\ab$,
the pairing identity \eqref{eq:Schubert-cycle-pairing-relation} implies
\begin{equation}\label{eq:phi-and-sigma-pairing}
\int_{\sigma_{\bb^*}}\phi_\ab 
= \bigl\la q_\ab,[\sigma_{\bb^*}]\bigr\ra = \delta^\bb_\ab\,.
\end{equation}
Consequently, no constant linear combination~$\phi = c^{\ab}\,\phi_\ab$ 
could possibly be a (weakly) positive $2p$-form unless~$c^\ab\ge0$ 
for all~$\ab$.  A result of Fulton and Lazarsfeld~\cite{MR85e:14021} shows
that this necessary condition is actually sufficient: 

\begin{theorem}[Fulton-Lazarsfeld]\label{thm:Fulton-Lazarsfeld}
For any~$\ab\in \sfP(m,n)$, the form~$\phi_\ab$ 
is positive.  
\end{theorem}

\begin{remark}
In \cite[Appendix A]{MR85e:14021},  Fulton and Lazarsfeld
give an explicit formula for~$\phi_\ab$ that makes this clear.
Since I will need their formula in what follows, I will 
sketch their proof.  

Alternatively, it follows from 
general results of Kostant~\cite[Corollary 6.15]{MR26:266} that 
there exist unique $\SU(n)$-invariant forms~$\phi_\ab$ on~$\Gr(m,n)$
that satisfy~\eqref{eq:phi-and-sigma-pairing} and that these uniquely
defined forms are necessarily positive.  However, the very
explicit Giambelli formula~\eqref{eq:Schub-form-def} for~$\phi_\ab$ 
requires a separate argument.
\end{remark}

\begin{proof}[Sketch of proof]
Fix~$\ab\in \sfP(m,n)$ with~$|\ab| = p>0$ and let~$S_p$
denote the symmetric group on~$[1,p] = \{1,\dots,p\}$.  
Recall from~\cite[Lecture 4]{MR93a:20069} (whose notation I will follow)
that one can associate to~$\ab$ an irreducible, unitary 
representation~$\rho_\ab:S_p \to\Un(V_\ab)$.
For example, $\ab=(p)$ corresponds
to the trivial representation of~$S_p$, 
while~$\ab=(1,\dots,1)$ corresponds
to the alternating representation, i.e., 
$\sigma\mapsto\sgn(\sigma)\in \{\pm 1\}$.  Let~$\chi_\ab:S_p\to\bbC$
be the corresponding character.

Define~$\Theta^a_b = \omega^a_i\w\ov{\omega^b_i}$
for~$a,b\in\{m{+}1,\dots,n\} = [m{+}1,n]$. According to 
Fulton and Lazarsfeld~\cite[(A.6)]{MR85e:14021}, the formula%
\footnote{The careful reader will notice 
a difference between equation~(A.6)
of~\cite{MR85e:14021} and \eqref{eq:Fulton-Lazarsfeld-formula}, namely
that it is the character of~$\ab'$ rather than that of~$\ab$ that
enters into \eqref{eq:Fulton-Lazarsfeld-formula}.  
This is caused by the fact
that the convention in~\cite{MR85e:14021} for associating a representation
to a partition differs from that of~\cite{MR93a:20069}, which is the one that
I am following in this article.}
\begin{equation}\label{eq:Fulton-Lazarsfeld-formula}
\pi^*_m(\phi_\ab) = \frac1{p!} \left(\frac\iC{2\pi}\right)^p
\sum_{\tau\in S_p}\ 
\sum_{\alpha\in[m{+}1,n]^p}
\chi_{\ab'}\bigl(\tau\bigr)\>
  \Theta^{\alpha_1}_{\alpha_{\tau(1)}}\w\dots\w
  \Theta^{\alpha_p}_{\alpha_{\tau(p)}}
\end{equation}
holds for any~$\ab\in \sfP(m,n)$.  Using manipulations similar 
to those in~\cite[Appendix~A]{MR85e:14021}, 
\eqref{eq:Fulton-Lazarsfeld-formula} can be rewritten
in the form
\begin{equation}\label{eq:visibly-positive}
\pi^*_m(\phi_\ab)  
= \frac{\iC^{p^2}}{(2\pi)^p} 
     \sum_{\alpha\in[m{+}1,n]^p}\ 
     \sum_{i\in[1,m]^p}
   \tr\left( \zeta^\alpha_{i}(\ab')
             \w \bigl(\zeta^\alpha_{i}(\ab')\bigr)^*\right)
\end{equation}
where, for~$i\in[1,m]^p$ and~$\alpha\in[m{+}1,n]^p$, I have set
\begin{equation}\label{eq: zeta elements}
\zeta^\alpha_{i}(\ab') 
= \frac1{p!} \sum_{\tau\in S_p} \rho_{\ab'}(\tau^{-1}) \>  
\omega^{\alpha_{\tau(1)}}_{i_1}\w\cdots\w\omega^{\alpha_{\tau(p)}}_{i_p}\,,
\end{equation}
so that~$\zeta^\alpha_{i}(\ab')$ is a $p$-form on~$\SU(n)$ with
values in~$\Aut(V_{\ab'})$.

Since~$\pi_m:\SU(n)\to\Gr(m,n)$ is a submersion, 
\eqref{eq:visibly-positive} shows that~$\phi_\ab$ is indeed positive.
\end{proof}

\subsection{Bundles generated by global sections}
\label{ssec: bundles generated by global sections}

Let~$F\to M$ be a holomorphic vector bundle of rank~$r$ 
over a compact complex manifold~$M$ and denote the
vector space of its global holomorphic sections by~$H^0(M,F)$.
This space is finite dimensional, with, say, 
dimension~$n=h^0(F)$.  Consider the \emph{evaluation mapping}
$$
\ev_{\!F} : M\times H^0(M,F)\to F
$$ 
defined by~$\ev_{\!F}(x,s) = s(x)$.  
If this is a surjective bundle mapping, then~$F$ 
is said to be \emph{generated by global sections}.  

Assuming that~$F$ is generated by global sections, 
let~$K\subset M\times H^0(M,F)$ be the kernel of~$\ev_{\!F}$.  
Then~$K$ is a holomorphic subbundle of rank~$m=h^0(F)-\rank(F) = n{-}r$.  
The holomorphic mapping~
$\kappa_F:M\to\Gr\bigl(m,H^0(M,F)\bigr)\simeq\Gr(m,n)$ 
defined by~$\kappa_F(x) = K_x$ satisfies~$F = \kappa_F^*(Q)$ 
where~$Q$, as usual, denotes the quotient bundle over~$\Gr(m,n)$
as defined in~\S\ref{ssec: canonical bundles over Gr(m,n)}.  

The Chern classes of~$F$ are given 
by~$c_a(F)=\kappa_F^*\bigl(c_a(Q)\bigr)$.
Generalizing this, for any partition~$\ab\in \sfP(m,n)$, 
one can define~$c_\ab(F)$ to be $\kappa_F^*(q_\ab)$.
Of course, each $c_\ab(F)$ can be written as a
polynomial in the usual Chern classes of~$F$:
$$
c_\ab(F) = 
\det \begin{vmatrix}
      \hfill c_{a_1}(F) & \hfill c_{a_1+1}(F) & \hdots & c_{a_1+m-1}(F)\\
      \hfill c_{a_2-1}(F) & \hfill c_{a_2}(F) & \hdots & c_{a_2+m-2}(F)\\
       \hdotsfor[2.0]4\\
      c_{a_m-m+1}(F) & c_{a_m-m+2}(F) & \hdots & \hfill c_{a_m}(F)\\
\end{vmatrix}.
$$
Thus, one takes this to be the definition of the Schur-Chern 
class~$c_\ab(F)$ even when~$F$ is not generated by global sections.

Theorem~\ref{thm:Fulton-Lazarsfeld} implies that, when~$F$
is generated by global sections, each Schur-Chern class~$c_\ab(F)$ 
is represented by a positive $\bigl(|\ab|,|\ab|\bigr)$-form, i.e., 
that~$c_\ab(F)$ is positive in the sense of~\S\ref{ssec:positivity}.  
This observation yields the following basic fact.

\begin{corollary}\label{cor: F generated by sections and ca = 0}
Suppose that~$M$ is compact and K\"ahler, and that~$F\to M$ is
a holomorphic bundle that is generated by global sections.
Then $c_\ab(F)\ge0$ and equality holds 
and only if $\kappa_F^*(\phi_\ab) = 0$.
\end{corollary}

\begin{proof}
Fix any positive definite Hermitian inner product on~$H^0(M,F)$ 
and define the corresponding invariant 
forms~$\phi_\ab$ on $\Gr(m,n)$.
If~$c_\ab(F) = 0$, then, since it is represented by~$\kappa_F^*(\phi_\ab)$,
which is positive by Theorem~\ref{thm:Fulton-Lazarsfeld}, 
and since~$M$ is K\"ahler,
the Generalized Wirtinger Inequality of~\S\ref{ssec:positivity} implies
that $\kappa_F^*(\phi_\ab)$ must vanish identically.
\end{proof}

Assume now that~$M$ is connected and that~$F\to M$ is
generated by its global sections.  Then~$\kappa_F(M)\subset\Gr(m,n)$ is
a irreducible algebraic variety of some dimension~$\dim_F(M)\le \dim M$.  
Since~$\phi_1$ is a K\"ahler form on~$\Gr(m,n)$, Wirtinger's theorem
implies that~$\dim_F(M)$ is the largest integer~$p\ge0$ 
so that~$\bigl(c_1(F)\bigr)^p\not=0$.  

One consequence of the Frobenius Formula~\cite[p. 49]{MR93a:20069} is
the identity
\begin{equation}\label{eq: sum of positive forms}
(\phi_1)^p = \sum_{|\ab| = p}
               (\dim V_\ab)\,\phi_{\ab}\,.
\end{equation}
In particular, $\dim_F(M)$ is the largest integer~$p$ for which
there exists an~$\ab\in \sfP(m,n)$ with~$|\ab|=p$ and~$c_\ab(F)\not=0$.

\begin{remark}[Relation with ampleness]
\label{rem: reln w ampleness}
Fulton and Lazarsfeld~\cite{MR85e:14021} prove that if~$F\to M$ 
is \emph{ample}%
\footnote{in the sense of Hartshorne, which is different from
Griffiths' notion of ample in~\cite{MR41:2717}, for example. }
then~$c_\ab(F)\not=0$ for all~$\ab\in \sfP(m,n)$ with~$|\ab|\le\dim M $.  
Their work was the
culmination of the efforts of several authors who had established
partial results along these lines relating the notion of ampleness
with that of positivity of various Chern classes.  
For a full discussion of the historical development, see~\cite{MR85e:14021}.
\end{remark}

In this article, I am going to be characterize
the `extremal cases' where~$F\to M$ is generated by its
global sections but~$c_\ab(F)$ vanishes for some~$\ab$ with~$|\ab|\le 3$.  
Of course, such a bundle is not ample if~$|\ab|\le \dim M$.

\begin{example}[When~$\dim_F(M)=0$ or~$1$]
Here are two particularly simple cases.  In each case,
I am assuming that~$M$ is compact, connected, and K\"ahler
and that~$F\to M$ is a holomorphic bundle that is generated by its
global sections.  
In particular,~$c_1(F)\ge0$, so~$\bigl(c_1(F)\bigr)^p\ge0$
for all~$p\ge0$.

First, if~$c_1(F)=0$, then~$\kappa_F^*(\phi_1)=0$, so~$\kappa_F(M)$ 
has dimension~$0$ and therefore is a single point.
Equivalently,~$K_x\subset \C{n}$ is independent of~$x\in M$.  Of course,
this implies that any section of~$F$ that vanishes at one point of~$M$
vanishes at all points of~$M$.  
Consequently,~$\ev_{\!F}: M\times H^0(F)\to F$
is an isomorphism, i.e.,~$F$ is trivial.

Second, suppose that~$c_1(F)\not=0$ but that~$\bigl(c_1(F)\bigr)^2=0$.
Then~$\kappa_F^*({\phi_1}^2)=0$, so~$\kappa_F(M)$ 
has dimension~$1$ and thus is an irreducible algebraic 
curve in~$\Gr(m,n)$.  Let~$C\to\kappa_F(M)$ be the (canonical) 
desingularization of~$\kappa_F(M)$ and let~$Q_C$ be the pullback to~$C$
of~$Q$ under the composition~$C\to\kappa_F(M)\subset\Gr(m,n)$.
Then there exists a unique `lifting'~$\kappa:M\to C$ of~$\kappa_F$
and it satisfies~$F = \kappa^*(Q_C)$.  Thus, the vanishing of~
$\bigl(c_1(F)\bigr)^2$ implies that~$F$ is the pullback of a bundle 
over a curve.  Conversely, it is obvious that if $Q_C\to C$ is any
bundle over a curve that is generated by its global sections, then
for any map~$\kappa:M\to C$, the bundle~$\kappa^*(Q_C)$ is generated
by its global sections and satisfies~$\bigl(c_1(F)\bigr)^2=0$.
\end{example}

Of course, this description generalizes to the cases where~$\dim_F(M)>1$, 
but in these cases there need not be a desingularization~$X\to\kappa_F(M)$
that allows a holomorphic lifting~$\kappa:M\to X$ of~$\kappa_F$.  Thus,
one can only say that~$F$ is the pullback of a bundle over a 
singular variety of dimension~$\dim_F(M)$.  It would be interesting to
know conditions implying that the singularities of the image~$\kappa_F(M)$
can be resolved in a manner compatible with the mapping~$\kappa_F$.

\subsection{The ideal~$\cI_\ab$}\label{ssec: ideal I-a}
By Corollary~\ref{cor: F generated by sections and ca = 0},
if~$F\to M$ is generated by global sections, then~$c_\ab(F)=0$ 
if and only if~$\phi_\ab$ vanishes on the tangent planes 
to~$\kappa_F(M)\subset\Gr(m,n)$ at the smooth points of~$\kappa_F(M)$.
Of course, when~$|\ab|<\dim_F(M)$, this vanishing puts nontrivial 
conditions on the image~$\kappa_F(M)$.  It is to the analysis
of these conditions that I now turn.

It follows from~\eqref{eq: sum of positive forms} that 
there is no complex $p$-plane on which all of 
the forms~$\phi_\ab$ with~$|\ab|=p$ vanish.  However, except
when~$\ab = (1)$ or~$(1)^*$ (i.e., the cases for which there 
is only one term in the sum), the locus~$Z(\phi_\ab)$ 
is nonempty:

\begin{corollary}\label{cor: integral elements of phi-a}
$Z(\phi_\ab)$ contains the 
$|\ab|$-planes~$E\subset T_V\Gr(m,n)\simeq V^\perp\ot V^*$ 
of type~$\bb^*$ for every~$\bb\in \sfP(m,n)$ 
with~$\bb\not=\ab$ and~$|\bb|=|\ab|$.
\end{corollary}

\begin{proof}
By~\eqref{eq:phi-and-sigma-pairing} and the positivity of~$\phi_\ab$, 
it follows that, when~$\bb\not=\ab$ and~$|\bb| = |\ab|$, 
the form~$\phi_{\ab}$ vanishes on the Schubert variety~$\sigma_{\bb^*}$.
As has already been seen, at each smooth point~$V\in\sigma_{\bb^*}$,
the tangent plane~$T_V\sigma_{\bb^*}\subset T_V\Gr(m,n)\simeq V^\perp
\ot V^*$ is of type~$\bb^*$ and so must belong to~$Z(\phi_\ab)$.

Conversely, every subspace~$E\subset T_V\Gr(m,n)\simeq V^\perp\ot V^*$
of type~$\bb^*$ is tangent to the Schubert cell~$W_{\bb^*}(F)$
for some flag~$F$ on~$\C{n}$.  Since~$W_{\bb^*}(F)$ is dense in the smooth
locus of~$\sigma_{\bb^*}(F)$, it follows that $E$ must belong 
to~$Z(\phi_\ab)$. 
\end{proof}

\begin{remark}
As will be seen (cf. Lemma~\ref{lem: integral elements of I3}), 
it is not generally true that every element 
of~$Z(\phi_\ab)$ is of type~$\bb^*$ for some~$\bb\in \sfP(m,n)$ with
$|\bb| = |\ab|$ and $\bb\not=\ab$.
\end{remark}

\begin{lemma}\label{lem:ideal-of-Schub-form}
Suppose that~$\ab\in \sfP(m,n)$ satisfies~$|\ab| = p$.  
Then~$Z(\phi_\ab)$ consists of the complex 
$p$-planes~$P\subset T_V\Gr(m,n)$ that are integral elements
of~$\cI_\ab$.  More generally,~$\phi_\ab$ vanishes on
a complex subspace~$E\subset T_V\Gr(m,n)$ if and only if
it is an integral element of~$\cI_\ab$.
\end{lemma}

\begin{proof}
It follows from equations~\eqref{eq:visibly-positive} 
and \eqref{eq: zeta elements}, together with the discussion in
\cite[\S6.1]{MR93a:20069} of Weyl's construction of the
of Schur functors (especially the Exercises~6.14 and 6.15),
that, when the $(p,p)$-form~$\phi_\ab$ is written locally in the form 
$$
\phi_\ab = \iC^{p^2}\,\sum_{k=1}^K\zeta_k\w\ov{\zeta_k}
$$
for some local~$(p,0)$-forms~$\zeta_1,\dots,\zeta_K$, these 
latter forms must be a local basis of the 
subspace~$\I_\ab = \bbS_\ab(S)\ot\bbS_{\ab'}(Q^*)\subset\L^p(S\ot Q^*)
= \L^p\bigl(T^*\Gr(m,n)\bigr)$.  The statements of the lemma 
follow immediately from this and Definition~\ref{def:ideal-a}.  
The representation-theoretic details are left to the reader.
\end{proof}

\begin{example}
Consider the case of $\ab = (2)$, for which~$\phi_{(2)} = \phi_2$, i.e.,
the form representing the second Chern class of the quotient bundle.
Since~$\ab' = (1,1)$, the representation~$\rho_{\ab'}$ has degree~$1$
and~$\rho_{\ab'}(\tau)$ is simply the sign of~$\tau\in S_2$.  This gives
\begin{equation}\label{eq:I-2-generators}
\zeta^{\alpha_1\alpha_2}_{i_1i_2}(\ab') 
= \frac12\left(\omega^{\alpha_1}_{i_1}\w\omega^{\alpha_2}_{i_2}
               - \omega^{\alpha_2}_{i_1}\w\omega^{\alpha_1}_{i_2}\right)
= \frac12\left(\omega^{\alpha_1}_{i_1}\w\omega^{\alpha_2}_{i_2}
               + \omega^{\alpha_1}_{i_2}\w\omega^{\alpha_2}_{i_1}\right).
\end{equation}
Note that this expression is skew-symmetric in~$\alpha_1,\alpha_2$
and symmetric in~$i_1,i_2$.  It then follows from 
Definition~\ref{def:ideal-a} that these 2-forms span the $\pi_m$-pullback
of the ideal~$\cI_{(2)}$, as is claimed by 
Lemma~\ref{lem:ideal-of-Schub-form}.
\end{example}

\begin{corollary}\label{cor:ideals-for-homology}
A subvariety~$V\subset\Gr(m,n)$ of dimension~$|\ab|$ satisfies
$[V] = r[\sigma_{\ab^*}]$ for some~$r\in\bbZ^+$ if and only if it is 
an integral variety of~$\cI_\bb$ for all~$\bb\in\sfP(m,n)$ with~$\bb\not=\ab$ 
and~$|\bb|=|\ab|$.  
\end{corollary}

\begin{proof}
It has already been noted that~$[V] = r[\sigma_{\ab^*}]$ for 
some~$r\in\bbZ^+$ if and only if~$\int_V \phi_\bb = 0$ for
all $\bb\in \sfP(m,n)$ with~$\bb\not=\ab$ and~$|\bb|=|\ab|$.  Since,
by Theorem~\ref{thm:Fulton-Lazarsfeld}, $\phi_\bb$ is positive, the
equation~$\int_V\phi_\bb = 0$ holds if and only if $\phi_\bb$
vanishes on~$V$.  In turn, by Lemma~\ref{lem:ideal-of-Schub-form},
this holds if and only if~$V$ is an integral manifold of~$\cI_\bb$
for all $\bb\in \sfP(m,n)$ with~$\bb\not=\ab$ and~$|\bb|=|\ab|$,
as claimed.
\end{proof}

\subsubsection{Ideal inclusions}
\label{sssec:ideal inclusions}
Equation~\eqref{eq: Pieri p} in cohomology implies the form equation
\begin{equation}
\label{eq: Pieri p for forms}
(\phi_1)^p\w\phi_\ab 
= \sum_{\substack{\bb\in\sfP(m,n) \\ |\bb|= p+|\ab|}} 
      \mu^\bb_\ab\,\phi_\bb\,,
\end{equation}
which leads to the following result, which characterizes the integral 
elements of~$\cI_\ab$ in terms of the~$Z(\phi_\bb)$ with~$\bb\ge\ab$. 

\begin{lemma}
\label{lem: ideal inclusions}
The following relationships hold between ideals and integral elements:
\begin{enumerate}
\item $\cI_\bb \subseteq \cI_\ab$  if and only if~$\ab \le \bb$.  
\item A complex subspace~$E\subset T_V\Gr(m,n)$ of dimension~$r+|\ab|$ 
is an integral element of~$\cI_\ab$ if and only if it is an integral element
of~$\cI_\bb$ for all~$\bb\in\sfP(m,n)$ with~$|\bb|=r+|\ab|$ and~$\bb\ge\ab$.
\item A complex subspace~$E\subset T_V\Gr(m,n)$ of type~$\bb^*$ 
is an integral element of~$\cI_\ab$ if and only if~$\bb\not\ge\ab$.
\end{enumerate}
\end{lemma}

\begin{proof}
Assertion~(1) follows from equation~\eqref{eq: Pieri p for forms}
together with equation~\eqref{eq: ideals wedge} and the fact 
that~$\cI_1=\Omega^{(*,0)}\bigl(\Gr(m,n)\bigr)$.

For assertion~(2), one direction is easy:  If $E$ is
an integral element of~$\cI_\ab$, then, by the first statement, 
$E$ is an integral element of~$\cI_\bb$ for all~$\bb$ with~$\bb \ge \ab$.

Conversely, suppose that~$E$ has dimension~$r+|\ab|$ and 
is an integral element of~$\cI_\bb$ for all~$\bb\in\sfP(m,n)$ with
$|\bb| = r + |\ab|$ and~$\bb\ge\ab$.  Then~\eqref{eq: Pieri p for forms}
implies that the form~$(\phi_1)^r\w\phi_\ab$ vanishes on~$E$.  Since
$\phi_1$ pulls back to~$E$ to be a strictly positive $(1,1)$-form,
the Generalized Wirtinger 
Inequality~\eqref{eq:Generalized Wirtinger Inequality} implies that
$\phi_\ab$ must vanish on~$E$ as well, i.e., that~$E$ is an integral
element of~$\cI_\ab$, as desired.

Finally, (3) now follows from~(2) 
and Corollary~\ref{cor: integral elements of phi-a}.
\end{proof}

\subsubsection{Some integral manifolds of $\cI_\ab$}
\label{sssec: integral manifolds of I a}
As has already been remarked, the ideals~$\cI_\ab$ are invariant
under the action of~$\SL(n,\bbC)$ and so the space of integral
elements of each~$\cI_\ab$ at any given point in~$\Gr(m,n)$ 
is essentially independent of the point.  Moreover, two 
subspaces~$A\subset Q_V\otimes V^* = T_V\Gr(m,n)$ and
$B\subset Q_W\otimes W^* = T_W\Gr(m,n)$ of the same type 
(see Remark~\ref{rem: A-cycles}) 
are either both integral elements of~$\cI_\ab$
or neither integral elements of~$\cI_\ab$.  

In particular, if~$E\subset Q_V\otimes V^* = T_V\Gr(m,n)$ is
an integral element of~$\cI_\ab$, consider an 
$E$-cycle~$\sigma_E\subset\Gr(m,n)$.  Since~$\sigma_E$ (which is
birational to a projective space) is quasi-homogeneous, it contains
a Zariski-open set~$\sigma^\circ_E$ in its smooth locus such that
the tangent space at each point of~$\sigma^\circ_E$ is a subspace
of type~$E$ and hence, in particular, an integral element of~$\cI_\ab$.
Since~$\sigma_E$ is irreducible, this implies that~$\sigma_E$ is
actually an integral manifold of~$\cI_\ab$.  This yields the 
following elementary but important result:

\begin{proposition}
\label{prop: integral E-cycles of I-a}
Every integral element~$E$ of~$\cI_\ab$ is tangent 
to an integral manifold of~$\cI_\ab$ that is an $E$-cycle~$\sigma_E$.
\end{proposition}

\begin{remark}[Non-uniqueness]
\label{rem: nonuniqueness of integral manifolds}
It is not generally true that all of the integral manifolds
of a given~$\cI_\ab$ (even the ones of maximal dimension) are
of the form~$\sigma_E$ for some integral element~$E$ of~$\cI_\ab$.
In fact, this seems to be very rare and several examples of its
failure will be seen in the next section.
\end{remark}

\begin{lemma}
\label{lem: max dim integral elements of I-a}
The maximum dimension for integral elements of~$\cI_\ab$ is
equal to the maximum value of~$|\bb|$ for~$\bb\in\sfP(m,n)$
that satisfy~$\bb\not\ge\ab$. 
\end{lemma}

\begin{proof}
Set~$r = \max\{\,|\bb|{-}|\ab|\,\mid\,\bb\not\ge\ab\}\ge0$ 
and suppose that~$\bb\in\sfP(m,n)$ is such that
$\bb\not\ge\ab$ and $|\bb| = r{+}|\ab|$.  
By Lemma~\ref{lem: ideal inclusions}(3), 
any subspace~$E\subset T_V\Gr(m,n)$ of type~$\bb^*$
is an integral element of~$\cI_\ab$.  Since the dimension of
such an~$E$ is~$|\bb|$, it follows that~$\cI_\ab$ has integral 
elements of dimension~$r{+}|\ab|$.  

It remains to show that $\cI_\ab$
has no integral elements of dimension~$r{+}1{+}|\ab|$.
By the defining property of~$r$, in the equation
\begin{equation}
(\phi_1)^{r+1}\w\phi_\ab 
= \sum_{\substack{\bb\in\sfP(m,n) \\ |\bb|= r+1+|\ab|}} 
      \mu^\bb_\ab\,\phi_\bb\,,
\end{equation}
all of the coefficients~$\mu^\bb_\ab$ with~$|\bb| = r{+}1{+}|\ab|$ are
positive.  It follows that if~$E$ were an integral element of~$\cI_\ab$
of dimension~$r{+}1{+}|\ab|$, then~$E$ would be an integral element 
of~$\cI_\bb$ for all~$\bb$ with~$|\bb| = r{+}1{+}|\ab|$.  
Since all such~$\bb$ satisfy $\bb\ge0$, Lemma~\ref{lem: ideal inclusions}(2),
would then imply that~$E$ was an integral element of~$\cI_{(0)}$, 
which is absurd, since~$\cI_{(0)}$ has no positive dimensional 
integral elements.
\end{proof}

\begin{remark}[Maximal vs. maximum dimension]
\label{rem: max vs max dim}
As will be seen during the computation of the 
integral elements of~$\cI_{(3)}$ below, it is not true that all
of the maximal integral elements of~$\cI_\ab$ have the maximum dimension
allowed by~Lemma~\ref{lem: max dim integral elements of I-a}.  Moreover,
it can also happen that there are integral elements of the maximum
dimension that are not of type~$\bb^*$ for any~$\bb$.  Thus, 
Lemma~\ref{lem: max dim integral elements of I-a}, while very useful, 
is still quite a long way from determining the space of integral
elements of~$\cI_\ab$.
\end{remark}

\begin{remark}[Explicit computation]
\label{rem: explicit computation}
It is actually quite easy to explicitly determine the maximum
dimension of integral elements of~$\cI_\ab$.  Let~$\ab = (a_1,\dots,a_m)$
and, for convenience, set~$a_{m+1}=0$.  For each~$q$ 
in the range~$1\le q\le m$ for which~$a_q > a_{q+1}$, consider
the partition~$\ab^q = (a_1^q,\dots,a_m^q)$ defined by the 
conditions~$a_i^q = n{-}m$ for all~$i<q$ and~$a_i^q = a_q-1$ for~$i\ge q$.  
Since~$a_q^q = a_q-1 < a_q$, it follows that~$\ab^q\not\ge\ab$. 
Any~$\bb\not=\ab^q$ that satisfies~$\bb\ge\ab^q$ also satisfies
$\bb\ge\ab$ and, moreover, these~$\ab^q$ (there are at most~$m$ of them)
are the maximal elements in~$\sfP(m,n)$ that are not greater than~$\ab$.
Thus, the maximal dimension of an integral element of~$\cI_\ab$ is
the maximum of~$|\ab^q|$ where~$a_q>a_{q+1}$.
\end{remark}

\subsubsection{Complementarity}
\label{sssec:complementarity}
Every~$V\in\Gr(m,n)$ has an orthogonal complement~$V^\perp\in\Gr(n{-}m,n)$ 
with respect to the standard Hermitian inner product.
There is an $\SU(n)$-equivariant identification
\begin{equation}
T_V\Gr(m,n)\simeq V^\perp\ot V^*
\end{equation}
for which the Hermitian metric on~$T_V\Gr(m,n)$ induced by~$\phi_1$ agrees
with the tensor product Hermitian metric induced by the Hermitian metrics 
on~$V$ and~$V^\perp$.  This identification will be used implicitly from 
now on.  

The assignment~$V\mapsto V^\perp$ induces an 
\emph{anti}-holomorphic isometry~$\perp:\Gr(m,n)\to\Gr(n{-}m,n)$.

It is not difficult to show that~$\perp^*(\phi_{\ab'}) 
= (-1)^{|\ab|} \phi_\ab$, so knowledge of the integral elements and
integral manifolds of~$\cI_\ab$ on~$\Gr(m,n)$ implies such 
information about~$\cI_{\ab'}$ on~$\Gr(n{-}m,n)$.
Specifically, 
\begin{equation}
\label{eq:complementary-integral-elements}
Z(\phi_{\ab'}) =\left\{\ {\perp_*}(\bar E) \mid E\in Z(\phi_\ab)\ \right\}
\end{equation}
and, by the definition of the ideals~$\cI_{\ab}$,
\begin{equation}\label{eq: complementary ideals}
{\perp^*}\bigl(\cI_{\ab}\bigr) = \ov{\cI_{\ab'}},
\end{equation} 
so that $\perp$ exchanges the integral manifolds of~$\cI_{\ab}$ 
on~$\Gr(m,n)$ with those of~$\cI_{\ab'}$ on~$\Gr(n{-}m,n)$.  

The relationship~\eqref{eq:complementary-integral-elements} 
substantially reduces the number of cases
one needs to treat in computing the integral elements of the various
$\cI_\ab$.  For example, the knowledge of~$Z(\phi_{(2)})$ for all the
cases where~$2\le m\le n{-}2$ implies the knowledge of~$Z(\phi_{(1,1)})$
for all cases where~$2\le m\le n{-}2$, as will be seen. 

\subsubsection{Duality}\label{sssec:duality}
On any oriented Riemannian $n$-manifold~$M$, the Hodge 
star 
\begin{equation}
*:\Omega^k(M)\to\Omega^{n-k}(M)
\end{equation} 
is defined in such a way that any oriented orthonormal basis~
$(e_1,\ldots,e_n)$ of~$T_xM$ and any~$\phi\in\Omega^k(M)$ satisfy
\begin{equation}\label{eq: Hodge dual values}
\phi(e_1,\ldots,e_k) = {*}\phi(e_{k+1},\ldots,e_n).
\end{equation}

For any oriented $k$-dimensional subspace~$E\subset T_xM$, let
$E^\perp\subset T_xM$ be its orthogonal complement, oriented so
that~$T_xM = E\oplus E^\perp$ as oriented vector spaces.

Harvey and Knapp show~\cite[Corollary 1.3(b)]{MR50:7573} that if~$\phi$ 
is a positive $(p,p)$-form on a K\"ahler manifold~$M$ of dimension~$m$, 
then ${*}\phi$ is a positive~$(m{-}p,m{-}p)$-form. 
Moreover,~\eqref{eq: Hodge dual values} implies
\begin{equation}\label{eq: zeros of Hodge duals}
Z({*}\phi) = \{\> E^\perp \mid E\in Z(\phi) \>\}.
\end{equation}
By Definition~\ref{def:n-of-type-a}, if~$E\subset T_V\Gr(m,n)$ 
is of type~$\ab$, then~$E^\perp$ is of type~$\ab^*$.  
It follows from this,~\eqref{eq: zeros of Hodge duals},
and Corollary~\ref{cor: integral elements of phi-a}
that~$*\phi_\ab$ vanishes on the subspaces of type~$\bb$ 
where~$\bb\not=\ab$ and~$|\bb| = |\ab|$.  
Consequently,~${*}\phi_\ab$ is some positive multiple of~$\phi_{\ab^*}$.

In particular, for all~$\ab\in \sfP(m,n)$,
\begin{equation}\label{eq:dual-integral-elements}
Z(\phi_{\ab^*}) = Z({*}\phi_\ab) 
  = \left\{\> E^\perp \mid E\in Z(\phi_\ab) \>\right\}.
\end{equation}

This identity reduces by a factor of two the task of computing the 
integral elements of the various~$\cI_\ab$.

\begin{remark}[The action of the Hodge star operator]
\label{rem: effect of star}
Although it will not be needed in this article, the reader may be curious 
about the multiplier in the relationship between~$*\phi_\ab$ 
and~$\phi_{\ab^*}$. This multiplier can be calculated easily by first using 
Wirtinger's theorem to note that, for each~$p$ 
in the range~$0\le p\le m(n{-}m)$, the expression~$\frac{1}{p!}{\phi_1}^p$ 
restricts to each complex $p$-plane to be the volume form.  This implies
$$
*\left(\frac{{\phi_1}^p}{p!}\right) 
= \frac{{\phi_1}^{m(n-m)-p}}{\bigl(m(n{-}m)-p\bigr)!}.
$$
Now, applying this equality to~\eqref{eq: sum of positive forms}
and using the fact that~$*\phi_\ab$ is a multiple of~$\phi_{\ab^*}$ 
yields the relation
$$
\frac{\dim V_\ab}{|\ab|!}\,{*}\phi_\ab 
= \frac{\dim V_{\ab^*}}{|\ab^*|!}\,\phi_{\ab^*}\,.
$$
The dimension of~$V_\ab$ is computed in~\cite[Lecture 4]{MR93a:20069}.  
The reader might also compare~\cite{kKhT99}, 
where the ideas of this calculation
are generalized to the other Hermitian symmetric spaces.

One other consequence of~\eqref{eq: sum of positive forms} that bears
mentioning is that, when 
combined with Wirtinger's theorem~\eqref{eq: Wirtinger theorem} 
and~\eqref{eq:phi-and-sigma-pairing},
it yields the useful formula
\begin{equation}\label{eq: volume of sigma-a}
\vol(\sigma_{\ab^*}) = \frac{\dim V_\ab}{|\ab|!}\,.
\end{equation}
\end{remark}

\subsection{Walters' differential systems}
\label{ssec: Walters diff sys}
The thesis of Maria Walters~\cite{mW97} is particularly focussed 
on the study of the subvarieties~$V\subset\Gr(m,n)$ that 
satisfy~$[V] = r[\sigma_{\ab^*}]$ for some~$\ab\in\sfP(m,n)$. To this
end, she defines two differential systems~\cite[\S5.1]{mW97}
and discusses some related rigidity questions.

\subsubsection{The two differential systems}
\label{sssec: the two differential systems} 
The first system~\cite[Definition~40]{mW97}, 
which she denotes~$\mathcal{R}_\ab$
and calls a \emph{Schur differential system}, 
is the intersection%
\footnote{This definition does not quite work when
$\ab = (1)$ or~$(1)^*$ because there are no~$Z(\phi_\bb)$
to intersect in this case.  In these two extreme cases,
we set~$\mathcal{R}_{(1)} = \Gr\bigl(1,T\Gr(m,n)\bigr)$
and $\mathcal{R}_{(1)^*} = \Gr\bigl(m(n{-}m)-1,T\Gr(m,n)\bigr)$.}
of the~$Z(\phi_{\bb})$ for all
$\bb\in\sfP(m,n)$ with~$\bb\not=\ab$ and~$|\bb|=|\ab|$.  
Thus, the (local) integrals of this system are
the subvarieties~$X\subset\Gr(m,n)$ of dimension~$|\ab|$ 
with the property that~$\phi_\bb$ vanishes when pulled back to~$X$
for all~$\bb\in\sfP(m,n)$ with~$\bb\not=\ab$ and~$|\bb|=|\ab|$.
By Corollary~\ref{cor:ideals-for-homology}, a closed 
subvariety~$X\subset\Gr(m,n)$ is an integral of the 
system~$\mathcal{R}_{\ab}$ if and only if $[V] = r[\sigma_{\ab^*}]$ for 
some~$r\in\bbZ^+$.  Note that~$\mathcal{R}_{\ab}$ is a closed 
subvariety of $\Gr\bigl(|\ab|,T\Gr(m,n)\bigr)$ that is invariant
under the natural action of~$\SL(n,\bbC)$, and that it may be singular
and/or disconnected.

The second~\cite[Definition~41]{mW97}, 
which she denotes~$\mathcal{B}_\ab$
and calls a \emph{Schubert differential system},
is more restrictive, being made up of the 
subspaces~$A\subset T_V\Gr(m,n) = \C{n}/V\otimes V^*$ 
of type~$\ab$ (see Definition~\ref{def:n-of-type-a}).
By Remark~\ref{rem: closure of type}, the system~$\mathcal{B}_\ab$
is a closed subvariety of~$\Gr\bigl(d(\ab),T\Gr(m,n)\bigr)$.  
In fact, it is homogeneous under the isometry group of~$\Gr(m,n)$
and hence is a smooth bundle over~$\Gr(m,n)$.  
Since the tangent spaces to a Schubert cell~$W_\ab(F)$ 
are of type~$\ab$, it follows that, at all of its smooth points, 
the tangent spaces to~$\sigma_\ab(F) = \overline{W_\ab(F)}$ 
are of type~$\ab$.  
Thus, a  subvariety~$X\subset\Gr(m,n)$ of codimension~$|\ab|$ is
an integral of~$\mathcal{B}_\ab$ if and only if, at each 
smooth point~$x\in X$, there exists some Schubert variety~$\sigma_\ab(F)$
passing through~$x$ and smooth there and so that~$T_xX = T_x\sigma_\ab(F)$.

\subsubsection{Inclusion relations}
\label{sssec: inclusion relations} 
The two systems are related by the inclusion~$\mathcal{B}_\ab
\subseteq \mathcal{R}_{\ab^*}$.  In some cases, equality holds,
such as for~$\ab = (p)^*$ and $((p)')^*$ when~$p>1$ 
(see Remarks~\ref{rem: schur rigidity quasi-rigidity} 
and~\ref{rem: B p* = R p}), 
but this appears to be rare.  Even in the simple case~$\ab = (1)$, 
the two are different as soon as~$2\le m\le n{-}m$.

Walters shows the difference between~$\mathcal{B}_{(2,1)}$
and~$\mathcal{R}_{(2,1)^*} = \mathcal{R}_{(2,1)}$ in~$\Gr(2,5)$ 
by exhibiting a three-dimensional subvariety~$X\subset\Gr(2,5)$ 
that is an integral of~$\mathcal{R}_{(2,1)^*}$~\cite[Example~2]{mW97}
but not an integral of~$\mathcal{B}_{(2,1)}$~\cite[Proposition~16]{mW97}.
This difference can be exhibited more directly by 
computing~$\mathcal{R}_{(2,1)^*}$ 
(see Lemma~\ref{lem: integral elements of I111 and I3}).

\begin{example}[When~$\mathcal{B}_\ab\not=\mathcal{R}_{\ab^*}$]
\label{ex: Ba not eq Ra*}
Walters' example is one of a general family. 
Let~$p$ and~$q$ be integers satisfying~$1\le p\le n{-}m$ and~$1\le q\le m$
and consider~$\ab = (p,\dots,p,p{-}1)$ where~$|\ab| = pq-1$ for some
$p,q\ge1$ and~$\bb = (p,\dots,p)$ where~$|\bb| = pq$.  Then~$\ab$
is the unique element of~$\sfP(m,n)$ satisfying~$|\ab| = |\bb|{-}1$
and~$\ab\le\bb$.  (In other words,~$\ab$ is
the unique predecessor of~$\bb$.)
It then follows from Lemma~\ref{lem: ideal inclusions}(3)
that a subspace~$E\subset T_V\Gr(m,n)$ of type~$\bb^*$ (and hence of
dimension~$pq$) is an integral element of~$\cI_\cb$ 
for all~$\cb\not=\ab$ with~$|\cb|=pq{-}1$.  
In particular, any hyperplane~$H\subset E$ is a $(pq{-}1)$-dimensional
integral element of of~$\cI_\cb$ 
for all~$\cb\not=\ab$ with~$|\cb|=pq{-}1$ and so, by definition, belongs
to~$\mathcal{R}_{\ab^*}$.  

When~$p$ and~$q$ are each at least~$2$, the general hyperplane in~$E$
is not of type~$\ab^*$, so~$\mathcal{B}_\ab\not=\mathcal{R}_{\ab^*}$ 
for such~$\ab$.  In fact,
the hyperplanes in~$E\simeq\C{p,q}$ break up into~$\min\{p,q\}$ orbits 
under the action of~$\GL(p,\bbC)\times\GL(q,\bbC)$, so $\mathcal{R}_{\ab^*}$
consists of at least~$\min\{p,q\}$ distinct~$\SL(n,\bbC)$-orbits in
this case.  

The case~$\ab = (2,1)$ shows that there can be other orbits 
in~$\mathcal{R}_{\ab^*}$ besides these `obvious' ones 
(see Remark~\ref{rem:I111+3 integral element closures and intersections}).
\end{example}

\begin{remark}[Connectedness of~$\mathcal{R}_{\ab^*}$]
\label{rem: connectedness of Rab*}
Since~$\mathcal{R}_{\ab^*}$ is invariant under the~$\SL(n,\bbC)$ action
on~$\Gr\bigl(|\ab^*|,T\Gr(m,n)\bigr)$,
it is a union of $\SL(n,\bbC)$-orbits.  
By Remark~\ref{rem: closure of type}, the only closed $\SL(n,\bbC)$-orbit 
in~$\mathcal{R}_{\ab^*}$ is~$\mathcal{B}_{\ab}$.  In particular,
the closure of any $\SL(n,\bbC)$-orbit in $\mathcal{R}_{\ab^*}$ contains
$\mathcal{B}_{\ab}$, so it follows that~$\mathcal{R}_{\ab^*}$ is connected
(though it may be, and often is, reducible).  
For a discussion of a specific case, 
see Remark~\ref{rem:I111+3 integral element closures and intersections}.
\end{remark}

\subsubsection{Rigidity questions}
\label{sssec: rigidity questions} 
Walters asks whether every (smooth) irreducible integral variety 
of~$\mathcal{R}_{\ab}$ is necessarily equal to (an open subset of) 
some Schubert cycle~$\sigma_{\ab^*}$ and shows that, for certain~$\ab$
the answer is `yes', while, for others, the answer is `no'.  Although
she does not introduce this terminology, in the cases where the
answer is `yes', one might describe this by saying
that~$\sigma_\ab$ is \emph{Schur rigid}.

\begin{example}[Schur non-rigidity]
\label{ex: N in Gr25}
Walters cites the classical example~\cite[Example~2]{mW97} 
of the (smooth) variety~$N\subset\Gr(2,5)$ consisting of the $2$-planes 
that are isotropic for a nondegenerate complex inner product on~$\C5$.
The dimension of~$N$ is~$3$ and~$[N] = 4[\sigma_{(2,1)}]$, so $N$
must be a solution of~$\mathcal{R}_{(2,1)^*}$. However,~$N$ is 
not a Schubert variety.  (It is not even a solution 
of~$\mathcal{B}_{(2,1)}$~\cite[Proposition 16]{mW97}.)

More generally, Schur rigidity fails for \emph{any}~$\ab$ 
for which~$\mathcal{B}_\ab\not=\mathcal{R}_{\ab^*}$ since,
if~$A\subset T_V\Gr(m,n)$ belongs to~$\mathcal{R}_{\ab^*}$
but not~$\mathcal{B}_\ab$, then the $A$-cycle~$\sigma_A$ will
be an integral variety of $\mathcal{R}_{\ab^*}$ that is not
a Schubert cycle~$\sigma_\ab$.  (See Remark~\ref{rem: A-cycles}.)
\end{example}

Walters also asks whether every (smooth) irreducible integral variety 
of~$\mathcal{B}_{\ab}$ is necessarily equal to (an open subset of) 
some Schubert cycle~$\sigma_{\ab}$ and, again, shows that, for 
certain~$\ab$ the answer is `yes', while, for others, 
the answer is `no'.    Again, although
she does not introduce this terminology, 
in the cases where the answer is `yes', 
one might describe this by saying that~$\sigma_\ab$
is~\emph{Schubert rigid}.

\begin{example}[Schubert rigidity]
\label{ex: walters Schubert rigid}
Walters shows~\cite[Theorem~8 and Corollary~5]{mW97} 
that when
\begin{enumerate}
\item $\ab = (p,\dots,p)^*$ for some~$p>1$ (except for~$\ab = (p)^*$), 
\item $\ab = (n-m)^*$, or
\item $\ab = ((m)')^*$, 
\end{enumerate}
then any local solution of~$\mathcal{B}_\ab$
is a Schubert cycle~$\sigma_\ab(F)$ for some flag~$F$.  

She does this as follows: 
First, she observes that, in all of the cases 
listed above, the Schubert cycle~$\sigma_\ab$ is smooth and, in fact,
homogeneous.  She then shows that, for an~$\ab$ from
one of the cases listed above, two Schubert cycles~$\sigma_\ab(F)$
and~$\sigma_\ab(F')$ that are tangent at some common point must coincide.
Finally, she shows that if~$W\subset\Gr(m,n)$ is a solution 
of~$\mathcal{B}_\ab$, then the `Gauss map', defined by sending 
each point~$x\in W$ to the (unique) Schubert cycle~$\sigma_\ab(x)$ passing
through~$x$ and having~$T_xW$ as its tangent space, must be constant.
\end{example}

\begin{example}[Schubert non-rigidity]
\label{ex: walters Schubert non-rigid}
By contrast, Walters provides 
examples~\cite[Proposition~17 and Example~3]{mW97}  
that show that, when
\begin{enumerate}
\item $\ab = (p)^*$ for~$p$ in the range $1\le p < n{-}m$,
\item $\ab = ((q)')^*$ for~$q$ in the range $1\le q < m$, or
\item $\ab = (2,1)$ where~$(m,n) = (2,5)$,
\end{enumerate}
there are solutions of~$\mathcal{B}_\ab$ that are not Schubert cycles.

While she does not give a complete classification of the solutions
of $\mathcal{B}_{(2,1)}$, she does 
show~\cite[Proposition~18]{mW97} that such solutions~$W\subset\Gr(2,5)$
are ruled.  For more on these solutions, 
see Remark~\ref{rem: solns of B 21*}.
\end{example}

\begin{remark}[Higher order rigidity]
\label{rem: higher order rigidity}
For general~$\ab\in\sfP(m,n)$ the cycle~$\sigma_\ab$ is singular
and it is also not true that two cycles~$\sigma_\ab(F)$ 
and~$\sigma_\ab(F')$ that are tangent at a common smooth 
point must be equal.  

The simplest example of this is $\ab=(1)$ when~$(m,n) = (2,4)$. 
A Schubert cycle~$\sigma_{(1)}\subset\Gr(2,4)$ is uniquely determined 
by a $2$-plane~$W\in\Gr(2,4)$, e.g., the cycle~$\sigma_{(1)}(W)$ is
simply the set of $2$-planes~$V\subset~\C{4}$ such that~$V\cap W\not=(0)$.
It follows that there is a $3$-parameter family of $\sigma_{(1)}$s
passing through a given~$V\in\Gr(2,4)$ and all but one of these, 
namely~$\sigma_{(1)}(V)$ itself, is smooth there.  However, there
is only a $2$-parameter family of subspaces of~$T_V\Gr(2,4)$ that
are of type~$(1)$.  Thus, there is a $1$-parameter family of 
$\sigma_{(1)}$s passing through~$V$ and having a given 
tangent plane there. 

This might seem to account for the non-rigidity of the solutions of
$\mathcal{B}_{(1)}$ in~$\Gr(2,4)$.  At least, it provides one place
where Walters' argument for rigidity (see 
Example~\ref{ex: walters Schubert rigid}) would fail in this case.

However, one should not immediately assume the failure of rigidity
based on this non-uniqueness alone:

\begin{example}[Second order rigidity]
\label{ex: 2nd order rigidity}  
Consider, the case of~$\ab=(2,2)$ when~$(m,n) = (3,6)$.
A Schubert cycle~$\sigma_{(2,2)}\subset\Gr(3,6)$ is uniquely determined 
by a $3$-plane~$W\in\Gr(3,6)$, e.g., the cycle~$\sigma_{(2,2)}(W)$ 
is simply the set of $3$-planes~$V\subset~\C{6}$ 
such that~$\dim(V\cap W)\ge2$.
It follows that there is a $5$-parameter family of $\sigma_{(2,2)}$s
passing through a given~$V\in\Gr(3,6)$ and all but one of these, 
namely~$\sigma_{(2,2)}(V)$ itself, is smooth there.  However, there
is only a $4$-parameter family of subspaces of~$T_V\Gr(3,6)$ that
are of type~$(2,2)$. Thus, there is a $1$-parameter family of 
$\sigma_{(2,2)}$s passing through~$V$ and having a given 
tangent plane there.  

Nevertheless, it turns out that any irreducible solution 
to~$\mathcal{B}_{(2,2)}$ in~$\Gr(3,6)$ is $\sigma_{(2,2)}(V)$ 
for some~$V\in\Gr(3,6)$.  The proof of this result depends on 
going to a second order Gauss map:  One shows that for any point~$x$
of a (nonsingular, local) solution~$W\subset\Gr(3,6)$ 
of~$\mathcal{B}_{(2,2)}$, there is a unique~$V\in\Gr(3,6)$ such 
that~$x$ is a smooth point of~$\sigma_{(2,2)}(V)$ and so that~$W$
and~$\sigma_{(2,2)}(V)$ osculate to order~$2$ at~$x$.  This defines
a `second order Gauss map' from~$W$ to~$\Gr(3,6)$ and consideration
of the structure equations for this Gauss map show that it is 
constant.  

In fact, this second order argument generalizes to prove Schubert 
rigidity in all of the cases~$\ab = (p,\dots,p)$ in~$\Gr(m,n)$ 
where $|\ab|=pq$ and where~$p$ and~$q$ satisfy~$2\le p\le n{-}m$
and~$2\le q\le m$. The argument is very much like the moving frame 
arguments for the last two cases in the proof of 
Proposition~\ref{prop: integral varieties of I3}.  
This is not accidental; see Remark~\ref{rem: a case of A-rigidity}.
\end{example}

It could well be that there are examples 
of~$\ab$ for which all irreducible solutions of~$\mathcal{B}_{\ab}$ 
are of the form~$\sigma_\ab$, but where the proof of such rigidity
requires consideration of a suitable `Gauss map' of order even greater
than~$2$.
\end{remark}

\begin{remark}[$A$-rigidity]
\label{rem: A-rigidity}
Generalizing the case of~$\mathcal{B}_\ab$, 
for any subspace~$A\subset\C{n-m,m}$ of dimension~$d$, 
one can consider the subset~$\mathcal{B}_A\subset\Gr\bigl(d,T\Gr(m,n)\bigr)$
consisting of the subspaces~$E\subset T_V\Gr(m,n)$ of type~$A$.  
Of course, $\mathcal{B}_A$ is a single~$\SL(n,\bbC)$-orbit 
in~$\Gr\bigl(d,T\Gr(m,n)\bigr)$, but it is not compact unless~$A$ has
type~$\ab$ for some~$\ab\in\sfP(m,n)$.  One can also pose the 
more general $A$-rigidity problem:  Is every connected 
solution of~$\mathcal{B}_A$ an open subset of some~$A$-cycle~$\sigma_A$?

As pointed out in Remark~\ref{rem: a case of A-rigidity}, there are
examples of non-Schubert~$A$ where this sort of `$A$-rigidity' does hold.

\end{remark}

\subsection{Integral element computations}
\label{ssec: integral element computations}
In this section, 
I will compute the space of integral elements of~$\cI_\ab$,
$\cI_{\ab^*}$, $\cI_{\ab'}$, 
and $\cI_{{\ab'}{}^*} = \cI_{{\ab^*}{}'}$ for
the first three nontrivial cases:  $\ab = (2)$, $(3)$, and $(2,1)$.

To simplify the notation, I will begin with some conventions:  
For any~$V\in\Gr(m,n)$, 
I will write~$Q_V$ for the quotient space~$\C{n}/V$
and abbreviate this to~$Q$ when there is no danger of confusion.  Also,
for a vector~$\zb\in\C{n}$, I will usually denote its class in~$Q_V$ 
by~$\dl\zb\dr_V$, abbreviated to~$\dl\zb\dr$ 
when there is no danger of confusion.

Once an element~$V\in\Gr(m,n)$ is fixed, I will consider only
unimodular bases~$\vb = (\vb_1,\dots,\vb_n)$ of~$\C{n}$ 
with the property that~$V$ is spanned by~$\vb_1,\dots,\vb_m$.  (These
bases will \emph{not} be assumed to be unitary.)  The dual basis 
of~$\bigl(\C{n}\bigr)^*$ will be denoted~$\vb^* = (\vb^1,\dots,\vb^n)$,
and the elements~$\vb^1,\dots,\vb^m$ will be regarded as a basis of~$V^*$
in the obvious way.   I will adopt the usual index 
ranges~$1\le i,j,k\le m < \alpha,\beta,\gamma\le n$.

Using the canonical isomorphism~$T_V\Gr(m,n)=Q\ot V^*$, the
identity map~$\eta:T_V\Gr(m,n)\to Q\ot V^*$ can be expanded in
the form
\begin{equation}
\eta = \dl\vb_\alpha\dr{\ot}\vb^i\,\eta^\alpha_i\,,
\end{equation}
so that~$\{\,\eta^\alpha_i\mid 1\le i\le m <\alpha\le n\,\}$ are a basis
for the $(1,0)$-forms on~$T_V\Gr(m,n)$.  This basis depends, of course,
on the choice of the basis~$\vb$, and it is important to understand 
this dependence.  

It is customary to write~$\eta = (\eta^\alpha_i)$ and to think of
it as having values in~$\C{n-m,m}$, so I will follow this convention.
If~$\tilde\vb = (\tilde\vb_1,\dots,\tilde\vb_n)$ is any other unimodular 
basis with the property that~$V$ is spanned 
by~$\tilde\vb_1,\dots,\tilde\vb_m$, then 
$\tilde\vb = \vb\,\ub$ where~$\ub$ lies in~$P_m\subset\SL(n,\bbC)$, i.e.,
\begin{equation}
\ub = \begin{pmatrix} A & C\\ 0 & B\end{pmatrix}
\end{equation}
where~$A$ lies in~$\GL(m,\bbC)$ and~$B$ lies in~$\GL(n{-}m,\bbC)$ and,
of course, they satisfy~$\det(A)\det(B)=1$.  It
is not difficult to compute that the corresponding matrix~$\tilde\eta$
satisfies
\begin{equation}
\tilde\eta = B^{-1}\,\eta\,A.
\end{equation}
Thus, the effect of allowable basis changes is to pre- and post-multiply
$\eta$ by invertible matrices.  

\subsubsection{Dimension and codimension~$2$}\label{sssec: dim and codim 2}
The first task is to determine the integral elements of~$\cI_{(2)}$ 
and~$\cI_{(1,1)}$.

It is simpler to first state a result 
that characterizes the \emph{maximal} 
integral elements of these ideals 
and then deduce the structure of the space 
of integral elements of any given dimension from the maximal list.

\begin{lemma}\label{lem: integral elements of I2}
The maximal integral elements 
of~$\cI_{(2)}$ in~$T_V\Gr(m,n)\simeq Q\ot V^*$ 
fall into two distinct classes:
\begin{enumerate}
\item The $m$-dimensional subspaces~$E = L{\ot}V^*$, where~$L\subset Q$
      is any line.
\item The $1$-dimensional subspaces~$E\subset Q{\ot}V^*$ that do not
      lie in any subspace of the first kind.
\end{enumerate}
\end{lemma}

\begin{proof}
Fix~$V\in\Gr(m,n)$ and consider any basis~$\vb = (\vb_1,\dots,\vb_n)$ 
of~$\C{n}$ with the property that~$V$ is spanned by~$\vb_1,\dots,\vb_m$.  
Let~$\vb^1,\dots,\vb^n$ be the dual basis of~$(\C{n})^*$. 
The identification~$\eta:T_V\Gr(m,n)\to Q\ot V^*$ can be
written in the form
\begin{equation}
\eta = \dl\vb_\alpha\dr{\ot}\vb^i\,\eta^\alpha_i\,,
\end{equation}
where~$\{\,\eta^\alpha_i\mid 1\le i\le m <\alpha\le n\,\}$ are a basis
for the $(1,0)$-forms on~$T_V\Gr(m,n)$.  

In terms of these $(1,0)$-forms, the $(2,0)$-forms
\begin{equation}
\theta^{\alpha_1\alpha_2}_{i_1i_2} 
= \frac12\left(\eta^{\alpha_1}_{i_1}\w\eta^{\alpha_2}_{i_2}
               - \eta^{\alpha_2}_{i_1}\w\eta^{\alpha_1}_{i_2}\right)
= -\theta^{\alpha_2\alpha_1}_{i_1i_2}  = \theta^{\alpha_1\alpha_2}_{i_2i_1} 
\end{equation}
with~$1\le i_1,i_2\le m < \alpha_1,\alpha_2\le n$ generate~$\cI_{(2)}$ 
on~$T_V\Gr(m,n)$ (see \eqref{eq:I-2-generators}).

A subspace~$E\subset T_V\Gr(m,n)$ of dimension~$d$ is defined by 
a set of~$m(n{-}m)-d$ independent linear relations among 
the~$\eta^\alpha_i$.  Let~$\xi^\alpha_i$ denote the restriction
of~$\eta^\alpha_i$ to~$E$, so that exactly~$d\ge 2$ of the~$\xi^\alpha_i$
are linearly independent.  The hypothesis that~$E$ be an integral
element of~$\cI_{(2)}$ is then just that 
\begin{equation}\label{eq:E-an-integral-of-I-2}
\xi^{\alpha_1}_{i_1}\w\xi^{\alpha_2}_{i_2}
               - \xi^{\alpha_2}_{i_1}\w\xi^{\alpha_1}_{i_2} = 0,
\end{equation}
so I assume these quadratic relations from now on.

The~$\eta^\alpha_i$ and (hence) the~$\xi^\alpha_i$ depend on the choice
of~$\vb$.   Choose the basis~$\vb$ so that the maximum number, 
say~$p$, of~$\{\,\xi^\alpha_1\mid \alpha>m\,\}$ are linearly independent. 
(I.e., so that the first `column' of~$\xi$ contains the maximal number of 
linearly independent $1$-forms.)
Note that~$p$ satisfies~$1\le p\le \min(d,n{-}m)$. 
By making an allowable basis change, I can assume that
$\xi^{m+1}_1,\ldots,\xi^{m+p}_1$ are linearly independent and that
$\xi^\alpha_1 = 0$ for~$\alpha > m+p$.  

Setting~$\alpha_1 = 1$, $\alpha_2 =2$, and~$i_1 = i_2 = 1$ 
in~\eqref{eq:E-an-integral-of-I-2} 
yields~$2\,\xi^{m+1}_1\w\xi^{m+2}_1 = 0$. 
Thus, it follows that~$p = 1$.  

All of the forms~$\{\>\xi^\alpha_i\mid 2\le i\le m,\>m{+}1<\alpha\le n\>\}$ 
must be multiples of~$\xi^{m+1}_1$, since, otherwise,
a new allowable basis~$\tilde\vb$ could be found that would result
in at least two independent forms among the 
corresponding~$\tilde\xi^\alpha_1$, contradicting the maximality of~$p$, 
which is equal to~$1$.

Since~$d>1$, there must be~$d{-}1>0$ forms 
among~$\{\xi^{m+1}_2,\dots,\xi^{m+1}_m\}$
that are linearly independent modulo~$\xi^{m+1}_1$.  
By making a basis change that fixes~$\vb_1$, 
I can assume that~$\{\,\xi^{m+1}_1,\dots,\xi^{m+1}_d\,\}$ 
are linearly independent,  but that ~$\xi^{m+1}_j=0$ when~$d<j\le m$.  

Since there cannot be two linearly independent forms 
among~$\{\>\xi^\alpha_i\mid \alpha>m\>\}$ for any~$1\le i\le m$, it 
follows that~$\xi^\alpha_i\w \xi^{m+1}_i = 0$ for~$1<i\le d$ 
and~$\alpha>m{+}1$, but it has already been shown that
$\xi^\alpha_i\w \xi^{m+1}_1 = 0$ for~$1< i\le d$ and~$\alpha>m{+}1$.
Since~$\xi^{m+1}_1\w\xi^{m+1}_i\not=0$ for~$1<i\le d$, 
it follows that~$\xi^\alpha_i = 0$ for~$1< i\le d$  and~$\alpha>m{+}1$.

Finally, when~$j$ satisfies~$d<j\le m$, the same 
argument that showed that~$\xi^\alpha_j$ is a multiple of~$\xi^{m+1}_1$
when~$\alpha>m{+}1$ shows that $\xi^\alpha_j$ 
is also a multiple of~$\xi^{m+1}_2$ when~$\alpha>m{+}1$.  
Of course, this implies that~$\xi^\alpha_j = 0$ when~$\alpha>m{+}1$. 

The result of all this vanishing is that
$$
\xi = \eta_{|_E} = \dl\vb_{m+1}\dr\ot
     \bigl(\vb^1\,\xi^{m+1}_1 + \dots + \vb^d\,\xi^{m+1}_d \bigr).
$$
Since~$\eta$ is the identity map, $\xi:E\to Q\ot V^*$ is
just inclusion.  In particular, $E$ is a subspace of~$L\ot V^*$
where~$L = \bbC{\cdot}\dl\vb_{m+1}\dr$, as desired.

For the converse, just note that, when~$E\subset\dl\vb_{m+1}\dr\ot V^*$,
it follows that $\xi^\alpha_i = 0$ when~$\alpha>m{+}1$.  
Since the left hand side of~\eqref{eq:E-an-integral-of-I-2} 
clearly vanishes when~$\alpha_1=\alpha_2$,
it follows that all of these expressions must vanish on~$E$.  Thus, $E$
is an integral element of~$\cI_{(2)}$.
\end{proof}

\begin{remark}[Non-involutivity of~$\cI_{(2)}$]
\label{rem: non-involutive I-2}
Note that~$\cI_{(2)}$ is trivial unless~$n{-}m\ge2$, so assume that
this holds.  Lemma~\ref{lem: integral elements of I2} implies 
that~$\cI_{(2)}$ is not involutive when~$m\ge 2$, 
since its generic integral element of dimension~$1$ does not lie in any 
integral element of dimension~$2$. 
However, each integral element of dimension~$2$ or more lies
in a unique integral element of dimension~$m$.

For~$2\le d\le m$, the space of $d$-dimensional integral elements 
of~$\cI_{(2)}$ in~$T_V\!\Gr(m,n)$ is the same as the set of subspaces 
of type~$(1,\dots,1)^*$ (where the sequence of~$1$s has length~$d$).
\end{remark}

\begin{corollary}\label{cor:Z(phi(2)) computation}
Every element of~$Z(\phi_{(2)})$ is of type~$(1,1)^*$.
\end{corollary}

In particular,  $Z(\phi_{(2)})$ is no larger than it is forced to be by 
Corollary~\ref{cor: integral elements of phi-a}.  The proof is immediate.

\begin{remark}[Walters' results when~$m=2$]
\label{rem: walters m=2 phi 2}
Although she does not remark on this explicitly, the case~$m=2$
of Corollary~\ref{cor:Z(phi(2)) computation} is contained 
implicitly in her proof of Theorem~$5$ of~\cite{mW97}.  Specifically,
her Claim~4.2.3 is equivalent to Corollary~\ref{cor:Z(phi(2)) computation}
in the case~$m=2$.
\end{remark}

\begin{lemma}\label{lem: integral elements of I11}
The maximal integral elements of~$\cI_{(1,1)}$ 
in~$T_V\Gr(m,n)\simeq Q\ot V^*$ 
fall into two distinct classes:
\begin{enumerate}
\item The $(n{-}m)$-dimensional 
      subspaces~$E = Q{\ot}L$, where~$L\subset V^*$ is any line.
\item The $1$-dimensional subspaces~$E\subset Q{\ot}V^*$ that do not
      lie in any subspace of the first kind.
\end{enumerate}
\end{lemma}

\begin{proof}
Apply Lemma~\ref{lem: integral elements of I2} and the complementarity
principle~(\S\ref{sssec:complementarity}).
\end{proof}

\begin{corollary}\label{cor:Z(phi(1,1)) computation}
Every element of~$Z(\phi_{(1,1)})$ is of type~$(2)^*$.
\end{corollary}

In particular,  $Z(\phi_{(1,1)})$ is no larger than it is forced to be by 
Corollary~\ref{cor: integral elements of phi-a}.  The proof is immediate.

\begin{remark}[Walters' results when~$m=2$]
\label{rem: walters m=2 phi 11}
Although she does not remark on this explicitly, the case~$m=2$
of Corollary~\ref{cor:Z(phi(1,1)) computation} is contained 
implicitly in her proof of Theorem~$6$ of~\cite{mW97}.
\end{remark}

Now, for the ideals~$\cI_{(2)^*}$ and~$\cI_{(1,1)^*}$, only the integral
elements of dimension $m(n{-}m)-2$ will be of interest, so I 
state the next two results for those cases only.

\begin{lemma}\label{lem: integral elements of I2*}
Suppose~$2\le m\le n{-}2$.  For any~$V\in \Gr(m,n)$, any codimension~$2$
subspace~$W\subset V^*$, and any hyperplane~$H\subset Q$,
the subspace
\begin{equation}\label{eq: integral elements of I2*}
E = (H{\ot}V^*)+(Q{\ot}W)\subset T_V\Gr(m,n)
\end{equation}
is an integral element of~$\cI_{(2)^*}$ of dimension~$m(n{-}m)-2$.
  
Conversely, if $E\subset T_V\Gr(m,n)$ is an integral element 
of~$\cI_{(2)^*}$ of dimension $m(n{-}m)-2$, 
there exist uniquely a codimension~$2$
subspace~$W\subset V^*$ and a hyperplane~$H\subset Q$
so that~$E$ is of the form~\eqref{eq: integral elements of I2*}.
\end{lemma}

\begin{proof}
Apply Lemma~\ref{lem: integral elements of I2} and the duality
principle~(\S\ref{sssec:duality}).
\end{proof}

\begin{corollary}\label{cor:Z(phi(2)*) computation}
Every element of~$Z(\phi_{(2)^*})$ is of type~$(1,1)$.
\end{corollary}

In particular,  $Z(\phi_{(2)^*})$ is no larger than it is forced to be by 
Corollary~\ref{cor: integral elements of phi-a}.  The proof is immediate.

\begin{lemma}\label{lem: integral elements of I11*}
Suppose~$2\le m\le n{-}2$.  For any~$V\in \Gr(m,n)$, any codimension~$2$
subspace~$W\subset Q$, and any hyperplane~$H\subset V^*$,
the subspace
\begin{equation}\label{eq: integral elements of I11*}
E = (W{\ot}V^*)+(Q{\ot}H)\subset T_V\Gr(m,n)
\end{equation}
is an integral element of~$\cI_{(1,1)^*}$ of dimension $m(n{-}m)-2$.  
Conversely, if $E\subset T_V\Gr(m,n)$ is an integral element 
of~$\cI_{(1,1)^*}$ of dimension~$m(n{-}m)-2$, 
there exist uniquely a codimension~$2$
subspace~$W\subset Q$ and a hyperplane~$H\subset V^*$
so that~$E$ is of the form~\eqref{eq: integral elements of I11*}
\end{lemma}

\begin{proof}
Apply Lemma~\ref{lem: integral elements of I11} and the duality
principle~(\S\ref{sssec:duality}).
\end{proof}

\begin{corollary}\label{cor:Z(phi(1,1)*) computation}
Every element of~$Z(\phi_{(1,1)^*})$ is of type~$(2)$.
\end{corollary}

In particular,  $Z(\phi_{(1,1)^*})$ is 
no larger than it is forced to be by 
Corollary~\ref{cor: integral elements of phi-a}.  
The proof is immediate.

\subsubsection{Dimension~$3$}\label{sssec: dim 3}
I will now treat the cases~$\ab = (3)$,~$(1,1,1)=(3)'$, and~$(2,1)
= (2,1)'$.  For these classes, the structure of the space of integral
elements of~$\cI_\ab$ is more complicated 
than it was for the classes of degree~$2$.
 
It is simpler to first state a result 
that characterizes the \emph{maximal} 
integral elements of these ideals and 
then deduce the structure of the space 
of integral elements of any given dimension from the maximal list.

\begin{remark}[Codimension~$3$]
\label{rem: codim 3}
By complementarity,
the calculations in this subsubsection also determine~$Z(\phi_\ab)$ 
when $\ab = (3)^*$, $(2,1)^*$, and~$(1,1,1)^*$.  
However, I will not actually use these results 
in later sections, so I will not remark on them explicitly.
\end{remark}

\begin{lemma}\label{lem: integral elements of I3}
The maximal integral elements of~$\cI_{(3)}$ 
in~$T_V\Gr(m,n)\simeq Q\ot V^*$ fall into four disjoint classes:
\begin{enumerate}
\item Any $2m$-dimensional subspace~$E = P{\ot}V^*$ where~$P\subset Q$
      is a subspace of dimension~$2$.
\item Any $(m{+}1)$-dimensional subspace~$E = L{\ot}V^*+\bbC{\cdot}R$
      where~$L\subset Q$ is a line and~$R\in Q{\ot}V^*$ is any element
      for which~$\bar R\in (Q/L)\ot V^*$ has rank at least~$2$.
\item Any $3$-dimensional subspace~$E$ that has a basis of the form
$$
(q_2{\ot}l_3-q_3{\ot}l_2,\ q_3{\ot}l_1-q_1{\ot}l_3,\
 q_1{\ot}l_2-q_2{\ot}l_1) 
$$
     where~$(q_1,q_2,q_3)$ and $(l_1,l_2,l_3)$ are each 
     linearly independent in~$Q$ and~$V^*$, respectively.%
     \footnote{This case only occurs when~$m\ge 3$.}
\item Any $2$-dimensional subspace~$E\subset Q\ot V^*$ that is not
      a subspace of an integral element of any of the first three kinds.
\end{enumerate}
\end{lemma}

\begin{remark}[Relations among the types]
\label{rem: integral elements of I3--trivial cases}
When~$n{-}m=2$, the ideal~$\cI_{(3)}$ is empty since~$\phi_3 = 0$.  
In this case, there is only the first type of maximal integral element, 
i.e., the whole tangent space.  
This case will be set aside as trivial in the 
discussion that follows. Also, I remind the
reader that~$m\ge 2$, so that~$2m>m{+}1>2$.  

Only the integral elements of the first type form a closed 
set in the appropriate Grassmannian.  
Indeed, these $2m$-dimensional integral
elements form a smooth variety~$X_1\subset\Gr(2m, Q{\ot}V^*)$ 
that is isomorphic to~$\Gr(2,Q)$.  

The closure of the set of integral elements of the second type 
is a (generally singular) variety~$X_2\subset\Gr(m{+}1, Q{\ot}V^*)$.
Let~$X'_2\subset X_2$ denote the integral elements of the second type.
The `extra' elements in the closure~$X_2$ 
are evidently $(m{+}1)$-dimensional
integral elements of~$\cI_{(3)}$ that lie in a (necessarily unique) 
$2m$-dimensional integral element of the first type.

The closure of the set of integral elements of the third type 
is a (generally singular) variety~$X_3\subset\Gr(3, Q{\ot}V^*)$.
Let~$X'_3\subset X_3$ denote the integral elements of the third type.  
The complement~$X_3\setminus X'_3$ 
can be written as a union~$Y_1\cup Y_2$, 
with~$Y_1$ consisting of integral elements that lie in a
$2m$-dimensional integral element of the first type
and $Y_2$ consisting of integral elements that lie in a 
$(m{+}1)$-dimensional integral element of the second type.  In general,
neither of these two varieties~$Y_i$ contains the other and the
intersection~$Y_1\cap Y_2$ is usually non-empty.

The integral elements of the fourth type form an open subset
of~$\Gr(2, Q{\ot}V^*)$, since, evidently, 
every $2$-plane~$E\subset Q{\ot}V^*$ is 
an integral element of~$\cI_{(3)}$ but, when~$n-m\ge 3$,
the generic $2$-plane~$E\subset Q{\ot}V^*$ does not lie a subspace
of any of the first three types.  In any case, these integral
elements are not of interest, since only integral varieties of~$\cI_{(3)}$ 
of dimension~$3$ or more will be considered in what follows.
\end{remark}

\begin{remark}[The structure of~$Z(\phi_{(3)})$]
\label{rem: 3-dim'l integral elements of I3}
Any $3$-dimensional integral element of~$\cI_{(3)}$ must lie in
a maximal integral element of one of the first three types, so this
affords a description of~$Z(\phi_{(3)})$.  One 
notices immediately is that~$Z(\phi_{(3)})$ contains many $3$-planes
in~$T_V\Gr(m,n)$ that are neither of type~$(2,1)^*$ 
nor of type~$(1,1,1)^*$. 
In fact, the set of subspaces of these types constitutes a
rather small part of~$Z(\phi_{(3)})$, which, for large~$m$ and~$n$
is the union of a large number of distinct~$\SL(n,\bbC)$-orbits.
This will make the analysis of the corresponding integral manifolds 
and varieties of~$\cI_{(3)}$ much more interesting 
than those of~$\cI_{(2)}$. 
\end{remark}

\begin{proof}
I will maintain the basic notation established during the proof of 
Lemma~\ref{lem: integral elements of I2}, especially the
identification~$\eta:T_V\Gr(m,n)\to Q\ot V^*$, which will be used
implicitly throughout the proof.   

Now, the ideal~$\cI_{(3)}$ is generated by the~$(3,0)$-forms
\begin{equation}
\theta^{\alpha_1\alpha_2\alpha_3}_{i_1i_2i_3} 
= \frac16\sum_{\tau\in S_3} \sgn(\tau)\,
\eta^{\alpha_{\tau(1)}}_{i_1}\w\eta^{\alpha_{\tau(2)}}_{i_2}
     \w\eta^{\alpha_{\tau(3)}}_{i_3}\,.
\end{equation}
Note that~$\theta^\alpha_i$ is skewsymmetric in its upper indices
and symmetric in its lower indices.  

As in the proof of Lemma~\ref{lem: integral elements of I2}, 
let~$E\subset T_V\Gr(m,n)$ be an integral element of~$\cI_{(3)}$
of dimension~$d$ and let~$\xi$ be the restriction of~$\eta$ to~$E$.
Then exactly~$d$ of the~$\xi^\alpha_i$ are linearly independent
and they satisfy the cubic relations
\begin{equation}\label{eq:I3 integral element relations}
0 = \sum_{\tau\in S_3} \sgn(\tau)\,
\xi^{\alpha_{\tau(1)}}_{i_1}\w\xi^{\alpha_{\tau(2)}}_{i_2}
     \w\xi^{\alpha_{\tau(3)}}_{i_3}\,.
\end{equation}
where~$\alpha_1<\alpha_2<\alpha_3$ and~$i_1\le i_2\le i_3$.

Before embarking on the classification, I first verify that each
of the four types of subspaces listed in the lemma are indeed integral
elements of~$\cI_{(3)}$.  

If~$E$ is of the first type, then
it is possible to choose the basis~$\vb$ so that~$\xi^\alpha_i=0$
for all~$\alpha>m{+}2$.  In other words~$\xi^\alpha_i$ is zero 
unless~$\alpha=m{+}1$ or~$m{+}2$.  Since the expression on the right
hand side of~\eqref{eq:I3 integral element relations} vanishes identically
unless~$\alpha_1$, $\alpha_2$, and~$\alpha_3$ are distinct, it follows
immediately that these expressions all vanish on~$E$, i.e., that~$E$ is
an integral element of~$\cI_{(3)}$.  

If~$E$ is of the second type, then it is possible to choose 
the basis~$\vb$ so that all of the~$\xi^\alpha_i$ with~$\alpha>m{+}1$
are multiples of a single $1$-form, say~$\psi$.  Again, since 
the expression on the right hand side 
of~\eqref{eq:I3 integral element relations} vanishes identically
unless~$\alpha_1$, $\alpha_2$, and~$\alpha_3$ are distinct, it follows
that every potentially nonzero term in any of these expressions 
contains a wedge product of two forms that are multiples of~$\psi$,
and hence must vanish.  Thus, all of these expressions vanish on~$E$, 
so that~$E$ is indeed an integral element of~$\cI_{(3)}$.  

If~$E$ is of the third type, then it is possible to choose the basis~$\vb$
so that~$\xi^{m+j}_i = 0$ unless $1\le i,j\le 3$ and so that~$\xi^{m+j}_i
=-\xi^{m+i}_j$, while~$\xi^{m+1}_2\w\xi^{m+2}_3\w\xi^{m+3}_1\not=0$.  It is
now not difficult to verify directly that all of the expressions on the
right hand side of~\eqref{eq:I3 integral element relations} vanish.

If~$E$ is of the fourth type, then it has dimension~$2$, so
any $3$-form on~$E$ is trivially zero.  Hence, all of the $2$-dimensional
subspaces~$E$ are integral elements of~$\cI_{(3)}$.  

Now, suppose that~$E\subset T_V\Gr(m,n)$ is an integral element.  There
is nothing to prove unless~$d = \dim E$ is at least~$3$, so assume this.
I am going to show that~$E$ necessarily lies in an integral element of
one of the first three types.  Since no integral element of one of these
types lies in an integral element of a different type, it will then follow
that they are all maximal.%
\footnote{I apologize in advance for the complexity of the argument
to follow.  Unfortunately, I have not been able to discover a simpler one.}

As before, choose the basis~$\vb$ so as to have the maximum number~$p$
of linearly independent~$\xi^\alpha_i$ in the first `column' and make
a basis change so that~$\xi^{m+1}_1,\dots,\xi^{m+p}_1$ are linearly
independent while~$\xi^\alpha_1 = 0$ for~$\alpha>m{+}p$.  Then the
argument made in the course of~Lemma~\ref{lem: integral elements of I2}
shows that all of the forms~$\xi^\alpha_i$ for~$\alpha>m{+}p$ 
must be linear combinations of~$\xi^{m+1}_1,\dots,\xi^{m+p}_1$ 
(or else the maximality of $p$ would be contradicted).

Now, setting~$\alpha_1=m{+}1$, $\alpha_2=m{+}2$, and~$\alpha_3 = m{+}3$ 
and $i_1=i_2=i_3=1$ in~\eqref{eq:I3 integral element relations} yields
$6\,\xi^{m+1}_1\w\xi^{m+2}_1\w\xi^{m+3}_1 = 0$. Thus, $p\le 2$.

First, suppose that~$p=1$.  Since~$\xi^\alpha_i$
is a multiple of~$\xi^{m+1}_1$ when~$\alpha>m{+}1$,
there must be at least~$d\ge 3$ linearly independent forms in the first
`row' of~$\xi$, i.e., among~$\{\xi^{m+1}_1,\dots,\xi^{m+1}_m\}$.
Write~$\xi^\alpha_i = R^\alpha_i\,\xi^{m+1}_1$ for~$\alpha>m{+}1$
and consider the $(m{+}1)$-dimensional subspace~$W$ of~$Q\ot V^*$ 
spanned by the $m$ elements $\eb^i = \dl\vb_{m+1}\dr{\ot}v^i$ 
and the element~$R = R^\alpha_i\,\dl\vb_{\alpha}\dr{\ot} v^i$ (note
that the sum only contains terms with~$\alpha>m{+}1$).  Then~$W$
contains~$E$.  If the rank of~$R$ is greater than~$1$, 
then~$W = \bbC{\cdot}\dl\vb_{m+1}\dr{\ot} V^* + \bbC{\cdot}R$,
so~$W$ is an integral element of~$\cI_{(3)}$ of the second
kind.  If the rank of~$R$ is less than or equal to~$1$, then~$W$ is a 
subspace  of~$P\ot V^*$ where~$P\subset Q$ is any $2$-dimensional subspace 
that contains~$\dl\vb_{m+1}\dr$ and the range of~$R$.   
Thus,~$E$ lies in an integral element of the first kind.  Either way, 
the assumption that~$p=1$ implies that~$E$ is a subspace 
of an integral element of one of the kinds listed in the lemma.

Next, suppose that~$p=2$, so that~$\xi^{m+1}_1\w\xi^{m+2}_1\not=0$, 
but~$\xi^\alpha_1=0$ for~$\alpha>m{+}2$.  Then~$\xi^\alpha_i\equiv 0
\mod~\xi^{m+1}_1,\xi^{m+2}_1$ for all~$\alpha>m{+}2$.  
Since~$d = \dim E\ge 3$, there must be at least one $1$-form 
in~$\{\>\xi^{m+1}_i,\xi^{m+2}_i\mid i>1\>\}$ that is nonzero 
modulo~$\xi^{m+1}_1,\xi^{m+2}_1$.  By making a basis change
in~$\vb_2,\dots,\vb_m$ and in~$\vb_{m+1},\vb_{m+2}$, I can assume 
that~$\xi^{m+1}_2\w\xi^{m+1}_1\w\xi^{m+2}_1\not=0$.  

Suppose, first, that it is possible to make such a basis change so that
the four $1$-forms~$\xi^{m+1}_1,\xi^{m+2}_1,\xi^{m+1}_2,\xi^{m+2}_2$ are
linearly independent.  Since any three elements in any column of~$\xi$
must be linearly dependent, it follows that~$\xi^\alpha_2\equiv 0
\mod~\xi^{m+1}_2,\xi^{m+2}_2$ for all~$\alpha>m{+}2$.  However,
it has already been shown that $\xi^\alpha_2\equiv 0
\mod~\xi^{m+1}_1,\xi^{m+2}_1$ for all~$\alpha>m{+}2$ and the 
linear independence of $\xi^{m+1}_1,\xi^{m+2}_1,\xi^{m+1}_2,\xi^{m+2}_2$
then implies that~$\xi^\alpha_2=0$ for all~$\alpha>m{+}2$.  Once this has
been established, the same argument that showed 
that~$\xi^\alpha_i\equiv0\mod~\xi^{m+1}_1,\xi^{m+2}_1$ 
for all~$\alpha>m{+}2$
can be applied to the second column of~$\xi$ to conclude that 
$\xi^\alpha_i\equiv0\mod~\xi^{m+1}_2,\xi^{m+2}_2$ for all~$\alpha>m{+}2$.  
Combining these two congruences yields that $\xi^\alpha_i=0$ 
for all~$\alpha>m{+}2$.   In other words, $E$ is a subspace of the span
of~$\{\,\dl\vb_{m+1}\dr{\ot}\vb^i,\dl\vb_{m+2}\dr{\ot}\vb^i
\mid 1\le i\le m\,\}$, i.e., $E\subset P\ot V^*$ where~$P\subset Q$
is the $2$-plane spanned by~$\dl\vb_{m+1}\dr$ and~$\dl\vb_{m+2}\dr$.  
Thus, $E$ lies inside an integral element of the first type.

Suppose, then, that for any choice of basis,
~$\xi^{m+1}_1,\xi^{m+2}_1,\xi^{m+1}_2,\xi^{m+2}_2$ are
linearly dependent.  
By making a basis change in~$\vb_{m+1},\vb_{m+2}$, 
I can assume that the linear dependence is 
that~$\xi^{m+2}_2\w\xi^{m+1}_1\w\xi^{m+2}_1=0$, 
i.e., that~$\xi^{m+2}_2=a_1\xi^{m+1}_1+a_2\xi^{m+2}_1$ 
for some~$a_1,a_2\in\bbC$.
By subtracting~$a_2$ times the first column from the second column,
I can assume that~$a_2=0$, so~$\xi^{m+2}_2=a\,\xi^{m+1}_1$ for
some~$a\in\bbC$.

On the other hand, adding~$t$ times the first column to the second column
and wedging together the first (i.e., top), second, 
and $\alpha$-th entries of the result gives
$$
(\xi^{m+1}_2+t\,\xi^{m+1}_1)
  \w(a\,\xi^{m+1}_1+t\,\xi^{m+2}_1)\w\xi^{\alpha}_2\,,
$$
which must vanish for all values of~$t$.  The ~$t^2$-coefficient is
$\xi^{m+1}_1\w\xi^{m+2}_1\w\xi^{\alpha}_2$, 
which is already known to vanish.
The $t$-coefficient is $\xi^{m+1}_2\w\xi^{m+2}_1\w\xi^{\alpha}_2$.  Since
this must vanish as well, it follows that~$\xi^{m+2}_1\w\xi^\alpha_2 = 0$ 
for~$\alpha>m{+}2$, so that there exist numbers~$b^\alpha\in\bbC$ for
$\alpha>m{+}2$ so that $\xi^\alpha_2 = b^\alpha\,\xi^ {m+2}_1$.

First, suppose that~$a\not=0$.  Then the vanishing of the constant
coefficient 
of the above expression yields~$\xi^{m+1}_2\w\xi^{m+1}_1\w\xi^\alpha_2=0$,
which, combined with~$\xi^{m+2}_1\w\xi^\alpha_2 = 0$ implies 
that~$\xi^\alpha_2=0$ for~$\alpha>m{+}2$.  Now, since the top two
entries of the second column of~$\xi$ are linearly independent, the same
argument as was applied to the first column applies to the second and,
indeed,
to any linear combination of the first and second.  In particular, it now
follows that, for all~$t$, 
$$
(\xi^{m+1}_2+t\,\xi^{m+1}_1)
 \w(a\,\xi^{m+1}_1+t\,\xi^{m+2}_1)\w\xi^{\alpha}_i=0
$$
for any~$i>2$ and~$\alpha>m{+}2$.  Using the fact that~$a$ is nonzero
and separating the terms out by $t$-degree then leads to the conclusion
that~$\xi^\alpha_2 = 0$ for all~$i$ and all~$\alpha>m{+}2$.  In other
words, the only nonzero entries of~$\xi$ are in the first two rows.  
Thus, $E$ is a subspace of~$P\ot V^*$ where~$P\subset Q$
is the $2$-plane spanned by~$\dl\vb_{m+1}\dr$ and~$\dl\vb_{m+2}\dr$.  
Thus, $E$ lies inside an integral element of the first type.

Thus, suppose, from now on, that~$a=0$, i.e., that~$\xi^{m+2}_2=0$.

If~$b^\alpha = 0$ for~$\alpha>m{+}2$, then all of the entries in the 
first two columns of~$\xi$ beyond the first two rows are zero.  
In particular, if I were to add~$t$ times the first column to the second, 
I would have a new second column whose only nonzero entries were the
top~$\xi^{m+1}_2+t\,\xi^{m+1}_1$ and the second entry~$t\,\xi^{m+2}_1$.
It would then follow 
that~$(\xi^{m+1}_2+t\,\xi^{m+1}_1)\w(t\,\xi^{m+2}_1)\w
\xi^{\alpha}_i = 0$ for all~$\alpha>m{+}2$ and for all~$t$.  Separating
out the powers of~$t$ in this expression, it would then follow that
$$
  0 = \xi^{m+1}_2\w\xi^{m+2}_1\w\xi^{\alpha}_i 
    = \xi^{m+1}_1\w\xi^{m+2}_1\w\xi^{\alpha}_i \,,
$$
so that~$\xi^{m+2}_1\w\xi^{\alpha}_i=0$ for all~$\alpha>m{+}2$. 
Thus, write $\xi^{\alpha}_i = R^\alpha_i\,\xi^{m+2}_1$ for~$\alpha>m{+}2$.  

If all of the~$R^\alpha_i$ vanish, then, again, $E$ is a subspace 
of~$P\ot V^*$ where~$P\subset Q$ is the $2$-plane 
spanned by~$\dl\vb_{m+1}\dr$ and~$\dl\vb_{m+2}\dr$, so again, 
$E$ lies inside an integral element of the first type.

If not all of the~$R^\alpha_i$ vanish, 
then there is some integer~$r\ge1$
that is the rank of the~$(n{-}m{-}2)$-by-$(m{-}2)$ 
matrix~$(R^\alpha_i)$.
By making a basis change 
in~$\vb_3,\dots,\vb_m$ and~$\vb_{m+3},\dots,\vb_n$,
I can assume that $R^{m+i}_i=1$ for~$3\le i\le r{+}2$ 
and that~$R^\alpha_i=0$
otherwise when~$\alpha>m{+}2$, so I do this.  

I want to show that~$\xi^{m+2}_i\w\xi^{m+2}_1 = 0$ for $i>2$.  
Suppose I can
do this, say, $\xi^{m+2}_i = R^{m+2}_i\,\xi^{m+2}_1$ for all~$i$.  Then
$E$ will be a subspace of the $(m{+}1)$-dimensional integral element 
of~$\cI_{(3)}$ that is spanned by~$\dl\vb_{m+1}\dr{\ot}v^i$ and the
element~$R = R^\alpha_i\,\dl\vb_{\alpha}\dr{\ot}v^i$, i.e., an integral
element of the second type of the lemma, 
and this subcase will be completed.

To prove this claim, first suppose that $3\le i\le r{+}2$, consider
the `column' obtained by first adding $t$ times the first column of~$\xi$
and $s$ times the second column of~$\xi$ to the $i$-th column of~$\xi$, 
and then wedging together the first (i.e., top), 
second, and $i$-th entries. This must vanish, so 
$$
0 = (\xi^{m+1}_i+t\,\xi^{m+1}_1+s\,\xi^{m+1}_2)
\w(\xi^{m+2}_i + t\,\xi^{m+2}_1)\w\xi^{m+2}_1\,.
$$
If this vanishes for all~$t$ and~$s$, then 
$$
  0 = \xi^{m+1}_2\w\xi^{m+2}_i\w\xi^{m+2}_1 
    = \xi^{m+1}_1\w\xi^{m+2}_i\w\xi^{m+2}_1 \,,
$$
so, by the linear independence of~$\xi^{m+1}_1,\xi^{m+2}_1,\xi^{m+1}_2$,
it follows that~$\xi^{m+2}_i\w\xi^{m+2}_1 = 0$, i.e., 
that~$\xi^{m+2}_i = R^{m+2}_i\,\xi^{m+2}_1$ when~$i\le r{+}2$.  
Next, suppose that~$i>r{+}2$.  Consider
the `column' obtained by first adding $t$ times the first column of~$\xi$
and $s$ times the second column of~$\xi$ and then the third
column to the $i$-th column of~$\xi$, and
then wedging together the first (i.e., top), second, and third entries.
This must vanish, so 
$$
0 = (\xi^{m+1}_i+t\,\xi^{m+1}_1+s\,\xi^{m+1}_2+\xi^{m+3}_1)
\w(\xi^{m+2}_i + t\,\xi^{m+2}_1+R^{m+2}_3\,\xi^{m+2}_1)\w\xi^{m+2}_1\,.
$$
Again, since this vanishes for all~$t$ and~$s$,
$$
  0 = \xi^{m+1}_2\w\xi^{m+2}_i\w\xi^{m+2}_1 
    = \xi^{m+1}_1\w\xi^{m+2}_i\w\xi^{m+2}_1 \,,
$$
so, by the linear independence of~$\xi^{m+1}_1,\xi^{m+2}_1,\xi^{m+1}_2$,
it follows that~$\xi^{m+2}_i\w\xi^{m+2}_1 = 0$, i.e., 
that~$\xi^{m+2}_i = R^{m+2}_i\,\xi^{m+2}_1$ when~$i> r{+}2$.   Thus, the
desired claim is established.

The only subcase left to treat now is when not all of the~$b^\alpha$
vanish, so assume this.  
By making a basis change in~$\vb_{m+3},\dots,\vb_n$,
I can assume that~$\xi^{m+3}_2 = \xi^{m+2}_1$, 
but that~$\xi^\alpha_2 = 0$
for~$\alpha>m{+}3$ (and $\alpha = m{+}2$, of course).  

The argument applied to the first column that showed that all
of the forms~$\xi^\alpha_i$ with~$\alpha>m{+}2$ 
must be linear combinations
of~$\xi^{m+1}_1,\xi^{m+2}_1$ can now be applied to the second column.
The result is that all of the forms~$\xi^\alpha_i$ 
with~$\alpha>m{+}3$ or $\alpha=m{+}2 $ must be linear combinations
of~$\xi^{m+1}_2,\xi^{m+2}_1$.  
Explicitly, there are constants~$R^\alpha_i$, $S_i$, $T_i$
when~$\alpha>m{+}1$ so that, when~$i>2$,
$$
\xi^\alpha_i =
\begin{cases}
R^{m+2}_i\,\xi^{m+2}_1 + S_i\,\xi^{m+1}_2 &\qquad \alpha = m{+}2,\\
R^{m+3}_i\,\xi^{m+2}_1 + T_i\,\xi^{m+1}_1 &\qquad \alpha = m{+}3,\\
R^{\alpha\phantom{+3}}_i\,\xi^{m+2}_1     &\qquad \alpha > m{+}3.
\end{cases}
$$
 
If~$S_i = T_i = 0$ for all~$i$, then~$E$ 
lies in the span of the elements~$\dl\vb_{m+1}\dr{\ot}\vb^i$ and
the element
$$
R = \dl\vb_{m+2}\dr{\ot}\vb^1 +  \dl\vb_{m+3}\dr{\ot}\vb^2
    + \sum_{\alpha>m{+}3,\>i} R^\alpha_i\dl\vb_{\alpha}\dr{\ot}\vb^i. 
$$
Consequently,~$E$ lies in an integral element of the second kind
listed in the lemma.

Thus, the subcase that remains to be treated is when not all of 
the~$S_i$ and~$T_i$ vanish, so assume this.  (Note, by the way, that this
subcase can only occur if~$m\ge 3$.)
By subtracting from the~$i$-th column $R^{m+2}_i$ times the first column
and $R^{m+3}_i$ times the second column (which is effected by an
appropriate basis change in~$\vb_1,\dots,\vb_m$), I can actually 
assume that~$R^{m+2}_i = R^{m+3}_i =0$, so I do this.  

I claim that~$S_i + T_i = 0$ for all~$i>2$.  To see this, note that,
adding to the $i$-th column~$t$ times the first column and~$s$ times
the second column and then wedging together the top three entries of the
resulting column gives
$$
0 = (\xi^{m+1}_i+t\,\xi^{m+1}_1+s\,\xi^{m+1}_2)
\w(S_i\,\xi^{m+1}_2 + t\,\xi^{m+2}_1)
\w(T_i\,\xi^{m+1}_1 + s\,\xi^{m+2}_1)\,.
$$
This must vanish for all values of~$s$ and~$t$. Expanding this out
and taking the~$st$ coefficient 
yields~$(S_i+T_i)\,\xi^{m+1}_1\w\xi^{m+1}_2
\w\xi^{m+2}_1$, so $S_i+T_i=0$ as claimed.  

Now, by making a basis change in~$\vb_3,\dots,\vb_m$, I can assume that
$S_i = T_i = 0$ for~$i>3$ while~$S_3 = -T_3 = 1$, so I do this. 

Now, I claim that $\xi^\alpha_i=0$ for~$\alpha>m{+}3$.  This has already
been established for~$i=1$ and~$2$.  If there were some~$\alpha>m{+}3$
for which~$\xi^\alpha_3= R^\alpha_3\,\xi^{m+2}_1\not=0$, then the
second, third, and $\alpha$-th entries of the third column would be 
linearly independent, contrary to hypothesis.  Thus, $\xi^\alpha_3 =0$
for~$\alpha>m{+}3$.   
Now, if $\xi^\alpha_i= R^\alpha_i\,\xi^{m+2}_1\not=0$,
for some~$i>3$ and~$\alpha>m{+}3$, then adding the $i$-th column of~$\xi$
to the third column will produce a column with three linearly independent
entries.  
Thus, $\xi^\alpha_i=0$ for~$\alpha>m{+}3$ and all~$i$, as claimed.

The entries of~$\xi$ that remain to be understood 
are $\{\xi^{m+1}_3,\dots,\xi^{m+1}_m\}$ (the remainder of the first row).
Since $0 = \xi^{m+1}_3\w\xi^{m+2}_3\w\xi^{m+3}_3
= -\xi^{m+1}_3\w\xi^{m+1}_2\w\xi^{m+1}_1$, there are
constants~$c^1_3,c^2_3\in\bbC$ so 
that~$\xi^{m+1}_3=c^1_3\,\xi^{m+1}_1+c^2_3\,\xi^{m+1}_2$.
Adding~$t$ times the first column and $s$ times the second column to the
third column and wedging the top three entries yields
\begin{align*}
0 &= \bigl((t+c^1_3)\,\xi^{m+1}_1+(s+c^2_3)\,\xi^{m+1}_2\bigr)
\w(\xi^{m+1}_2 + t\,\xi^{m+2}_1)
\w(-\xi^{m+1}_1 + s\,\xi^{m+2}_1)\\
&= (c^1_3\,s-c^2_3\,t)\,\xi^{m+1}_1\w\xi^{m+1}_2\w\xi^{m+2}_1\,.
\end{align*}
Since this must vanish for all~$s$ and~$t$, this gives~$c^1_3=c^2_3=0$.
Thus~$\xi^{m+1}_3=0$.

For~$i>3$, adding the $i$-th column to the third column has the effect
of replacing~$\xi^{m+1}_3$ by~$\xi^{m+1}_i$ in the upper left hand
$3$-by-$3$ minor.  The above argument can then be repeated to conclude
that~$\xi^{m+1}_i = 0$ as well.  

Now, exchanging the first and third rows and then multiplying
the top row by~$-1$ yields a~$\xi$ whose upper left hand $3$-by-$3$
minor is of the form
$$
\begin{pmatrix} 
0 &-\psi^3 & \phantom{-}\psi^2\\
\phantom{-}\psi^3& 0 & -\psi^1\\
-\psi^2 & \phantom{-}\psi^1 & 0
\end{pmatrix},
\qquad\qquad (\psi^1\w\psi^2\w\psi^3\not=0)
$$
while all of the other entries of~$\xi$ vanish. 
Thus~$E$ has dimension~$3$
and has the third type listed in the lemma.

Finally, it has been shown that every integral element of~$\cI_{(3)}$
of dimension at least~$3$ lies in either a ($2m$-dimensional) integral
element
of the first type, a ($(m{+}1)$-dimensional) 
integral element of the second 
type, or a ($3$-dimensional) integral element of the third type.
It only remains to observe that none of the
integral elements of the second type 
lie in an integral element of the first 
type, and none of the integral elements of the third
type lie in an integral element of either of the first two types.
Thus, the first three types listed in the statement of the lemma are each 
maximal.  The only integral elements not accounted for are the maximal ones 
of dimension at most~$2$.  Since every $2$-dimensional subspace is an
integral
element, the ones that do not lie in a subspace of any of the first three 
types must be maximal.   The classification is now complete.
\end{proof}

\begin{lemma}\label{lem: integral elements of I111}
The maximal integral elements of~$\cI_{(1,1,1)}$ 
in~$T_V\Gr(m,n)\simeq Q\ot V^*$ fall into four disjoint classes:
\begin{enumerate}
\item Any $2(n{-}m)$-dimensional subspace~$E = Q\ot P$ 
       where~$P\subset V^*$ is a subspace of dimension~$2$.
\item Any $(n{-}m{+}1)$-dimensional subspace~$E = Q\ot L + \bbC{\cdot} R$
      where $L\subset V^*$ is a line and~$R\in Q\ot V^*$ is any element
      for which~$\bar R\in Q\ot V^*/L$ has rank at least~$2$.
\item Any $3$-dimensional subspace~$E$ that has a basis of the form
$$
(q_2{\ot}l_3-q_3{\ot}l_2,\ q_3{\ot}l_1-q_1{\ot}l_3,\
 q_1{\ot}l_2-q_2{\ot}l_1) 
$$
     where~$(q_1,q_2,q_3)$ and $(l_1,l_2,l_3)$ 
     are each linearly independent in~$Q$ and~$V^*$, respectively.%
     \footnote{This case only exists when~$n{-}m\ge 3$.}
\item Any $2$-dimensional subspace~$E\subset Q\ot V^*$ that is not
      a subspace of an integral element of any of the first three kinds.
\end{enumerate}
\end{lemma}

\begin{proof}
Apply Lemma~\ref{lem: integral elements of I3} and the complementarity
principle~(\S\ref{sssec:complementarity}).
\end{proof}

Before going on to study the case~$\ab = (2,1) = (2,1)'$, I state 
the following `combined' result.  It will be needed in the next section.

\begin{lemma}\label{lem: integral elements of I111 and I3}
The maximal integral elements of~$\cI_{(1,1,1)}\cup\cI_{(3)}$ 
in~$T_V\Gr(m,n)\simeq Q\ot V^*$ fall into five disjoint types:
\begin{enumerate}
\item Any $4$-dimensional subspace~$E = W\ot P$ 
    where~$W\subset Q$ and~$P\subset V^*$ are subspaces of dimension~$2$.
\item Any $3$-dimensional subspace~$E$ that has a basis of the form
$$
(q_2{\ot}v^3-q_3{\ot}v^2,\ q_3{\ot}v^1-q_1{\ot}v^3,
\ q_1{\ot}v^2-q_2{\ot}v^1)
$$
   where~$(q_1,q_2,q_3)$ and $(v^1,v^2,v^3)$ are each 
   linearly independent in~$Q$ and~$V^*$, respectively.
\item Any $3$-dimensional subspace~$E$ that has a basis of the form
$$
(q_2{\ot}v^3,\ -q_1{\ot}v^3,\ q_1{\ot}v^2-q_2{\ot}v^1) 
$$
      where~$(q_1,q_2)$ and $(v^1,v^2,v^3)$ are each 
     linearly independent in~$Q$ and~$V^*$, respectively.
\item Any $3$-dimensional subspace~$E$ that has a basis of the form
$$
(-q_3{\ot}v^2,\ q_3{\ot}v^1,\ q_1{\ot}v^2-q_2{\ot}v^1) 
$$
      where~$(q_1,q_2,q_3)$ and $(v^1,v^2)$ are each linearly independent 
     in~$Q$ and~$V^*$, respectively.
\item Any $2$-dimensional subspace~$E\subset Q\ot V^*$ that is not
      a subspace of an integral element of any of the first four kinds.
\end{enumerate}
\end{lemma}

\begin{proof}
Combine Lemmas~\ref{lem: integral elements of I3} 
and~\ref{lem: integral elements of I111}.
\end{proof}

\begin{remark}[$\mathcal{R}_{(2,1)}=Z(\phi_{(1,1,1)})\cap Z(\phi_{(3)})$]
\label{rem:I111+3 integral element closures and intersections}
Lemma~\ref{lem: integral elements of I111 and I3} allows a rather
complete description of the three dimensional integral elements
of~$\cI_{(1,1,1)} \cup \cI_{(3)}$, which is what Maria Walters calls
the Schubert differential system~$\mathcal{R}_{(2,1)}$. 
(See \S\ref{sssec: the two differential systems}.)

The set~$Y_1\subset\Gr(4,Q{\ot}V^*)$ 
of $4$-dimensional integral elements of~$\cI_{(1,1,1)} \cup \cI_{(3)}$
of the first type is a variety isomorphic to~$\Gr(2,Q)\times\Gr(2,V^*)$.  
The set~$X_1\subset\Gr(3,Q{\ot}V^*)$ consisting of $3$-dimensional 
subspaces~$E$ lying in a $4$-dimensional integral element~$E^+$, 
i.e., an element of~$Y_1$, is also a closed submanifold.  
In fact, because the extension $E\mapsto E^+$ is unique, this defines
an algebraic submersion~$X_1\to Y_1$ whose fiber over~$E^+\in V_1$
is simply~$\Gr(3,E^+)\simeq\bbP^3$.  From this, one can show that~$X_1$ 
is a smooth submanifold of~$\Gr(3,Q{\ot}V^*)$ of dimension~$2n{-}5$.

The space~$X_1$ contains the subvariety~$B_{(2,1)^*}$, 
consisting of the subspaces of type~$(2,1)^*$, as a hypersurface.  
Both $B_{(2,1)^*}$ and its complement~$X_1' = X_1\setminus B_{(2,1)^*}$ 
are single $\bigl(\GL(Q)\times\GL(V)\bigr)$-orbits in~$\Gr(3,Q{\ot}V^*)$.

For~$i = 2$, $3$, or~$4$, let ~$X_i\subset\Gr(3,Q{\ot}V^*)$ 
denote the closure of the set~$X'_i$ of $3$-dimensional 
integral elements of type~$(i)$ in the list 
of Lemma~\ref{lem: integral elements of I111 and I3}.  Each of
the spaces~$X_i'$ is a single $\bigl(\GL(Q)\times\GL(V)\bigr)$-orbit 
in~$\Gr(3,Q{\ot}V^*)$ and is not closed.

If $m\ge3$ and~$n{-}m\ge3$, the sets~$X_2$, $X_3$, and~$X_4$ 
are nonempty and they have dimensions~$3n{-}10$,
$2n{+}m{-}8$ and~$3n{-}m{-}8$, respectively. 
Furthermore,~$X_2\setminus X'_2= X_3\cup X_4$.
One can show that $X_{34}=X_3\cap X_4=X_1\cap X_3=X_1\cap X_4$  
is the complement of~$X'_3$ in~$X_3$ and of~$X'_4$ in~$X_4$. 
In fact, this intersection is simply~$B_{(2,1)^*}$.

If~$m\ge 3$ but~$n{-}m=2$, then $X_2$ and~$X_4$ are empty, 
but $X_3$ is nonempty and has dimension~$3m{-}4$.  
The intersection~$X_1\cap X_3$ is~$B_{(2,1)^*}$.

If~$m=2$ but~$n{-}m=3$, then $X_2$ and~$X_3$ are empty, 
but $X_4$ is nonempty and has dimension~$2m{+}1$.  
The intersection~$X_1\cap X_4$ is~$B_{(2,1)^*}$.

Of course, if~$m=n{-}m=2$, then~$X_2$, $X_3$, and $X_4$ are 
empty.  But in this case, $\cI_{(1,1,1)}$ and~$\cI_{(3)}$ are
both trivial ideals and every $3$-plane is an integral element.

Note that, except in this last (trivial) case, the space 
$R_{(2,1)}\subset\Gr(3,Q\otimes V^*)$ has two irreducible components 
and that they intersect in the locus~$B_{(2,1)^*}$.  
This closure information can be displayed in a diagram
\begin{equation}
\label{eq: R21 closure diagram}
\begin{matrix}
    &        & X_3'  &\\
    &\nearrow&           &\searrow\\
X_2'&    &\longrightarrow&     &\mathcal{B}_{(2,1)^*}&\longleftarrow&X_1'\\
    &\searrow&           &\nearrow\\
    &        & X_4' &
\end{matrix}
\end{equation}
where each of the five entries is a single 
$\bigl(\GL(Q)\times\GL(V)\bigr)$-orbit
and each arrow points from a given orbit to an orbit in its closure. 

I will have more to say about the geometry of these integral 
elements in the next section.
\end{remark}

\begin{lemma}\label{lem: integral elements of I21}
The maximal integral elements of~$\cI_{(2,1)}$ 
in~$T_V\Gr(m,n)\simeq Q\ot V^*$ fall into three disjoint classes:
\begin{enumerate}
\item Any $m$-dimensional subspace~$E = L{\ot}V^*$ where~$L\subset Q$
      is a line.
\item Any $(n{-}m)$-dimensional subspace~$E = Q{\ot}L$
      where~$L\subset V^*$ is a line.
\item Any $2$-dimensional subspace~$E\subset Q\ot V^*$ that is not
      a subspace of an integral element of either of the first two kinds.
\end{enumerate}
\end{lemma}

\begin{proof}
Again, the notation established in the previous proof-analyses 
will be maintained.  The first difference is that the ideal~$\cI_{(2,1)}$
is generated by $(3,0)$-forms of the form
\begin{equation}
\theta_{i_1i_2i_3}(c) 
= \sum_{\alpha\in[m{+}1,n]^3} c_{\alpha_1\alpha_2\alpha_3}\,
\eta^{\alpha_1}_{i_1}\w\eta^{\alpha_2}_{i_2}\w\eta^{\alpha_3}_{i_3}\,.
\end{equation}
where~$c:[m{+}1,n]^3\to\bbC$ satisfies the relations
\begin{equation}\label{eq:c in S(2,1)}
c_{\alpha_1\alpha_2\alpha_3} = c_{\alpha_2\alpha_1\alpha_3},
\qquad\qquad
c_{\alpha_1\alpha_2\alpha_3}
+c_{\alpha_2\alpha_3\alpha_1}+c_{\alpha_3\alpha_1\alpha_2} = 0.
\end{equation}
(Essentially,~$c$ is the general element of~$\bbS_{(2,1)}(\C{n-m})$.)
Note that~$\theta_i(c)$ satisfies
\begin{equation}\label{eq: theta(c) relations}
\theta_{i_1i_2i_3}(c) = -\theta_{i_2i_1i_3}(c),
\qquad\qquad
\theta_{i_1i_2i_3}(c)+\theta_{i_2i_3i_1}(c)+\theta_{i_3i_1i_2}(c) = 0.
\end{equation}

As in the proofs of Lemmas~\ref{lem: integral elements of I2}
and~\ref{lem: integral elements of I3}, 
let~$E\subset T_V\Gr(m,n)$ be an integral element of~$\cI_{(2,1)}$
of dimension~$d$ and let~$\xi$ be the restriction of~$\eta$ to~$E$.
Then exactly~$d$ of the~$\xi^\alpha_i$ are linearly independent
and they satisfy the cubic relations
\begin{equation}\label{eq:I21 integral element relations}
0 = \sum_{\alpha\in[m{+}1,n]^3} c_{\alpha_1\alpha_2\alpha_3}\,
\xi^{\alpha_1}_{i_1}\w\xi^{\alpha_2}_{i_2}\w\xi^{\alpha_3}_{i_3}\,.
\end{equation}
for all~$c$ that satisfy the relations~\eqref{eq:c in S(2,1)}.

Before going on to the classification, it is a good idea to verify
that the subspaces listed in the statement of the lemma are indeed
integral elements of~$\cI_{(2,1)}$.  

If~$E = L{\ot}V^*$ for some
line~$L\subset Q$, then it is possible to choose the basis~$\vb$ so
that~$L$ is spanned by~$\dl\vb_{m+1}\dr$.  
In this case, $\xi^\alpha_i = 0$
for all~$\alpha>m{+}1$.  Since the relations~\eqref{eq:c in S(2,1)} imply
that~$c_{\alpha\alpha\alpha}=0$ for all~$\alpha$, it follows that
the right hand side of~\eqref{eq:I21 integral element relations} must
vanish identically for all~$c$ satisfying~\eqref{eq:c in S(2,1)}.  
Thus,~$L{\ot}V^*$ is an integral element of~$\cI_{(2,1)}$.

If~$E = Q{\ot}L$ for some line~$L\subset V^*$, then it is possible to
choose the basis~$\vb$ so that~$L$ is spanned by~$\vb^1$.  
In this case,
$\xi^\alpha_i=0$ for all~$i>1$.  Thus, the right hand side 
of~\eqref{eq:I21 integral element relations} 
vanishes unless~$i_1=i_2=i_3=1$.
However, the remaining expression
$$
 \sum_{\alpha\in[m{+}1,n]^3} c_{\alpha_1\alpha_2\alpha_3}\,
\xi^{\alpha_1}_{1}\w\xi^{\alpha_2}_{1}\w\xi^{\alpha_3}_{1}
$$
vanishes because~$c_{\alpha_1\alpha_2\alpha_3}
=c_{\alpha_2\alpha_1\alpha_3}$.
Thus~$Q{\ot}L$ is an integral element of~$\cI_{(2,1)}$.

Now, on to the classification.
Fix~$\alpha$, and~$\beta$ satisfying~$m < \alpha \not= \beta \le n$.  
Let~$c:[m{+}1,n]^3\to\bbC$ satisfy~$c_{\alpha\alpha\beta} = 2$ 
while~$c_{\alpha\beta\alpha}=c_{\beta\alpha\alpha} = -1$ 
and suppose further 
that~$c_{\alpha_1\alpha_2\alpha_3}=0$ except in these three cases.  
Then~$c$ satisfies~\eqref{eq:c in S(2,1)}.  
The relation~\eqref{eq:I21 integral element relations} specializes 
in this case to
$$
0 =  2\,\xi^{\alpha}_{i_1}\w\xi^{\alpha}_{i_2}\w\xi^{\beta}_{i_3}
- \xi^{\alpha}_{i_1}\w\xi^{\beta}_{i_2}\w\xi^{\alpha}_{i_3}
- \xi^{\beta}_{i_1}\w\xi^{\alpha}_{i_2}\w\xi^{\alpha}_{i_3}\,.
$$
Now, setting~$i_1 = i$ and~$i_2 = i_3 = j$, 
this relation reduces to the simple relation
\begin{equation}\label{eq: simple I21 relation}
0 =  3\,\xi^{\alpha}_{i}\w\xi^{\alpha}_{j}\w\xi^{\beta}_{j}\,.
\end{equation}
Thus,~\eqref{eq: simple I21 relation}
holds whenever~$1\le i, j\le m < \alpha,\beta \le n$.

The relation~\eqref{eq: simple I21 relation} must hold on~$E$ and, 
moreover, because the condition of
being an integral element of~$\cI_{(2,1)}$ is unaffected by the choice
of basis~$\vb$, it follows that these relations among triples of matrix
entries in any $2$-by-$2$ minor of~$\xi$ must continue to hold when
$\xi$ is pre- or post-multiplied by any matrices.  This device will
be very helpful in what follows.

Now suppose that~$E\subset T_V\Gr(m,n)$ is an integral element 
of~$\cI_{(2,1)}$ of dimension~$d\ge 3$.  (Unless $d\ge 3$, there
is nothing to prove.)

Suppose that the basis~$\vb$ has been chosen so as to have the
maximum number~$p$ of linearly independent forms 
in the first column of~$\xi$
and that, moreover, it has been arranged that~$\xi^\alpha_1 = 0$ 
for~$\alpha>m{+}p$.  Then all of the 1-forms~$\xi^\alpha_i$ 
with~$\alpha>m{+}p$ must be linear combinations 
of~$\xi^{m+1}_1,\dots,\xi^{m+p}_1$ (otherwise, the maximality of~$p$
would be contradicted).   

Suppose, first, that~$p>1$.  Then, for~$j>1$ 
and any~$\alpha$ and~$\beta$ 
satisfying~$m{+}1\le \alpha<\beta\le m{+}p$,
the relation~$\xi^{\alpha}_{j}\w\xi^{\alpha}_{1}\w\xi^{\beta}_{1} = 0$
shows that~$\xi^\alpha_j$ is a linear combination 
of~$\xi^{\alpha}_1$ and $\xi^{\beta}_1$, so it follows 
that~$\xi^{m+1}_1,\dots,\xi^{m+p}_1$ 
is actually a basis for the $1$-forms
on~$E$.  In other words,~$p = d\ge 3$.  
Thus, choosing any~$\alpha$, $\beta$,
and~$\gamma$ satisfying $m{+}1\le \alpha<\beta<\gamma\le m{+}p$ and~$j$
satisfying~$1<j\le m$, the relations
$$
  \xi^{\alpha}_{j}\w\xi^{\alpha}_{1}\w\xi^{\beta}_{1} 
 =\xi^{\alpha}_{j}\w\xi^{\alpha}_{1}\w\xi^{\gamma}_{1} = 0
$$
and the independence 
of~$\{\xi^{\alpha}_{1},\xi^{\beta}_{1},\xi^\gamma_1\}$
imply that~$\xi^{\alpha}_{j}\w\xi^{\alpha}_{1} = 0$.  
In other words, there are constants~$R^\alpha_j$ so that 
$\xi^{\alpha}_{j}= R^\alpha_j\,\xi^{\alpha}_{1}$
when~$m{+}1\le\alpha\le m{+}p$.

Making a basis change in~$\vb_1,\dots,\vb_m$ (which has the 
effect of post-multiplying~$\xi$ by an invertible $m$-by-$m$ matrix),
I can assume that~$R^{m+1}_j=0$ for~$j>1$.   
I claim that, for this choice
of basis,~$R^\alpha_j=0$ whenever~$1<j\le m<\alpha\le m{+}p$.
To see this, fix any~$j$ and~$q$ satisfying $1<j\le m< m{+}q\le m{+}p$.  
Add the $q$-th row of~$\xi$ to the first row, 
resulting in a new matrix~$\tilde\xi$.
Choose an~$r\not=1,q$ with~$r<p$ and wedge together the $(1,1)$, $(1,j)$,
and~$(r,1)$ entries of~$\tilde\xi$, obtaining
$$
 0 = \tilde\xi^{m+1}_j\w\tilde\xi^{m+1}_1\w\tilde\xi^{m+r}_1 
   = R^{m+q}_j\,\xi^{m+q}_1\w\xi^{m+1}_1\w\tilde\xi^{m+r}_1 \,.
$$
However~$\{\xi^{m+q}_1,\xi^{m+1}_1,\tilde\xi^{m+r}_1\}$ are linearly
independent, so~$R^{m+q}_j=0$.
Thus,~$\xi^\alpha_j=0$ for~$1<j\le m<\alpha\le m{+}p$, as desired.

Now, if any~$\xi^\alpha_j$ with~$1<j\le m$ and~$m{+}p<\alpha\le n$ 
were nonzero, I could add the row that it appears on to, say,
the top row, and get a new~$\xi$ that still satisfies all of the 
hypotheses so far but has a nonzero entry  on the top row after the
first column.  Since I have just shown that this is impossible, it
follows that~$\xi^\alpha_j = 0$ for all~$j>1$.  

Of course, this implies that~$E\subset Q\ot \bbC{\cdot}\vb^1$, so 
that~$E$ lies in an integral element of the second kind.

Now suppose, instead, that~$p=1$.  Then, by the first part of the
argument, there must be $(d{-}1)$ $1$-forms among 
the~$\{\xi^{m+1}_2,\dots,\xi^{m+1}_m\}$ that are linearly independent 
modulo~$\xi^{m+1}$.  By a change of basis in~$\vb$, I can assume that
$\{\xi^{m+1}_1,\dots,\xi^{m+1}_d\}$ are linearly independent and that
$\xi^{m+1}_j=0$ for~$j>d$.  Recall that, by hypothesis,~$d\ge 3$.  

Now, I claim that ~$\xi^\alpha_j=0$ for all~$\alpha>m{+}1$.
To see this, first note that, when~$1<i\not=j\le d$, the relation
$\xi^{m+1}_i\w\xi^{m+1}_j\w\xi^\alpha_j=0$ implies that~$\xi^\alpha_j$
is a linear combination of~$\{\xi^{m+1}_i,\xi^{m+1}_j\}$
for~$\alpha>{m+1}$ and~$1<j\le d$.  However, the maximality of~$p$ 
has already shown that~$\xi^\alpha_j\w\xi^{m+1}_1=0$.  
Thus,~$\xi^\alpha_j=0$ when~$j\le d$.  

If~$\xi^\alpha_j$ were nonzero for some~$j>d$, then adding the $j$-th
column of~$\xi$ to the second column would produced a~$\tilde\xi$ that
still satisfied the~$\cI_{(2,1)}$ relations, but had a nonzero entry 
in the second column other than the top entry.  It has just been shown,
though, that this is impossible.  Thus, ~$\xi^\alpha_j=0$ 
whenever~$\alpha>m{+}1$.  

Of course, this implies that~$E\subset 
\bbC{\cdot}\dl\vb_{m+1}\dr\ot V^*$, so that~$E$ 
lies in an integral element of the first kind.  

Thus, the argument has shown that any integral element of~$\cI_{(2,1)}$
of dimension~$3$ or more lies in an integral element of one of the first 
two types listed in the statement of the lemma, as desired.
\end{proof}

\begin{corollary}\label{cor: Z(phi(2,1))}
Every $3$-dimensional integral element 
of~$\cI_{(2,1)}$ is a subspace that is either of type~$(3)^*$ or of
type~$(1,1,1)^*$. 
\end{corollary}

\begin{remark}[Non-involutivity of~$\cI_{(2,1)}$]
\label{rem: non-involutivity of I 21}
By Lemma~\ref{lem: integral elements of I21}, the ideal ~$\cI_{(2,1)}$ 
has no $3$-dimensional integral elements when~$m=n{-}m=2$, 
a fact that could have been seen directly 
from~\eqref{eq:dual-integral-elements} since~$(2,1)^* = (1)$ 
in this case.

On the other hand, Corollary~\ref{cor: Z(phi(2,1))} implies that
the set~$Z\bigl(\phi_{(2,1)}\bigr)$
has two components when both~$m$ and~$n{-}m$ are at least~$3$.
Each of these components is smooth and closed.

Since the intersection of an integral element of~$\cI_{(2,1)}$ 
of the first kind listed in Lemma~\ref{lem: integral elements of I21}
with an integral element of the second kind 
is at most of dimension~$1$, it follows that every integral 
element of~$\cI_{(2,1)}$ of 
dimension~$2$ or more lies in a unique maximal integral element. 

Finally,~$\cI_{(2,1)}$ is not involutive for integral
manifolds of dimension~$3$ or more since the generic $2$-plane
does not lie in any $3$-dimensional integral element.
\end{remark}

\subsubsection{General remarks on higher degrees}
\label{sssec: higher degs}

As the reader will have noticed, the analysis of the integral elements
of~$\cI_{(3)}$ was considerably more difficult than the analysis
of the integral elements of~$\cI_{(2)}$ and also more difficult than
the analysis of the integral elements of~$\cI_{(2,1)}$.  In this
last subsection of this section, I will collect together a few remarks
about the calculations in general.

First of all, calculation of the integral elements 
of~$\cI_{\ab}$ for general~$\ab$ seems to be difficult.  Of course, 
Corollary~\ref{cor: integral elements of phi-a} provides an important 
`lower bound' for these integral elements, but the example of~$\cI_{(3)}$ 
shows that this can be very far from a complete description.

In general, when~$|\ab| = p$, recall that~$\cI_{\ab}$ is generated 
by the $\GL(n,\bbC)$-invariant subspace
$\bbS_\ab(V^*)\ot\bbS_{\ab'}(Q)\subset\L^{p,0}\bigl(T^*\Gr(m,n)\bigr)$.
Given this, is not difficult to show 
that~$\pi_m^*(\cI_\ab)\subset \Omega^p\bigl(\SU(n)\bigr)$ 
is generated by the forms
\begin{equation}\label{eq: generators of I-a}
\theta_{i_1\dots i_p}(c) 
= \sum_{\alpha\in[m{+}1,n]^p} c_{\alpha_1\dots\alpha_p}\,
\omega^{\alpha_1}_{i_1}\w\dots\w\omega^{\alpha_p}_{i_p}\,,
\end{equation}
where~$c:[m{+}1,n]^p\to\bbC$ ranges over the elements 
of~$\bbS_{\ab'}(\C{n-m})$.  Moreover, the forms~$\theta_{i_1\dots i_p}(c)$
for fixed~$c$ have the same $i$-index symmetry as the general element
of~$\bbS_\ab(\C{m})$.

Consequently, in the notation for integral elements~$E\subset Q\ot V^*$
that has been employed in this section, one sees that the corresponding
matrix~$\xi$ must satisfy the relations
\begin{equation}\label{eq: relations on xi for I-a}
0= \sum_{\alpha\in[m{+}1,n]^p} c_{\alpha_1\dots\alpha_p}\,
  \xi^{\alpha_1}_{i_1}\w\dots\w\xi^{\alpha_p}_{i_p}\,,
\end{equation}
where~$c:[m{+}1,n]^p\to\bbC$ ranges over the elements 
of~$\bbS_{\ab'}(\C{n-m})$.  Unfortunately, these relations appear 
to be rather difficult to understand directly except in the simplest
cases.  

\begin{example}\label{ex: I(p) generators}
Consider~$\ab = (p)$ (where~$p\le n{-}m$, of course).
Since~$\bbS_{(p)}(\C{m}) = S^p(\C{m})$ while~$\bbS_{(p)'}(\C{n-m})
= \L^p(\C{n-m})$, the above relations are equivalent to
\begin{equation}\label{eq: relations on xi for I(p)}
0= \sum_{\tau\in S_p} \sgn(\tau)\,
  \xi^{\alpha_{\tau(1)}}_{i_1}\w\dots\w\xi^{\alpha_{\tau(p)}}_{i_p}\,,
\end{equation}
whenever~$1\le i_1\le\dots\le i_p\le m <\alpha_1<\dots<\alpha_p\le n$.
Now, the expression on the right is symmetric in~$i_1,\dots,i_p$,
which motivates considering the general linear 
combination~$\xi^\alpha(t) = t^i\xi^\alpha_i$ and rewriting the above
relation in the form
\begin{equation}\label{eq: highest weight relations on xi for I(p) }
0= \xi^{\alpha_1}(t)\w\dots\w\xi^{\alpha_p}(t)\,.
\end{equation}
Thus, the relations~\eqref{eq: relations on xi for I(p)} are equivalent 
to the condition that the wedge product of any $p$ of the entries of any 
linear combination of the columns of~$\xi$ should vanish.  

Note that this last formulation is precisely the condition stated by 
Griffiths and Harris as~\cite[(4.6)]{MR81k:53004} for the vanishing of the 
Chern form~$c_p$, as it should be.  However, a glance at the analysis
of the integral elements of~$\cI_{(3)}$ shows that this formulation,
though an important first step, is still very far from a determination 
of the integral elements of~$\cI_{(3)}$.
\end{example} 

One strategy for studying the 
relations~\eqref{eq: relations on xi for I-a}
is to choose~$c$ to be a highest weight vector for the representation~
$\bbS_{\ab'}(\C{n-m})$ and/or to combine the relations so as to reflect
a highest weight vector for~$\bbS_{\ab}(\C{m})$.  This usually gives
the simplest relations.  

For example, the formulation~
\eqref{eq: highest weight relations on xi for I(p) } is nothing but
considering the relation for the orbit of a highest weight vector in
$\bbS_{\ab}(\C{m})= S^p(\C{m})$.  Similarly, the 
relation~\eqref{eq: simple I21 relation} that was so fundamental to
the analysis of integral elements of~$\cI_{(2,1)}$ is merely the
relation corresponding to a highest weight vector.  

Since an irreducible representation is spanned by the orbit of its
highest weight vector, 
all the relations~\eqref{eq: relations on xi for I(p)}
will be generated by starting with a highest weight relation and 
considering all the relations it implies after pre- and post-multiplying
$\xi$ by arbitrary invertible matrices of the appropriate size.  This
was essentially the strategy I used in constructing the proofs of 
the various Lemmas in this section.

Finally, since, among all the 
representations~$\bbS_\ab(V)$ with~$|\ab|=p$,
the ones with~$\ab = (p)$ and~$\ab = (p)'$ are 
generally the lowest dimensional, the ideals~$\cI_{(p)}$ and~$\cI_{(p)'}$
are usually the smallest in size.  For that reason, one might expect that
their integral elements would display a greater variety than the
integral elements for~$\cI_\ab$ with other~$\ab$ of the same degree.
This expectation was born out in the degree~$3$ case, since the analysis
for~$\cI_{(3)}$ and~$\cI_{(1,1,1)}$ was considerably more complicated than
the analysis for~$\cI_{(2,1)}$, an ideal with approximately four times
as many generators as either of the other two.

The reader experienced with exterior differential systems will know to take
this sort of `dimension count' with a grain of salt, since it is usually
more subtle algebraic features than the rank of an ideal that
play the major role in determining the integral elements and integral
manifolds.  However, this `dimension count' does seem to correspond 
somewhat to the complexity of the analysis in each case, so I offer it
as an observation to the interested reader.

\begin{figure}
\includegraphics[width=\linewidth]{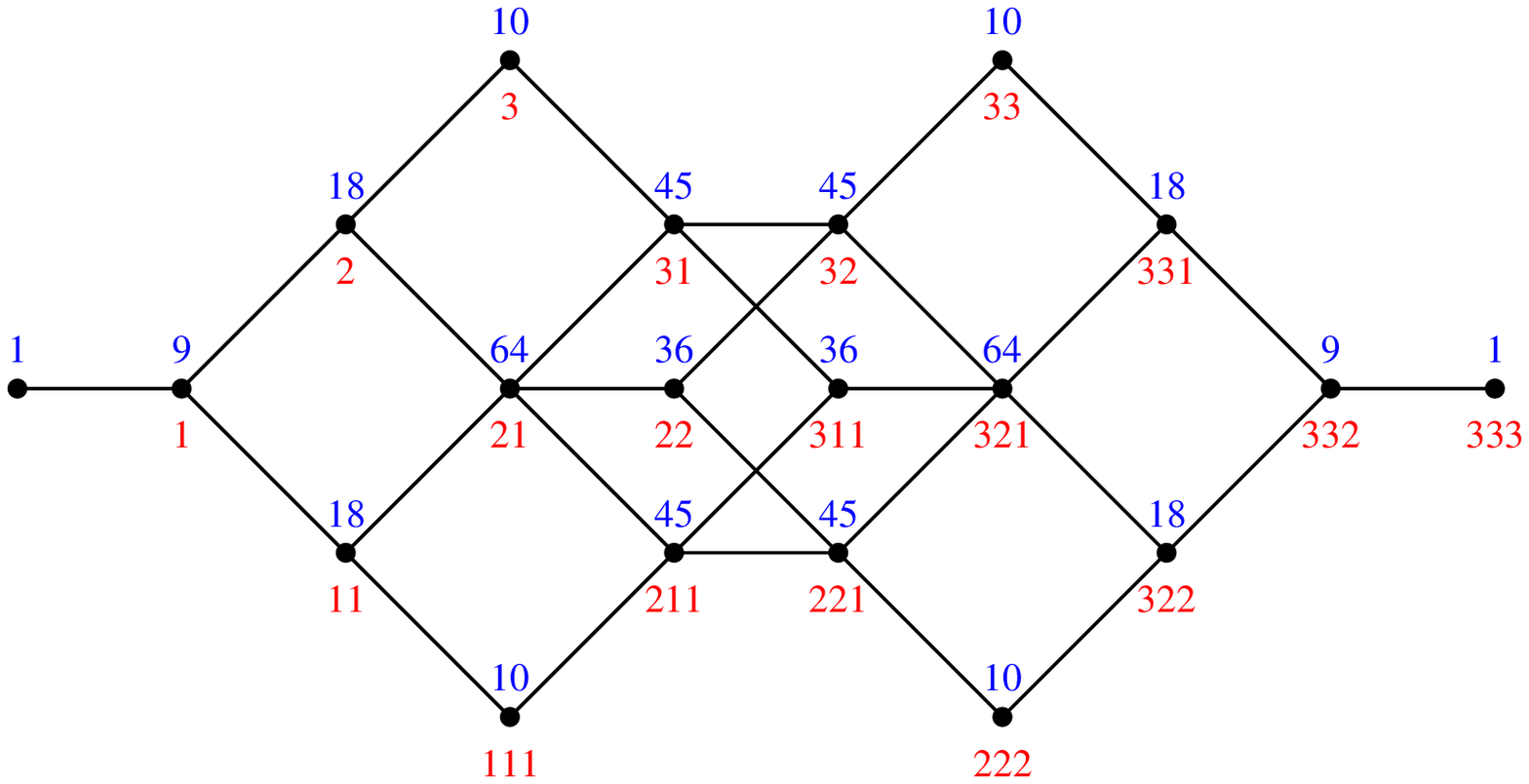}
\caption[The ideal poset for~$\Gr(3,6)$]
{\label{fig:Gr36poset} The ideal poset for~$\Gr(3,6)$.
The lower label on each node is the~$\ab\in \sfP(3,6)$ associated to the 
node and the upper label is the dimension of the corresponding
subspace of $\L^{*,0}(\eum)$.}
\end{figure}

\begin{remark}[The Hasse diagram of~$\sfP(3,6)$]
\label{rem: Hasse for G36}
The Hasse diagram%
\footnote{For an explanation of how this diagram encodes the
structure of the poset, see~\S\ref{ssec:hss-ideal-poset}.}
for the ideal poset of the 
Grassmannian~$\Gr(3,6)$ is drawn in Figure~\ref{fig:Gr36poset}.
Each node is labeled below according to its 
partition~$\ab$ and labeled above according to the rank of
$\bbS_{\ab'}(Q)\ot\bbS_{\ab}(S^*)$.  Inspection of this figure 
may clarify some of the relationships discussed in this section.
I certainly found it helpful.
\end{remark}

\section[Extremal Cycles in Grassmannians]{Extremal Cycles in Grassmannians}
\label{sec:extremals-in-grass}

\subsection{Ideals of degree~$2$}\label{ssec:2-cycles}
In this section, I am going to analyze the integral varieties of
$\cI_{(2)}$ and~$\cI_{(1,1)}$.  The main application will be
to giving a complete description of the the effective cycles of 
dimension~$2$ or more whose homology classes are either of the 
form~$r[\sigma_{(1,\dots,1)^*}]$
(see Theorem~\ref{thm: homologies to sigma(1,...,1)*}) 
or~$r[\sigma_{(p)^*}]$ 
(see Theorem~\ref{thm: homologies to sigma(p)*}). 
Essentially, these results state that such cycles are represented 
only by subvarieties of projective spaces `in disguise'.

\subsubsection{Integrals of~$\cI_{(2)}$}
\label{sssec:I2integrals}
First, a characterization of the irreducible integral 
varieties of~$\cI_{(2)}$.

\begin{proposition}\label{prop: integral varieties of I2}
For every~$A\in\Gr(m{+}1,n)$, 
the submanifold~$\Gr(m,A)\subset\Gr(m,n)$ 
is an integral manifold of~$\cI_{(2)}$.

Conversely if~$X\subset\Gr(m,n)$ is an integral variety of~$\cI_{(2)}$
that is irreducible and of dimension~$d\ge 2$,   
then $X\subset\Gr(m,A)$ for some unique~$A\in\Gr(m{+}1,n)$.
\end{proposition}

\begin{proof}
It is immediate that, for any~$A\in\Gr(m{+}1,n)$, 
the submanifold~$\Gr(m,A)$ is an integral manifold of~$\cI_{(2)}$.

Suppose that~$X\subset\Gr(m,n)$ satisfies the stated hypotheses
and let~$X^\circ\subset X$ denote the smooth part of~$X$, which 
is connected since~$X$ is irreducible~\cite[p.~21]{MR80b:14001}. 
By hypothesis, ~$X^\circ$ is an integral manifold of~$\cI_{(2)}$, 
i.e., its tangent planes are integral elements of~$\cI_{(2)}$.

Since~$d\ge 2$,  Lemma~\ref{lem: integral elements of I2} implies that
for every~$V\in X^\circ$ there is a $d$-plane~$P_V\subset V^*$
and a line~$L_V\subset Q_V = \C{n}/V$ so that~
\begin{equation}\label{eq:11-tangents}
T_VX = L_V\ot P_V\subset Q_V\ot V^* = T_V\Gr(m,n).
\end{equation}

Now consider the set~$F\subset X^\circ\times \SL(n,\bbC)$ 
consisting of the
set of pairs~$(V,\vs)$ so that~$V\in X^\circ$ and 
\begin{enumerate}
\item $\vs_1,\dots,\vs_m$ spans~$V$;
\item $\vs_1,\dots,\vs_{m-d}$ spans the annihilator of~$P_V$; and
\item $\dl\vs_{m+1}\dr$ spans~$L_V\subset Q_V$.
\end{enumerate}
Then~$F$ is a holomorphic $G$-bundle over~$X^\circ$ 
where~$G\subset\SL(n,\bbC)$ is the parabolic subgroup that stabilizes
$\C{m-d}$, $\C{m}$, and $\C{m+1}$, i.e.,
$$
G = P_{m-d}\cap P_m\cap P_{m+1}.
$$
 From now on, all computations will
take place on~$F$ or subbundles of~$F$.  Since~~$X^\circ$ and $G$ are
connected, $F$ is connected as well.

Consider the structure equations
\begin{equation}
\d \vs_A = \vs_B\,\omega^B_A\,,\qquad\qquad
\d \omega^A_B = -\omega^A_C\w\omega^C_B\,.
\end{equation} 
By \eqref{eq:11-tangents} and the definition of~$F$, 
it follows that~$\omega^\alpha_i = 0$ for all 
pairs~$(i,\alpha)$ satisfying 
either~$1\le i\le m{-}d$ and $\alpha = m{+}1$
or~$1\le i \le m$ and~$m{+}1 < \alpha \le n$.  On the other hand,
$\omega^{m+1}_{m-d+1}\w\dots\w\omega^{m+1}_m\not=0$.  

Consequently, taking $\alpha>m{+}1$ and $j$ satisfying ~$m{-}d<j\le m$
and computing exterior derivatives 
via the structure equations yields
\begin{equation}
0 = \d \omega^{\alpha}_{j} 
 = -\omega^{\alpha}_A\w\omega^A_{j}
= -\omega^{\alpha}_{m+1}\w\omega^{m+1}_{j}\,.
\end{equation} 
Since~$\omega^{m+1}_{m-d+1}\w\dots\w\omega^{m+1}_m\not=0$ 
and~$d\ge 2$, this implies~$\omega^\alpha_{m+1}=0$ when~$\alpha>m{+}1$.  

Since $\omega^\alpha_i=0$ for all pairs~$(i,\alpha)$ 
satisfying~$1\le i\le m{+}1< \alpha\le n$, it follows
that
\begin{equation}
\d \vs_1 \equiv \dots \equiv \d \vs_{m+1}
  \equiv 0\mod \vs_1,\dots,\vs_{m+1}\,.
\end{equation}
Thus, the span of the $\C{n}$-valued functions~$\vs_1,\dots,\vs_{m+1}$
is locally constant on~$F$.  Since~$F$ is connected, this span is 
constant.  Let~$A\subset\C{n}$ be this span.  By construction~$A$
contains~$V$ for all~$V\in X^\circ$.  Thus,~$X^\circ$ lies in~$\Gr(m,A)$.  
Since~$X^\circ$ is dense in~$X$, it follows
that~$X$ itself lies in~$\Gr(m,A)\simeq\bbP(A^*)$, as claimed.
\end{proof}

Proposition~\ref{prop: integral varieties of I2} 
has some interesting consequences.

\begin{theorem}\label{thm: homologies to sigma(1,...,1)*}
Suppose that~$1<p\le m$ and let~$\ab = (1,\dots,1)$ with~$|\ab| = p$.
Let~$X\subset\Gr(m,n)$ be an irreducible $p$-dimensional variety 
that satisfies $[X]=r\,[\sigma_{\ab^*}]$ for some~$r\in\bbZ^+$.
Then there exists a unique $A\in\Gr(m{+}1,n)$ so that
$$
X\subset\Gr(m,A)\quad\left(\simeq\bbP(A^*)\simeq\bbP^m\right)
$$
and~$r$ is the degree of~$X$ as a variety in~$\bbP^m$.

Conversely, for~$A\in\Gr(m{+}1,n)$, 
any subvariety~$X\subset \Gr(m,A)\subset\Gr(m,n)$
of dimension~$p$ and degree~$r$ satisfies~$[X] = r\,[\sigma_{\ab^*}]$.
\end{theorem}

\begin{proof}
Since~$[X] = r\,[\sigma_{\ab^*}]$, 
it follows that~$\phi_\bb$ vanishes on~$X^\circ$, the smooth part of~$X$, 
for all~$\bb\in\sfP(m,n)$ with~$|\bb| = p$ and~$\bb\not=\ab$.  
Consider the positive~$(p,p)$-form~$(\phi_1)^{p-2}\w\phi_2$.  
(This is where the hypothesis~$p>1$ is used.)
By Pieri's formula~\eqref{eq: Pieri p for forms}
$$
(\phi_1)^{p-2}\w\phi_{(2)} 
= \sum_{\substack{\bb\in\sfP(m,n) \\ |\bb|= p}} 
      \mu^\bb_{(2)}\,\phi_\bb\,.
$$
Since~$\ab \not\ge (2)$, it follows that~$\mu^\ab_{(2)}=0$, so 
every term on the right hand side of the above equation vanishes
on~$X^\circ$.  Thus $(\phi_1)^{p-2}\w\phi_2$ vanishes
on~$X^\circ$ as well.  Since~$\phi_1$ defines a K\"ahler
form on~$X^\circ$, the generalized Wirtinger 
inequality~\eqref{eq:Generalized Wirtinger Inequality} implies that~$\phi_2$
must vanish on~$X^\circ$.  In particular,~$X$ is an integral manifold
of~$\cI_{(2)}$.  Since~$X$ is irreducible and of dimension~$p>1$, 
Proposition~\ref{prop: integral varieties of I2}
applies.  The statements about degree now follow immediately.
\end{proof}

This result has an interesting consequence of its own:

\begin{theorem}\label{thm: homologies to sigma(m)'*}
Suppose that~$m>1$ and let~$\ab = (1,\dots,1)$, where~$|\ab| = m$.
Let~$X\subset\Gr(m,n)$ be an irreducible $m$-dimensional variety 
that satisfies $[X]=r\,[\sigma_{\ab^*}]$ for some~$r\in\bbZ^+$.
Then $r=1$ and there exists a unique $A\in\Gr(m{+}1,n)$ so that
$X=\Gr(m,A)$.  
{\upshape(}In particular,~$X$ is a Schubert 
variety~$\sigma_{\ab^*}$.{\upshape)}
\end{theorem}

\begin{remark}[Walters' results]
\label{rem: Walters results 1}
When~$m=2$, Theorem~\ref{thm: homologies to sigma(m)'*} was proved
by Walters~\cite[Theorem~5 and Corollary~3]{mW97}.  
Her proof relies on a local
computation in coordinates that, essentially, computes the $2$-dimensional 
integral elements of~$\cI_{(2)}$ in the case~$m=2$ 
(see \cite[Claim 4.2.3]{mW97}) and then, using this, proves a coordinate 
version of~Proposition~\ref{prop: integral varieties of I2} in this
case.
\end{remark}

\begin{remark}[Schur rigidity and quasi-rigidity]
\label{rem: schur rigidity quasi-rigidity}
As was mentioned already in \S\ref{sssec: rigidity questions}, 
Walters showed that~$\ab = ((m)')^*$, has Schubert rigidity, i.e.,
that any (local) solution of~$\mathcal{B}_\ab$ is a Schubert variety
$\sigma_{\ab}$.  Theorem~\ref{thm: homologies to sigma(m)'*} shows
that this~$\ab$ even has Schur rigidity.  
This is not surprising, though, because, as the proof of 
Theorem~\ref{thm: homologies to sigma(1,...,1)*} makes clear, 
$\mathcal{B}_{\ab} = \mathcal{R}_{\ab^*}$ for all~$\ab$ of the
form~$((p)')^*$ when~$p>1$.

Walters also showed that, when~$p<m$, the type~$\ab = ((p)')^*$ 
does not have Schubert rigidity.  However, she did not classify 
the solutions to~$\mathcal{B}_{\ab}$ in this case.  
Theorem~\ref{thm: homologies to sigma(1,...,1)*} 
does this classification in the range~$1<p<m$, showing that 
such an irreducible variety is a $p$-dimensional subvariety 
of a projective space~$\Gr(m,A) \simeq \bbP^m$ for some~$A\in\Gr(m{+}1,n)$.

Since~$\mathcal{B}_{(1)^*}\not=\mathcal{R}_{(1)}$, 
the solutions of~$\mathcal{B}_{(1)^*}$ (the omitted case) are 
of a different nature.  However, it is not difficult to show that 
an irreducible curve in~$\Gr(m,n)$ that is a solution 
of~$\mathcal{B}_{(1)^*}$ can be described as follows:  
Fix a subspace~$W\subset \C{n}$ of dimension~$k<m$ 
and let~$C\subset \bbP\bigl(\C{n}/W\bigr)\simeq\bbP^{n-k-1}$ 
be a projective curve that spans a projective subspace of dimension 
at least~$m{-}k$.  Let~$C^{[m-k-1]}\subset\Gr\bigl(m{-}k,\C{n}/W\bigr)$
be the $(m{-}k{-}1)$-th osculating curve of~$C$ and consider its
inclusion into~$\Gr(m,n)$ by the canonical injection
$\Gr\bigl(m{-}k,\C{n}/W\bigr)\hookrightarrow \Gr(m,n).$
This image curve is a solution to $\mathcal{B}_{(1)^*}$ and every
irreducible solution is of this form.
\end{remark}

\subsubsection{Bundles with $c_2 = 0$}
\label{sssec: c2=0}

Proposition~\ref{prop: integral varieties of I2} can also be applied
to characterize bundles generated by their sections but with vanishing
second Chern class.

\begin{theorem}\label{thm: bundles with c2 = 0}
Let~$M$ be a connected compact K\"ahler manifold and let~$F\to M$ be 
a holomorphic vector bundle that is generated by its sections.  
If~$c_2(F) = 0$, then either there exists a holomorphic 
splitting~$F = L\oplus T$ where~$L$
is a line bundle and~$T$ is trivial or else there exists an algebraic
curve~$C$, a bundle~$F'\to C$ that is generated by its sections,
and a holomorphic mapping~$\kappa:M\to C$ so that~$F = \kappa^*(F')$.
\end{theorem}

\begin{proof}
Since~$F$ is generated by its sections, there exists an~$n>0$ 
and a surjective holomorphic bundle map~$\phi: M\times\C{n}\to F$.  
Let~$m\le n$ be the rank of the kernel bundle~$K\subset M\times\C{n}$.  
The mapping~$\kappa:M\to\Gr(m,n)$ defined by~$\kappa(x)=K_x$ then
has the property that~$F=\kappa^*(Q)$.  Moreover, 
$c_2(F)=\kappa^*\bigl(c_2(Q)\bigr) = \kappa^*(q_2)$.  

Since~$[\phi_2]=c_2(Q) = q_2$ and~$\phi_2$ is a positive~$(2,2)$-form,
it follows that~$\kappa^*(\phi_2)$ is a positive~$(2,2)$-form on~$M$
that represents~$c_2(F)$.  Since~$M$ is compact and K\"ahler, 
the hypothesis~$c_2(F)=0$ implies that the representing positive 
form~$\kappa^*(\phi_2)$ must also be zero 
by Corollary~\ref{cor: F generated by sections and ca = 0}.  
Equivalently, $\kappa(M)\subset\Gr(m,n)$ is an integral variety 
of~$\cI_{(2)}$.   Since~$M$ is connected, $\kappa(M)$ is irreducible.
 Now there are three cases: 

If the dimension of~$\kappa(M)$ is equal to~$0$, then $\kappa$ is 
constant and $F = \kappa^*(Q)$ is trivial.  This falls into both
of the two cases allowed by the proposition.

If the dimension of~$\kappa(M)$ is equal to~$1$,
then $\kappa(M)$ is an algebraic curve~$C\subset\Gr(m,n)$.  
Replace~$C$ by its normalization if necessary and 
define~$F'\to C$ to be the pullback to~$C$ of the bundle~$Q\to\Gr(m,n)$.
Then~$\kappa^*(F') = F$.

If the dimension of~$\kappa(M)$ is greater than~$1$, 
then Proposition~\ref{prop: integral varieties of I2} 
implies that there exists an~$(m{+}1)$-plane~$A\subset\C{n}$ so 
that~$\kappa(M)\subset\Gr(m,A)$.  Let~$B\subset\C{n}$ be a subspace
of dimension~$n{-}m{-}1$ that is a complement to~$A$.  The bundle~$Q$
restricted to~$\Gr(m,A)$ splits as a sum~$L'\oplus T'$ where~$L'$ is
a line bundle and~$T'$ is trivial.  Explicitly, for~$V\in\Gr(m,A)$,
$L'_V = A/V$ and~$T'_V = B$.  
Setting~$L=\kappa^*(L')$ and~$T=\kappa^*(T')$
yields the desired splitting.
\end{proof}

\subsubsection{Integrals of~$\cI_{(1,1)}$}
\label{sssec:I11integrals}
Now, all of these results from the analysis of~$\cI_{(2)}$ can
be translated by complementarity into corresponding results from
the analysis of~$\cI_{(1,1)}$:

\begin{proposition}\label{prop: integral varieties of I11}
Suppose that~$n{-}m\ge2$ and
let~$V\subset\Gr(m,n)$ be an irreducible subvariety 
of dimension~$d>1$ that is an integral manifold of~$\cI_{(1,1)}$. 
Then there is a unique $A\in\Gr(m{-}1,n)$ so that
$$
V\subset [A,\C{n}]_m \quad\left(\ =\Gr(1,\C{n}/A)\simeq\bbP^{n-m}\right).
$$
Conversely, for every~$A\in\Gr(m{-}1,n)$, the submanifold 
$[A,\C{n}]_m\subset\Gr(m,n)$ is an integral manifold of~$\cI_{(1,1)}$.
\end{proposition}

\begin{proof}
Apply complementarity to the proof of 
Proposition~\ref{prop: integral varieties of I2}.
\end{proof}

This proposition can be applied to characterize the cycles representing
a number of extremal classes.

\begin{theorem}\label{thm: homologies to sigma(p)*}
Suppose that~$1<p\le n{-}m$ and let~$V\subset\Gr(m,n)$ 
be an irreducible $p$-dimensional variety 
that satisfies $[V]=r\,[\sigma_{(p)^*}]$ for some~$r\in\bbZ^+$.
Then there exists a unique $A\in\Gr(m{-}1,n)$ so that
$$
V\subset[A,\C{n}]_m\quad\left(\simeq\bbP(\C{n}/A)\simeq\bbP^{n-m}\right)
$$
and~$r$ is the degree of~$V$ as a variety in~$\bbP^{n-m}$.
Conversely, for each~$A\in\Gr(m{-}1,n)$, any $p$-dimensional
variety~$V\subset [A,\C{n}]_m$ of degree~$r$
satisfies~$[V] = r\,[\sigma_{(p)^*}]$.
\end{theorem}

\begin{proof}
Apply complementarity to the proof of
Theorem~\ref{thm: homologies to sigma(1,...,1)*}.
\end{proof}

\begin{remark}[Walters' results]
\label{rem: Walters results 2}
These results when~$p=2$ and~$m=2$ are also (partially) to be found
in the work of Walters.  
\end{remark}

\begin{remark}[Another Schur-Schubert coincidence]
\label{rem: B p* = R p}
Note that Theorem~\ref{thm: homologies to sigma(p)*} also implies that
$\mathcal{B}_\ab = \mathcal{R}_{\ab^*}$ when~$\ab = (p)^*$ for~$p>1$.
However, this could have been proved directly, using the algebraic ideas 
that went into the proof.
\end{remark}

\subsubsection{Bundles with ${c_{1}}^2-c_2 = 0$}
\label{sssec: c11=0}

Proposition~\ref{prop: integral varieties of I11} can also be applied
to characterize bundles generated by their sections but with vanishing
Schur-Chern polynomial~$c_{(1,1)}$.

\begin{theorem}\label{thm: bundles with c2 = c1*c1}
Let~$M$ be a connected compact K\"ahler manifold and let~$F\to M$ be 
a holomorphic vector bundle of rank~$r$ that is generated by 
its sections.  

Then~$c_1(F)^2-c_2(F)\ge0$ and, 
if equality holds, either $F = (M{\times}\C{r+1})/L$ for some line 
bundle~$L\subset M{\times}\C{r+1}$ 
or else there exists an algebraic
curve~$C$, a bundle~$F'\to C$ that is generated by its sections,
and a holomorphic mapping~$\kappa:M\to C$ so that~$F=\kappa^*(F')$.
\end{theorem}

\begin{proof}
Let~$H^0(F)$ be the space of global sections of~$F$, a vector
space of dimension~$n = h^0(F)$.  Let~$\ev_{\!F}: M\times H^0(F)\to F$ 
be the evaluation mapping, which, by assumption,
is surjective, so that~$n\ge r$.  The kernel~$K\subset M\times H^0(F)$
is then a subbundle of rank~$m=n{-}r$ and can be used to define a 
mapping~$\kappa_F:M\to\Gr\bigl(m,H^0(F)\bigr)$ that satisfies~
$\kappa_F^*(Q) = F$.  Consequently,
$$
c_1(F)^2-c_2(F) = c_{(1,1)}(F) = \left[\kappa^*(\phi_{(1,1)})\right].
$$
Thus, the inequality~$c_1(F)^2-c_2(F)\ge 0$ follows directly from
Corollary~\ref{cor: F generated by sections and ca = 0}.  Moreover,
if equality holds, $\kappa_F(M)\subset\Gr\bigl(m,H^0(F)\bigr)$ must
be an integral variety of~$\cI_{(1,1)}$.  
As in the proof of Theorem~\ref{thm: bundles with c2 = 0}, 
there are now three cases:

If $\kappa_F(M)$ is a single point, then~$F$ is trivial.  

If $\kappa_F(M)$ is a curve, 
let~$C\to\kappa_F(M)\subset\Gr\bigl(m,H^0(F)\bigr)$ 
be its normalization and let~$F'$ be the pullback of~$Q$ to~$C$.

If $\kappa_F(M)$ has dimension greater than~$1$, then 
Proposition~\ref{prop: integral varieties of I11} implies that there
is an $(m{-}1)$-plane~$A\subset H^0(F)$ 
so that~$\kappa_F(M)\subset\bigl[A,H^0(F)\bigr]_m$.  However, this 
implies that~$M\times A$ is a subset of~$K$.  In other words,~$A$
consists of the global sections of~$F$ that vanish at all points of~$M$.
Of course, this implies that~$A = (0)$, i.e., that~$m=1$, so that~$
H^0(F) = n = r{+}m = r{+}1$, which is what needed to be shown.
\end{proof}

\subsection{Ideals of codegree~$2$}\label{ssec:2-cocycles}
In this subsection, 
I will analyze the maximal dimension integral varieties 
of the ideals~$\cI_{(1,1)^*}$ and~$\cI_{(2)^*}$, namely their integral
varieties of codimension~$2$ in~$\Gr(m,n)$.  This is not really of 
interest unless both ideals are nontrivial, so I will assume that
$2\le m\le n{-}2$.  In this case, $[\sigma_{(2)}]$ and~$[\sigma_{(1,1)}]$ 
are linearly independent and give a basis 
of~$H_{m(n-m)-4}\bigl(\Gr(m,n),\bbZ\bigr)$.
The goal of this subsection is 
to give a description of the effective cycles 
representing these classes, so I state the main results in those terms,
rather than directly in terms of the integral varieties of the ideals.

One reason for interest in these results is that the comparison
of the nonsmoothability of these varieties with  
the nonsmoothability results 
of Hartshorne, Rees, and Thomas~\cite{MR50:9870}.

\subsubsection{Integrals of~$\cI_{(1,1)^*}$}
\label{sssec:I11starintegrals}
Consider the class~$[\sigma_{(2)}]$.
Recall that~$\sigma_{(2)}\subset\Gr(m,n)$ consists of the $m$-planes
$E\subset\C{n}$ that meet~$\C{n-m-1}$ in at least a line. Equivalently,
this is the same as requiring that~$\bbP E\subset\bbP^{n-1}$ 
meet a fixed~$\bbP^{n-m-2}\subset\bbP^{n-1}$.  

\begin{theorem}\label{thm:codim2-sigma2}
For any algebraic variety $A\subset\Gr(1,n)$ 
of codimension~$m{+}1$ and degree~$r$, the
variety
$$
\Psi_m(A) = \bigl\{\,E\in\Gr(m,n)\mid \bbP E \cap A \not=\emptyset\,\}
$$ 
is of codimension~$2$ in~$\Gr(m,n)$ and 
satisfies~$\bigl[\Psi_m(A)\bigr] = r\,[\sigma_{(2)}]$.

Conversely, if~$V\subset\Gr(m,n)$ is an algebraic variety 
of codimension~$2$ that satisfies $[V] = r\,[\sigma_{(2)}]$ 
for some~$r\in\bbZ^+$, then $V = \Psi_m(A)$ for some unique
algebraic variety $A\subset\Gr(1,n)$ 
of codimension~$m{+}1$ and degree~$r$.
\end{theorem}

\begin{proof}
Let~$A\subset\Gr(1,n)\simeq\bbP^{n-1}$ be an algebraic variety of
codimension~$m{+}1$ and degree~$r$.  
It is immediate that~$\Psi_m(A)$ is of codimension~$2$.  Moreover,
a simple local calculation shows that, at its smooth points, its tangent
spaces are integral elements of~$\cI_{(1,1)^*}$ 
(cf. Corollary~\ref{cor:Z(phi(1,1)*) computation}).  Thus, $\phi_{(1,1)^*}$
vanishes on~$\Psi_m(A)$, which implies that  
$\bigl[\Psi_m(A)\bigr] = r'\,[\sigma_{(2)}]$ for some~$r'$.  Since~$r'=1$
when~$r=1$, it follows easily that~$r'=r$ in all cases.

Now, for the converse statement, it clearly suffices to prove the
characterization when~$V$ is irreducible, so assume this.  Thus
suppose that~$V\subset\Gr(m,n)$ satisfies the stated hypotheses
and is irreducible. Let~$V^\circ\subset V$ denote the smooth part of~$V$, 
which is connected, since~$V$ is irreducible.  
Since $[V] = r\,[\sigma_{(2)}]$, this smooth part~$V^\circ$
must be an integral manifold of~$\cI_{(1,1)^*}$, i.e., 
its tangent planes must be integral elements of~$\cI_{(1,1)^*}$.
Thus, by Lemma~\ref{lem: integral elements of I11*}, 
for every~$E\in V^\circ$ there exists a line~$L_E\subset E$
and a codimension $2$ subspace~$R_E\subset Q_E = \C{n}/E$ so that~
\begin{equation}\label{eq:co-2tangents}
T_EV = \bigl(R_E{\ot}E^*\bigr) + \bigl(Q_E{\ot}L_E^\perp\bigr)
   \subset Q_E\ot E^* = T_E\Gr(m,n).
\end{equation}
Consider the set~$F\subset V^\circ\times \GL(n,\bbC)$ consisting of the
set of pairs~$(E,\vs)$ so that 
\begin{enumerate}
\item $\vs_1$ spans~$L_E$;
\item $\vs_1,\dots,\vs_m$ spans~$E$; and
\item $\dl\vs_{m+1}\dr_E,\dots,\dl\vs_{n-2}\dr_E$ spans~$R_E\subset Q_E$.
\end{enumerate}
Then~$F$ is a holomorphic $G$-bundle over~$V^\circ$ 
where~$G\subset\GL(n,\bbC)$ is the parabolic subgroup that stabilizes
$\C{1}$, $\C{m}$, and $\C{n-2}$.  From now on, all computations will
take place on~$F$ or subbundles of~$F$.  Of course, since~$G$ is
connected and since~$V^\circ$ is connected, it follows that~$F$ is
connected as well.

Consider the structure equations
\begin{equation}
\begin{split}
\d \vs_A &= \vs_B\,\omega^B_A\,,\\
\d \omega^A_B &= -\omega^A_C\w\omega^C_B\,.\\
\end{split}
\end{equation} 
By \eqref{eq:co-2tangents} and the definition of~$F$, it follows 
that~$\omega^{n-1}_1 =\omega^{n}_1 = 0$ and, moreover, that these two
relations are the only linear relations among the~$m(n{-}m)$ 
forms~$\omega^a_i$ with~$1\le i\le m< a\le n$.  Note that these
latter forms generate the module of forms that are semibasic
for the fibration~$\pi:F\to V^\circ$.

Computing exterior derivatives via 
the structure equations yields
\begin{equation}
\begin{split}
0 = \d \omega^{n-1}_1 
   &= { }-\omega^{n-1}_A\w\omega^A_1
    = -\sum_{j=2}^m \omega^{n-1}_j\w\omega^j_1
      -\sum_{a = m+1}^{n-2} \omega^{n-1}_a\w\omega^a_1\,,\\
0 = \d \omega^{n\phantom{-1}}_1 
   &= { } -\omega^{n\phantom{-1}}_A\w\omega^A_1
    = -\sum_{j=2}^m \omega^{n\phantom{-1}}_j\w\omega^j_1
      -\sum_{a = m+1}^{n} \omega^{n\phantom{-1}}_a\w\omega^a_1\,.\\
\end{split}
\end{equation} 
Reducing these equations modulo $\{\omega^{m+1}_1,\ldots,\omega^{n-2}_1\}$
yields
\begin{equation}
\sum_{j=2}^m \omega^{n-1}_j\w\omega^j_1
\equiv \sum_{j=2}^m \omega^{n}_j\w\omega^j_1
\equiv 0 \mod \omega^{m+1}_1,\ldots,\omega^{n-2}_1\,.
\end{equation}
Since~$\{\omega^{n-1}_j,\omega^n_j\mid 2\le j\le m\}$
are linearly independent modulo~$\{\omega^{m+1}_1,\ldots,\omega^{n-2}_1\}$,
these equations imply that
\begin{equation}
\omega^2_1\equiv \dots \equiv \omega^k_1\equiv 0 \mod
\omega^{m+1}_1,\ldots,\omega^{n-2}_1\,.
\end{equation}
Thus, there exist functions~$\{\,s^j_a\mid 2\le j\le m<a\le n{-}2\,\}$
on~$F$ so that~$\omega^j_1 = s^j_a\,\omega^a_1$.  The structure equations
now imply that
\begin{equation}
d\vs_1 \equiv \sum_{a=m+1}^{n-2} (\vs_a + s^j_a\,\vs_j)\,\omega^a_1
\mod\vs_1\,.
\end{equation}
Consequently, the map~$[\vs_1]:F\to \Gr(1,\C{n}) = \bbP^{n-1}$
is a holomorphic map of constant rank~$n{-}m{-}2$.  Moreover, since
the forms~$\{\omega^{m+1}_1,\ldots,\omega^{n-2}_1\}$ are $\pi$-semibasic,
and~$\pi$ has connected fibers, it follows that there is a 
well-defined holomorphic map~$\lambda:V^\circ\to\bbP^{n-1}$
of constant rank~$n{-}m{-}2$ that satisfies~$\lambda(E) = L_E$, 
i.e., $[\vs_1] = \lambda\circ\pi$.

By dimension count, the fibers of~$\lambda$ have dimension~$m(n{-}m){-}2
-(n{-}m{-}2) = (m{-}1)(n{-}m)$.   Moreover, by construction, 
for each line~$L = \lambda(E)$, the fiber~$\lambda^{-1}(L)$ is
embedded as a submanifold of the sub-Grassmannian~$[L,\C{n}]_m
\subset\Gr(m,n)$, which also has dimension~$(m{-}1)(n{-}m)$.  It follows
that $V^\circ\cap[L,\C{n}]_m$ is an open subset
of~$[L,\C{n}]_m$.  Since~$V$ is algebraic, 
it follows that~$V$ must actually 
contain~$[L,\C{n}]_m$ for all~$L\in\lambda(V^\circ)$.

Let~$B\subset~\Gr(1,\C{n})$ be the set of lines~$L$
for which~$[L,\C{n}]_m$ lies in~$V$.  Then~$B$
is evidently a variety that has at least one component~$B'$ of 
dimension~$d\ge n{-}m{-}2$.  Since, in particular,~$V$ must contain all of
the~$m$-planes that meet~$B'$, it follows that the dimension of~$B'$
cannot be more than~$n{-}m{-}2$ (otherwise, it would impose at most one 
condition for a $\bbP^{m-1}\subset\bbP^{n-1}$ to meet~$B'$).  
By connectedness, the image~$\lambda(V^\circ)$ must lie in a component
of~$B$, say~$A\subset B$ that has the maximum possible dimension, 
namely~$n{-}m{-}2$.  

Let~$W\subset V$ denote the set of $m$-planes~$E$ 
satisfying~$\bbP(E)\cap A\not=\emptyset$.  Then, since~$A$ has 
dimension~$n{-}m{-}2$, it follows that~$W$ has the same dimension as~$V$.
Since~$V$ is irreducible,~$W=V$.  Now, if~$B$ had any other 
component~$A'\not=A$ of dimension~$n{-}m{-}2$, then the corresponding~$W'$
would also satisfy~$W'=V$ and, consequently, $W'=W$.  However, this would
imply that every $\bbP^{m-1}$ that meets~$A$ must also meet~$A'$ and
vice versa.  But if~$A\not=A'$, this is absurd.  Thus,~$A$ is unique.

Finally, the equation~$r = \deg(A)$ and the converse follow by the 
Schubert calculus and calculation, respectively.
\end{proof}

\begin{remark}[Integral varieties of~$\cI_{(1,1)^*}$]
\label{rem: integral varieties of I11star}
The reader may have expected Theorem~\ref{thm:codim2-sigma2}
to have been stated in terms of a characterization of the irreducible
integral varieties of~$\cI_{(1,1)^*}$ of codimension~$2$.  There is,
of course, such a characterization and it follows the more-or-less expected 
lines except for one caveat:  What is true (and the above proof
can be easily adapted so as to prove it) is that every 
codimension~$2$, irreducible integral variety~$V$ of~$\cI_{(1,1)^*}$ is
\emph{locally} of the form~$\Psi_m(A)$ for some irreducible 
subvariety~$A\subset \Gr(1,n)$ of dimension~$n{-}m{-}2$.  In this
sense, the codimension~$2$ integral manifolds of~$\cI_{(1,1)^*}$
depend (in Cartan's sense)
on $m{+}1$ functions of~$n{-}m{-}2$ variables (in the
holomorphic category).  However, without some hypotheses on the
`finiteness' of the variety~$V$, 
the mapping~$\lambda:V^\circ\to\Gr(1,n)$ can be very far from proper,
so that the image is not a variety, even locally.  It is for this reason
that I incorporated compactness into the statement of
Theorem~\ref{thm:codim2-sigma2}.
\end{remark}

\begin{example}[Rigidity of~$\sigma_{(2)}$]
\label{ex:sigma2-must-be-Schubert}
Theorem~\ref{thm:codim2-sigma2} shows that any~$V\subset\Gr(m,n)$
of codimension~$2$ that satisfies~$[V] = [\sigma_{(2)}]$ must actually
\emph{be} a Schubert cycle, a very strong form of rigidity.
Of course, when~$(m,n)=(2,4)$, this result is classical~\cite{gFano}
(and, in any case, already follows from the previous results on 
extremal subvarieties of dimension~$2$).
\end{example}

\begin{example}[Nonsmoothability]
\label{ex:sigma2-must-be-singular}
Note that the variety~$\Psi_m(A)\subset\Gr(m,n)$ will be
singular as soon as $A$ is not a single point.  Thus,
Theorem~\ref{thm:codim2-sigma2} implies that if~$n>m{+}2$, then
any~$V\subset\Gr(m,n)$ of codimension~$2$ 
satisfying~$[V] = r[\sigma_{(2)}]$
with~$r>0$ must necessarily be singular.  

In~\cite{MR50:9870}, the authors
use results of Thom~\cite{MR15:890a} to prove that, for any~$r>0$, 
there are integers~$m$ and~$n$ for which
the class~$r[\sigma_{(2)}]$ in~$H_*\bigl(\Gr(m,n),\bbZ\bigr)$ cannot be 
represented by a smooth manifold.  
While Thom's results also show that, for every~$m$ and~$n>m{+}2$,
there exists an~$r>0$ for which~$r[\sigma_{(2)}]$ 
is representable by a smooth manifold,
Theorem~\ref{thm:codim2-sigma2}
implies that, when~$n>m{+}2$, the class~$r[\sigma_{(2)}]$ is not 
representable by a smooth algebraic variety for \emph{any}~$r>0$.

If~$n=m{+}2$, the cycle~$\sigma_{(2)}=[\C{1},\C{m+2}]_m$ 
is isomorphic to~$\Gr(m{-}1,m{+}1)$ and hence is smooth.  
Theorem~\ref{thm:codim2-sigma2} 
implies that if~$V\subset\Gr(m,m{+}2)$ 
satisfies~$[V] = r[\sigma_{(2)}]$, 
then
$$
V = [L_1,\C{m+2}]_m \cup \dots \cup [L_r,\C{m+2}]_m
$$
for some lines~$L_1,\dots,L_r\in\Gr(1,\C{m+1}) = \bbP^{m+1}$.  
In particular, $V$ is singular when~$r>1$.
\end{example}

\subsubsection{Integrals of~$\cI_{(2)^*}$}
\label{sssec:I2starintegrals}
Recall that~$\sigma_{(1,1)}$ consists of the $m$-planes
$E\subset\C{n}$ that meet~$\C{n-m+1}$ in at least a $2$-plane.  
Now, $E$ satisfies this condition when~$E+\C{n-m+1}$
has dimension at most~$n{-}1$, i.e., when~$E$ lies in a hyperplane 
containing~$\C{n-m+1}$.  Thus, another way of describing~$\sigma_{(1,1)}$
is as the $m$-planes that lie in one of the hyperplanes in~$\C{n}$
that contain~$\C{n-m+1}$.  The set of such hyperplanes forms
a~$\bbP^{m-2}$ in~$\Gr(n{-}1,n)\simeq\bbP^{n-1}$.  
This description of~$\sigma_{(1,1)}$ motivates the following result.

\begin{theorem}\label{thm:codim2-sigma11}
For any algebraic variety $A\subset\Gr(n{-}1,n)$ of dimension~$m{-}2$ 
and degree~$r$, the subvariety~$V = \Sigma_m(A)$ defined by
$$
\Sigma_m(A) = \cup_{H\in A}\Gr(m,H) 
$$
has codimension~$2$ in~$\Gr(m,n)$ 
and satisfies~$[V] = r\,[\sigma_{(1,1)}]$.

Conversely, if~$V\subset\Gr(m,n)$ is 
an algebraic variety of codimension~$2$
that satisfies $[V] = r\,[\sigma_{(1,1)}]$ for some~$r\in\bbZ^+$,
then there exists a subvariety $A\subset\Gr(n{-}1,n)$ 
of dimension~$m{-}2$ and degree~$r$ so that~$V =\Sigma_m(A)$.
\end{theorem}

\begin{proof}
Let~$A\subset\Gr(n{-}1,n)\simeq\bbP^{n-1}$ be an algebraic variety
of dimension~$m{-}2$ and degree~$r$.  A simple local calculation verifies
that~$\Sigma_m(A)$ as defined in the proposition has codimension~$2$
in~$\Gr(m,n)$ and that its tangent space at a smooth point is 
of type~$(1,1)$.  Consequently 
(cf. Corollary~\ref{cor:Z(phi(2)*) computation}), it follows that
$\bigl[\Sigma_m(A)\bigr] = r'\bigl[\sigma_{(1,1)}\bigr]$ for some
$r'>0$.  Since~$r=1$ if and only if~$A$ is a linear~$\bbP^{m-2}$ and 
since, in this case,~$\Sigma_m(A)$ is a Schubert cycle of type~$(1,1)$,
it follows that~$r'=1$ in this case.  It now follows easily by a degree
argument that~$r'=r$ in all cases.

Now, to prove, the converse, it clearly suffices to treat the case 
in which~$V\subset\Gr(m,n)$ is irreducible, so assume this. Thus,
suppose that~$V\subset\Gr(m,n)$ satisfies the stated hypotheses
and let~$V^\circ\subset V$ denote the smooth part of~$V$, which 
is connected, since~$V$ is irreducible.  
Since $[V] = r\,[\sigma_{(1,1)}]$, this smooth part~$V^\circ$
must be an integral manifold of~$\cI_{(2)^*}$, i.e., its tangent
planes must be integral elements of~$\cI_{(2)^*}$.
Thus, by Lemma~\ref{lem: integral elements of I2*}, 
for every~$E\in V^\circ$ there exists a 2-plane~$P_E\subset E$
and a hyperplane~$R_E\subset Q_E = \C{n}/E$ so that~
\begin{equation}\label{eq:codim2-tangents-11}
T_EV = \bigl(R_E{\ot}E^*\bigr)+ \bigl(Q_E{\ot}P_E^\perp\bigr)
   \subset Q_E\ot E^* = T_E\Gr(m,n).
\end{equation}
Consider the set~$F\subset V^\circ\times \GL(n,\bbC)$ consisting of the
set of pairs~$(E,\vs)$ so that 
\begin{enumerate}
\item $\vs_1,\vs_2$ spans~$P_E$;
\item $\vs_1,\dots,\vs_m$ spans~$E$; and
\item $\dl\vs_{m+1}\dr_E,\dots,\dl\vs_{n-1}\dr_E$ spans~$R_E\subset Q_E$.
\end{enumerate}
Then~$F$ is a holomorphic $G$-bundle over~$V^\circ$ 
where~$G\subset\GL(n,\bbC)$ is the parabolic subgroup that stabilizes
$\C{2}$, $\C{m}$, and $\C{n-1}$.  From now on, all computations will
take place on~$F$ or subbundles of~$F$.  Of course, since~$G$ is
connected and since~$V^\circ$ is connected, it follows that~$F$ is
connected as well.

Consider the structure equations
\begin{equation}
\begin{split}
\d \vs_A &= \vs_B\,\omega^B_A\,,\\
\d \omega^A_B &= -\omega^A_C\w\omega^C_B\,.\\
\end{split}
\end{equation} 
By \eqref{eq:codim2-tangents-11} and the definition of~$F$, it follows 
that~$\omega^{n}_1 =\omega^{n}_2 = 0$ and, moreover, that these two
relations are the only linear relations among the~$m(n{-}m)$ 
forms~$\omega^a_i$ with~$1\le i\le m< a \le n$.  Note that these
latter forms generate the module of forms that are semibasic
for the fibration~$\pi:F\to V^\circ$.

Computing exterior derivatives via 
the structure equations yields
\begin{equation}
\begin{split}
0 = \d \omega^{n}_1 
   &= { }-\omega^{n}_A\w\omega^A_1
    = -\sum_{j=3}^m \omega^{n}_j\w\omega^j_1
      -\sum_{a = m+1}^{n-1} \omega^{n}_a\w\omega^a_1\,,\\
0 = \d \omega^{n}_2 
   &= { } -\omega^{n}_A\w\omega^A_2
    = -\sum_{j=3}^m \omega^{n}_j\w\omega^j_2
      -\sum_{a = m+1}^{n-1} \omega^{n}_a\w\omega^a_2\,.\\
\end{split}
\end{equation} 
Reducing these equations modulo $\{\omega^{n}_3,\ldots,\omega^{n}_m\}$
yields
\begin{equation}
\sum_{a = m+1}^{n-1} \omega^{n}_a\w\omega^a_1
\equiv \sum_{a = m+1}^{n-1} \omega^{n}_a\w\omega^a_2
\equiv 0 \mod \omega^{n}_3,\ldots,\omega^{n}_m\,.
\end{equation}
Since~$\{\omega^{a}_1,\omega^a_2\mid m{+}1\le a\le n{-}1\}$
are linearly independent modulo~$\{\omega^{n}_3,\ldots,\omega^{n}_m\}$,
these equations imply that
\begin{equation}
\omega^n_{m+1}\equiv \dots \equiv \omega^n_{n-1}
              \equiv 0 \mod \omega^{n}_3,\ldots,\omega^{n}_m\,.
\end{equation}
Thus, there exist functions~$\{\,s_j^a\mid 3\le j\le m<a\le n{-}1\,\}$
on~$F$ so that~$\omega^n_a = s^j_a\,\omega^n_j$.  The structure equations
now imply that for all~$A<n$, 
\begin{equation}
d\vs_a \equiv 0 \mod \vs_1,\dots,\vs_{n-1},
                     \omega^{n}_3,\ldots,\omega^{n}_m\,.
\end{equation}
Consequently, the map~$[\vs_1\w\dots\w\vs_{n-1}]:F\to \Gr(n{-}1,\C{n}) 
\simeq \bbP^{n-1}$
is a holomorphic map of constant rank~$m{-}2$.  Moreover, since
the forms~$\{\omega^{n}_3,\ldots,\omega^{n}_m\}$ are $\pi$-semibasic,
and~$\pi$ has connected fibers, it follows that there is a 
well-defined holomorphic map~$\xi:V^\circ\to\Gr(n{-}1,\C{n}) $
of constant rank~$m{-}2$ that 
satisfies~$\xi(E) = H_E$, where~$H_E/E = R_E$.
In particular, $[\vs_1\w\dots\w\vs_{n-1}] = \xi\circ\pi$.

The same sort of argument as was made 
in Theorem~\ref{thm:codim2-sigma2}
now shows that there is an irreducible 
variety~$A\subset \Gr(n{-}1,\C{n}) $
of dimension~$m{-}2$ such that~$A$ is the closure of~$\xi(V^\circ)$
and, moreover, that~$V$ consists exactly 
of the union of the~$\Gr(m,H)$
for~$H\in A$.  Details are left to the reader.

Finally, the equation~$r = \deg(A)$ and the converse follow by the 
Schubert calculus and calculation, respectively.
\end{proof}

\begin{remark}[Integral varieties of~$\cI_{(2)^*}$]
\label{rem: integral varieties of I2star}
What was said before in Remark~\ref{rem: integral varieties of I11star}
about the characterization of the local integrals varieties 
of~$\cI_{(1,1)^*}$ applies also to the characterization of the local
integral varieties of~$\cI_{(2)^*}$.  Namely,  every 
codimension~$2$, irreducible integral variety~$V$ of~$\cI_{(2)^*}$ is
\emph{locally} of the form~$\Sigma_m(A)$ for some irreducible 
subvariety~$A\subset \Gr(n{-}1,n)$ of dimension~$m{-}2$.  In this
sense, the codimension~$2$ integral manifolds of~$\cI_{(2)^*}$
depend (in Cartan's sense) on $n{-}m{+}1$ functions of~$m{-}2$ 
variables (in the holomorphic category). 
\end{remark}

\subsection{Ideals of degree~$3$}\label{ssec: extremal 3-cycles}
In this subsection, I turn to the much more interesting case
of the ideals~$\cI_{(3)}$, $\cI_{(2,1)}$, and~$\cI_{(1,1,1)}$
and their application to describing the cycles representing
certain homology classes in~$\Gr(m,n)$.  Of course, the first
and the last of these three are linked by complementarity, so 
there are really only two cases to consider in depth.

\subsubsection{Integrals of~$\cI_{(2,1)}$}
\label{sssec: I21 integral varieties}
It turns out (for reasons that stem from the discussion in
\S\ref{sssec: higher degs}) that the analysis of the integral varieties 
of~$\cI_{(3)}$ and~$\cI_{(1,1,1)}$ is much more difficult than 
the analysis of the integral varieties of~$\cI_{(2,1)}$.  
Thus, I will start with this `middle' case.

\begin{proposition}\label{prop: integral varieties of I21}
For any~$W_-\in\Gr(m{-}1,n)$ and~$W_+\in\Gr(m{+}1,n)$,
the submanifolds~$[W_-,\C{n}]_m\ (\simeq\bbP^{n-m})$
and~$\Gr(m,W_+)\ (\simeq\bbP^m)$ in~$\Gr(m,n)$
are integral manifolds of~$\cI_{(2,1)}$. 

Conversely, for any irreducible subvariety~$V\subset\Gr(m,n)$
of dimension~$d\ge 3$ that is an integral variety of~$\cI_{(2,1)}$,
there exists either a fixed~$(m{+}1)$-plane~$W_+\subset \C{n}$ so 
that~$V\subset\Gr(m,W_+)$
or an $(m{-}1)$-plane~$W_-\subset\C{n}$ 
so that~$V\subset [W_-,\C{n}]_m$.
\end{proposition}

\begin{proof}
By Lemma~\ref{lem: integral elements of I21}, 
every integral element of $\cI_{(2,1)}$ of dimension at least~$3$ 
is either an integral element of~$\cI_{(2)}$ or of~$\cI_{(1,1)}$.  
Moreover,
these two ideals have no integral elements of dimension~$2$ or more 
in common.  Thus, any irreducible integral variety of~$\cI_{(2,1)}$ is
either an integral variety of~$\cI_{(2)}$ or of~$\cI_{(1,1)}$.  Now
apply either Proposition~\ref{prop: integral varieties of I2} 
or Proposition~\ref{prop: integral varieties of I11}, as appropriate.  
\end{proof}

Proposition~\ref{prop: integral varieties of I21} 
has an application to the rigidity of certain
extremal cycles:

\begin{theorem}\label{thm: cycles killed by phi-21}
If~$X\subset \Gr(m,n)$ is an irreducible variety of dimension~$d\ge3$ 
satisfying~$[X] = r[\sigma_{(d)^*}]+s[\sigma_{{(d)'}^*}]$, 
then either
\begin{enumerate}
\item $s=0$ and~$X\subset[A,\C{n}]_m$ for 
some~$(m{-}1)$-plane~$A\subset\C{n}$, or else 
\item $r=0$ and~$X\subset\Gr(m,A)$ 
for some~$(m{+}1)$-plane~$A\subset\C{n}$.
\end{enumerate}
\end{theorem}

\begin{proof}
In each case, the homological assumption implies that~$\phi_{(2,1)}$
vanishes on~$V$ and, hence, that~$V$ is
an integral variety of~$\cI_{(2,1)}$. 
Now apply Proposition~\ref{prop: integral varieties of I21}.
\end{proof}

\begin{remark}\label{rem: cycles killed by phi-21}
The significance of Theorem~\ref{thm: cycles killed by phi-21}
is that it shows how rigid the effective cycles are on an entire
$2$-dimensional `face' of $H^+_{2d}\bigl(\Gr(m,n),\bbZ\bigr)$, namely, 
the semigroup spanned by~$[\sigma_{(d)^*}]$ and~$[\sigma_{{(d)'}^*}]$.
In fact, any effective $d$-cycle in~$\Gr(m,n)$ whose homology class
lies on this `face' is a union of $d$-cycles whose homology classes
lie on the two bounding extremal rays 
generated (individually) by~$[\sigma_{(d)^*}]$ and~$[\sigma_{{(d)'}^*}]$.
\end{remark}

\subsubsection{Bundles with~$c_1c_2-c_3 = 0$}
\label{sssec: c21=0}
Proposition~\ref{prop: integral varieties of I21} can be 
applied to characterize the bundles generated by sections
that satisfy~$c_1c_2 = c_3$.

\begin{theorem}\label{thm: bundles with c3 = c1*c2}
Let~$M$ be a connected compact K\"ahler manifold and let~$F\to M$ be 
a holomorphic vector bundle of rank~$r$ that is generated by 
its sections.  

If~$c_1(F)c_2(F)-c_3(F)=0$,
then one {\upshape(}or more{\upshape)} of the following is true:
\begin{enumerate} 
\item $F = (M{\times}\C{r+1})/L$ for some 
           line bundle~$L\subset M{\times}\C{r+1}$.
\item $F = L\oplus (M{\times}\C{r-1})$ for some line bundle~$L$
           that is generated by its sections.
\item There is an algebraic curve~$C$, a bundle~$F'\to C$ 
      that is generated by its sections, and a holomorphic 
      mapping~$\kappa:M\to C$ so that~$F=\kappa^*(F')$.
\item There is a {\upshape(}possibly singular{\upshape)}
      algebraic surface~$S$, a bundle~$F'\to S$, and a 
      holomorphic mapping~$\kappa:M\to S$ so that~$F = \kappa^*(F')$.
\end{enumerate}
\end{theorem}

\begin{proof}
Let~$H^0(F)$ be the space of global sections of~$F$, a vector
space of dimension~$n = h^0(F)$.  Let~$\ev_{\!F}: M\times H^0(F)\to F$ 
be the evaluation mapping, which, by assumption,
is surjective, so that~$n\ge r$.  The kernel~$K\subset M\times H^0(F)$
is then a subbundle of rank~$m=n{-}r$ and can be used to define a 
mapping~$\kappa_F:M\to\Gr\bigl(m,H^0(F)\bigr)$ that satisfies~
$\kappa_F^*(Q) = F$.  Consequently,
$$
c_1(F)c_2(F)-c_3(F) = c_{(2,1)}(F) = \left[\kappa^*(\phi_{(2,1)})\right].
$$
Thus, the inequality~$c_1(F)c_2(F)-c_3(F)\ge 0$ follows directly from
Theorem~\ref{thm:Fulton-Lazarsfeld}.  

If equality holds, $\kappa_F(M)\subset\Gr\bigl(m,H^0(F)\bigr)$ must
be an integral variety of~$\cI_{(2,1)}$.  There are now five cases:

If $\kappa_F(M)$ is a single point, then~$F$ is trivial.  

If $\kappa_F(M)$ is a curve, 
let~$C\to\kappa_F(M)\subset\Gr\bigl(m,H^0(F)\bigr)$ 
be its normalization and let~$F'$ be the pullback of~$Q$ to~$C$.

If $\kappa_F(M)$ is a (possibly singular)
surface~$S\subset\Gr(m,n)$, let~$F'$ be the pullback of~$Q$ to~$S$.

If $\kappa_F(M)$ has dimension greater than~$2$, then 
Proposition~\ref{prop: integral varieties of I21} implies that 
$\kappa_F(M)$ is either an integral variety of~$\cI_{(2)}$, in which
case~$c_2(F)=0$, so that Theorem~\ref{thm: bundles with c2 = 0} applies,
or else of~$\cI_{(1,1)}$, in which case~$c_{(1,1)}(F)=0$, so that 
Theorem~\ref{thm: bundles with c2 = c1*c1} applies.
\end{proof}

\subsubsection{Integrals of~$\cI_{(3)}$ and~$\cI_{(1,1,1)}$}
\label{sssec: I3 and I111 integral varieties}
Now I turn to the more difficult problem of classifying the integral 
varieties of~$\cI_{(3)}$ (and, by complementarity,~$\cI_{(1,1,1)}$).  
To avoid trivial cases in which~$\cI_{(3)} = (0)$, 
assume that~$n\ge m{+}3$.

\begin{proposition}\label{prop: integral varieties of I3}
The following are integral varieties of~$\cI_{(3)}$ in~$\Gr(m,n)$:
\begin{enumerate}
\item For any~$A\in\Gr(m{+}2,\C{n})$, 
      the $2m$-dimensional submanifold~$\Gr(m,A)$.
\item For any curve~$C\subset\Gr(m{+}1,n)$, the $(m{+}1)$-dimensional 
      subvariety
$$
\Sigma_m(C) = \cup_{B\in C}\Gr(m,B).
$$
\item For any pair~$W_-\in \Gr(m{-}2,n)$ and~$W_+\in\Gr(m{+}3,n)$ with
      $W_-\subset W_+$ and any nondegenerate quadratic form
      $G$ on~$W_+/W_-\simeq\C{5}$, the $3$-dimensional submanifold
      $N_G(W_-,W_+)\subset\Gr(m,n)$ that consists of 
      the $m$-planes~$W$ in $[W_-,W_+]_m$ for which~$W/W_-$
      is $G$-isotropic in~$W_+/W_-$.
\item For any pair~$W_-\in \Gr(m{-}3,n)$ and~$W_+\in\Gr(m{+}3,n)$ with
      $W_-\subset W_+$ and any nondegenerate quadratic form
      $G$ on~$W_+/W_-\simeq\C{6}$, the $3$-dimensional submanifold
      $N_G(W_-,W_+)\subset\Gr(m,n)$ that consists of 
      the $m$-planes~$W$ in $[W_-,W_+]_m$ for which~$W/W_-$
      is $G$-isotropic in~$W_+/W_-$.
\item Any subvariety~$V\subset\Gr(m,n)$ of dimension at most~$2$.
\end{enumerate}
Moreover, any irreducible algebraic integral variety 
of~$\cI_{(3)}$ is a subvariety of an algebraic integral variety of
one of the five listed types.
\end{proposition}

\begin{remark}[The `exceptional' integrals] 
The submanifold~$N_G(W_-,W_+)$ defined in Item~$3$ is isomorphic 
to~$\bbP^3$.  The submanifold~$N_G(W_-,W_+)$ defined in Item~$4$
is isomorphic to the disjoint union of two copies of~$\bbP^3$. 
See~\S\ref{ssec: isotropic Gr} for more information.
\end{remark}

\begin{remark}[The structure of the proof] 
The following proof is rather complex, but I am at a loss as to how
to simplify it.  Perhaps, though, this complexity is unavoidable
in view of the complexity of the resulting classification.  
However, an overview without details may be of some use to the
reader, so I will give it here.

Roughly speaking, the strategy of the proof will be to identify 
the types of maximal integral manifolds of~$\cI_{(3)}$ with the
types of maximal integral elements of~$\cI_{(3)}$ as listed in
Lemma~\ref{lem: integral elements of I3}.  As the reader will see,
the correspondence is not perfect, but this will at least serve
as a guide to organizing the proof.

The first step is to restrict to the smooth locus~$X^\circ$
of an irreducible integral variety~$X$ of~$\cI_{(3)}$ 
of dimension~$d\ge 3$
and introduce an integer~$\delta(V)$ for~$V\in X^\circ$ 
that is the dimension of the smallest subspace~$Q'_V\subset Q_V$ for 
which~$T_VX\subset T_V\Gr(m,n)$ lies in~$Q'_V\ot V^*$.  The
set~$X^\star\subset X^\circ$ on which~$\delta$ attains its
maximum value~$\delta(X)\le n{-}m$ is a Zariski open subset 
of~$X$ and~$Q'\to X^\star$ is a holomorphic bundle over~$X^\star$.

The proof is then broken up into cases according to~$\delta(X)$.  
When~$\delta(X)=1$ (an easy case),~$X$ actually lies in~$\Gr(m,B)$ 
for some~$B\in \Gr(m{+}1,n)$, so that~$X$ falls into the first category
of the proposition.

The case~$\delta(X)=2$ turns out to be the most complicated, as there
are several subcases and some of these even have their own subsubcases.
The basic idea is that if the tangent spaces are sufficiently `free'
in an appropriate sense, then one can show that~$X$ must lie 
in~$\Gr(m,A)$ for some~$A\in\Gr(m{+}2,n)$.  However, there is one
`degenerate' subsubcase in which a lack of `freeness' allows the
tangent spaces to~$X$ to vary in such a way that~$X$ can only be
shown to lie in a curve of~$\Gr(m,m{+}1)$s, as described in the
second category of the proposition.

The case~$\delta(X)=3$ is, in some ways, the most interesting. 
The only possibility for the tangent spaces
of~$X^\circ$ are the integral elements of~$\cI_{(3)}$ 
that are of the second and third types listed in 
Lemma~\ref{lem: integral elements of I3}.  
When the tangent spaces are of the second type, one can find
a Zariski-open~$X^\bullet\subset X$ (that lies in~$X^\star$)
and define a canonical~$A:X^\bullet\to\Gr(m{+}1,n)$ whose differential 
has rank is at least equal to~$1$.  In the case that the rank of~$dA$
is identically~$1$, it is not difficult to show that~$X$ must belong
to the second category of the proposition (in fact,~$C$ is the closure
of the image~$A(X^\bullet)$).  When the rank of~$dA$ is greater 
than~$1$, one shows that~$\dim X = 3$ and then a (rather
involved) moving frame analysis shows that~$X$ belongs 
to the third category of the proposition.
When the tangent spaces are of the third type listed in 
Lemma~\ref{lem: integral elements of I3}, 
then an analysis via the moving frame (also rather involved)
shows that~$X$ must belong to the fourth category of the proposition.

When~$\delta(X)>3$, the only possibility for the tangent spaces
of~$X^\circ$ are the integral elements of~$\cI_{(3)}$ 
that are of the second type listed in 
Lemma~\ref{lem: integral elements of I3}.  The moving frame
analysis in this case is straightforward, with the result 
that~$X$ belongs to the second category of the proposition.
\end{remark}

\begin{proof} Verifying that each of the types listed is
indeed an integral manifold of~$\cI_{(3)}$ is relatively straightforward.
Simple calculations via the moving frame show that the tangent spaces
to these subvarieties at their smooth points are integral elements 
of~$\cI_{(3)}$.  Alternatively, the calculations below will provide a
direct proof.

Thus, let~$X\subset\Gr(m,n)$ be an irreducible integral variety 
of~$\cI_{(3)}$.  If $\dim X \le 2$, there is nothing to prove, so assume
that~$\dim X = d\ge 3$.  

Let~$X^\circ\subset X$ be the smooth part of~$V$.
Then for each~$V\in X^\circ$, the subspace~$T_VX\subset T_V\Gr(m,n)
=Q_V{\ot} V^*$ is
an integral element of~$\cI_{(3)}$ of dimension at least~$3$.  
For~$V\in X^\circ$, let~$Q'_V\subset Q_V = \C{n}/V$ be the smallest 
subspace for which $T_VX$ is contained in~$Q'_V\ot V^*$.  
The function~$\delta:X^\circ\to \bbZ^+$ 
defined by~$\delta(V)=\dim Q'_V\ge1$
is equal to its maximum value, say~$\delta(X)$, 
on a subset~$X^\star\subset X^\circ$ that is open and dense in~$X^\circ$
and connected.  (For any~$q$, the set of~$V\in X^\circ$
for which~$\delta(V)< q$ is easily seen to be 
an analytic subvariety of~$X^\circ$.)

Suppose, first, that~$\delta(X) = 1$. 
Then~$Q'_V{\ot}V^*$ is an integral 
element of~$\cI_{(2)}$ for all~$V\in X^\star$, so $X^\star$ and, hence,
$X^\circ$ and~$X$ are integral varieties of~$\cI_{(2)}$.  By 
Proposition~\ref{prop: integral varieties of I2}, 
there exists a~$B\in \Gr(m{+}1,n)$ so that~$X\subset\Gr(m,B)$. 
Choosing~$A\in[B,\C{n}]_{m+2}$, it follows that~$X\subset\Gr(m,A)$.
Thus,~$X$ lies in an integral manifold of the 
first category listed in the proposition.

Suppose, second, that~$\delta(X)=2$.  In this case, I claim that,
for each~$V\in X^\star$, the space~$T_VX$ lies in a unique maximal
integral element, namely~$Q'_V\ot V^*$.  To see this, first
note that~$Q'_V\ot V^*$ is a maximal integral element of the 
first type.  Now, since it has dimension at least~$3$, $T_VX$ 
does not lie in an integral element
of the fourth type listed in Lemma~\ref{lem: integral elements of I3}.
Also, it cannot lie in an integral element of the third type, because
these integral elements have dimension~$3$, which would force~$T_VX$
to equal a maximal integral element of the third type, but these
integral elements do not lie in any integral element of the first type.
Finally, suppose~$T_VX$ were to lie in an integral element of the 
second type, say~$T_VX\subset L{\ot}V^* + \bbC{\cdot}R$,
where~$L\subset Q_V$ is a line and~$R\in Q_V{\ot} V^*$ has the 
property that~$\bar R\in Q_V/L \ot V^*$ has rank at least~$2$.  Now,
is easy to see that the only subspaces of~$L{\ot}V^* + \bbC{\cdot}R$
that have dimension at least~$3$ and that lie in a subspace of the 
form~$P\ot V^*$ where~$P\subset Q$ has dimension at most~$2$ are
the subspaces of~$L{\ot}V^*$.  Consequently, if~$T_VX$ were to
lie in~$L{\ot}V^* + \bbC{\cdot}R$, then it would follow that
$T_VX$ lies in~$L{\ot}V^*$.  But this would violate the assumption
that~$\delta(V)=2$.  Thus, the claim has been established.

Because of the evident uniqueness of~$Q'_V$ for~$V\in X^\star$, 
the family of vector spaces~$Q'\to X^\star$ is a holomorphic
subbundle of~$Q\to X^\star$ of rank~$2$.  Now I need to introduce 
another invariant.  For~$V\in X^\star$, say that a subspace~$S\subset V$ 
is \emph{free} if the composition
$$
\rho_S: T_VX\hookrightarrow Q'_V{\ot}V^*
\longrightarrow Q'_V\ot S^*
$$
is surjective.  For each~$p\le m$, the set of non-free subspaces 
of~$V$ of dimension~$p$ is an algebraic subset of~$\Gr(p,V)$.
Let~$\sigma(V)\le m$ be the dimension of the 
largest free subspace of~$V$ and let~$s$ be the 
maximum of~$\sigma(V)$ for~$V\in X^\star$.  Let~$X^\bullet\subset X^\star$
be subset consisting of those~$V\in X^\star$ for which $\sigma(V)=s$.
Since the complement of~$X^\bullet$ in~$X^\star$ is evidently a
proper analytic subvariety of~$X^\star$, it follows that~$X^\bullet$ is
open and dense in~$X^\star$ and is connected.  

Now, I claim that~$s\ge 1$.  This follows by elementary linear
algebra from the assumptions~$\dim T_VX\ge 3$ and~$\delta(X)\ge 2$,
so I will leave this to the reader.  

Suppose first that~$s\ge 2$ and let~$F\subset X^\bullet\times\GL(n,\bbC)$ 
be the set of pairs~$(V,\vs)$ that satisfy the conditions
\begin{enumerate}
\item $\vs_1,\dots,\vs_m$ span $V$,
\item $\vs_{m-s+1},\dots,\vs_m$ span a free subspace of~$V$, and
\item $\dl\vs_{m+1}\dr,\ \dl\vs_{m+2}\dr$ spans~$Q'_V$.
\end{enumerate}
Then~$F$ is connected, as follows from the connectedness of~$X^\bullet$.

Consider the usual structure equations:
$$
\d\vs_A = \vs_B\omega^B_A\qquad\qquad 
\d\omega^A_B = -\omega^A_C\w\omega^C_B\,.
$$
The conditions defining~$F$ imply that~$\omega^\alpha_i = 0$ 
whenever~$1\le i\le m$ and~$\alpha>m{+}2$.  Moreover, the freeness
assumption implies that the~$2s$ entries of the matrix~
$$
\begin{pmatrix} 
\omega^{m+1}_{m-s+1}&\dots&\omega^{m+1}_m\\
\noalign{\vskip2pt}\\ 
\omega^{m+2}_{m-s+1}&\dots&\omega^{m+2}_m
\end{pmatrix}
$$
are linearly independent on~$F$.  Now, when~$m{-}s{+}1\le i\le m$ 
and~$\alpha>m{+}2$, the structure equations combined with the
stated vanishing of forms give
$$
0 = \d\omega^\alpha_i = -\omega^\alpha_A\w\omega^A_i 
 = {} - \omega^\alpha_{m+1}\w\omega^{m+1}_i 
      - \omega^\alpha_{m+2}\w\omega^{m+2}_i\,.
$$
This implies, because of the stated linear independence, that 
$$
\omega^\alpha_{m+1}\equiv \omega^\alpha_{m+2}\equiv 0\mod
   \omega^{m+1}_i,\omega^{m+2}_i
$$
for each~$m{-}s{+}1\le i\le m$.  Since~$s\ge2$ by hypothesis, the stated
linear independence implies that 
$\omega^\alpha_{m+1}=\omega^\alpha_{m+2}=0$, for all~$\alpha>m{+}2$.

In turn, this implies that
$$
\d\vs_1\equiv\cdots\equiv\d\vs_{m+2}\equiv 0 \mod \vs_1,\cdots,\vs_{m+2}\,.
$$
In other words, the $(m{+}2)$-plane spanned by~$\vs_1,\cdots,\vs_{m+2}$
is locally constant on~$F$.  Since~$F$ is connected, this implies that
there is a fixed~$A\in\Gr(m{+}2,n)$ that is spanned 
by~$\vs_1,\cdots,\vs_{m+2}$ at all points of~$F$.  By construction,
this implies that~$A$ contains~$V$ for all~$V\in X^\bullet$.  Of course,
this implies that~$X^\bullet$ lies in~$\Gr(m,A)$, and hence that~$X$
lies in~$\Gr(m,A)$.  

Thus, suppose instead that~$s = 1$.   For each~$V\in X^\bullet$, the
set of lines~$L\subset V$ that are free is the complement in~$\Gr(1,V)$
of an algebraic subset and is therefore open, dense, and connected.
By hypothesis, for any $2$-plane~$P\subset V$, the induced
mapping~$\rho_P:T_VX\to Q'_V{\ot}P^*$ is not surjective.  
I claim that, for~$P$ outside a closed algebraic set in~$\Gr(2,V)$, 
the mapping~$\rho_P$ has rank~$3$.  Certainly, the set of~$P$ for
which the rank of~$\rho_P$ is at most~$2$ is an algebraic subvariety
of~$\Gr(2,V)$, so it suffices to show that it is not everything.
However, this follows by linear algebra since~$\dim T_VX\ge 3$ 
and~$\delta(V)\ge2$.  I leave details to the reader.  

Say that a~$P\in\Gr(2,V)$ is \emph{semi-free} 
if the rank of~$\rho_P$ is~$3$.  When~$P$ is semi-free, the
annihilator of~$\rho_P(T_VX)\subset Q'_V{\ot}P^*$ is a line 
in~$(Q'_V)^*{\ot}P$.  There are two subcases now to consider.  
The first is when the tensor rank of a generator of this line 
is generically equal to~$2$.  The second is when the tensor rank
of a generator of this line is equal to~$1$ everywhere.

Consider the first subcase and let~$X^\diamond\subset X^\bullet$ 
be the open, dense, connected subset consisting of those~$V\in X^\bullet$ 
for which there exist~$P\in\Gr(2,V)$ so that the rank of a generator
of the annihilator of~$\rho_P(T_VX)$ in~$(Q'_V)^*{\ot}P$ is equal to~$2$.
Let~$F^\diamond\subset X^\diamond\times\GL(n,\bbC)$ 
be the set of pairs~$(V,\vs)$ that satisfy the conditions
\begin{enumerate}
\item $\vs_1,\dots,\vs_m$ span $V$,
\item $\dl\vs_{m+1}\dr,\ \dl \vs_{m+2}\dr$ spans~$Q'_V$.
\item $\vs_{m-1},\ \vs_m$ span a semi-free plane~$P\subset V$, 
      and, moreover,
      the annihilator of~$\rho_P(T_VX)$ is spanned by 
      $\vs^{m+2}{\ot}\vs_{m-1}+\vs^{m+1}{\ot}\vs_{m}\in (Q'_V)^*{\ot}P$.
\end{enumerate}
Then~$F^\diamond$ is connected, as follows from 
the connectedness of~$X^\diamond$.

Consider the usual structure equations:
$$
\d\vs_A = \vs_B\omega^B_A\qquad\qquad 
\d\omega^A_B = -\omega^A_C\w\omega^C_B\,.
$$
The conditions that define~$F^\diamond$ give~$\omega^\alpha_i = 0$ 
whenever~$1\le i\le m$ and~$\alpha>m{+}2$.  Moreover, the semi-freeness
and the assumption about the annihilator imply that the four $1$-forms~
$$
\{\omega^{m+1}_{m-1},\omega^{m+1}_m,\omega^{m+2}_{m-1},\omega^{m+2}_m\}
$$
satisfy exactly one linear relation, 
$\omega^{m+2}_{m-1}+\omega^{m+1}_m = 0$,
and are otherwise linearly independent on~$F^\diamond$.  

Just as in the case~$s\ge2$, when~$i= m{-}1$ or~$m$ and~$\alpha>m{+}2$, 
the structure equations combined with the stated vanishing of forms give
$$
0 = \d\omega^\alpha_i = -\omega^\alpha_A\w\omega^A_i 
 = {} - \omega^\alpha_{m+1}\w\omega^{m+1}_i 
      - \omega^\alpha_{m+2}\w\omega^{m+2}_i\,.
$$
This implies, because of the stated linear independence, that 
$$
\omega^\alpha_{m+1}\equiv \omega^\alpha_{m+2}\equiv 0\mod
   \omega^{m+1}_i,\omega^{m+2}_i
$$
when~$i= m{-}1$ or~$m$.  Now, however, because the span 
of~$\{\omega^{m+1}_{m-1},\omega^{m+2}_{m-1}\}$ intersects the span
of~$\{\omega^{m+1}_{m},\omega^{m+2}_{m}\}$ in the multiples 
of~$\omega^{m+1}_m$, these congruences only imply that 
$$
\omega^\alpha_{m+1}\equiv\omega^\alpha_{m+2}\equiv0\mod \omega^{m+1}_m
$$
for all~$\alpha>m{+}2$.  Setting~$\omega^\alpha_j 
= R^\alpha_j\,\omega^{m+1}_m$ for~$j = m{+}1$ 
and~$m{+}2$ and~$\alpha>m{+}2$
and substituting this back into the relation
$$
0 = - \omega^\alpha_{m+1}\w\omega^{m+1}_i 
    - \omega^\alpha_{m+2}\w\omega^{m+2}_i
$$
for~$i = m{-}1$ and~$m$ shows that~$R^\alpha_j=0$ 
for $j = m{+}1$ and~$m{+}2$ and~$\alpha>m{+}2$.  
Thus, $\omega^\alpha_j=0$
for~$\alpha$ and~$j$ with these ranges, 
just as in the~$s\ge 2$ case.

In turn, this implies that
$$
\d\vs_1\equiv\cdots\equiv\d\vs_{m+2} 
       \equiv 0 \mod \vs_1,\cdots,\vs_{m+2}\,.
$$
In other words, the $(m{+}2)$-plane spanned by~$\vs_1,\cdots,\vs_{m+2}$
is locally constant on~$F$.  Since~$F$ is connected, this implies that
there is a fixed~$A\in\Gr(m{+}2,n)$ that is spanned 
by~$\vs_1,\cdots,\vs_{m+2}$ at all points of~$F$.  By construction,
this implies that~$A$ contains~$V$ for all~$V\in X^\bullet$.  Of course,
this implies that~$X^\bullet$ lies in~$\Gr(m,A)$, and hence that~$X$
lies in~$\Gr(m,A)$.  This finishes the first subcase of~$s=1$.

Consider the second subcase, in which the rank of the annihilator
of~$\rho_P(T_VX)$ in~$(Q'_V)^*{\ot}P$ is equal to~$1$ for 
all semi-free~$P\in\Gr(2,V)$ and~$V\in X^\bullet$.

Let~$F^\bullet\subset X^\bullet\times\GL(n,\bbC)$ 
be the set of pairs~$(V,\vs)$ that satisfy the conditions
\begin{enumerate}
\item $\vs_1,\dots,\vs_m$ span $V$,
\item $\dl\vs_{m+1}\dr,\ \dl \vs_{m+2}\dr$ spans~$Q'_V$.
\item $\vs_{m-1},\ \vs_m$ span a semi-free plane~$P\subset V$, 
      and, moreover, the annihilator of~$\rho_P(T_VX)$ is spanned 
      by~$\vs^{m+2}{\ot}\vs_{m-1}\in (Q'_V)^*{\ot}P$.
\end{enumerate}
Then~$F^\bullet$ is connected, as follows from 
the connectedness of~$X^\bullet$.

Consider the usual structure equations:
$$
\d\vs_A = \vs_B\,\omega^B_A\qquad\qquad 
\d\omega^A_B = -\omega^A_C\w\omega^C_B\,.
$$
The conditions defining~$F^\bullet$ imply that~$\omega^\alpha_i = 0$ 
whenever~$1\le i\le m$ and~$\alpha>m{+}2$.  Moreover, the semi-freeness
and the assumption about the annihilator imply that the three $1$-forms~
$\{\omega^{m+1}_{m-1},\omega^{m+1}_m,\omega^{m+2}_m\}$
are linearly independent, while $\omega^{m+2}_{m-1} = 0$.

When~$\alpha>m{+}2$, the structure equations combined with 
the stated vanishing of forms give
$$
0 = \d\omega^\alpha_{m-1} = -\omega^\alpha_A\w\omega^A_{m-1} 
 = {} - \omega^\alpha_{m+1}\w\omega^{m+1}_{m-1}\,.
$$
so~$\omega^\alpha_{m+1} = R^\alpha\,\omega^{m+1}_{m-1}$ 
for some functions~$R^\alpha$.  The structure equations then give
\begin{equation*}
\begin{split}
0 = \d\omega^\alpha_m = -\omega^\alpha_A\w\omega^A_m 
 &= {} - \omega^\alpha_{m+1}\w\omega^{m+1}_m 
      - \omega^\alpha_{m+2}\w\omega^{m+2}_m\\
 &= {} - R^\alpha\omega^{m+1}_{m-1}\w\omega^{m+1}_m 
      - \omega^\alpha_{m+2}\w\omega^{m+2}_m\,,
\end{split}
\end{equation*}
which implies, first, that~$R^\alpha = 0$ for all~$\alpha>m{+}2$
and then that there must exist functions~$S^\alpha$ so that
$\omega^\alpha_{m+2} = S^\alpha\,\omega^{m+2}_{m}$.

If all of the~$S^\alpha$ vanish identically, then~$\omega^\alpha_{m+1}
=\omega^\alpha_{m+2} = 0$ when~$\alpha>m{+}2$, which implies, as before, 
that
$$
\d\vs_1\equiv\cdots\equiv\d\vs_{m+2}\equiv 0 \mod \vs_1,\cdots,\vs_{m+2}\,.
$$
In other words, the $(m{+}2)$-plane spanned by~$\vs_1,\cdots,\vs_{m+2}$
is locally constant on~$F^\bullet$.  Since~$F^\bullet$ is connected, 
this implies that there is a fixed~$A\in\Gr(m{+}2,n)$ that is spanned 
by~$\vs_1,\cdots,\vs_{m+2}$ at all points of~$F$.  By construction,
this implies that~$A$ contains~$V$ for all~$V\in X^\bullet$.  Of course,
this implies that~$X^\bullet$ lies in~$\Gr(m,A)$, and hence that~$X$
lies in~$\Gr(m,A)$.

Suppose, instead, that the~$S^\alpha$ do not vanish identically.
The set~$F^\diamond\subset F^\bullet$ on which at least one of 
the~$S^\alpha$ is nonzero is the complement of an analytic subvariety
of~$F^\bullet$ and hence is open and dense in~$F^\bullet$ and connected.
Its image~$X^\diamond\subset X^\bullet$ is easily seen to be the 
complement of a proper analytic subvariety in~$X^\bullet$, so~$X^\diamond$
is open and dense in~$X^\bullet$ and is connected also.  

Now, I claim that for~$1\le i<m{-}1$, there exist functions~$U_i$ 
on~$F^\diamond$ so that~$\omega^{m+2}_i = U_i\,\omega^{m+2}_m$.  
To see this, differentiate the relation~$\omega^\alpha_i = 0$ 
for~$\alpha>m{+}2$ and~$i<m{-}1$, which now yields
$$
0 = -\omega^\alpha_A\w\omega^A_i 
= -\omega^\alpha_{m+2}\w\omega^{m+2}_i
= -S^\alpha\,\omega^{m+2}_m\w \omega^{m+2}_i\,.
$$
Since not all of the~$S^\alpha$ vanish, this implies the desired
relations.  Now that this has been established, 
the fact that there must be 
at least~$d{-}1\ge2$ one-forms 
among~$\{\omega^{m+1}_1,\dots\omega^{m+1}_m\}$
that are independent modulo~$\omega^{m+2}_m$ 
shows that~$d\le m{+}1$ 
and that~$\dl\vs_{m+1}\dr$ 
spans the unique line~$L_V\subset Q'_V$ with the 
property that~$T_VX$ meets~$L_V{\ot}V^*$ 
in a subspace of dimension~$d{-}1$. 

Write~$T_VX\cap (L_V{\ot}V^*) = L_V{\ot}P_V$ where~$P_V\subset V^*$ has 
dimension~$d{-}1$.   

Consider the subset~$F'\subset F^\diamond$ consisting 
of the~$(V,\vs)\in F^\diamond$ for which~$(P_V)^\perp\subset V$ is 
spanned by the~$\vb_i$ for which~$i\le m{-}d{+}1$.  The above arguments
show that~$F'\to X^\diamond$ is a $G$-bundle where
$$
G' = P_{m-d+1}\cap P_m\cap P_{m+1}\cap P_{m+2}\subset\GL(n,\bbC),
$$
and so is connected. Consequently, there are well-defined 
mappings~$A:X^\diamond\to\Gr(m{+}1,n)$ and~$B:X^\diamond\to\Gr(m{+}2,n)$
with the property that, for all~$V\in X^\diamond$ and~$(V,\vs)\in F'$,
the span of~$\{\vs_1,\dots,\vs_{m+1}\}$ is~$A(V)$ and 
the span of~$\{\vs_1,\dots,\vs_{m+2}\}$ is~$B(V)$.

In addition to the relations already found, 
the relations~$\omega^{m+1}_i\equiv 0 \mod\omega^{m+2}_m$ hold for
$i\le m{-}d{+}1$ while the $d$ 
one-forms~$\{\omega^{m+1}_{m-d+2},\dots,\omega^{m+1}_m,\omega^{m+2}_m\}$
are linearly independent and generate the semibasic forms
for the map~$F'\to X^\diamond$.

When~$\alpha>m{+}2$, the structure equations and the stated and derived 
vanishing (including~$\omega^\alpha_{m+1}=0$) give
\begin{equation*}
\begin{split}
S^\alpha\,\d\omega^{m+2}_{m}
&\equiv d\bigl(S^\alpha\,\omega^{m+2}_m\bigr)\mod \omega^{m+2}_m\\
&\equiv d\omega^\alpha_{m+2} \equiv -\omega^\alpha_A\w\omega^A_{m+2}\\
&\equiv {}-\omega^\alpha_{m+2}\w\omega^{m+2}_{m+2}
          -\sum_{\beta>m{+}2}\omega^\alpha_{\beta}\w\omega^{\beta}_{m+2}\\
&\equiv {} - S^\alpha\,\omega^{m+2}_{m}\w\omega^{m+2}_{m+2}
          - \sum_{\beta>m{+}2}\omega^\alpha_{\beta}
        \w (S^\beta\,\omega^{m+2}_{m})\\
&\equiv 0 \mod \omega^{m+2}_m\,.
\end{split}
\end{equation*}
Since not all of the~$S^\alpha$ vanish, it follows that
\begin{equation}\label{eq: omega m+2 m is integrable}
\d\omega^{m+2}_{m} \equiv 0 \mod \omega^{m+2}_m
\end{equation}
on~$F'$.   Thus~$F'$ is foliated by hypersurfaces that are the 
leaves of~$\omega^{m+2}_m=0$.  Since~$\omega^{m+2}_m$ is semibasic
for~$F'\to X^\diamond$ and since the fibers of this submersion are
connected, this foliation pushes down to
define a codimension~$1$ foliation~$\cF$ of~$X^\diamond$.  The tangent
space to the $\cF$-leaf through~$V\in X^\diamond$ is simply 
the unique rank~$1$ subspace of~$T_VX\subset Q'_V\ot V^*$ 
of dimension~$d{-}1$, namely~$L_V\ot P_V$.

Now, differentiating the 
relation~$\omega^{m+2}_{m-1}=0$ yields
\begin{equation*}
\begin{split}
0 &= -\omega^{m+2}_A\w\omega^A_{m-1} \\
&= {} - U_i\,\omega^{m+2}_m\w\omega^i_{m-1} 
      - \omega^{m+2}_m\w\omega^{m}_{m-1} 
      - \omega^{m+2}_{m+1}\w\omega^{m+1}_{m-1}\\
&\equiv - \omega^{m+2}_{m+1}\w\omega^{m+1}_{m-1} \mod\omega^{m+2}_m\,,
\end{split}
\end{equation*}
so~$\omega^{m+2}_{m+1}\equiv V\,\omega^{m+1}_{m-1}
\mod\omega^{m+2}_m$ for some function~$V$ on~$F'$.  
On the other hand, applying the structure
equations to expand the relation~\eqref{eq: omega m+2 m is integrable}
yields the relation
$$
- \omega^{m+2}_{m+1}\w\omega^{m+1}_{m} \equiv 0 \mod \omega^{m+2}_m\,
$$
which implies that~$V = 0$.  Thus, $\omega^{m+2}_{m+1} 
= U_{m+1}\,\omega^{m+2}_{m}$ for some function~$U_{m+1}$ on~$F'$.

Thus, setting~$U_m=1$ for notational consistency, the identities derived
so far imply that, for~$1\le k\le m{+}1$, 
$$
\d \vs_k \equiv\vs_{m+2}\>\bigl(U_k\,\omega^{m+2}_m\bigr)
            \mod \vs_1,\dots,\vs_{m+1}\,.
$$
Consequently, the differential of the mapping~$A:X^\diamond\to\Gr(m{+}1,n)$ 
has rank equal to~$1$ everywhere and each fiber is a union of leaves
of the foliation~$\cF$.  In other words, if~$W = A(V)$, then~$\Gr(m,W)$
intersects~$X^\diamond$ in a subvariety of dimension~$d{-}1$.  

Now suppose that~$X$ is algebraic (as well as irreducible).  
Then~$X^\diamond$ is the complement of a proper algebraic 
subvariety~$Z\subset X$ and it is not difficult to see 
that~$A:X^\diamond\to\Gr(m{+}1,n)$ is the restriction of a 
rational%
\footnote{The rationality of this map follows
from the fact that the space~$A(V)$ can be found as the kernel
of a linear map constructed from the second fundamental form
of~$X$ in~$\Gr(m,n)$.  The (routine) details are left to the reader.}
map of~$X$ into~$\Gr(m{+}1,n)$ whose 
indeterminacy locus is contained in~$Z$.
Since the rank of the differential of~$A$ is
equal to~$1$ on~$X^\diamond = X\setminus Z$, an elementary argument
shows that there is an irreducible algebraic curve~$C\subset\Gr(m{+}1,n)$
so that the graph of $A$ over~$X^\diamond$ is contained in~$X\times C$.
In particular, the closure of this graph is contained in~$X\times C$.

Now consider the subvariety~$\Sigma_m(C)\subset\Gr(m,n)$ 
that is the union of the~$\Gr(m,B)$ for~$B\in C$.   
At its smooth points, the tangent spaces to~$\Sigma_m(C)$ are
integral elements of~$\cI_{(3)}$, so~$\Sigma_m(C)$ is an integral
variety of~$\cI_{(3)}$.  Now, for~$V\in X^\diamond$, the
intersection~$X^\diamond\cap\Gr\bigl(m,A(V)\bigr)\subset\Sigma_m(C)$
has dimension~$d{-}1$ and contains~$V$, so it follows that~$X^\diamond$
itself is contained in~$\Sigma_m(C)$.  Since~$X^\diamond$ is Zariski
dense in~$X$, it follows that~$X$ is contained in~$\Sigma_m(C)$
as well, so the case~$\delta(V)=2$ is finally completed. 

Third, suppose that~$\delta(X)\ge3$.  Each tangent space~$T_VX$ 
for~$V\in X^\star$ must then lie in a maximal integral element 
of~$\cI_{(3)}$ of either the second or third types listed in 
Lemma~\ref{lem: integral elements of I3}.  
Now, none of the vectors in a maximal 
integral element of the third type is of tensor rank one 
when regarded as an element of~$Q_V{\ot}V^*$, 
while the vectors of tensor rank one in a maximal integral element 
of the second type form a canonical subspace of codimension~$1$.  
Since the dimension of~$T_VX$ is at least~$3$, it follows that~$T_VX$
lies in a unique integral element of exactly one of these two types,
depending on whether or not~$T_VX$ contains any vectors of tensor rank one.
Consequently, there are two possibilities:  
Either there is an open, dense, connected
subset of~$X^\star$ consisting of those~$V\in X^\star$ for which~$T_VX$
lies in a maximal integral element of the second type, or else, there
is an open dense, connected subset of~$X^\star$ consisting of those~$V
\in X^\star$ for which~$T_VX$ is an integral element of the third type.
I will now treat these two cases in turn.

Thus, suppose first that~$X^\bullet$ is an open, dense, connected
subset of~$X^\star$ with the property that~$T_VX$ lies in an
integral element of the second type for all~$V\in X^\bullet$.  
In particular,
there exists a line~$L_V\subset Q_V$ so that~$T_VX\cap L_V{\ot}V^*$
has dimension $d{-}1$ for all~$V\in X^\bullet$.
As before, let~$P_V\subset V^*$ be the subspace of dimension~$d{-}1$
so that~$T_VX\cap L_V{\ot}V^* = L_V\ot P_V$.  The uniqueness of
this line~$L_V$ implies that the family of lines~$L\to X^\bullet$
is a holomorphic line subbundle of~$Q\to X^\bullet$, 
while the family of subspaces~$P\to X^\bullet$ 
is a holomorphic subbundle of~$S^*\to X^\bullet$.  
For each~$V\in X^\bullet$, let~$A(V)\in\Gr(m{+}1,n)$ be the subspace
that satisfies~$A(V)/V = L_V$.  The rank of the differential 
of~$A:X^\bullet\to\Gr(m{+}1,n)$ is at least~$1$ everywhere 
since the kernel of~$dA$ lies in~$T_VX\cap L_V{\ot}V^*$.

Suppose first that that~$A:X^\bullet\to\Gr(m{+}1,n)$ 
is a holomorphic map whose differential
has rank equal to~$1$ everywhere on~$X^\bullet$.  
Then, again, just
as in the concluding subsubcase of the~$\delta(X)=2$ argument, 
$X^\bullet$ is foliated in codimension~$1$ by leaves of the 
form~$X^\bullet\cap\Gr\bigl(m,A(V)\bigr)$ for~$V\in X^\bullet$.
Since~$X$ is an irreducible algebraic variety, it is not difficult to
show that~$A$ is a rational map from~$X$ to~$\Gr(m{+}1,n)$ and, again,
it follows, just as in the previous argument, 
that there is an irreducible algebraic curve~$C\subset\Gr(m{+}1,n)$ 
with the property that~$X^\bullet\times C$ contains the graph of~$A$.  
In particular,~$X$ is a subvariety of~$\Sigma_m(C)$ and therefore
belongs to the second category of the proposition.

Thus, suppose that the rank of~$dA$ is sometimes greater than~$1$.
I am going to show that this implies that~$d=\delta(X)=3$ and then 
that~$X$ necessarily belongs to the third category of the proposition.
Let~$X^\diamond\subset X^\bullet$ be the Zariski open subset on 
which the rank of~$dA$ reaches its maximum and 
let~$F^\diamond\subset X^\diamond\times\SL(n,\bbC)$ denote the
set of pairs~$(V,\vs)$ that satisfy
\begin{enumerate}
\item $\vs_1,\dots,\vs_m$ spans~$V$,
\item $\vs_1,\dots,\vs_{m-d+1}$ spans~$(P_V)^\perp\subset V$,
\item $\dl\vs_{m+1}\dr$ spans~$L_V$.
\end{enumerate}
Then~$F^\diamond$ is a~$G$-bundle over~$X^\diamond$ where
$G = P_{m-d+1}\cap P_{m}\cap P_{m+1}$, so $F^\diamond$ is connected.
Consider the structure equations as usual.

By the construction of~$F^\diamond$, the forms $\omega^\alpha_i$ 
with~$i\le m$ and~$\alpha>m{+}1$ together with the forms $\omega^{m+1}_i$ 
with~$i\le m{-}d{+}1$ are pairwise linearly dependent.  
Choose%
\footnote{This choice is not canonical, of course, but this will not
matter.  The reader who wants a canonical construction at this point
is free to consider instead the $\C{*}$-bundle over~$F^\diamond$ on
which such an $\omega_0$ can be canonically defined.}
a $1$-form~$\omega_0$ (which will be unique up to multiples) so that
there exist~$R^\alpha_i$ and~$S_i$ so that
$\omega^\alpha_i = R^\alpha_i\,\omega_0$ when~$i\le m$ and~$\alpha>m{+}1$
while~$\omega^{m+1}_i = S_i\,\omega_0$ when~$i\le m{-}d{+}1$.

By construction, the forms~$\{\omega_0,\omega^{m+1}_{m-d+2},\dots,
\omega^{m+1}_m\}$ are a basis for the $X^\diamond$-semibasic $1$-forms
on~$F^\diamond$.  Moreover, because $T_VX$ is an integral element of
the second type in Lemma~\ref{lem: integral elements of I3}, the
rank of the $(n{-}m{-}1)$-by-$m$ matrix~$R = (R^\alpha_i)$ is
at least equal to~$2$ everywhere.  

Now, I claim that~$\omega_0$ cannot be integrable.  Indeed, suppose
that~$d\omega_0\equiv0\mod\omega_0$.  Then, taking~$\alpha>m{+}1$
and~$i>m{-}d{+}1$, expanding out the structure equation~$\d\omega^\alpha_i
= -\omega^\alpha_C\w\omega^C_i$, using the 
congruences $\omega^\beta_i\equiv0\mod\omega_0$ 
when~$i\le m$ and~$\beta>m{+}1$, and reducing modulo~$\omega_0$ yields
the relation
$$
0 \equiv -\omega^\alpha_{m+1}\w\omega^{m+1}_i\mod\omega_0\,.
$$
Since~$d\ge 3$, these congruences, together with the linear independence
of the $1$-forms~$\{\omega_0,\omega^{m+1}_{m-d+2},\dots,\omega^{m+1}_m\}$, 
imply that $\omega^\alpha_{m+1}\equiv 0\mod\omega_0$.  Thus, set
$\omega^\alpha_{m+1} = R^\alpha_{m+1}\,\omega_0$.  Then the structure
equations so far imply the relations
$$
\d \vs_j \equiv \sum_{\alpha=m+2}^n  \vs_\alpha\>R^\alpha_j\,\omega_0
\mod \vs_1,\dots,\vs_{m+1}
$$
for~$1\le j\le m{+}1$ and this implies that~$A$ has rank~$1$ everywhere,
contrary to hypothesis.  Thus, $\omega_0$ is not integrable, as claimed.

Next, I claim that~$R^\alpha_i=0$ for~$i\le m{-}d{+}1$ and~$\alpha>m{+}1$.
This follows because when the pair~$(i,\alpha)$ satisfy these
restrictions, expanding the structure equation~$\d\omega^\alpha_i
= -\omega^\alpha_C\w\omega^C_i$, using the congruences
$\omega^\beta_i\equiv0\mod\omega_0$ 
when~$i\le m$ and~$\beta>m{+}1$, 
and reducing modulo~$\omega_0$ yields $R^\alpha_i\,\d\omega_0\equiv 0
\mod\omega_0$, which implies~$R^\alpha_i = 0$.  

Expanding the structure equation for~$d\omega^\alpha_i$
when~$\alpha>m{+}1$ and~$i>m{-}d{+}1$ and reducing modulo~$\omega_0$
yields
$$
R^\alpha_i\,\d\omega_0 
   \equiv -\omega^\alpha_{m+1}\w\omega^{m+1}_i\mod\omega_0\,.
$$
Since the matrix~$R$ has rank at least~$2$ everywhere,
it follows that~$d\omega_0$ must be decomposable modulo~$\omega_0$.
Moreover, since $R^\alpha_j=0$ for~$j\le m{-}d{+}1$, it follows 
that~$d\omega_0\w\omega^{m+1}_i\w\omega_0 = 0$
for at least two distinct values of~$i\in\{m{-}d{+}2,\dots, m\}$.  
If~$R$ were to have rank equal to~$3$
at any point, then it would have rank~$3$ on a dense open set
and this would force $d\omega_0\w\omega^{m+1}_i\w\omega_0 = 0$
to hold for at least three distinct values 
of~$i\in\{m{-}d{+}2,\dots, m\}$, which would, in turn, force
$d\omega_0\w\omega_0 = 0$ to hold,
contradicting the nonintegrability of~$\omega_0$.  
Thus,~$R$ has rank equal to~$2$ at all points.

Because~$R$ has constant rank equal to~$2$, it is now
possible to define a sub-bundle of~$F^\diamond$ on which~$R$
is normalized to some particular normal form.  The choice of
this normal form is not important for the structure of the
argument, but a judicious choice (made with the desired
end result in mind, I must confess) that simplifies the notation
is to normalize so that~$\omega^{m+2}_m = -\omega^{m+3}_{m-1}$
and so that~$\omega^\alpha_i = 0$ for all pairs~$(i,\alpha)$
satisfying~$i\le m$ and~$\alpha>m{+1}$ except~$(i,\alpha) = 
(m,m{+}2)$ and~$(m{-}1,m{+}3)$.  The subset~$F'\subset F^\diamond$
on which this holds is easily seen 
to be a $G'$-subbundle over~$X^\diamond$
with a connected structure group~$G'\subset G$ that will be made
explicit later on in the argument when it will be useful to do so.
I will now use~$\omega_0$ to stand for~$\omega^{m+3}_{m-1}$.  

Now, I claim that~$d = 3$.  Suppose, instead that~$d>3$.
Then, $\omega^{m+1}_{m-2}\w\omega_0\not=0$ while 
$\omega^{m+2}_{m-2} = \omega^{m+3}_{m-2} = 0$.  Differentiating
these two equations and reducing modulo~$\omega_0$ then yields
$$
-\omega^{m+2}_{m+1}\w\omega^{m+1}_{m-2} 
\equiv -\omega^{m+3}_{m+1}\w\omega^{m+1}_{m-2}
\equiv 0 \mod\omega_0.
$$
Of course, this implies that there exist 
functions~$T^{m+2}$ and $T^{m+3}$
so that~$\omega^{m+2}_{m+1}
\equiv T^{m+2}\,\omega^{m+1}_{m-2} \mod\omega_0$
and~$\omega^{m+3}_{m+1}\equiv T^{m+3}\,\omega^{m+1}_{m-2}\mod\omega_0$.
Substituting these relations into the equations~
$R^\alpha_i\,\d\omega_0 
    \equiv -\omega^\alpha_{m+1}\w\omega^{m+1}_i\mod\omega_0$
for~$(i,\alpha) = (m,m{+}2)$ and~$(m{-}1,m{+}3)$, then yields
$$
\d\omega_0 \equiv -T^{m+3}\,\omega^{m+1}_{m-2}\w\omega^{m+1}_{m-1}
          \equiv -T^{m+2}\,\omega^{m+1}_{m-2}\w\omega^{m+1}_{m}
\mod \omega_0\,.
$$
Since~$\{\omega^{m+1}_{m-2},\omega^{m+1}_{m-1},
\omega^{m+1}_{m},\omega_0\}$
are linearly independent by hypothesis, this is impossible unless~$T^{m+2}
=T^{m+3} = 0$, but this vanishing would make~$d\omega_0$ integrable.  
Thus, $d=3$, as claimed.  

Since~$\omega^\alpha_{m-1} = \omega^\alpha_{m} = 0$ for~$\alpha>m{+}3$,
differentiating these equations and reducing modulo~$\omega_0$ yields
$$
\omega^\alpha_{m+1}\w \omega^{m+1}_{m-1} 
\equiv \omega^\alpha_{m+1}\w \omega^{m+1}_{m} 
\equiv 0 \mod\omega_0\,,
$$
from which it follows that $\omega^\alpha_{m+1}\equiv 0 \mod\omega_0$
for all~$\alpha>m{+}3$.  Thus, there exist~$T^\alpha$ for~$\alpha>m{+}3$
so that ~$\omega^\alpha_{m+1}=T^\alpha\,\omega_0$
for all~$\alpha>m{+}3$.

Differentiating the relations~$\omega^{m+2}_{m-1}=\omega^{m+3}_{m} = 0$
and reducing modulo~$\omega_0$ implies
$$
\omega^{m+2}_{m+1}\w\omega^{m+1}_{m-1}
\equiv \omega^{m+3}_{m+1}\w\omega^{m+1}_{m}
\equiv 0 \mod\omega_0,
$$
so there exist~$a$, $b$, ~$T^{m+2}$, and~$T^{m+3}$  
so that~$\omega^{m+2}_{m+1} =  a\,\omega^{m+1}_{m-1} + T^{m+2}\,\omega_0$ 
and $\omega^{m+3}_{m+1} = b\,\omega^{m+1}_{m} + T^{m+3}\,\omega_0$.  
Substituting this into
the derivatives of the equations~$\omega^{m+3}_{m-1} = \omega_0$
and~$\omega^{m+2}_{m} = -\omega_0$ and reducing modulo~$\omega_0$
yields
$$
\d\omega_0 \equiv -b\,\omega^{m+1}_m\w\omega^{m+1}_{m-1}
\equiv a\,\omega^{m+1}_{m-1}\w\omega^{m+1}_{m} \mod\omega_0\,,
$$
so it follows that~$a = b$.  The structure equations now imply
$$
\left.
\begin{aligned}
\d \vs_j&\equiv 0 \qquad\qquad\qquad \text{(when $j<m-1$)}\\
\d \vs_{m-1}&\equiv \phantom{ { } - { } } \vs_{m+3}\,\omega_0\\
\d \vs_{m\phantom{+0}}&\equiv { } - \vs_{m+2}\,\omega_0\\
\d \vs_{m+1}
&\equiv a\,(\vs_{m+2}\,\omega^{m+1}_{m-1}+\vs_{m+3}\,\omega^{m+1}_{m})
+ \sum_{\alpha=m+2}^n \vs_\alpha\, T^\alpha\,\omega_0
\end{aligned}\right\}
\mod \vs_1,\dots,\vs_{m+1}\,.
$$
Since the rank of~$dA$ is
greater than~$1$, it follows that~$a$ cannot vanish and, hence, that
the rank of~$dA$ is identically equal to~$3$.

To save writing, introduce the 
abbreviations~$\omega^{m+1}_{m-1} = \omega_1$ 
and~$\omega^{m+1}_m = \omega_2$.  Thus, for example, $\d\omega_0
\equiv a\,\omega_1\w\omega_2\mod\omega_0$.
Moreover,~$\omega^{m+2}_{m+1} \equiv a\,\omega_1\mod\omega_0$
and~$\omega^{m+3}_{m+1}\equiv a\,\omega_2\mod\omega_0$.  

I am now going to reduce to the case~$m=2$.
(If $m=2$ already, these next two paragraphs are unnecessary.)  
Fix~$i<m{-}1$.  
Differentiating the identities~$\omega^{m+2}_i = \omega^{m+3}_i = 0$
and using the structure equations yields
$$
\omega_0\w\omega^m_i - a S_i\,\omega_1\w\omega_0
 = -\omega_0\w\omega^{m-1}_i - a S_i\,\omega_2\w\omega_0 = 0,
$$
i.e., $\omega^m_i\equiv -aS_i\omega_1\mod\omega_0$ 
and~$\omega^{m-1}_i\equiv aS_i\omega_2\mod\omega_0$.
On the other hand, differentiating the 
equation~$\omega^{m+1}_i = S_i\,\omega_0$ and reducing modulo~$\omega_0$
yields
$$
-\omega_1\w\omega^{m-1}_i-\omega_2\w\omega^m_i
\equiv aS_i\,\omega_1\w\omega_2\mod\omega_0\,.
$$
In view of the congruences for~$\omega^m_i$ and~$\omega^{m-1}_i$,
this gives~$aS_i\,\omega_1\w\omega_2\w\omega_0 = 0$.  

Consequently, $S_i = 0$ for all~$i<m{-}1$.  
Moreover, there exist~$S^{m-1}_i$ and~$S^m_i$ 
so that~$\omega^{m-1}_i = S^{m-1}_i\,\omega_0$
and~$\omega^m_i = S^m_i\,\omega_0$ for~$i<m{-}1$.  However, I now
claim that~$S^{m-1}_i=S^m_i=0$.  This follows since, if I now
differentiate the relations $\omega^{m-1}_i = S^{m-1}_i\,\omega_0$
and~$\omega^m_i = S^m_i\,\omega_0$ and reduce modulo~$\omega_0$, the
result is
$$
0 \equiv S^m_i\,\d \omega_0 \equiv S^{m-1}_i\,\d\omega_0
\mod\omega_0\,,
$$
from which the claim follows, since~$\d\omega_0\not\equiv0\mod\omega_0$. 
This vanishing, in turn, now implies the congruences
$$
\d\vs_1\equiv\dots\equiv \d\vs_{m-2}\equiv 0 \mod \vs_1,\dots,\vs_{m-2}\,.
$$
In other words, the $(m{-}2)$-plane~$(P_V)^\perp\subset V$ is locally
constant on~$X^\diamond$ and hence, by connectedness constant.  It follows
that~$X^\diamond$ lies inside $[W_-,\C{n}]_m$ 
for some~$W_-\in\Gr(m{-}2,n)$.
Thus, it clearly suffices to take~$m=2$ for the rest of the
argument, so I will do so.

I am now going to reduce further to the case~$n = m{+}3 = 5$.  
(I remind the reader that~$m=2$ in what follows.  If~$n=5$ already,
then these next two paragraphs are unnecessary.)  This argument is
essentially the `dual' of the argument just given.  Differentiating
the relations~$\omega^\alpha_1 = \omega^\alpha_2 = 0$ for~$\alpha>5$ 
and using the structure equations yields
$$
-T^\alpha\,\omega_0\w\omega_1 -\omega^\alpha_5\w\omega_0
 = -T^\alpha\,\omega_0\w\omega_2 + \omega^\alpha_4\w\omega_0 = 0,
$$
so that $\omega^\alpha_4\equiv -T^\alpha\,\omega_2\mod\omega_0$
and~$\omega^\alpha_5\equiv T^\alpha\,\omega_1\mod\omega_0$.
Now differentiating the relation~$\omega^\alpha_3 = T^\alpha\,\omega_0$,
using the structure equations, and reducing modulo~$\omega_0$ yields
$$
T^\alpha\,\d\omega_0 
\equiv -\omega^\alpha_4\w\omega^4_3-\omega^\alpha_5\w\omega^5_3 
\mod\omega_0,.
$$
Using the known congruences, this yields
$$
aT^\alpha\,\omega_1\w\omega_2 
\equiv aT^\alpha\,\omega_2\w\omega_1 - aT^\alpha\,\omega_1\w\omega_2
\mod\omega_0\,,
$$
which, of course, implies that~$T^\alpha = 0$.  

Now that~$T^\alpha = 0$, 
it follows that~$\omega^\alpha_4= T^\alpha_4\,\omega_0$
and~$\omega^\alpha_5= T^\alpha_5\,\omega_0$ 
for some functions~$T^\alpha_4$
and~$T^\alpha_5$.  Again, differentiating these relations and then
reducing modulo~$\omega_0$  implies the relations
$$
0 \equiv T^\alpha_4\,\d\omega_0 \equiv T^\alpha_5\,\d\omega_0
\mod\omega_0\,,
$$
which, of course, implies that~$T^\alpha_4 = T^\alpha_5 = 0$ 
for all~$\alpha>5$.  
This vanishing, in turn, now implies the congruences
$$
\d\vs_1\equiv\dots\equiv \d\vs_{5}\equiv 0 \mod \vs_1,\dots,\vs_{5}\,.
$$
In other words, the $5$-plane spanned by~$\vs_1,\dots,\vs_5$ is locally,
and, hence, globally constant on~$X^\diamond$.  Let~$W_+\in\Gr(5,n)$
be this constant $5$-plane.  Then~$X^\diamond$ and, hence,~$X$ are
contained in~$\Gr(2,W_+)$.  Thus, as claimed, it suffices to take~$n=5$
for the rest of the argument, so I will do so.

At this point, I am going to switch over to the standard language of
the moving frame and assume that the reader is familiar with it.  
(The reader who is not might consult~\cite{MR81k:53004}.  Of course,
such a reader probably could not have followed the argument to this
point anyway.) 

It is a good idea to take stock of the problem.
At this moment, $X\subset\Gr(2,5)$ is an algebraic variety of
dimension~$3$ that contains a Zariski-open subset~$X^\diamond$ 
over which there exists
a `moving frame'~$\vs: X^\diamond\to\SL(5,\bbC)$
satisfying the condition that~$V\in X^\diamond$ is spanned by
~$\vs_1(V)$ and~$\vs_2(V)$ as well as the structure equations
\begin{align*}
d\vs_1 &\equiv \vs_3\,\omega_1 + \vs_5\,\omega_0 &&\mod\vs_1,\vs_2\,,\\
d\vs_2 &\equiv \vs_3\,\omega_2 - \vs_4\,\omega_0 &&\mod\vs_1,\vs_2\,,\\
\noalign{\vskip3pt}
d\vs_3 &\equiv 
          \vs_4\,(a\,\omega_1+T^4\,\omega_0) 
        + \vs_5\,(a\,\omega_2+T^5\,\omega_0)&&\mod\vs_1,\vs_2,\vs_3\,.
\end{align*}
where~$\omega_0,\omega_1,\omega_2$ are linearly independent 
and~$a\not=0$.  

My goal now is to show that the existence of 
such a frame field implies that~$X^\diamond$ is an open subset
of the isotropic Grassmannian associated to some nondegenerate quadratic
form on~$\C{5}$.  This will have to be done in a series of steps.%
\footnote{The reason that the following argument is somewhat complicated
can be seen as follows:  
It is not difficult to see that the condition
that~$X^\diamond\subset\Gr(2,5)$ 
have such a coframing constitutes
a set of four first-order PDE for~$X^\diamond$ 
as a submanifold of~$\Gr(2,5)$
(plus an open condition on the second 
derivatives to ensure that~$a\not=0$).
Since the codimension of~$X^\diamond$ in~$\Gr(2,5)$ is~$3$, 
this means that this system of equations is overdetermined, but by
only one equation.  
Not surprisingly, this system is not involutive.  
Moreover, it only goes 
into involution after several cycles of 
prolongation and torsion reduction.
Any proof of the claimed rigidity that works locally will have to 
reproduce this calculation  in some form, so it cannot be too simple. 
The proof in the text has been designed to get to the classification
as quickly as possible, and explicit discussion of the exterior
differential systems analysis that inspired it has been suppressed.
I apologize if this makes the proof seem unmotivated.}
First, to simplify the argument, note that it suffices to prove this
in the case where~$X^\diamond$ is a non-singular, connected, and simply 
connected $3$-dimensional complex submanifold 
in~$\Gr(2,5)$ that possesses
such a frame field, so assume that~$X^\diamond$ has these properties.

Replacing~$\vs_3$ by~$\vs_3 - T^5\,\vs_1 + T^4\,\vs_2$
yields a new frame for which the corresponding~$T^4$ and~$T^5$ are
zero.   Thus, without loss of generality, I can assume that~$T^4=T^5=0$.
Using the simple-connectivity of~$X^\diamond$, write~$a = t^5$
for some function~$t$ on~$X^\diamond$.
Replacing the given frame~$\vs = (\vs_1,\dots,\vs_5)$
by~$(t^{-2}\vs_1,t^{-2}\vs_2,t^{-2}\vs_3,t^{3}\vs_4,t^{3}\vs_5)$ 
yields a new unimodular frame for which the corresponding~$a$ 
is equal to~$1$.  Thus, again, without loss of generality, 
I can further assume that~$a=1$.

Thus, I will say that a frame field~$\vs:X^\diamond\to\SL(5,\bbC)$ is 
\emph{$0$-adapted} to~$X^\diamond$ if
the map~$[\vs_1\w\vs_2]:X^\diamond\to\Gr(2,5)$
is the inclusion~$X^\diamond\hookrightarrow\Gr(2,5)$ 
and, moreover, $\vs$ satisfies
\begin{equation}
\begin{aligned}
\d\vs_1 &\equiv \vs_3\,\omega_1+\vs_5\,\omega_0 &&\mod\vs_1,\vs_2\,,\\
\d\vs_2 &\equiv \vs_3\,\omega_2-\vs_4\,\omega_0 &&\mod\vs_1,\vs_2\,,\\
\noalign{\vskip3pt}
\d\vs_3 &\equiv \vs_4\,\omega_1+\vs_5\,\omega_2 &&\mod\vs_1,\vs_2,\vs_3\,.
\end{aligned}
\end{equation}
for some independent $1$-forms~$\omega_0,\omega_1,\omega_2$ 
on~$X^\diamond$.  

By standard methods in the theory of moving frames, 
one sees that the $0$-adapted frame fields over~$X^\diamond$
are the sections of a principal $G_0$-bundle~$F_0\subset 
X^\diamond\times\SL(5,\bbC)$ where~$G_0\subset\SL(5,\bbC)$ is a
$10$-dimensional Lie subgroup whose 
Lie algebra~$\eug_0\subset\eusl(5,\bbC)$
is the set of matrices of the form
$$
\begin{pmatrix}
x^1_1&x^1_2&x^1_3& x^1_4& x^1_5\\
x^2_1&x^2_2&x^2_3& x^2_4& x^2_5\\
  0  &  0  &  0  & x^1_3& x^2_3\\
  0  &  0  &  0  &-x^1_1&-x^2_1\\
  0  &  0  &  0  &-x^1_2&-x^2_2 
\end{pmatrix}.
$$
The group~$G_0$ is not connected because of the usual complication
caused by the fact that~$\SL(5,\bbC)$ has a nontrivial, finite center.
Instead, $G_0$ is equal to the product of its identity component
(which is determined by the Lie algebra~$\eug_0$) and the elements
of the form~$\varepsilon\I_5$ where~$\varepsilon^5=1$.  Thus,~$G_0$
has five components.  It follows either that~$F_0$ is connected
or else that it has $5$ components.

Following the standard method of the moving frame, for any $0$-adapted
frame field~$\vs:X^\diamond\to\SL(5,\bbC)$,
let~$\d\vs_a = \vs_b\,\omega^b_a$, where~$\omega^a_a = 0$
and $\d\omega^a_b = -\omega^a_c\w\omega^c_b$.  By hypothesis, the
relations
\begin{equation}\label{eq: Gr25 0-frame relations}
\omega^4_1 = \omega^5_2 = \omega^5_1+\omega^4_2 = \omega^4_3-\omega^3_1
= \omega^5_3 - \omega^3_2 = 0
\end{equation}
hold.  (I will continue to use the abbreviations $\omega^5_1 = \omega_0$,
$\omega^3_1 = \omega_1$, and $\omega^3_2 = \omega_2$.)   Taking
the exterior derivatives of these relations,  applying the
structure equations, and then collecting terms and applying Cartan's
Lemma shows that there exist functions~$t^a_b$ so that
$$
\begin{aligned}
\omega^4_4 &= -\omega^1_1 + t^4_4\,\omega_0\,,\\
\omega^4_5 &= -\omega^2_1 + t^4_5\,\omega_0\,,\\
\omega^5_4 &= -\omega^1_2 + t^5_4\,\omega_0\,,\\
\omega^5_5 &= -\omega^2_2 + t^5_5\,\omega_0\,,
\end{aligned}
\qquad\text{and}\qquad
\begin{aligned}
\omega^3_4 &= \omega^1_3 + t^5_4\,\omega_1 + (3t^5_5+2t^4_4)\,\omega_2 
                         + t^3_4\,\omega_0\,,\\
\omega^3_5 &= \omega^2_3 - t^4_5\,\omega_2 - (2t^5_5+3t^4_4)\,\omega_1 
                         + t^3_5\,\omega_0\,.
\end{aligned}
$$
In particular~$\omega^3_3 = -(t^4_4+t^5_5)\,\omega_0$.  
Since~$d\omega_0\equiv\omega_1\w\omega_2\mod\omega_0$, applying the
structure equation for $\d\omega^3_3=-\omega^3_a\w\omega^a_3$ and
using the above relations yields 
$$
-(t^4_4+t^5_5)\,\d\omega_0\equiv 5(t^4_4+t^5_5)\,\omega_1\w\omega_2
\mod\omega_0
$$
which implies~$t^4_4+t^5_5=0$.  In particular~$\omega^3_3=0$, so
going back to its structure equation yields
$$
0 = \d\omega^3_3 = -\omega^3_a\w\omega^a_3 
= (t^3_4\,\omega_1+t^3_5\,\omega_2)\w\omega_0\,,
$$
so it follows that~$t^3_4 = t^3_5=0$ also.  

Now, computing how the~$t^a_b$ vary under a change of $0$-adapted
frame (a detail that can be safely left to the reader), one sees that
by adding the appropriate multiples of~$\vs_1$ and~$\vs_2$ to~$\vs_4$
and~$\vs_5$, one can construct a $0$-adapted frame for 
which~$t^4_4 = t^5_4 = t^4_5 = t^5_5 = 0$.  
I will say that a $0$-adapted
frame that satisfies this additional property is \emph{$1$-adapted}.

Again, the usual methods show that the $1$-adapted frame fields are the
sections of a principal~$G_1$-bundle~$F_1\subset F_0$ 
where~$G_1\subset G_0$
is the $7$-dimensional Lie subgroup whose 
Lie algebra~$\eug_1\subset\eug_0$
is the space of matrices of the form
$$
\begin{pmatrix}
x^1_1&x^1_2&x^1_3&   0  & x^1_5\\
x^2_1&x^2_2&x^2_3&-x^1_5&   0  \\
  0  &  0  &  0  & x^1_3& x^2_3\\
  0  &  0  &  0  &-x^1_1&-x^2_1\\
  0  &  0  &  0  &-x^1_2&-x^2_2
\end{pmatrix}
$$
and that is generated by its identity component together with
the elements of the form~$\varepsilon\I_5$ where~$\varepsilon^5=1$.
(Thus,~$G_1$, like $G_0$, has five components.)

For a $1$-adapted coframe field, in addition to the relations
\eqref{eq: Gr25 0-frame relations}, there are now relations
\begin{equation}\label{eq: Gr25 1-frame relations}
\omega^4_4+\omega^1_1 = \omega^4_5+\omega^2_1 
= \omega^5_4+\omega^1_2 = \omega^5_5+\omega^2_2
= \omega^3_4-\omega^1_3 = \omega^3_5 - \omega^2_3 = 0.
\end{equation}
Taking the exterior derivatives of these six relations,  applying the
structure equations, and then collecting terms and applying Cartan's
Lemma (keeping in mind that~$\omega_0\w\omega_1\w\omega_2\not=0$) 
implies the further relations
\begin{equation}\label{eq: Gr25 final relations}
\omega^1_4 = \omega^2_5 = \omega^1_5+\omega^2_4 = 0.
\end{equation}
The relations \eqref{eq: Gr25 0-frame relations},
\eqref{eq: Gr25 1-frame relations}, and \eqref{eq: Gr25 final relations}
combine to imply that the matrix~$\omega= (\omega^a_b)$
satisfies~${}^t(Q\omega) = -Q\omega$, where
$$
Q = {}^tQ = 
\begin{pmatrix}
  0  &  0  &  0  &  1  &  0  \\
  0  &  0  &  0  &  0  &  1  \\
  0  &  0  & -1  &  0  &  0  \\
  1  &  0  &  0  &  0  &  0  \\
  0  &  1  &  0  &  0  &  0  
\end{pmatrix}.
$$ 
I.e., $\omega$ takes values
in the subspace~$\euso(Q)\subset\eusl(5,\bbC)$ that
is the Lie algebra of the group~$\SO(Q)\subset\SL(5,\bbC)$
of matrices~$A$ that satisfy~${}^t\!A\,Q\,A = Q$ and~$\det(A)=1$.
Of course, $\SO(Q)$ is isomorphic to~$\SO(5,\bbC)$.

Since the~$\omega^a_b$ are linearly independent except for
the relations \eqref{eq: Gr25 0-frame relations},
\eqref{eq: Gr25 1-frame relations}, and \eqref{eq: Gr25 final relations},
it follows that the projection~$\vs:F_1\to\SL(5,\bbC)$ 
immerses each component of~$F_1$ into a single left coset 
of~$\SO(Q)$. Moreover, each component of~$\vs(F_1)$ is open in such
a coset.  Since the~$\bbZ_5$ that forms the center of~$\SL(5,\bbC)$
does not lie in~$\SO(Q)$ but does lie in~$G_1$, it follows that 
the image of~$F_1$ actually maps into five distinct left cosets 
of~$\SO(Q)$ and so must consist of five components instead of one.
In particular, the inverse image of one of these 
components is a component~$F^\circ_1\subset F_1$ that is 
a $G^\circ_1$-bundle over~$X^\circ$.  I now restrict all forms and
functions to this component~$F^\circ_1$.

Since~$\vs:F^\circ_1\to\SL(5,\bbC)$ is an open immersion into a
single left coset of~$\SO(Q)$, it now follows that there 
exists a (unique) non-degenerate inner 
product~$\la,\ra$ on~$\C5$ with the property that~ 
$\la \vs_a,\vs_b\ra = Q_{ab}$.  Thus, the $2$-plane
$[\vs_1(x)\w\vs_2(x)] = x\in X^\diamond$ is $\la,\ra$-isotropic 
for all~$x\in X^\diamond$.  Since~$X^\diamond$ and the $\la,\ra$-isotropic
Grassmannian are both $3$-dimensional submanifolds of~$\Gr(2,5)$, it
follows that~$X^\diamond$ is an open subset of the $\la,\ra$-isotropic
Grassmannian, as was to be proved.

As already mentioned, this implies that any $3$-dimensional irreducible 
algebraic variety~$X\subset\Gr(2,5)$ 
that contains a Zariski-open subset~$X^\diamond$
that supports a $0$-adapted frame field must actually be the 
$\la,\ra$-isotropic Grassmannian for some non-degenerate inner product
$\la,\ra$ on~$\C5$.  Thus, at last, this subcase is finished; 
such varieties fall into the third category of the proposition.

\emph{Finally}, all that remains is to address the last subcase, 
that of an
irreducible $3$-dimensional subvariety~$X\subset\Gr(m,n)$ that contains
a Zariski-open subset~$X^\diamond$ that is smooth and whose tangent
spaces are integral elements of the third
type listed in Lemma~\ref{lem: integral elements of I3}.  My goal is
to prove that such an~$X$ is necessarily a component of the variety
listed in the fourth category of the proposition.

The first task is to reduce to the case~$(m,n) = (3,6)$.  This
will be reminiscent of the previous argument's reduction to~$\Gr(2,5)$.
Let~$F^\diamond\subset X^\diamond\times\SL(n,\bbC)$
consist of the pairs~$(V,\vs)$  that satisfy the following conditions:
\begin{enumerate}
\item $\vs_1,\dots,\vs_m$ spans~$V$, and
\item the tangent space~$T_VX^\diamond\subset Q_V\ot V^*$ is spanned
by the three vectors
$$
\begin{matrix}
 \dl\vs_{m+2}\dr{\ot}\vs^{m-2}-\dl\vs_{m+1}\dr{\ot}\vs^{m-1}\,,\\ 
\noalign{\vskip2pt}
 \dl\vs_{m+1}\dr{\ot}\vs^{m\phantom{-0}}-\dl\vs_{m+3}\dr{\ot}\vs^{m-2}\,,\\
\noalign{\vskip2pt}
 \dl\vs_{m+3}\dr{\ot}\vs^{m-1}-\dl\vs_{m+2}\dr{\ot}\vs^{m\phantom{-0}}\,.
\end{matrix}
$$
\end{enumerate}
(The indexing is a little unfortunate, but, for consistency with
the conventions I have used so far, it is unavoidable.)  
Now, $F^\diamond\to X^\diamond$ is a submersion and is a 
principal~$G$-bundle where~$G\subset\SL(n,\bbC)$ is a closed subgroup 
of $P_{m-3}\cap P_m\cap P_{m+3}$ of codimension~$9$.  I will not need
the full definition of~$G$ right now, so I postpone this.

As usual, the structure equations hold on~$F^\diamond$.  By construction,
$\omega^\alpha_i = 0$ if either~$i<m{-}2$ or $\alpha>m{+}3$ 
while the matrix
$$
\begin{pmatrix}
\omega^{m+1}_{m-2}&\omega^{m+1}_{m-1}&\omega^{m+1}_{m\phantom{-0}}\\
\noalign{\vskip2pt}
\omega^{m+2}_{m-2}&\omega^{m+2}_{m-1}&\omega^{m+2}_{m\phantom{-0}}\\
\noalign{\vskip2pt}
\omega^{m+3}_{m-2}&\omega^{m+3}_{m-1}&\omega^{m+3}_{m\phantom{-0}}
\end{pmatrix}
$$
is skew-symmetric.  
Moreover, introducing $1$-forms~$\omega_1$, $\omega_2$,
and~$\omega_3$ by the equations
$$
\begin{pmatrix}
\omega^{m+1}_{m-2}&\omega^{m+1}_{m-1}&\omega^{m+1}_{m\phantom{-0}}\\
\noalign{\vskip2pt}
\omega^{m+2}_{m-2}&\omega^{m+2}_{m-1}&\omega^{m+2}_{m\phantom{-0}}\\
\noalign{\vskip2pt}
\omega^{m+3}_{m-2}&\omega^{m+3}_{m-1}&\omega^{m+3}_{m\phantom{-0}}
\end{pmatrix}
= 
\begin{pmatrix}
\>0&-\omega_3&\phantom{-}\omega_2\\
\phantom{-}\omega_3&\>0&-\omega_1\\
-\omega_2&\phantom{-}\omega_1&\>0
\end{pmatrix},
$$
one has~$\omega_1\w\omega_2\w\omega_3\not=0$.  

Suppose~$m>3$ and fix~$i<m{-}2$.  Differentiating the 
equations~$\omega^{m+1}_i=\omega^{m+2}_i=\omega^{m+3}_i=0$ and
using the structure equations yields the relations
$$
\begin{pmatrix}
\>0&-\omega_3&\phantom{-}\omega_2\\
\phantom{-}\omega_3&\>0&-\omega_1\\
-\omega_2&\phantom{-}\omega_1&\>0
\end{pmatrix}
\w
\begin{pmatrix}
\omega^{m-2}_i\\ 
\noalign{\vskip2pt}
\omega^{m-1}_i\\ 
\noalign{\vskip2pt}
\omega^{m\phantom{-0}}_i
\end{pmatrix} = 0.
$$
The linear independence of $\omega_1,\omega_2,\omega_3$ then
implies that $\omega^{m-2}_i = \omega^{m-1}_i = \omega^{m}_i = 0$.
This vanishing implies the relations
$$
\d\vs_1\equiv\dots\equiv \d\vs_{m-3}\equiv 0
\mod \vs_1,\dots,\vs_{m-3}\,.
$$
In other words, the $(m{-}3)$-plane spanned by~$\vs_1,\dots,\vs_{m-3}$
is locally constant.  Set $W_- = [\vs_1\w\dots\w\vs_{m-3}]$
for some fixed~$W_-\in\Gr(m{-}3,n)$.  
Since~$X^\diamond$ is connected, it follows that $X^\diamond$ 
(and hence~$X$) must lie 
in~$[W_-,\C{n}]_m\simeq\Gr\bigl(3,\C{n}/W_-\bigr)$. 

Thus, it suffices to analyze the case~$m=3$, 
so I will assume this from now on.

Suppose $n>6$ (remember that~$m=3$ now) and fix~$\alpha>6$.  
Differentiating the 
equations~$\omega^\alpha_1=\omega^\alpha_2=\omega^\alpha_3=0$ 
and using the structure equations yields the relations
$$
\begin{pmatrix} \omega^\alpha_4&\omega^\alpha_5&\omega^\alpha_6
\end{pmatrix}
\w
\begin{pmatrix}
\>0&-\omega_3&\phantom{-}\omega_2\\
\phantom{-}\omega_3&\>0&-\omega_1\\
-\omega_2&\phantom{-}\omega_1&\>0
\end{pmatrix} = 0.
$$
Again, the linear independence of $\omega_1,\omega_2,\omega_3$ then
implies that $\omega^\alpha_4 = \omega^\alpha_5 = \omega^\alpha_6= 0$.
This vanishing implies the relations
$$
\d\vs_1\equiv\dots\equiv \d\vs_{6}\equiv 0
\mod \vs_1,\dots,\vs_{6}\,.
$$
In other words, the $6$-plane spanned by~$\vs_1,\dots,\vs_{6}$
is locally constant.  Set $W_+ = [\vs_1\w\dots\w\vs_{6}]$
for some fixed~$W_+\in\Gr(6,n)$.  
Since~$X^\diamond$ is connected, it follows that $X^\diamond$ 
(and hence~$X$) must lie in~$\Gr\bigl(3,W_+\bigr)$.  
Thus, it suffices to analyze the case~$n=6$, 
so I will assume this from now on.

At this point, it is worthwhile to make the group~$G\subset\SL(6,\bbC)$
explicit.  The usual calculation shows that this is a Lie subgroup of
matrices whose Lie algebra is the space~$\eug_0\subset\eusl(6,\bbC)$ 
of matrices of the form
$$
\begin{pmatrix} a & b \\ 0 & -{}^ta \end{pmatrix}
$$
where~$a$ and~$b$ are arbitrary $3$-by-$3$ matrices.  The group~$G$
is not connected, but is generated by its identity component and
the center of~$\SL(6,\bbC)$, a cyclic group of order~$6$ that 
consists of the matrices of the form~$\varepsilon\I_6$
where~$\varepsilon^6=1$.  Consequently, $G$ actually has $3$
components (the $\bbZ_2$-subgroup~$\left\{\pm\I_6\right\}$
already lies in the identity component of~$G$).

From this point on, the argument is much like the argument for the
integral manifolds  of the third category, so I will just indicate
the steps without explicitly writing out the details.  

The first step is to note that the following six relations hold
on~$F^\diamond$:
\begin{equation}\label{eq: Gr36 0-frame equations}
\omega^{i+3}_j + \omega^{j+3}_i = 0.
\end{equation}
Differentiating these relations, applying the structure equations, 
and applying Cartan's Lemma 
shows that there exist functions~$s^{ij}=s^{ji}$
on~$F^\diamond$ so that the relations
$$
\omega^{i+3}_{j+3} = -\omega^j_i + s^{jk}\,\omega_l -  s^{jl}\,\omega_k
$$
hold for any~$1\le i,j\le 3$ 
and where $(i,k,l)$ is an even permutation 
of~$(1,2,3)$.  (It is important to remember 
that~$\omega^1_1+\dots+\omega^6_6 = 0$, 
since this relation figures into
these calculations.)

The six equations~$s^{ij}=0$ define a principal~$G_1$ 
subbundle~$F_1\subset F^\diamond$ where~$G_1\subset G$ is the
subgroup whose Lie algebra~$\eug_1$ consists of the matrices
of the form
$$
\begin{pmatrix} a & b \\ 0 & -{}^ta \end{pmatrix}
$$
where~$a$ and~$b$ are $3$-by-$3$ matrices with~$b = -{}^tb$.  
Again,~$G_1$ is generated by its identity component and the (finite)
center of~$\SL(6,\bbC)$ and so has three components.

On~$F_1$, in addition to the six equations~\eqref{eq: Gr36 0-frame
equations},
the nine equations
\begin{equation}\label{eq: Gr36 1-frame equations}
\omega^{i+3}_{j+3} +\omega^j_i= 0
\end{equation}
also hold for all~$1\le i,j\le 3$. 
Differentiating these relations, applying the structure equations, 
and applying Cartan's Lemma shows that 
\begin{equation}\label{eq: Gr36 final equations}
\omega^{i}_{j+3} +\omega^j_{i+3}= 0.
\end{equation}   

In view of the relations~\eqref{eq: Gr36 0-frame equations}, 
\eqref{eq: Gr36 1-frame equations}, 
and \eqref{eq: Gr36 final equations}, 
it follows that $\omega = (\omega^a_b)$ 
satisfies~${}^t(Q\omega) = -Q\omega$
where
$$
Q = {}^tQ = \begin{pmatrix} 0_3 & \I_3\\ \I_3&0_3 \end{pmatrix}.
$$
The~$\omega^a_b$ are otherwise linearly independent,
so the map~$\vs:F_1\to\SL(6,\bbC)$ immerses 
each component of~$F_1$ as an open subset of a left coset 
of
$$
\SO(Q)=\left\{\,A\in\SL(6,\bbC) \mid {}^t\!A\,Q\,A = Q\,\right\},
$$
a subgroup isomorphic to~$\SO(6,\bbC)$.  Since~$\SO(Q)$ does not
contain the full center of~$\SL(6,\bbC)$ while~$G_1$ does, it follows
that the image~$\vs(F_1)$ lies in three distinct left cosets of~$\SO(Q)$
and that~$F_1$ must therefore consist of three distinct components.
Let~$F^\circ_1\subset F_1$ be one of these three components and
restrict all forms and functions to~$F^\circ_1$ henceforth.

As in the argument for the third case, it follows that there exists
a nondegenerate quadratic form~$\la,\ra$ on~$\C{6}$ 
with the property that~$X^\diamond\subset\Gr(3,6)$
lies in the $3$-dimensional submanifold 
of $\la,\ra$-isotropic $3$-planes
in~$\Gr(3,6)$.  Thus~$X^\diamond$ (and hence~$X$) must be an open 
subset of one of the two components of this isotropic Grassmannian.
Since~$X$ was assumed to be irreducible and algebraic, 
it follows that~$X$ must actually be one of these components, i.e., 
it must belong to the fourth category of the proposition.

At last, the proof of the proposition is complete. 
\end{proof}

\begin{remark}[A case of $A$-rigidity]
\label{rem: a case of A-rigidity}
The alert reader may have noticed that 
Proposition~\ref{prop: integral varieties of I3} contains
a proof of the rigidity of~$A$-cycles for a certain 
$3$-dimensional subspace~$A$ of~$T_V\Gr(m,n)$.  

In fact, what the analysis of the last case shows is that
if~$A$ is the third type of integral element listed 
in Lemma~\ref{lem: integral elements of I3}, 
then any irreducible $3$-dimensional 
subvariety~$W\subset\Gr(m,n)$ whose tangent spaces at all
smooth points are of type~$A$ is an open subset of 
one of the components of~$N_G(W_-,W_+)\simeq\bbP^3\cup\bbP^3$
as described in the fourth case 
of Proposition~\ref{prop: integral varieties of I3}.

One interesting consequence of this analysis is that~$\sigma_A$ 
(the only irreducible $A$-cycle) is nonsingular and, 
in fact, homogeneous.

Another interesting feature of this proof is that it can serve
as an example of how the moving frame approach can be used to 
prove the sort of higher order rigidity results that were mentioned 
in Remark~\ref{rem: higher order rigidity}, 
particularly Example~\ref{ex: 2nd order rigidity}.  
Note the pattern of the proof:  
\begin{enumerate}
\item Use the hypothesis that the tangent space 
to the submanifold has type~$A$ to derive 
equations~\eqref{eq: Gr36 0-frame equations}.
\item Differentiate equations~\eqref{eq: Gr36 0-frame equations}
to derive the conditions for second order osculation and use
them to make a second order frame adaptation to arrive at
equations~\eqref{eq: Gr36 1-frame equations}.  (This essentially
amounts to defining the second order Gauss mapping.)
\item Differentiate equations~\eqref{eq: Gr36 1-frame equations}
to derive the conditions for third order osculation and conclude
that this third (and higher) order osculation is automatic, which
is equivalent to~\eqref{eq: Gr36 final equations}. (This essentially
amounts to showing that the second order Gauss mapping is constant.)
\end{enumerate}
Essentially this same pattern is repeated in the proofs of the
claims of Example~\ref{ex: 2nd order rigidity}.
\end{remark}

Fortunately, the long proof 
of Proposition~\ref{prop: integral varieties of I3}
pays off double.  To avoid triviality, assume that~$m\ge 3$.

\begin{proposition}\label{prop: integral varieties of I111}
The following are integral varieties of~$\cI_{(1,1,1)}$ in~$\Gr(m,n)$:
\begin{enumerate}
\item For any~$A\in\Gr(m{-}2,\C{n})$, 
      the $2(n{-}m)$-dimensional 
      submanifold
$$
[A,\C{n}]_m\simeq\Gr\bigl(2,\C{n}/A\bigr).
$$
\item For any curve~$C\subset\Gr(m{-}1,n)$, the $(n{-}m{+}1)$-dimensional 
      subvariety
$$
\Psi_m(C) = \cup_{B\in C}[B,\C{n}]_m.
$$
\item For any pair~$W_-\in\Gr(m{-}3,n)$ and~$W_+\in\Gr(m{+}2,n)$ with
      $W_-\subset W_+$ and any nondegenerate quadratic form
      $G$ on~$W_+/W_-\simeq\C{5}$, the $3$-dimensional submanifold
      $N^\perp_G(W_-,W_+)\subset\Gr(m,n)$ that consists of 
      the $m$-planes~$W$ in $[W_-,W_+]_m$ for which~$(W/W_-)^\perp$
      is $G$-isotropic in~$W_+/W_-$.
\item For any pair~$W_-\in \Gr(m{-}3,n)$ and~$W_+\in\Gr(m{+}3,n)$ with
      $W_-\subset W_+$ and any nondegenerate quadratic form
      $G$ on~$W_+/W_-\simeq\C{6}$, the $3$-dimensional submanifold
      $N_G(W_-,W_+)\subset\Gr(m,n)$ that consists of 
      the $m$-planes~$W$ in $[W_-,W_+]_m$ for which~$W/W_-$
      is $G$-isotropic in~$W_+/W_-$.
\item Any subvariety~$V\subset\Gr(m,n)$ of dimension at most~$2$.
\end{enumerate}
Moreover, any irreducible algebraic integral variety 
of~$\cI_{(1,1,1)}$ is a subvariety of an algebraic integral variety of
one of the five listed types.
\end{proposition}

\begin{proof}
Combine complementarity and 
Proposition~\ref{prop: integral varieties of I3}.
\end{proof}

Propositions~\ref{prop: integral varieties of I21}, 
\ref{prop: integral varieties of I3}, 
and \ref{prop: integral varieties of I111}
combine to provide an effective means of analyzing the irreducible
varieties in~$\Gr(m,n)$ whose homology classes are linear combinations
of the~$[\sigma_{\ab^*}]$ where $\ab$ satisfies~$a_1\le 2$ (i.e., the
integral varieties of~$\cI_{(3)}$)
or the irreducible varieties in~$\Gr(m,n)$ whose homology classes 
are linear combinations of the~$[\sigma_{\ab^*}]$ where $\ab$ 
satisfies~$a_3=0$ (i.e., the integral varieties of~$\cI_{(1,1,1)}$).  
In practice, though, there are combinatorial difficulties 
when the dimension of such a cycle is such that there are many possible
choices for~$\ab$.  On the other hand, at the extremes, the descriptions
are fairly simple:

\begin{theorem}
\label{thm: rigid (n-m,n-m)-star and (2,...,2)-star}
Assume~$2\le m\le n{-}2$.  Then any $X\in\cZ_{2(n-m)}\bigl(\Gr(m,n)\bigr)$ 
that satisfies~$[X] = r[\sigma_{(n-m,n-m)^*}]$ is of the form
$$
X = [B_1,\C{n}]_m+\dots+[B_r,\C{n}]_m
$$
for some~$B_1,\dots,B_r\in\Gr(m{-}2,n)$.
Moreover, if $\ab = (2,2,\dots,2)$ has~$|\ab|=2m$ 
(i.e., the length of~$\ab$ is~$m$), 
then any~$X\in\cZ_{m(n-m-2)}\bigl(\Gr(m,n)\bigr)$ 
that satisfies~$[X]=r[\sigma_{\ab^*}]$ is of the 
form
$$
X = \Gr(m,B_1)+\dots+ \Gr(m,B_r)
$$
for some~$B_1,\dots,B_r\in\Gr(m+2,n)$.
\end{theorem}

\begin{proof}
Any variety $X\subset\Gr(m,n)$ of pure dimension~$2(n{-}m)$
that satisfies~$[X] = r[\sigma_{(n-m,n-m)^*}]$ is necessarily 
an integral variety of~$\cI_{(1,1,1)}$.  
Proposition~\ref{prop: integral varieties of I111} implies that any
such irreducible variety must, for dimension reasons, 
fall into the first category listed there.  
Thus, if~$X$ is irreducible, then $X = [B,\C{n}]_m$ 
for some~$B\in\Gr(m{-}2,n)$.  Since the ray generated 
by~$[\sigma_{(n-m,n-m)^*}]$ is extremal, 
this implies the first rigidity statement of the theorem.

Similarly if~$\ab = (2,2,\dots,2)$ has~$|\ab|=2m$ 
(i.e., the length of~$\ab$ is~$m$), 
then any variety~$X\in\Gr(m,n)$ of pure dimension~$2m$ 
that satisfies~$[X]=r[\sigma_{\ab^*}]$ 
is an integral variety of ~$\cI_{(3)}$. 
Proposition~\ref{prop: integral varieties of I3} implies that any
such irreducible variety of dimension~$2m$ must fall into the 
first category listed there.
Thus, if~$X$ is irreducible, then~$X=\Gr(m,B)$ 
for some~$B\in\Gr(m+2,n)$.  Since the ray generated 
by~$[\sigma_{\ab^*}]$ is extremal, 
this implies the second rigidity statement of the theorem.
\end{proof}

Before stating the next theorem, I need to introduce some
constructions of certain $3$-folds in~$\Gr(m,n)$ that are
based on curves.

\begin{example}[A chordal $3$-fold]
\label{ex: chordal 3-fold}
This construction depends on a pair of curves.
Let~$W\in\Gr(m,n)$ be fixed, let 
$$
\alpha\subset\bigl[W,\C{n}\bigr]_{m+1}
\ \bigl(\ \simeq\bbP(Q_W)\simeq\bbP^{n-m-1}\bigr)
$$
be an irreducible curve of degree~$a$, and let~
$$
\beta\subset\Gr(m{-}1,W)\ \bigl(\ \simeq\bbP(W^*)\simeq\bbP^{m-1}\bigr)
$$
be an irreducible curve of degree~$b$.  
Define~$\Sigma_m(\beta,\alpha)\subset\Gr(m,n)$ 
to be the union of the lines~$[B,A]_m\subset\Gr(m,n)$ 
with~$B\in\beta$ and~$A\in\alpha$.  Thus,~$\Sigma_m(\beta,\alpha)$
is the image of a $\bbP^1$-bundle over the surface~$\beta\times\alpha$ 
and so has dimension~$3$.  

Straightforward calculation shows that
$$
\bigl[\Sigma_m(\beta,\alpha)\bigr] = ab\,[\sigma_{(2,1)^*}].
$$
When~$a = b = 1$, the curves~$\alpha$ and~$\beta$ are, of course,
of the form~$\alpha = [W,W_+]_{m+1}$ 
and~$\beta = [W_-,W]_{m-1}$ for some~$W_+\in[W,\C{n}]_{m+2}$ 
and~$W_-\in\Gr(m{-}2,W)$.  
Thus, the cycle~$\Sigma_m(\beta,\alpha)$ is a
Schubert cycle~$\sigma_{(2,1)^*}$ in some~$[W_-,W_+]_m\simeq\Gr(2,4)$.  
However, when~$ab>1$, the variety~$\Sigma_m(\beta,\alpha)$ does not
lie in any such~$[W_-,W_+]_m$.  

Note that all of the lines~$[B,A]_m\subset\Sigma_m(\beta,\alpha)$ 
pass through~$W$, which is thus a singular point 
of~$\Sigma_m(\beta,\alpha)$.  Thus, none
of these `chordal $3$-folds' are smooth.
\end{example}

\begin{example}[Suspension and extension $3$-folds]
\label{ex: suspension and extension 3-folds}
The next two constructions are `complementary' to one
another in the sense of~\S\ref{sssec:complementarity}.

For the first construction, fix a plane~$W_-\in\Gr(m{-}2,n)$
and an algebraic curve~$\alpha\subset[W_-,\C{n}]_{m+1}\simeq
\Gr(3,\C{n}/W_-)$.  Let
\begin{equation}
\label{eq: sigma-m alpha}
\Sigma_m(W_-,\alpha) =\bigcup_{E\in \alpha} [W_-,E]_m\,.
\end{equation}
Thus,~$\Sigma_m(W_-,\alpha)$ is a curve of~$\bbP^2$s.  It
will be called the \emph{suspension} of~$\alpha$ relative to~$W_-$.

For the second construction, fix a plane~$W_+\in\Gr(m{+}2,n)$ 
and an algebraic curve~$\beta\subset\Gr(m{-}1,W_+)$.  Let
\begin{equation}
\label{eq: psi-m beta}
\Psi_m(\beta,W_+) =\bigcup_{E\in \beta} [E,W_+]_m\,.
\end{equation}
Thus,~$\Psi_m(\beta,W_+)$ is a curve of~$\bbP^2$s.  It
will be called the \emph{extension} of~$\beta$ relative to~$W_+$.
\end{example}

\begin{theorem}\label{thm: homologies to sigma21-star}
An irreducible variety~$X\subset\Gr(m,n)$  of dimension~$3$
satisfies~$[X] = r[\sigma_{(2,1)^*}]$ for some~$r>0$ 
if and only if one of the following holds:
\begin{enumerate}
\item $X$ is an irreducible hypersurface in~$[W_-,W_+]_m$,
      where $W_+\in\Gr(m{+}2,n)$ contains~$W_-\in\Gr(m{-}2,n)$.
\item $X$ is $\Sigma_m(\beta,\alpha)$ for some~$W\in\Gr(m,n)$
      and a pair of irreducible curves 
      $\alpha\subset [W,\C{n}]_{m+1}$ and~$\beta\subset\Gr(m{-}1,W)$.
\item $X$ is $\Sigma_m(W_-,\alpha)$ for some~$W_-\in\Gr(m{-}2,n)$
      and some irreducible curve 
      $\alpha\subset [W_-,\C{n}]_{m+1}$.  
\item $X$ is $\Psi_m(\beta,W_+)$ for some~$W_+\in\Gr(m{+}2,n)$
      and some irreducible curve 
      $\beta\subset \Gr(m{-}1,W_+)$.  
\item $X$ is $N_G(W_-,W_+)$, where $W_+\in\Gr(m{+}3,n)$ 
      contains~$W_-\in\Gr(m{-}2,n)$ and~$G$ is a nondegenerate inner
      product on~$W_+/W_-$.
\item $X$ is $N^\perp_G(W_-,W_+)$, where $W_+\in\Gr(m{+}2,n)$ 
      contains~$W_-\in\Gr(m{-}3,n)$ and~$G$ is a nondegenerate inner
      product on~$W_+/W_-$.
\item $X$ is a component of~$N_G(W_-,W_+)$, 
      where $W_+\in\Gr(m{+}3,n)$ contains $W_-\in\Gr(m{-}3,n)$ 
      and~$G$ is a nondegenerate inner product on~$W_+/W_-$.  
\end{enumerate}
\end{theorem}

\begin{proof}
The proof is very similar to the proofs of 
Proposition~\ref{prop: integral varieties of I3} and
Proposition~\ref{prop: integral varieties of I111}, 
so I will only sketch the argument.

Such an~$X$, is, of course, an integral variety of both~$\cI_{(3)}$
and~$\cI_{(1,1,1)}$.  Conversely, any $3$-dimensional integral variety
of both of these ideals is homologous to~$r[\sigma_{(2,1)^*}]$ 
for some~$r\ge1$.

Lemma~\ref{lem: integral elements of I111 and I3} 
describes the three-dimensional integral elements 
of~$\cI_{(3)}\cup\cI_{(1,1,1)}$, pointing out that they form five
distinct orbits under~$\SL(n,\bbC)$, which, in
Remark~\ref{rem:I111+3 integral element closures and intersections},
are denoted~$X'_1$, $X'_2$, $X'_3$, $X'_4$ and~$\mathcal{B}_{(2,1)^*}$.

Let~$X^\circ\subset X$ be the complement of the singular locus of~$X$.
Each of the tangent spaces to~$X^\circ$ lies in one of the five orbits
and there is one of these five orbits for which the 
set~$X^\bullet\subset X^\circ$ consisting of the points whose tangent
spaces lie in that orbit is a non-empty Zariski open set in~$X$.  The
argument now breaks into five cases.

If the tangent spaces at the points of~$X^\bullet$ lie in~$X'_1$, then
one can apply a moving frame argument to show that~$X$ must fall into
the first category of Theorem~\ref{thm: homologies to sigma21-star}.
Conversely, any subvariety~$X$ that falls into this category is an integral
of~$\cI_{(3)}\cup\cI_{(1,1,1)}$, so it must be homologous to some 
multiple of~$\sigma_{(2,1)^*}$.

If the tangent spaces at the points of~$X^\bullet$ lie in~$X'_2$, 
then the final part of the proof 
of Proposition~\ref{prop: integral varieties of I3} shows that~$X$
falls into the last category 
of Theorem~\ref{thm: homologies to sigma21-star}.  Conversely, 
since the tangent spaces to~$N_G(W_-,W_+)$ are integral elements
of~$\cI_{(3)}\cup\cI_{(1,1,1)}$, it must be homologous to a 
multiple of~$\sigma_{(2,1)^*}$.

If the tangent spaces at the points of~$X^\bullet$ lie in~$X'_3$, 
then arguments similar to those
of Proposition~\ref{prop: integral varieties of I3} 
show that~$X$ falls into either the fourth category or the sixth category 
of Theorem~\ref{thm: homologies to sigma21-star}.  Conversely, 
varieties in these two categories have their tangent spaces at 
generic points of type~$X'_3$ or~$\mathcal{B}_{(2,1)^*}$, 
so they are integrals of~$\cI_{(3)}\cup\cI_{(1,1,1)}$.

If the tangent spaces at the points of~$X^\bullet$ lie in~$X'_4$, 
then arguments similar to those
of Proposition~\ref{prop: integral varieties of I3} 
show that~$X$ falls into either the third category or the fifth category 
of Theorem~\ref{thm: homologies to sigma21-star}.  Conversely, 
varieties in these two categories have their tangent spaces at 
generic points of type~$X'_4$ or~$\mathcal{B}_{(2,1)^*}$, 
so they are integrals of~$\cI_{(3)}\cup\cI_{(1,1,1)}$.

Finally, if the tangent spaces at the points of~$X^\bullet$ lie 
in~$\mathcal{B}_{(2,1)^*}$, then arguments similar to those
of Proposition~\ref{prop: integral varieties of I3} 
show that~$X$ falls into one of the first four categories
of Theorem~\ref{thm: homologies to sigma21-star}.  Conversely, 
varieties in these four categories have their tangent spaces at 
generic points fall into one of~$X'_1$, $X'_2$, $X'_3$, $X'_4$,
or~$\mathcal{B}_{(2,1)^*}$, so they are integrals 
of~$\cI_{(3)}\cup\cI_{(1,1,1)}$.
\end{proof}

\begin{remark}[The seven types of solutions]
\label{rem: seven types of solns}
It is probably worth remarking that none of the seven types listed
in~Theorem~\ref{thm: homologies to sigma21-star} is contained in one
of the other types.  However, there is some overlap among the first
four types when the curves~$\alpha$ and/or~$\beta$ have low 
degree.  Otherwise, there is no overlap.
\end{remark}

\begin{remark}[Solutions of~$\mathcal{B}_{(2,1)^*}$]
\label{rem: solns of B 21*}
In Walters' terminology, Theorem~\ref{thm: homologies to sigma21-star}
gives a classification of the irreducible solutions 
of~$\mathcal{R}_{(2,1)}$.  Of course, since none of the last three
types of~Theorem~\ref{thm: homologies to sigma21-star} are solutions 
of~$\mathcal{B}_{(2,1)^*}$, one sees immediately how much more
restrictive the differential system~$\mathcal{B}_{(2,1)^*}$ is than
$\mathcal{R}_{(2,1)}$ is.  (Walters herself pointed out that
the varieties of type~$(5)$ in~$\Gr(2,5)$ are solutions 
of~$\mathcal{R}_{(2,1)}$ that are not solutions of~$\mathcal{B}_{(2,1)^*}$.
See~\cite[Example~$2$, Proposition~$16$]{mW97}.)

It remains to point out exactly which of the solutions 
of~$\mathcal{R}_{(2,1)}$ in the first four categories are also
solutions of~$\mathcal{B}_{(2,1)^*}$.  This is not difficult to do.

First, all of the varieties of type~$(2)$ are solutions 
of~$\mathcal{B}_{(2,1)^*}$.

Second, a variety~$\Sigma_m(W_-,\alpha)$ (i.e., of type~$(3)$) 
is a solution of~$\mathcal{B}_{(2,1)^*}$ if and only if~$\alpha$
is a rank~$1$ curve%
\footnote{A curve~$\gamma\subset\Gr(m,n)$ is of \emph{rank~$r$}
if, at the generic point~$V\in\gamma$, the line~$T_V\gamma$
is spanned by a rank~$r$ element of $Q_V\otimes V^* = \Hom(V,Q_V)$.}
in~$\Gr(3,\C{n}/W_-)$, i.e., if and only if~$\alpha$ is a solution
of~$\mathcal{B}_{(1)^*}$ 
(see Remark~\ref{rem: schur rigidity quasi-rigidity}).

Third, a variety~$\Psi_m(\beta,W_+)$ (i.e., of type~$(4)$) 
is a solution of~$\mathcal{B}_{(2,1)^*}$ if and only if~$\beta$
is a rank~$1$ curve
in~$\Gr(m{-}1,W_+)$, i.e., if and only if~$\beta$ is a solution
of~$\mathcal{B}_{(1)^*}$ 
(see Remark~\ref{rem: schur rigidity quasi-rigidity}).

Finally, the most interesting case is that of a variety of type~$(1)$,
i.e., a hypersurface~$X\subset[W_-,W_+]_m\simeq\Gr(2,4)$.  Clearly,
the generic such hypersurface is not a solution of~$\mathcal{B}_{(2,1)^*}$.
However, it is not difficult to describe the ones that are.  Without
loss of generality, I can take~$m=2$ and~$n=4$, so that~$W_-=(0)$
and~$W_+=\C{4}$.  Then it turns out that if~$X\subset\Gr(2,4)$
is a solution of~$\mathcal{B}_{(2,1)^*}$, then either it is 
accounted for by one of the first three constructions already
mentioned or else~$X$ consists of the family of~$\bbP^1$s 
tangent to a (possibly singular) surface~$S\subset\bbP^3$.  (If
$S$ itself is ruled, then this case is already accounted for by 
a previous construction.)

Note that this classification provides an alternate proof of 
Walters' result~\cite[Proposition~18]{mW97} that any solution 
of~$\mathcal{B}_{(2,1)^*}$ in~$\Gr(2,5)$ is ruled.  The present
classification is somewhat more general, since it holds for
all Grassmannians.
\end{remark}

\begin{example}[$\Gr(2,5)$]
\label{ex: Gr25 analysis}
The considerations in this section do not, by any means, give 
a complete analysis of all the extremal classes in the Grassmannians.
However, the cases treated do suffice to treat all of the cases
that appear in~$\Gr(2,5)$. The Hasse diagram for the ideal poset 
for~$\Gr(2,5)$ is drawn in Figure~\ref{fig:Gr25poset}.  

\begin{figure}
\includegraphics[width=3.5in]{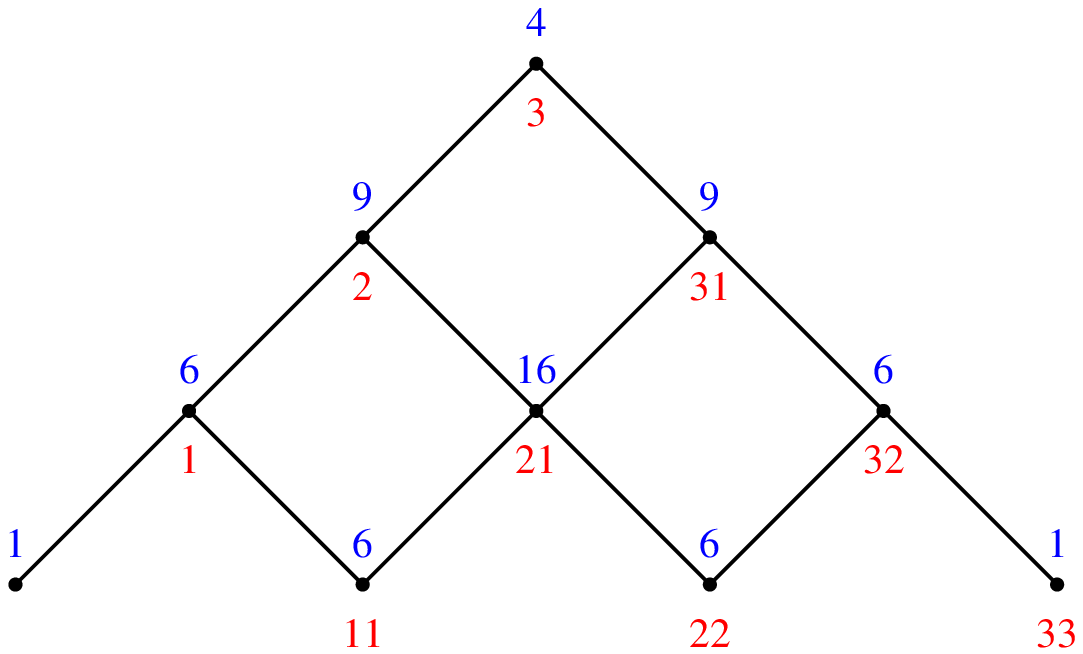}
\caption[The ideal poset for~$\Gr(2,5)$]
{\label{fig:Gr25poset} The ideal poset for~$\Gr(2,5)$.
The lower label on each node is the~$\ab\in \sfP(2,5)$ associated 
to the node and the upper label is the dimension of the corresponding
subspace of $\L^{*,0}(\eum)$.}
\end{figure}

From the labels indicating the 
`size' of each ideal at its generation level, one can see, 
for example, why the ideal~$\cI_{(3)}$, being smaller,
is less restrictive than the ideal~$\cI_{(2,1)}$.  This (partially)
explains why the effective $3$-cycles with homology
class~$r[\sigma_{(3)}]$ 
(which, by Proposition~\ref{prop: integral varieties of I21}, are 
 of the form~$[B_1,\C5]_2+\dots+[B_r,\C5]_2$ for some
lines~$B_1,\dots,B_r\in\Gr(1,5)$) are
more rigid than the $3$-cycles homologous 
to~$[\sigma_{(2,1)}]$.

Also, the relative sizes of the ideals explains, to some 
extent, why the effective $2$-cycles in~$\Gr(2,5)$
homologous to $r[\sigma_{(1,1)^*}]=r[\sigma_{(2,2)}]$ are
of the form~$\Gr(2,A_1)+\dots+\Gr(2,A_r)$ 
for some $A_1,\dots,A_r\in\Gr(3,5)$,
while the effective $2$-cycles in~$\Gr(2,5)$
homologous to $r[\sigma_{(2)*}]=r[\sigma_{(3,1)}]$ are
of the form~$S_1+\dots+S_k$ where each~$S_i$ is a surface 
of degree~$r_i$ in~$[L_i,\C5]_2\simeq\bbP^3$ for some~$L_1,\dots,L_k
\in \Gr(1,5)$ and where~$r_1+\dots+r_k = r$.
\end{example}

\subsubsection{Bundles with $c_3=0$}
\label{sssec: c3=0}

To finish this section, I give an application of 
Proposition~\ref{prop: integral varieties of I3} 
to the characterization of bundles generated by their sections
that satisfy~$c_3=0$.  Before stating the result, I need to introduce 
a particularly interesting $3$-plane bundle over~$\bbP^3=\Gr(1,4)$. 

\begin{example}[Two presentations of a $3$-plane bundle over~$\bbP^3$]
\label{ex: two presentations of L2Q}
The bundle is simply~$\L^2(Q)\to\bbP^3$, where~$Q\to\bbP^3$
is the canonical quotient bundle.  Note that there is
an exact sequence
\begin{equation}
\label{eq: C6 to L2Q presentation} 
0\longrightarrow S\w Q
 \longrightarrow \bbP^3\times\bigl(\L^2(\C4)\bigr)
 \longrightarrow \L^2(Q)
 \longrightarrow 0,
\end{equation}
where~$S\w Q$ is the image of~$S{\ot}Q$ in the trivial bundle
$\bbP^3\times\bigl(\L^2(\C4)\bigr)$ induced by wedge product.
Straightforward calculation verifies 
that~$c\bigl(\L^2(Q)\bigr) = 1+2u+2u^2$
where~$c(S) = 1 - u$.  In particular, $c_3\bigl(\L^2(Q)\bigr) =0$,
even though~$\L^2(Q)$ is obviously generated by its global sections,
since it is a quotient of a trivial bundle of rank~$6$.

For use in Theorem~\ref{thm: bundles with c3 = 0}, I want
to remark on the `Gauss mapping' induced by this presentation 
of~$\L^2(Q)$ as a quotient of a trivial bundle of rank~$6$.  
Recall that 
the wedge product~$\L^2(\C4)\times \L^2(\C4)\to \L^4(\C4)\simeq\bbC$
defines a nondegenerate symmetric quadratic form~$G$ on~$\L^2(\C4)\simeq\C6$
that is invariant under the action of~$\SL(4,\bbC)$. (In fact, this is 
the basis of the `exceptional isomorphism'~$\SL(4,\bbC)\simeq\Spin(6,\bbC)$.)
The $3$-dimensional subspaces of the form~$L\w\C4\subset \L^2(\C4)$
for~$L\in\Gr(1,4)=\bbP^3$ are isotropic for this inner product and,
in fact, this defines an embedding~$\lambda:\bbP^3\to\Gr(3,6)$ 
whose image is one component of~$N_G(0,\C6)$.  By its very definition, 
this~$\lambda$ is the Gauss mapping induced by the 
presentation~\eqref{eq: C6 to L2Q presentation}.  
In fact, this provides, via Proposition~\ref{prop: integral varieties of I3},
another proof that~$c_3\bigl(\L^2(Q)\bigr)=0$.

The bundle~$\L^2(Q)$ can also be presented as a quotient of a trivial bundle 
of rank~$5$ and this representation of~$\L^2(Q)$ will also be important
in Theorem~\ref{thm: bundles with c3 = 0}.

Fix a symplectic structure~$\Omega$ on~$\C4$, 
i.e., a nondegenerate element of~$\L^2(\C4)^*$.  
(Since these are all equivalent up to isomorphism,
it does not matter which one.)  Let~$\L^2_0(\C4)\subset \L^2(\C4)$
denote the $5$-dimensional subspace that is annihilated by~$\Omega$
and let~$(S\w Q)_0\subset S\w Q$ be the $2$-plane bundle over~$\bbP^3$
that is the intersection of~$S\w Q$ 
with~$\bbP^3\times\bigl(\L^2_0(\C4)\bigr)$.  
Define the $3$-plane bundle~$J\to\bbP^3$ by the exact sequence
\begin{equation}
\label{eq: C5 to L2Q presentation} 
0\longrightarrow (S\w Q)_0
 \longrightarrow \bbP^3\times\bigl(\L^2_0(\C4)\bigr)
 \longrightarrow J
 \longrightarrow 0.
\end{equation}

The inclusion of~$\L^2_0(\C4)$ into~$\L^2(\C4)$ and the definition
of~$(S\w Q)_0\subset S\w Q$ imply that~$J$ is isomorphic 
to~$\L^2(Q)$, so it may appear that~$J$ is `redundant'. However,
it is important to note that~$J$ is `equivariant' with respect 
to~$\Symp(2,\bbC)\subset\SL(4,\bbC)$ while~$\L^2(Q)$ is `equivariant'
with respect to the full group~$\SL(4,\bbC)$.   

Now, the quadratic form~$G$ on~$\L^2(\C4)$ restricts to~$\L^2_0(\C4)\simeq\C5$
to be a nondegenerate quadratic form, which I will continue to denote by~$G$.
The $2$-dimensional subspaces of the form~$(L\w\C4)_0\subset\L^2_0(\C4)$
for~$L\in\Gr(1,4)=\bbP^3$ are isotropic for this inner product and,
in fact, this defines an embedding~$\lambda_0:\bbP^3\to\Gr(2,5)$ 
whose image is~$N_G(0,\C5)$.  By its very definition, 
this~$\lambda_0$ is the Gauss mapping induced by the 
presentation~\eqref{eq: C5 to L2Q presentation}.  In fact, 
this provides, via Proposition~\ref{prop: integral varieties of I3},
another proof that~$c_3(J)=0$.
\end{example} 

With this discussion in place, I can now state the following 
corollary of Proposition~\ref{prop: integral varieties of I3}.

\begin{theorem}\label{thm: bundles with c3 = 0}
Let~$M$ be a compact K\"ahler manifold and let~$F\to M$ be a bundle
that is generated by its sections.  Then~$c_3(F)\ge0$ and, if
equality holds, then one of the following four possibilities holds:
\begin{enumerate}
\item $F=E\oplus T$ where~$E$ has rank~$2$ and~$T$ is trivial.
\item There is a line bundle~$L\subset F$ so that the quotient 
bundle~$F/L$ is pulled back from a curve~$C$ by some holomorphic 
map~$\alpha:M\to C$.
\item $F = \kappa^*\bigl(\L^2(Q)\bigr)\oplus T$ for a trivial bundle~$T$
and a holomorphic map~$\kappa:M\to\bbP^3$.
\item $F = \kappa^*(F')$ where $F'\to X$ is a holomorphic bundle over
some {\upshape(}possibly singular{\upshape)} 
surface~$X$ and~$\kappa:M\to X$
is holomorphic.
\end{enumerate}
Conversely, if any of these conditions holds, then~$c_3(F)=0$.
\end{theorem}

\begin{proof}
The inequality~$c_3(F)\ge0$ is, of course, immediate
from Corollary~\ref{cor: F generated by sections and ca = 0}.
Moreover, by Lemma~\ref{lem:ideal-of-Schub-form},
if~$c_3(F)=0$, then~$\kappa_F(M)\subset\Gr\bigl(m,H^0(F)\bigr)$ 
is an integral variety of~$\cI_{(3)}$,
where~$m=h^0(F)-\rank(F)$.  Now apply 
Proposition~\ref{prop: integral varieties of I3}
and interpret each of the possible cases.  
Only two cases require any comment: 

One of these is the second category of
integral varieties of Proposition~\ref{prop: integral varieties of I3}, 
i.e., the case in which~$\kappa_F(M)$ has dimension at least~$3$
and lies in an integral variety of the form~$\Sigma_m(C)$ for some 
curve~$C\subset\Gr\bigl(m{+}1,H^0(F)\bigr)$.  In such a case, examination
of the proof-analysis for Proposition~\ref{prop: integral varieties of I3}
shows that the rational mapping~$A:\kappa_F(M)\to C$ that was
defined in this case has the property that~$A{\circ}\kappa_F:M\to C$
is actually well-defined globally on~$M$.  Then, for every~$x\in M$,
$\kappa_F(x)$ is a hyperplane in~$A{\circ}\kappa_F(x)$.  Since~$F_x$
is canonically isomorphic to~$H^0(F)/\kappa_F(x)$, the line bundle~$L$
is then defined by~$L_x = A{\circ}\kappa_F(x)/\kappa_F(x)$.  
The quotient~$F_x/L_x$ is canonically isomorphic 
to~$H^0(F)/A{\circ}\kappa_F(x)$, but this latter space depends only on
the point~$A{\circ}\kappa_F(x)\in C$, so the quotient bundle~$F/L$ is
necessarily pulled back from~$C$, as claimed.

The other is the third and fourth category of integral varieties 
of Proposition~\ref{prop: integral varieties of I3}:
 
Suppose first, that~$\kappa_F(M)$ has dimension~$3$ 
and is an integral variety of the form~$N_G(W_-,W_+)$ 
for some~$W_-\in\Gr(m{-}2,n)$ and $W_+\in\Gr(m{+}3,n)$ 
with~$W_-\subset W_+$ and~$G$ 
is a nondegenerate inner product on~$W_+/W_-\simeq\C5$.
By the very definition of~$\kappa_F$, the subspace~$W_-$ must 
be zero and~$W_+$ must have dimension~$5$.  By the discussion 
in the second half of Example~\ref{ex: two presentations of L2Q},
it follows that, after identifying~$N_G(0,W_+)$ with~$\bbP^3$
via the embedding~$\lambda_0$, the bundle~$F$ can be written in the
form~$F = \kappa^*\bigl(\L^2(Q)\bigr)\oplus T$ 
where~$\kappa = {\lambda_0}^{-1}\circ\kappa_F$. 

Last, suppose that $\kappa_F(M)$ has dimension~$3$
and is a component of~$N_G(W_-,W_+)$ 
for some~$W_-\in\Gr(m{-}3,n)$ and $W_+\in\Gr(m{+}3,n)$ 
with~$W_-\subset W_+$ and~$G$ 
is a nondegenerate inner product on~$W_+/W_-\simeq\C6$.
By the very definition of~$\kappa_F$, the subspace~$W_-$ must 
be zero and~$W_+$ must have dimension~$6$.  By the discussion 
in the first half of Example~\ref{ex: two presentations of L2Q},
it follows that, after identifying~$N_G(0,W_+)$ with~$\bbP^3$
via the embedding~$\lambda$, the bundle~$F$ can be written in the
form~$F = \kappa^*\bigl(\L^2(Q)\bigr)\oplus T$ 
where~$\kappa = {\lambda}^{-1}\circ\kappa_F$. 
\end{proof}

Of course, there is an analog of Theorem~\ref{thm: bundles with c3 = 0} 
for bundles~$F$ generated by sections that satisfy~$c_{(1,1,1)}(F)=0$.
It can be deduced from Theorem~\ref{thm: bundles with c3 = 0} by 
applying the complementarity principle, but this is a task that
can be left to the interested reader.

\section[Other Symmetric Spaces]
{Extremal Cycles in Other Hermitian Symmetric Spaces }
\label{sec:extremals-in-others}

\subsection{Generalities}\label{ssec:hss-generalities}
Any irreducible Hermitian symmetric space~$M$ of compact type
can be written in the form~$M = U/K$ where $U$ is compact
and simple with Lie algebra~$\euu$ and~$K\subset U$
is a compact subgroup with Lie algebra~$\euk\subset\euu$ 
that is the fixed subgroup of an involution of~$U$
and that has a central subgroup~$T\subset K$ of dimension~$1$.  

\subsubsection{Classification}\label{ssec:hss-classification}
The list of irreducible Hermitian symmetric spaces of compact type 
is well-known~\cite[p.~518]{MR80k:53081}:
\begin{enumerate}
\item $\Gr(m,n) = \SU(n)/S\bigl(\Un(m)\times\Un(n{-}m)\bigr)$, 
     the complex Grassmannians;
\item $Q_n = \SO(n{+}2)/\bigl(\SO(2)\times\SO(n)\bigr)$, 
     the complex $n$-quadric;
\item $N^+_n = \SO(2n)/\Un(n)$, the space of isotropic%
\footnote{I.e., totally null for the inner product, 
           which explains the `N' in the notation.}
      $n$-planes 
      of positive chirality%
\footnote{For an explanation of this term, 
see~\S\ref{ssec:quadrics-as-symmetric}. }
in~$\C{2n}$ endowed with a nondegenerate 
      inner product and orientation;
\item $L_n = \Symp(n)/\Un(n)$, the space of Lagrangian $n$-planes 
      in~$\C{2n}$ endowed with a symplectic form;
\item $\E_6/\bigl(S^1{\cdot}\Spin(10)\bigr)$, 
      the singular locus of the
      projectivization of the null cone 
      of Cartan's $\E_6^\bbC$-invariant 
      cubic form on~$\C{27}$~\cite[pp.~142--143]{eC1894}; and
\item $\E_7/\bigl(S^1{\cdot}\E_6\bigr)$, 
      the second singular locus of the
      the projectivization of the null cone 
      of Cartan's $\E_7^\bbC$-invariant
      quartic form on~$\C{56}$~\cite[pp.~143--144]{eC1894}.
\end{enumerate}

\subsubsection{Positive forms}
\label{ssec:hss-positive-forms}

For the rest of this subsection, $M=U/K$ will represent 
one of the members of the above list, with~$U$ and~$K$ as indicated.
Set~$\euu = \euk+\eum$ where~$\eum = \euk^\perp\subset\euu$.  
Then~$\Ad(K)$ preserves this splitting of~$\euu$, so that~$\eum$
is naturally a~$K$-representation.  This representation is 
irreducible and complex, i.e., there exists complex structure on~$\eum$
that commutes with the action of~$K$.  The negative of the Killing form 
on~$\euu$ restricted to~$\eum$ is a positive definite inner product
that is compatible with this complex structure.

The projection~$U\to U/K = M$ defines a canonical isomorphism
between~$\eum$ and~$T_{eK}M$ and there is a unique K\"ahler structure
on~$M$ that is~$U$-invariant and that agrees with the complex structure
and inner product on~$\eum$ under this identification. 

Recall that the ring~$\Omega^*(M)^U$ of $U$-invariant 
forms on~$M$ consists entirely of closed forms and that the induced map 
to deRham cohomology~$\Omega^*(M)^U\to H^*(M,\bbR)$ is an isomorphism.
The mapping~$\Omega^*(M)^U\to\L^*(\eum^*)^K$ defined by 
evaluating~$\phi\in\Omega^*(M)^U$ at~$eK$ is also an isomorphism.

According to Borel and Hirzebruch~\cite[\S14.10]{MR21:1586},
all of the cohomology of~$M$ is of type~$(p,p)$.  Thus, each cohomology
class of~$M$ is represented by a unique~$U$-invariant $(p,p)$-form.

Let~$\cH^{p,p}(M)$ denote the real-valued $U$-invariant~$(p,p)$-forms 
on~$M$, and let
$$
\cH^{p,p}_+(M)\subset\cH^{p,p}(M)
$$ 
denote the closed, 
convex cone of positive $U$-invariant $(p,p)$-forms on~$M$.  
This cone has nonempty interior since
the~$p$-th power of the K\"ahler form obviously lies in its interior.  
In particular, a basis for~$\cH^{p,p}(M)$ can be chosen from among the 
extremal rays of~$\cH^{p,p}_+(M)$.%
\footnote{This relies on the fact that any closed, convex cone 
is the convex hull of its extremal rays.} 

Suppose that~$\phi\not=0$ lies on an extremal ray of $\cH^{p,p}_+(M)$.
Since~$\phi$ is positive, it can be written (locally) in the form
\begin{equation}\label{eq:phi-positive-expression}
\phi = \iC^{p^2}\sum_{k=1}^{r}\zeta_k\w\ov{\zeta_k}
\end{equation}
for some~$\zeta_1,\dots,\zeta_r\in\Omega^{p,0}$ whose (complex)
span~$\I_\phi\subset \L^{p,0}(M)$ is independent of this representation
and so is globally defined.  Moreover, since~$M$ is a symmetric space,
all the $U$-invariant forms are parallel with respect to the Levi-Civita
connection, so it follows easily that~$\I_\phi$ must be a $U$-invariant, 
parallel subbundle of~$\L^{p,0}(M)$ with respect to the Levi-Civita 
connection. Since the Levi-Civita connection is torsion-free and respects 
the holomorphic structure on~$M$, 
it follows that~$\I_\phi$ must actually be a 
holomorphic subbundle of~$\L^{p,0}(M)$, with a holomorphic sheaf of 
sections~$\cI_\phi$.

On the other hand, suppose that~$\I\subset \L^{p,0}(M)$ 
is a minimal $U$-invariant (and hence parallel) complex subbundle 
of $\L^{p,0}(M)$.  Then~$\I$ is holomorphic and has a $U$-invariant
Hermitian inner product, which is unique up to a constant multiple.%
\footnote{Once one fixes a $U$-invariant K\"ahler form on~$M$ 
(which is, itself, unique up to a constant multiple), this determines
a canonical choice of Hermitian inner 
product on each of the~$\L^{p,0}(M)$.}
If~$r$ is the rank of~$\I$, and~$\zeta_1,\ldots,\zeta_r$ is a local 
unitary basis of~$I$, then defining~$\phi_{\I}$ to be the right hand
side of~\eqref{eq:phi-positive-expression}, one sees that~$\phi_{\I}$ is
independent of the choice of unitary basis of~$\I$, so that~$\phi_{\I}$ is
globally, defined, positive, and $U$-invariant.  
Hence~$\phi_{\I}$ belongs
to~$\cH^{p,p}_+(M)$.    

Evidently, $\phi_{\I}$ will lie on the boundary of~$\cH^{p,p}_+(M)$ as long
as~$\I$ is a proper subbundle of~$\L^{p,0}(M)$ and, moreover, it will be
an extreme point of the boundary if and only if~$\I$ is minimal.  (If
$\I$ is not minimal, it can be written in the form~$\I = \I' \oplus \I''$
for some orthogonal invariant subbundles~$\I',\I''\subset \I$, 
in which case $\phi_{\I} = \phi_{\I'} + \phi_{\I''}$, so that $\phi_{\I}$
is not extremal.  Conversely, suppose~$\I$ is minimal but
that~$\phi$ is not extremal, i.e., that there 
exist~$\phi',\phi''\in\cH^{p,p}_+(M)$ that are not multiples of 
each other so that~$\phi = \phi' + \phi''$.  Then~$\phi'$ and $\phi''$
will be associated to nonzero parallel subbundles~$\I',\I''\subset\I$,
which, since~$\I$ is minimal by hypothesis, must be equal to~$\I$ itself.
The~$U$-invariance of~$\phi'$ and~$\phi''$ implies that they each
define a $U$-invariant Hermitian inner product on~$\I$ and the assumption
that they are not multiples of each other implies that these two Hermitian
inner products are not multiples of each other.  However, this would imply
that there exists a nontrivial $U$-invariant splitting of~$\I$, contrary
to hypothesis.  Thus, $\phi$ must have been extremal after all.)

This argument establishes a one-to-one correspondence between the 
extremal rays of the cone~$\cH^{p,p}_+(M)$ and the minimal $U$-invariant
subbundles of~$\L^{p,0}(M)$.  The $U$-invariant subbundles of~$\L^{p,0}(M)$ 
are, in turn, in one-to-one correspondence with the $K$-invariant
subspaces of~$\L^{p,0}(\eum)$.

\subsubsection{Kostant's description}
\label{ssec:hss-kostant-descrip}

By a theorem of Kostant~\cite[Corollary 8.2]{MR26:265},
the representation of~$K$ on~$\L^{*,0}(\eum)$ is multiplicity-free,
so, for each~$p$, there are only a finite number of minimal $K$-invariant 
subspaces of~$\L^{p,0}(\eum)$ and hence only a finite number of extremal 
rays in~$\cH^{p,p}_+(M)$.  Moreover, the generators of these rays are 
evidently linearly independent, implying that the base of the 
cone~$\cH^{p,p}_+(M)$ is a simplex of dimension~$h^{p,p}(M)-1$.

In fact, Kostant proved that there is a generalized Schubert cell 
decomposition of~$M$.  This will be useful in what follows, so I 
will now describe some of his results in~\cite{MR26:265,MR26:266}.  Let~$d$
be the dimension of~$M$ (as a complex manifold) and
let~$\sfP(M)$ be a set that indexes the minimal $K$-invariant subspaces
of~$\L^{*,0}(\eum)$.  Write
\begin{equation}\label{eq: indexing K-subspaces}
\L^{*,0}(\eum) = \oplus_{\ab \in \sfP(M)} \I_\ab\,.
\end{equation}
Define $|\ab|\in\bbZ^+$ so that $\I_\ab$ is a subspace 
of~$\L^{|\ab|,0}(\eum)$.  Then~$|\ab|\le d$ and, for every~$\ab\in \sfP(M)$
there is a unique~$\ab^*\in \sfP(M)$ satisfying $|\ab|+|\ab^*|=d$
and having the property that wedge product induces a nonzero (and hence
nondegenerate) pairing~$\I_\ab\times\I_{\ab^*}\to\L^{d,0}(\eum)$.  
(Of course, $\L^{d,0}(\eum)$, being the top exterior power of~$\eum$ 
as a complex vector space, has dimension~$1$.)

In \S\S6--8 of~\cite{MR26:265}, Kostant shows that there 
is decomposition 
\begin{equation}\label{eq: hss Schubert cells}
M = \cup_{\ab \in \sfP(M)} W_\ab\,,
\end{equation}
where~$W_{\ab}$ is biholomorphic to~$\C{d-|\ab|}$.  
Moreover, regarding~$M$
as a complex homogeneous space~$M = G/P$ where~$G$ is the 
connected complex Lie group whose maximal compact is~$U$ and~$P$
is a (maximal) parabolic subgroup of~$G$, the cell~$W_{\ab}$ can be
written as the orbit of a nilpotent subgroup~$N_{\ab}\subset G$ that
is transverse to~$P$ and of dimension~$d{-}|\ab|$.  

Let~$\sigma_\ab\subset M$ be the closure of~$W_\ab$.  
Then~$\sigma_\ab$ is an irreducible algebraic subvariety of~$M$.
The classes~$[\sigma_\ab]$ give a basis for~$H_*(M,\bbZ)$ and,
moreover, generate the semigroup~$H_*^+(M,\bbZ)$.  These varieties are
known as the (generalized) Schubert varieties of~$M$.

Kostant shows~\cite[Corollary 6.15]{MR26:266} 
that there exists a form~$\phi_\ab\in\cH^{|\ab|,|\ab|}_+(M)$
in the extremal ray associated to the subspace~$\I_\ab$ with
the property that
\begin{equation}\label{eq: hss Schur forms}
\int_{\sigma_{\bb^*}}\phi_\ab = \delta^\bb_\ab\,.
\end{equation}
Because of the positivity of~$\phi_\ab$, it follows that~$\phi_\ab$
vanishes identically on~$\sigma_{\bb^*}$ for all~$\bb\in \sfP(M)$ 
with~$|\bb| = |\ab|$ and~$\bb\not=\ab$.

Just as in the case of Grassmannians, the tangent spaces to
the $G$-images of the cell~$W_\ab$ define subspaces~$E\subset T_xM$
that are said to be of type~$\ab$.  Again, 
because of~\eqref{eq: hss Schur forms} it follows
that~$Z(\phi_\ab)$ contains all of the spaces of type~$\bb^*$
where~$\bb\in \sfP(M)$ satisfies~$|\bb| = |\ab|$ but~$\bb\not=\ab$.
Moreover, it is not difficult to show from Kostant's definitions
that~$E\subset T_xM$ is of type~$\ab$ if and only if~$E^\perp\subset T_xM$ 
is of type~$\ab^*$.  In turn, this implies that the extremal 
positive form~$*\phi_\ab$ must be a (positive) multiple of~$\phi_{\ab^*}$.

The subspace~$\L^{1,0}(\eum) = \eum$ is, of course, irreducible under~$K$
and, conforming to the case of Grassmannians, I will denote 
the~$\ab\in \sfP(M)$ for which~$\I_\ab = \L^{1,0}(\eum)$ by~$\ab=1$.
Then~$\phi_1$ defines a $U$-invariant K\"ahler metric on~$M$ and
there is an integer~$\mu^\ab>0$ for each~$\ab\in \sfP(M)$ so that,
when~$1\le p\le d$,
\begin{equation}\label{eq: hss sum of Schur forms}
{\phi_1}^p = \sum_{\{\ab\in \sfP(M) \mid\, |\ab| = p\}} \mu^\ab\,\phi_\ab\,.
\end{equation}
By Wirtinger's theorem and~\eqref{eq: hss Schur forms},
it follows that~$\mu^\ab = |\ab|!\,\vol(\sigma_{\ab^*})$, 
where the volume is computed with respect to the metric~$\phi_1$.
Of course, \eqref{eq: hss Schur forms} 
and~\eqref{eq: hss sum of Schur forms} imply that 
\begin{equation}
*\phi_\ab = \frac{|\ab|!\,\mu^{\ab^*}}{|\ab^*|!\,\mu^\ab}\>\phi_{\ab^*}
\end{equation}
The integers~$\mu^\ab$ can be calculated by representation theoretic means.
For explicit formulae and more details than are given
here, see~\cite{kKhT99}.  

Also, the fundamental reference~\cite{MR55:2941} contains an explicit
computation of~$\sfP(M)$ in each of the classical cases.  The reader
may find this helpful from time to time, though I will not need 
the full details in what follows.  

\subsubsection{The ideal poset}
\label{ssec:hss-ideal-poset}

The set~$\sfP(M)$ has a natural poset structure.  The partial 
ordering is defined by the condition that~$\ab\le\bb$
if and only if~$\I_\bb$ lies in the subspace~$\I_\ab\cdot\L^{*,0}(\eum)$.   
The initial element, usually denoted by~$\ab = 0$,
is the one such that~$\I_0 = \L^{0,0}(\eum)$, and
the element~$1\in \sfP(M)$ satisfies~$0\le 1\le \ab$ for all~$\ab\not=0$.
This poset is sometimes called the \emph{Bruhat poset}
associated with~$M$ (see \cite{kKhT99}).  

Obviously, $\ab\le \bb$ implies that~$|\ab|\le |\bb|$.  
Moreover, if $\ab\le \bb$, 
there exists a \emph{chain}~$(\ab_p,\ab_{p+1},\dots,\ab_q)$
where~
$$
\ab = \ab_p\le \ab_{p+1}\cdots\le \ab_{q-1}\le \ab_q = \bb,
$$
with~$|\ab_k| = k$ for ~$p\le k\le q$. It is shown in~\cite{kKhT99} 
that~$\mu^\ab$ is the number of  distinct chains from~$0$ to~$\ab$. 

The Pieri formula generalizes to
\begin{equation}\label{eq: hss Pieri}
{\phi_1}^r\phi_\ab 
= \sum_{\{\bb\in \sfP(M) \mid\, |\bb| = r+|\ab|\}} \mu^\bb_\ab\,\phi_\bb\,.
\end{equation}
where~$\mu^\bb_\ab$ is a nonnegative integer that is positive
if and only if~$\bb\ge\ab$.

The analog of Lemma~\ref{lem: ideal inclusions} 
holds for the general Hermitian symmetric space:  
A subspace of type~$\bb^*$ is an integral element of~$\cI_\ab$ 
if and only if~$\bb\not\ge\ab$.
Moreover, the maximum value of~$|\bb|$ for~$\bb$ satisfying~$\bb\not\ge\ab$
is also the maximum dimension of integral elements of~$\cI_\ab$.
(Of course, this does not generally provide a 
classification of the maximal integral elements of~$\cI_\ab$,
which seems to be a hard problem in general.)

\begin{remark}[Hasse diagrams]
\label{rem: Hasse diagrams}
There are several figures in this article that illustrate the structure
of the ideal poset for various Hermitian symmetric spaces by 
drawing the associated \emph{Hasse diagram}.  The 
convention followed is that the elements in the poset are represented
by the nodes of a graph.  The horizontal placement of the node 
corresponding to~$\ab$ is determined by~$|\ab|$, with this coordinate
increasing from left to right.%
\footnote{This is slightly
nonconventional; usually a Hasse diagram is drawn so that relative
order is indicated by relative \emph{height}.  However, in the interests
of saving space, I have reoriented the diagrams as indicated.}
An edge is drawn between
two nodes~$\ab$ and $\bb$ when~$\ab\le\bb$ and~$|\bb| = |\ab| + 1$.  
The vertical placement of the nodes in the graph is less algorithmic
but is chosen to minimize the number of different slopes of the edges
and the number of crossings of edges.
\end{remark}

\subsection{Quadrics}\label{ssec:quadrics-as-symmetric}
Let~$(,)$ be the standard complex inner product on~$\C{n+2}$
and let~$Q_n\subset\bbP^{n+1}$ be the space of null lines for
this inner product, i.e.,~$[v]$ lies in~$Q_n$ for~$v\not=0$
in~$\C{n+2}$ if and only if~$(v,v)=0$.  Then~$Q_n$ is a compact
complex manifold of dimension~$n$.  It can also be regarded as
an Hermitian symmetric space:
\begin{equation}
Q_n = \frac{\SO(n{+}2)}{\SO(2)\times\SO(n)}
\end{equation}
and so carries an $\SO(n{+}2)$-invariant K\"ahler structure,~$\omega$.
(Explicitly, the isomorphism takes an 
oriented $2$-plane~$P\subset\bbR^{n+2}$
to the line~$[v_1 - \iC\,v_2]\in Q_n$, where~$(v_1,v_2)$ is any oriented,
orthonormal basis of~$P$.)

\subsubsection{Topology}
\label{sssec: quadric topology }
When~$n$ is odd, $H_{2p}(Q_n,\bbZ) \simeq \bbZ$ for~$0\le p\le n$ and 
$\omega^p$ determines a generator~$a_p$ of~$H^+_{2p}(Q_n,\bbZ)$.

When~$n=2m$, one still has~$H_{2p}(Q_{2m},\bbZ) \simeq \bbZ$ 
for~$0\le p < m$ and~$m < p\le 2m$, 
but~$H_{2m}(Q_{2m},\bbZ) \simeq \bbZ^2$. 
A pair of generators of~$H^+_{2m}(Q_{2m},\bbZ)$ can be described as follows: 

Let~$\eb_1,\dots,\eb_{2m+2}$ be the standard basis of~$\C{2m+2}$.
Let~$V\subset\C{2m+2}$ be any maximal isotropic plane, 
let $\ub_1,\dots,\ub_{m+1}$
be a basis for~$V$, and let~$\vb^1,\dots,\vb^{m+1}\in \C{2m+2}$ be
chosen so that~$(\ub_i,\vb^j) = \delta^j_i$.  It is easy to show that
\begin{equation}
\ub_1\w\dots\w\ub_{m+1}\w\vb^1\w\dots\w\vb^{m+1}
= \pm\,\iC^{m+1}\,\eb_1\w\dots\w\eb_{2m+2}
\end{equation}
and that the sign~$\pm$ does not depend on the choices of~$\ub$ or~$\vb$,
but only on~$V$ itself.  (This is a manifestation of the fact 
that~$\Or(n{+}2,\C{})$ has two components.)
One says that~$V$ has \emph{positive chirality} 
or \emph{negative chirality} 
according to this sign.  

Let~$V_+$ and $V_-$ be two maximal isotropic planes in~$\C{2m+2}$
of positive and negative chiralities, respectively.  Then their 
projectivizations~$P_\pm = \bbP(V_\pm)\subset Q_m$ give generators 
for~$H^+_{2m}(Q_{2m},\bbZ)$.  

\subsubsection{The ideals~$\cI_\pm$}
\label{sssec: ideals I+- on Q}

Consider the
representation of~$K = \SO(2){\times}\SO(2m)$ on the
space~$\eum\simeq\C{2m}$.
This representation is seen to act as follows:  The factor~$\SO(2)$
acts as the unitary multiples of~$\I_{2m}$, i.e., 
as~$\textrm{e}^{i\theta}\I_{2m}$.
The factor~$\SO(2m)$ acts on~$\C{2m}$  by regarding~$\C{2m}$ 
as~$\bbC\ot\bbR^{2m}$ and letting~$\SO(2m)$ act on the~$\bbR^{2m}$ 
factor.  This action is irreducible as soon as~$m>1$, which I
assume from now on.  In this representation~$K = S^1{\cdot}\SO(2m)$ is a
maximal compact subgroup of the complex subgroup~$\C{*}{\cdot}\SO(2m,\bbC)$,
which certainly acts irreducibly on~$\C{2m}$.  According to 
\cite[Theorem 19.2]{MR93a:20069}, $\SO(2m,\bbC)$ acts irreducibly on
each of the (complex) exterior powers~$\L^p(\C{2m})$ for~$p<m$ while
$\L^m(\C{2m})$ is the direct sum of two irreducible 
subspaces~$\L^m_\pm(\C{2m})$.  Using duality,~$\L^{m+p}(\C{2m})\simeq 
\L^{m-p}(\C{2m})$ for~$p\ge0$, so the only reducible exterior power
is the middle one.

Since~$K$ is a maximal compact 
in~$\C{*}{\cdot}\SO(2m,\bbC)$, it follows without difficulty that 
each of the representations~$\L^p(\C{2m})$ for~$p\not=m$ 
and~$\L^m_\pm(\C{2m})$ are irreducible as complex representations of~$K$.  
Moreover, as representations of~$K$, the space~$\L^p(\C{2m})$ is isomorphic 
to~$\L^{p,0}(\C{2m})$.  Thus, $\L^{p,0}(\C{2m})$ is irreducible
for~$p\not=m$, 
and has two inequivalent irreducible summands for~$p=m$.

Now, corresponding to the irreducible 
summands~$\L^m_\pm\bigl((\C{2m})^*\bigr)$ in~$\L^m\bigl((\C{2m})^*\bigr)$,
there are two $\SO(2m{+}2)$-invariant holomorphic 
subbundles~$\I_\pm\subset\L^{m,0}(Q_{2m})$ and, according to the general
results of~\S\ref{ssec:hss-positive-forms}, 
two corresponding $\SO(2m{+}2)$-invariant positive~$(m,m)$-forms, 
say~$\phi_\pm$.  They can be normalized by requiring
that~$\omega^m = \phi_-+\phi_+$, so I do this. These two forms are 
linearly independent and so must span~$\cH^{m,m}(Q_{2m})$, 
which has dimension~$2$.  These forms lie on the extremal rays
of the convex cone~$H^{m,m}_+(Q_{2m},\bbZ)$.  

The sections of the two subbundles~$\I_\pm$ 
generate holomorphic ideals~$\cI_\pm$ on~$Q_{2m}$
and it is the integral manifolds of these that are of interest.
 
\subsubsection{Integral elements}
\label{sssec: integral elements on Q}

It will be useful to identify the 
spaces~$\L^m_\pm(\C{2m})$ more explicitly.
If~$\eb_1,\dots,\eb_{2m}$ is an oriented 
orthonormal basis of~$\C{2m}$, there
is a unique linear map~$*:\L^{m}(\C{2m}) \to \L^{m}(\C{2m})$ 
that satisfies
\begin{equation}
*\bigl(\eb_{i_1}\w\dots\eb_{i_{m}}\bigr) = \eb_{j_1}\w\dots\eb_{j_{m}}
\end{equation}
whenever~$(i_1,\dots,i_{m},j_1,\dots,j_{m})$ is an even permutation 
of~$(1,\dots,2m)$.  This map does not depend on the choice of basis, but
only on the orientation and inner product, i.e., it commutes with the
action of~$\SO(2m,\bbC)$.   By its definition, $*$ satisfies~$** 
= (-1)^{m^2} = (-1)^m$.  In fact, $*$ has two eigenvalues,
namely~$\pm\iC^m$.  In order to simplify some of the
statements appearing below, I will take~$\L^m_+(\C{2m})$ to be
the~$\iC^{-m}$ eigenspace of~$*$ and take~$\L^m_-(\C{2m})$ to be
the $-\iC^{-m}$ eigenspace of~$*$.  Thus, $\L^m_\pm(\C{2m})$
is spanned by the vectors
\begin{equation}
\eb_{i_1}\w\dots\eb_{i_{m}} \pm  \iC^m\,\eb_{j_1}\w\dots\eb_{j_{m}}
\end{equation}
where~$(i_1,\dots,i_{m},j_1,\dots,j_{m})$ is an even permutation 
of~$(1,\dots,2m)$.  

Of course, there is a corresponding $\SO(2m,\bbC)$-invariant 
decomposition of the exterior forms of degree~$m$ on~$\C{2m}$.  With
these definitions, I can now state the following lemma.

\begin{lemma}\label{lem:isotropic-ideals}
All of the forms in~$\L^m_-\bigl((\C{2m})^*\bigr)$
vanish on a given $m$-plane~$E\subset\C{2m}$ if and only if $E$ 
is isotropic with positive chirality.

All of the forms in~$\L^m_+\bigl((\C{2m})^*\bigr)$
vanish on a given $m$-plane~$E\subset\C{2m}$ if and only if $E$ 
is isotropic with negative chirality.
\end{lemma}

\begin{proof}
I will first show that if~$E$ is not isotropic
then one can construct forms in each of~$\L^m_\pm\bigl((\C{2m})^*\bigr)$
that do not vanish on~$E$.  To begin, suppose that the inner product
is nondegenerate on~$E$.  Then there exists an oriented orthonormal basis
$\vb_1,\dots,\vb_{2m}$ of~$\C{2m}$, with dual basis~$\vb^1,\dots,\vb^{2m}$ 
of~$(\C{2m})^*$, so that~$E$ is spanned by the vectors
$\vb_1,\dots,\vb_{m}$.  Then neither of the forms
$$
\psi_\pm = \vb^{1}\w\vb^{2}\w\dots\w\vb^{m} 
\pm \iC^m\,\vb^{m+1}\w\vb^{m+2}\w\dots\w\vb^{2m}
\quad \in\quad \L^m_\pm\bigl((\C{2m})^*\bigr)
$$
vanishes on~$E$.  If $E$ has nullity~$p<m$, then one can choose the
basis~$\vb_1,\dots,\vb_{2m}$ as above so that~$E$ is spanned by
$$
\vb_1-\iC\,\vb_{m+1},\ \dots,\ \vb_p-\iC\,\vb_{m+p},
\ \vb_{p+1},\ \dots\ \vb_m\,,
$$
and, again, both of~$\psi_\pm$ are nonzero on~$E$.

Now, suppose that~$E$ is isotropic and has positive
chirality.  Then there is an oriented orthonormal basis
$\vb_1,\dots,\vb_{2m}$ of~$\C{2m}$ so that~$E$ is spanned by the vectors
$$
\vb_1-\iC\,\vb_{m+1},\ \vb_2-\iC\,\vb_{m+2},
      \ \dots,\ \vb_m-\iC\,\vb_{2m}\,.
$$
Straightforward computation now shows that, 
when $(i_1,\dots,i_m,j_1,\dots,j_m)$ is any even permutation 
of~$(1,\dots,2m)$, the $m$-form
$$
\vb^{i_1}\w\vb^{i_2}\w\dots\w\vb^{i_m} 
- \iC^m\,\vb^{j_1}\w\vb^{j_2}\w\dots\w\vb^{j_m}
\quad \in \L^m_-\bigl((\C{2m})^*\bigr)
$$
vanishes on~$E$.  Since such $m$-forms 
span~$\L^m_-\bigl((\C{2m})^*\bigr)$,
all the forms in~$\L^m_-\bigl((\C{2m})^*\bigr)$ vanish on all of the 
isotropic planes of positive chirality.  Since~$\L^m\bigl((\C{2m})^*\bigr)$
is the direct sum of the spaces~$\L^m_\pm\bigl((\C{2m})^*\bigr)$, not all
of the forms in~$\L^m_+\bigl((\C{2m})^*\bigr)$ can vanish on~$E$.  By
applying an orientation reversing isometry, it follows that not all
of the elements of~$\L^m_-\bigl((\C{2m})^*\bigr)$ vanish on any given
$m$-plane of negative chirality.  This proves the first statement in
the lemma.

The proof of the second statement is similar.
\end{proof}

It follows from Lemma~\ref{lem:isotropic-ideals}
that~$Z(\phi_-)$ at~$[v]\in Q_{2m}$ 
consists of the isotropic $m$-dimensional 
subspaces~$E\subset T_{[v]}Q_{2m}$
of positive chirality while $Z(\phi_+)$ at~$[v]\in Q_{2m}$ 
consists of the isotropic $m$-dimensional 
subspaces~$E\subset T_{[v]}Q_{2m}$
of negative chirality.

\subsubsection{Integral varieties}
\label{sssec: integral varieties on Q}
The computation of the integral elements of~$\cI_{\pm}$ showing
that they are the maximal isotropic subspaces of~$T_xQ$ now
allows a classification of the integral manifolds of these two
ideals. 

\begin{proposition}\label{prop:isotropic-rigidity}
Any connected $m$-dimensional complex submanifold~$S\subset Q_{2m}$ 
whose tangent plane at each point is isotropic is an open subset of 
the projectivization of an isotropic $(m{+}1)$-dimensional subspace 
of~$\C{2m+2}$.  
\end{proposition}

\begin{proof}
Let~$H\simeq\Or(2m{+}2,\bbC)$ be the subgroup of~$\GL(2m{+}2,\bbC)$
consisting of the matrices~$\us$ that satisfy
\begin{equation}\label{eq:defining-H}
{}^t\us\begin{pmatrix}0_{m+1}&\I_{m+1}\\ \I_{m+1}&0_{m+1}\end{pmatrix} \us
= \begin{pmatrix} 0_{m+1} & \I_{m+1}\\ \I_{m+1} & 0_{m+1}\end{pmatrix}.
\end{equation}
Then $H$ acts on~$\C{2m+2}$ and induces a transitive action 
on~$Q_{2m}\subset\bbP^{2m+1}$.  Let~$\us:H\to\GL(2m{+}2,\bbC)$
denote the inclusion and write
$$
\us = (\us_0\quad\dots\quad \us_{m} \quad \us^0\quad\dots\quad \us^{m})
$$
where~$\us_a,\us^a:H\to\C{2m+2}$ are regarded as (holomorphic) mappings.
Then~$[\us_0]:H\to Q_{2m}$ is a holomorphic principal fiber bundle 
over~$Q_{2m}$.   Moreover, the 
map~$[\us_0\w\us_1\w\dots\w\us_m]:H\to\Gr(m,2m)$ 
makes~$H$ into a holomorphic
fiber bundle over~$N^+_m\cup N^-_m$, i.e., 
the set of isotropic $m$-planes
in~$\C{2m}$.

In accordance with the moving frame, write the structure
equations as
\begin{equation}\label{eq:H-structure-equations}
\d\us = \d (\us_a\ \us^a) = (\us_b\ \us^b)
\begin{pmatrix}\alpha^b_a & \gamma^{ba}\\ 
                \beta_{ba}& -\alpha^a_b \end{pmatrix}
= \us\,\theta
\end{equation}
where~$\beta_{ba} = -\beta_{ab}$ and~$\gamma^{ba}=-\gamma^{ab}$.  
(These relations follow in the usual way from the exterior 
derivative of~\eqref{eq:defining-H}.)  Moreover, the 
structure equation $d\theta = -\theta\w\theta$ holds since~$\theta 
= \us^{-1}\,\d \us$.

Now suppose that~$S\subset Q_{2m}$ is an $m$-dimensional complex 
submanifold with the property that all of its tangent planes are 
isotropic and let~$F\subset S \times H$ denote the set of pairs~$(x,\us)$ 
that satisfy
\begin{enumerate}
\item $[\us_0] = x\in S\subset Q_{2m}$; and
\item The projectivized isotropic $m$-plane~$[\us_0\w\us_1\w\dots\w\us_m]$ 
is tangent to~$S$ at~$x$.
\end{enumerate}
Then~$F\to S$ is a holomorphic fiber bundle over~$S$.  We now consider
the functions and forms on~$S\times H$ to be pulled back to~$F$ in the
usual way of the moving frame.  Then, by construction,
$$
d\us_0 = \us_0\,\alpha^0_0+\us_1\,\alpha^1_0+\dots+\us_m\,\alpha^m_0\,
$$
i.e., $\beta_{0a} = 0$ for~$1\le a\le m$, while the fact that~$F\to S$ is
a submersion implies that~$\alpha^1_0\w\dots\w \alpha^m_0\not=0$.  

When~$1\le a\le m$, the structure equations imply
\begin{equation}
0 = d\beta_{0a} = - \beta_{0b}\w\alpha^b_a + \alpha^b_0\w\beta_{ba} 
= \sum_{b=1}^{m}\alpha^b_0\w\beta_{ba}.
\end{equation}
Since $\alpha^1_0\w\dots\w \alpha^m_0\not=0$,
Cartan's Lemma implies that there exist functions~$B_{abc}=B_{acb}$
so that~$\beta_{ab} = B_{abc}\,\alpha^c_0$.  
However, since~$\beta_{ab} = -\beta_{ba}$, it follows that~$B_{abc}=0$,
i.e.,~$\beta_{ab} = 0$.  

The structure equations thus imply that~$[\us_0\w\us_1\w\dots\w\us_m]$
is locally constant, i.e., that, locally,~$S$ is tangent to the
projectivization of a fixed isotropic~$(m{+}1)$-plane in~$\C{2m+2}$.
Since~$S$ is connected, smooth, and holomorphic, it follows that~$S$
must be everywhere tangent to this linear~$\bbP^m$, as desired.
\end{proof}

Finally, this local rigidity statement yields the desired global
rigidity statements:

\begin{theorem}\label{thm: isotropic rigidity}
Suppose that $S\subset Q_{2m}$ is an $m$-dimensional 
subvariety~that satisfies~$[S] = r[P_+]$ 
{\upshape(}respectively, $[S]=r[P_-]${\upshape)} for some~$r>0$.  
Then $S$ is the union of $r$ linear isotropic $\bbP^m$s in~$Q_{2m}$ 
of positive {\upshape(}respectively, negative{\upshape)} chirality.  
\end{theorem}

\begin{proof}
Each of~$\phi_\pm$ is a positive $(m,m)$-form and the above analysis
shows that
$$
\int_{P_-}\phi_+ = \int_{P_+}\phi_+\- = 0.
$$
Since~$\omega^m = \phi_+ + \phi_-$, it follows that $[P_+]$ and~$[P_-]$
generate~$H^+_{2m}(Q_{2m},\bbZ)$.  

Thus, any $m$-dimensional subvariety $S\subset Q_{2m}$ that 
satisfies~$[S] = r[P_+]$ must be a union of irreducible subvarieties 
whose homology classes are multiples of~$[P_+]$.  Thus, I may assume that 
$S$ is irreducible and satisfies~$S=r[P_+]$.  In particular, $\phi_-$
vanishes on~$S$.  By Lemma~\ref{lem:isotropic-ideals}, it follows that, 
on the smooth locus of~$S$, each of its tangent spaces is an 
isotropic~$m$-plane.  By Proposition~\ref{prop:isotropic-rigidity}, 
it follows that~$S$ must contain an isotropic~$\bbP^m$ in~$Q_{2m}$.
Since~$S$ is irreducible, it follows that~$S$ must itself be
such a plane.  In particular,~$r=1$.  Since $S$ is homologous to~$P_+$, 
it must have positive chirality.

The argument when~$[S] = r[P_-]$ is essentially the same.
\end{proof}

\begin{remark}\label{rem:quadric-intersection-pairing}
It follows from the proposition in~\cite[p.~735]{MR80b:14001} that,
when~$m$ is odd, the intersection pairing on~$H_{2m}(Q_{2m},\bbZ)$ 
satisfies
$$
[P_\pm]\cap [P_\pm] = 0,\qquad\text{and}\qquad [P_-]\cap [P_+] = 1,
$$
while, when~$m$ is even, the pairing is
$$
[P_\pm]\cap [P_\pm] = 1,\qquad\text{and}\qquad [P_-]\cap [P_+] = 0.
$$
In particular, when~$m$ is even, any two linear~$\bbP^m$s of the 
same chirality must intersect.  Consequently, the classes~$r[P_\pm]$
with~$r>1$ contain only singular cycles (and single linear $\bbP^m$s of
multiplicity~$r$).
\end{remark}

\begin{remark}
When~$m=2$, the exceptional isomorphism~$\SU(4) = \Spin(6)$
leads to the isomorphism of symmetric spaces~$Q_4 = \Gr(2,4)$.
Thus, the results of this section for~$Q_4$ have already been covered
in the treatment of the Grassmannians.
\end{remark}

\subsection{Isotropic Grassmannians}
\label{ssec: isotropic Gr}

As in the previous section, fix the standard inner product on~$\C{2m}$
and consider the set~$N^+_m\subset\Gr(m,2m)$ consisting of the
isotropic $m$-planes of positive chirality.  This is a compact
manifold of complex dimension~$\frac12m(m{-}1)$ that is
homogeneous under the action of the group~$\SO(2m,\bbC)$.  
The maximal compact subgroup~$\SO(2m)\subset\SO(2m,\bbC)$ also
acts transitively on~$N^+_m$, with stabilizer isomorphic to~$\Un(m)$.
Thus, 
\begin{equation}
N^+_m = \frac{\SO(2m)}{\Un(m)}\,,
\end{equation}
which exhibits~$N^+_m$ as one of the classical Hermitian symmetric
spaces.%
\footnote{In~\cite[Section 16]{MR21:1586}, the notation~$F_m$ is 
used for this symmetric space.  In~\cite{MR93a:20069}, this variety is 
called the \emph{spinor variety} of $\SO(2m,\bbC)$.}

\begin{remark}[An exceptional case]
When~$m=4$, there is the `exceptional isomorphism' (due to 
triality)~\cite[p.~519--520]{MR80k:53081}
\begin{equation}
N^+_4 = \frac{\SO(8)}{\Un(4)}
= \frac{\SO(8)}{\SO(2)\times\SO(6)} = Q_{6}\,,
\end{equation}
so this case has already been treated and the rigidity of the
extremal $3$-cycles has already been established.  Thus,
I will assume for the rest of this subsection that $m\ge 4$ and, 
whenever it is convenient, that~$m>4$.
\end{remark}

\subsubsection{Topology}
\label{ssec: isotropic topology}

In~\cite{MR21:1586}, the Poincar\'e polynomial of $N^+_m$ is found to be
\begin{equation}
p(N^+_m,t) = (1+t^2)(1+t^4)\dots(1+t^{2m-2}) 
           = 1 + t^2 + t^4 + 2\,t^6 + \dots,
\end{equation}
so $6$ is the lowest degree in which the rank of a homology group is
greater than~$1$ and this only happens when~$m\ge4$.  

As defined,~$N^+_m$ is a submanifold of~$\Gr(m,2m)$, and so inherits
bundles~$S$ and~$Q$ by pullback.  Since~$V\in N^+_m$ is a maximal
isotropic subspace, the inner product induces an 
isomorphism~$Q_V\simeq V^*$, so these bundles satisfy~$S^* = Q$.
  
It will be important to understand the tangent space to~$N^+_m$ at 
a general point~$V\in N^+_m$.  Now, at~$V\in N^+_m$,  
isomorphism~
$$
T_V\Gr(m,2m)\simeq Q_V\ot V^*\simeq V^*\ot V^* = S^2(V^*)\oplus\L^2(V^*)
$$
is canonical.  Under this isomorphism, the tangent space at~$V$ 
to~$N^+_m$ corresponds to the subspace~$\L^2(V^*)$.   In other words,
$TN^+_m \simeq \L^2(Q) \simeq \L^2(S^*)$. 

\begin{figure}
\includegraphics[width=\linewidth]{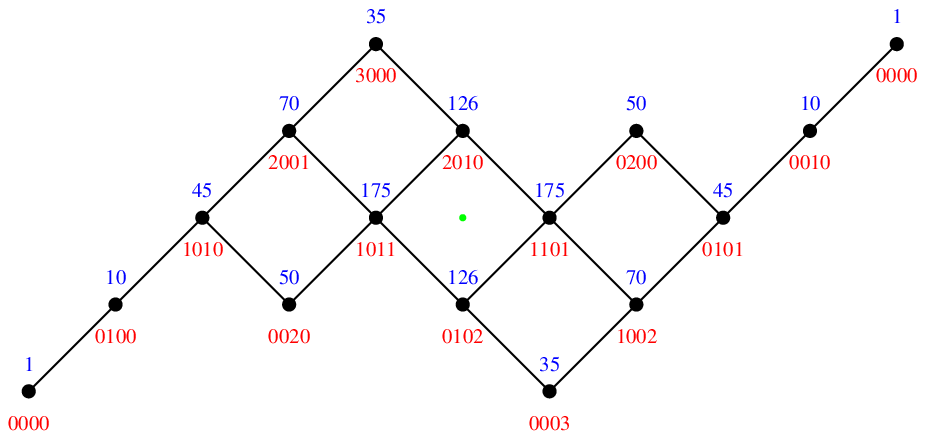}
\caption[The ideal poset for~$N_5^+=\SO(10)/\Un(5)$]
{\label{fig:SO10poset} The ideal poset for~$N_5^+=\SO(10)/\Un(5)$.
The upper label on each node is the dimension of the corresponding 
subrepresentation of~$\L^{*,0}\bigl(\euso(10)/\euu(5)\bigr)$ 
and the lower label is its highest weight 
as a representation of~$\SU(5)$. }
\end{figure}

More detail about the topology and Schubert cell decomposition of~$N^+_m$
can be found in~\cite{MR55:2941}.  Complete information about 
the irreducible constituents of the exterior powers of the cotangent
bundle of~$N^+_m$ and the consequent structure of its
ideal poset is collected in a convenient form in~\cite{kKhT99}. 
The corresponding Hasse diagram for the case~$m=5$ is drawn 
in Figure~\ref{fig:SO10poset}.

\subsubsection{Ideals of degree~$3$}
\label{ssec: isotropic ideals deg 3}
I have not analyzed the minimal ideals in all dimensions
for the isotropic Grassmannian, so I will confine myself to studying 
the cases in the first interesting degree, that of ideals of
degree or codegree equal to~$3$.

The first task is to describe the irreducible decomposition
of~$\L^{3,0}(\eum)$ under the action of~$K = \Un(m)$.  Fortunately,
this is relatively easy.  The above description of the tangent
space to~$N^+_m$ implies that, as a
representation of~$\Un(m)$, the space~$\eum$ is isomorphic to the 
representation~$\L^2(\C{m}) = \bbS_{(1,1)}(\C{m})$ 
associated to the standard representation of~$\Un(m)$ on~$\C{m}$.  
Then a little work with  multiplicity formulae from~\cite{MR93a:20069} 
shows that
\begin{equation}\label{eq: L3 decomp of TN}
\L^3\bigl(\L^2(\C{m})\bigr) 
\simeq \L^3\bigl(\bbS_{(1,1)}(\C{m})\bigr)
\simeq \bbS_{(2,2,2)}(\C{m})\oplus \bbS_{(3,1,1,1)}(\C{m}).
\end{equation}
These latter two representations are irreducible and their 
dimensions are given by~\cite[Theorem~6.3 or Exercise~6.4]{MR93a:20069} as:
\begin{equation}
\begin{split}
\dim\bbS_{(2,2,2)}(\C{m})  &=\frac{m^2(m{-}1)^2(m{-}2)(m{+}1)}{144},\\
\dim\bbS_{(3,1,1,1)}(\C{m})&= \frac{m(m^2{-}1)(m^2{-}4)(m{-}3)}{72}.
\end{split}
\end{equation}

Let~$\cI_{(2,2,2)}$ and $\cI_{(3,1,1,1)}$, respectively, 
denote the exterior differential systems on~$N^+_m$ generated in 
degree~$3$ by the sections of $\bbS_{(2,2,2)}(S)\subset\L^3(T^*N^+_m)$ 
and $\bbS_{(3,1,1,1)}(S)\subset\L^3(T^*N^+_m)$.

\subsubsection{Integral elements}
\label{ssec: isotropic ideals deg 3 int elem}

The following linear algebra lemma identifies
the integral elements of dimension three or more for each of the two
$\SL(m,\bbC)$-invariant subspaces of~$\L^3\bigl(\L^2((\C{m})^*)\bigr)$. 

\begin{lemma}\label{lem:integrals-of-ideals-in-L3N+m}
Any subspace $E\subset\L^2(\C{m})$ of dimension~$3$ or more on which
all of the elements of~$\bbS_{(2,2,2)}\bigl((\C{m})^*\bigr)$ vanish
is of the form~$E = L\w V$ where~$L\subset\C{m}$ is a line 
and~$V\subset\C{m}$ is a subspace containing~$L$ whose dimension is
one more than that of~$E$.

Any subspace $E\subset\L^2(\C{m})$ of dimension~$3$ or more on which
all of the elements of~$\bbS_{(3,1,1,1)}\bigl((\C{m})^*\bigr)$ vanish
has dimension~$3$ and is of the form~$E = \L^2(W)$ 
where~$W\subset\C{m}$ is a subspace of dimension~$3$.
\end{lemma}

\begin{proof}
Let~$\eta: \L^2(\C{m})\to \L^2(\C{m})$ be the identity map.
For any basis~$\vb_1,\dots,\vb_m$ of~$\C{m}$, 
write~$\eta = \frac12\eta^{ab}\vb_a\w\vb_b$, where~$\eta^{ab} 
= -\eta^{ba}$ are 1-forms on~$\L^2(\C{m})$.  
Note that~$\{\eta^{ab}\mid a<b\}$ is a basis for the dual space.  

Now,~$\bbS_{(2,2,2)}(\C{m})$ occurs as a constituent
of~$S^2(\C{m})^{\ot 3}$, but~$\bbS_{(3,1,1,1)}(\C{m})$ 
does not.  Consequently, the $3$-forms of the form
\begin{equation}\label{eq:222-generators}
\psi(X,Y,Z) = X_{i_1i_2}Y_{i_3i_4}Z_{i_5i_6}\,
     \eta^{i_1i_3}\w \eta^{i_2i_5}\w \eta^{i_4i_6}\,,
\end{equation}
when~$X$, $Y$, and $Z$ are symmetric in their indices,
must lie in~$\bbS_{(2,2,2)}\bigl((\C{m})^*\bigr)
\subset \L^3\bigl(\L^2((\C{m})^*)\bigr)$.  Taking,
as a particular example, $X_{11} = Y_{22} = Z_{33} = 1$ and 
all other~$X_{ij}$, $Y_{ij}$ and~$Z_{ij}$ equal to zero yields
$\psi(X, Y, Z) = \eta^{12}\w\eta^{13}\w\eta^{23}\not=0$, so
the span of the~$\psi(X,Y,Z)$ is nontrivial.  Since
this span is invariant under~$\GL(m,\bbC)$ and since 
$\bbS_{(2,2,2)}\bigl((\C{m})^*\bigr)$ is irreducible, 
this span must be all of this subspace.  Thus,
$\bbS_{(2,2,2)}\bigl((\C{m})^*\bigr)$ is the span of the
$3$-forms of the form~\eqref{eq:222-generators}.

Now suppose that~$E\subset\L^2(\C{m})$ 
is a subspace of dimension~$d\ge3$
on which all of the forms 
in~$\bbS_{(2,2,2)}\bigl((\C{m})^*\bigr)$ vanish.
My goal is to show that there is a unique line~$L\subset\C{m}$ so 
that~$E\subset L\w\C{m}$ and, conversely, that all of these
$3$-forms vanish on~$L\w\C{m}$ for any line~$L\subset\C{m}$.

The converse assertion is easy, so let me do this first.  Since all of
the conditions are invariant under the action of~$\GL(m,\bbC)$, 
it suffices to prove this for the case that~$L = \bbC\,\vb_1$.  
In this case,~$E = \bbC\,\vb_1\w\C{m}$ has dimension~$m{-}1$ and is
defined by the equations~$\eta^{ab} = 0$ when~$1<a<b$, so suppose that
all of these $1$-forms have been set to zero.  Then for~$\psi$ of the
form~\eqref{eq:222-generators}, all of the terms vanish unless
exactly one entry of each of the pairs~$(i_1,i_3)$, $(i_2,i_5)$, and
$(i_4,i_6)$ is equal to~$1$.  Moreover, the entries not equal to~$1$
in these pairs must all be distinct.  

Clearly, it suffices to treat the case where all of these entries 
are drawn from the set~$\{1,2,3,4\}$.  There are eight possible ways
to assign the value of~$1$ to one element of each of the pairs 
$(i_1,i_3)$, $(i_2,i_5)$, and~$(i_4,i_6)$.  In six of those ways, 
two of these indices will enter the same~$X$, $Y$, or~$Z$ coefficient 
and the corresponding sub-sum will vanish.  For example, when
$i_1 = i_2= i_4 = 1$, the part of the sum in~\eqref{eq:222-generators} 
corresponding to this choice is the sub-sum 
$X_{11}Y_{i_31}Z_{i_5i_6}\,\eta^{1i_3}\w\eta^{1i_5}\w \eta^{1i_6}$,
which vanishes, since~$Z_{i_5i_6}$ is symmetric 
and~$\eta^{1i_5}\w \eta^{1i_6}$ is skewsymmetric in the pair~${i_5,i_6}$.   
The two exceptional configurations are~$i_1=i_4=i_5=1$ and~$i_2=i_3=i_6=1$,
but these two sub-sums cancel:
\begin{multline*}
    X_{i_11}Y_{1i_4}Z_{i_51}\,\eta^{i_11}\w \eta^{1i_5}\w \eta^{i_41}
  + X_{1i_2}Y_{i_31}Z_{1i_6}\,\eta^{1i_3}\w \eta^{i_21}\w \eta^{1i_6}\\
= \quad X_{1i_1}Y_{1i_4}Z_{1i_5}\,\eta^{1i_1}\w \eta^{1i_5}\w \eta^{1i_4} 
  - X_{1i_2}Y_{1i_3}Z_{1i_6}\,\eta^{1i_3}\w \eta^{1i_2}\w \eta^{1i_6}= 0.
\end{multline*}
Thus, all the forms in~$\bbS_{(2,2,2)}\bigl((\C{m})^*\bigr)$ 
vanish on $E = L\w\C{m}$, as desired.

Now suppose that~$E\subset\L^2(\C{m})$ is an integral element 
of~$\bbS_{(2,2,2)}\bigl((\C{m})^*\bigr)$ and that $\dim E\ge 3$.   
Then~$E$ will be defined by some set of linear relations among 
the 1-forms~$\eta^{ab}$.  (By hypothesis, however, at least three
of the~$\eta^{ab}$ are linearly independent on~$E$.)
My goal is to show that one can choose the basis~$\vb_1,\dots,\vb_m$ 
so that these relations include~$\eta^{ab} = 0$ for~$a,b>1$.  

Suppose that the basis~$\vb_1,\dots,\vb_m$ has been chosen so that the 
maximum number, say~$p-1\ge0$, of the forms~$\eta^{12},\dots,\eta^{1m}$ 
are linearly independent on~$E$.  Clearly,~$2\le p\le \dim E + 1$.
By making a change of basis in~$\vb_2,\dots,\vb_m$, I can 
assume that~$\eta^{12}\w\dots\w\eta^{1p}\not=0$ on~$E$ but 
that~$\eta^{1a} = 0$ on~$E$ for~$a>p$.     

I claim that the maximality property implies that~$\eta^{ab}\equiv 0 
\mod \eta^{12},\dots,\eta^{1p}$ whenever~$a>p$.  
To see this, let~$\lambda_2,\dots,\lambda_m$ be parameters 
and consider the basis~$\vb_1,\vb^*_{2},\dots,\vb^*_m$ 
defined by~$\vb^*_a = \vb_a - \lambda_a\,\vb_1$ for~$a>1$.  Then
$$
\eta = \frac12\,\eta^{ab}\,\vb_a\w\vb_b 
=\sum_{1<a} (\eta^{1a}+\lambda_b\eta^{ba})\,\vb_1\w\vb^*_a
+ \frac12\sum_{1<a,b} \eta^{ab}\,\vb^*_a\w\vb^*_b\,.
$$
Suppose that there exist~$q>p$ and~$r>1$ so that~
$\eta^{12}\w\dots\w\eta^{1p}\w\eta^{qr}\not=0$.  Then set~$\lambda_a = 0$
for~$a\not=r$ and~$\lambda_r = t$ and consider the expansion
\begin{equation*}
(\eta^{12}{+}t\eta^{r2})\w\dots\w(\eta^{1p}{+}t\eta^{rp})\w
(\eta^{1q}{+}t\eta^{rq})
= t\,\eta^{12}{\w}\dots{\w}\eta^{1p}{\w}\eta^{rq} + O({t}^2)
\end{equation*}
Clearly, there will be a nonempty open set of values for~$t$
for which the left hand side of this equation will be nonzero, thus
contradicting the maximality of~$p$.  

It follows immediately that~$p\ge 3$. 
Now, so far, no use has been made of the hypothesis that~$E$ be an
integral element of~$\bbS_{(2,2,2)}\bigl((\C{m})^*\bigr)$.  Its first
use is to show that~$p\ge 4$.  This follows because,
as has already been noted, one of the elements 
of~$\bbS_{(2,2,2)}\bigl((\C{m})^*\bigr)$
is~$\eta^{12}\w\eta^{13}\w\eta^{23}$.
Since~$\eta^{12}\w\eta^{13}\not=0$, it follows that~$\eta^{23}$ is a
linear combination of~$\eta^{12}$ and~$\eta^{13}$.  If it were true that
$p = 3$, then all of the forms~$\eta^{ab}$ with~$a>3$ would also be linear
combinations of~$\eta^{12}$ and~$\eta^{13}$, so there could not be three
linearly independent $1$-forms among the~$\eta^{ab}$.  Thus, $p\ge 4$, as
claimed.  

Now, the same argument that showed that~$\eta^{12}\w\eta^{13}\w\eta^{23}$
is in $\bbS_{(2,2,2)}\bigl((\C{m})^*\bigr)$ shows that 
$\eta^{1a}\w\eta^{1b}\w\eta^{ab}$ is in 
$\bbS_{(2,2,2)}\bigl((\C{m})^*\bigr)$
for all~$1<a<b$.  In particular, it follows that $\eta^{ab}\equiv 0
\mod~\eta^{12},\dots,\eta^{1p}$ for all~$a,b\le p$.  Combined with
the previous argument, showing that $\eta^{ab}\equiv 0
\mod~\eta^{12},\dots,\eta^{1p}$ when either~$a$ or~$b$ is greater than~$p$,
this shows that $\eta^{12},\dots,\eta^{1p}$ must actually be a basis
for the 1-forms on~$E$, i.e., $p = \dim E + 1$.  

Now, the fact that~$\eta^{1a}\w\eta^{1b}\w\eta^{ab} = 0$ on~$E$ for~$a,b\le
p$
combined with the skew-symmetry $\eta^{ab} = -\eta^{ba}$
implies that there exist unique numbers~$A_{ab}$ for~$1<a\not=b\le p$ 
so that
$$
\eta^{ab} = A^{ab}\eta^{1a} - A^{ba}\eta^{1b}.
$$
Moreover, for any distinct~$a,b,c$ satisfying~$1<a,b,c\le p$, the 
formula~\eqref{eq:222-generators} shows that the form
$$
\eta^{1a}\w\eta^{1b}\w\eta^{ac} + \eta^{1a}\w\eta^{1c}\w\eta^{ab} 
$$
lies in~$\bbS_{(2,2,2)}\bigl((\C{m})^*\bigr)$, so the fact that this
vanishes on~$E$ implies that $A^{ca} = A^{ba}$.  It follows that there
are constants~$A^a$ so that~$A^{ba} = A^a$ for all~$b\not=a$.  
Consequently, $\eta^{ab} = A^{b}\eta^{1a} - A^{a}\eta^{1b}$ for~
$1<a,b\le p$.  Thus, by replacing~$\vb_1$
by~$\vb_1 + A^b\,\vb_b$, I get a new basis in which~$\eta^{ab}=0$ holds
on~$E$ for~$1<a,b\le p$, so I assume this from now on.

Next, taking~$a,b,c$ satisfying~$1<a<b\le p<c$, the above form
simplifies on~$E$ to~$\eta^{1a}\w\eta^{1b}\w\eta^{ac}$.  Since this
must vanish, it follows that~$\eta^{ac}\equiv 0 \mod \eta^{1a},\eta^{1b}$
for all such triples. Since~$p\ge 4$, this implies 
$\eta^{ac}\equiv 0 \mod \eta^{1a}$ whenever~$1<a\le p<c$. Thus, 
set~$\eta^{ac} = B^{ac}\eta^{1a}$ for some quantities~$B^{ac}$.  

Now, observe that the sum
$$
 \eta^{1a}\w\eta^{1b}\w\eta^{dc} 
+\eta^{1d}\w\eta^{1b}\w\eta^{ac}
+\eta^{1a}\w\eta^{1c}\w\eta^{db} 
+\eta^{1d}\w\eta^{1c}\w\eta^{ab}
$$  
is in~$\bbS_{(2,2,2)}\bigl((\C{m})^*\bigr)$, and take~$1<a<d\le p$
and~$1<b\le p<c$ with~$b$ not equal to either~$a$ or~$d$. (This last
is possible since~$p\ge4$.) The last two terms in the sum vanish 
since~$\eta^{1c}=0$ and the first two terms simplify to~$(B^{dc}-B^{ac})\,
\eta^{1a}\w\eta^{1b}\w\eta^{1d}$.  Since this must vanish, it follows
that~$B^{ac} = B^c$ for some constants~$B^c$ when~$c>p$.  
Thus,~$\eta^{ac}=B^c\eta^{1a}$. It follows that, 
by replacing~$\vb_1$ by~$\vb_1+B^c\vb_c$, 
I can arrange that~$\eta^{ac} = 0$
whenever~$1<a\le p<c$, so assume this from now on.

Finally, go back to the above sum and assume~$1<a<b\le p<d<c$. Then
all of the terms except the first are zero.  The vanishing of the first
term~$\eta^{1a}\w\eta^{1b}\w\eta^{dc}$ implies that~$\eta^{dc}$ 
is a linear combination of $\eta^{1a}$ and~$\eta^{1b}$ for any distinct
pair of indices~$a$ and~$b$ satisfying~$1<a<b\le p$.  Since~$p\ge 4$, this
forces~$\eta^{dc} = 0$.  

Thus~$E$ satisfies the relations~$\eta^{ab} = 0$ for $1<a,b\le m$ 
and for $a=1$ and~$b>p = \dim E + 1$.  
It follows that~$\vb_1\w\vb_2,\dots,
\vb_1\w\vb_p$ is a basis for~$E$, as desired.

Now, I will compute the integral elements 
of~$\bbS_{(3,1,1,1)}\bigl((\C{m})^*\bigr)$.  My goal is
to show that any integral element~$E$ of dimension~$3$ or more is actually
of the form~$E = \L^2(W)$ for some (unique) $3$-dimensional 
subspace~$W\subset\C{m}$.  

First, I must describe a set of generators 
of~$\bbS_{(3,1,1,1)}\bigl((\C{m})^*\bigr)$.  Now, it is easy to calculate
that $\bbS_{(3,1,1,1)}\bigl((\C{m})^*\bigr)$ occurs as a constituent of
$\L^4\bigl((\C{m})^*\bigr)\ot S^2\bigl((\C{m})^*\bigr)$, but that
$\bbS_{(2,2,2)}\bigl((\C{m})^*\bigr)$ does not.  Consequently,
every $3$-form of the form
\begin{equation}\label{eq:3111-generators}
\psi(X,Y) = X_{i_1i_2i_3i_4}Y_{j_1j_2}\,
     \eta^{i_1i_2}\w \eta^{i_3j_1}\w \eta^{i_4j_2}
\end{equation}
when~$X$ is skewsymmetric in its indices and ~$Y$ is symmetric
in its indices, must actually 
lie in~$\bbS_{(3,1,1,1)}\bigl((\C{m})^*\bigr)$.
Taking, as particular examples,~$Y_{ij} = \delta^1_i\delta^1_j$ and
letting~$X_{ijkl}$ be zero unless~$\{i,j,k,l\} = \{1,2,3,4\}$
while~$X_{1234} = 1$, gives~$\psi(X,Y) 
= 24\,\eta^{12}\w\eta^{13}\w\eta^{14}\not=0$, 
so it follows that the~$\psi(X,Y)$ span a nontrivial subspace 
of~$\bbS_{(3,1,1,1)}\bigl((\C{m})^*\bigr)$.  Since this subspace is
evidently invariant under~$\GL(m,\bbC)$ and 
since~$\bbS_{(3,1,1,1)}\bigl((\C{m})^*\bigr)$ is irreducible, it 
follows that the forms~$\psi(X,Y)$ of the form~\eqref{eq:3111-generators}
must span~$\bbS_{(3,1,1,1)}\bigl((\C{m})^*\bigr)$.

Now suppose that~$E\subset\L^2(\C{m})$ is an integral element of
$\bbS_{(3,1,1,1)}\bigl((\C{m})^*\bigr)$ whose dimension is at least~$3$.
Just as in the first part of the argument, 
choose a basis~$\vb_1,\dots,\vb_m$
of~$\C{m}$ so that the maximum number, say,~$p{-}1$, 
of~$\{\eta^{12},\dots,\eta^{1m}\}$ are linearly independent on~$E$ and
so that~$\eta^{1a}=0$ for~$a>p$.  The argument given in the first half of
the proof shows that~$p\ge 3$, but the fact 
that~$\eta^{12}\w\eta^{13}\w\eta^{14}$ lies 
in~$\bbS_{(3,1,1,1)}\bigl((\C{m})^*\bigr)$ implies that~$p<4$.  
Thus, $p=3$.  As before, the fact that~$p=3$ is maximal implies that
$\eta^{ab}\equiv0\mod\eta^{12},\eta^{13}$ whenever~$a>3$.  In particular,
all the 1-forms on~$E$ must be linear combinations of~$\eta^{12}$,
$\eta^{13}$, and~$\eta^{23}$.  Consequently, $\dim E\le 3$, but since
$\dim E\ge 3$ by hypothesis, $\dim E = 3$.  Moreover, since there
must be at least three linearly independent forms on~$E$, it follows
that~$\eta^{12}\w\eta^{13}\w\eta^{23}\not=0$. 

Now, the same argument as showed that there cannot be more than two
independent forms among the~$\eta^{1a}$ shows that there cannot
be more than two independent forms 
among the~$\eta^{2a}$ or the~$\eta^{3a}$.
In particular, for~$a>3$, the $1$-form~$\eta^{2a}$ must be a linear
combination of~$\eta^{21}$ and~$\eta^{23}$.  However, I have already
shown that it must also be a linear combination of~$\eta^{12}$
and~$\eta^{13}$.
Consequently~$\eta^{2a}$ must simply be a multiple of~$\eta^{21}$,
say~$\eta^{2a} = A^{2a}\,\eta^{21}$ for~$a>3$. 
Similarly,~$\eta^{3a}$ must simply be a multiple of~$\eta^{31}$,
say~$\eta^{3a} = A^{3a}\,\eta^{31}$ for~$a>3$.

Now, replacing~$\vb_1$ by~$\vb_1+\sum_{a>3}A^{2a}\vb_a$ produces a new
basis for which~$\eta^{2a} = 0$ for~$a>3$, so assume that this has been 
done.  Now, for each~$a>3$, consider~$\psi(X,Y)$ as 
in~\eqref{eq:3111-generators} where~$X_{123a}=1$ while~$X_{ijkl}=0$
unless~$\{i,j,k,l\} = \{1,2,3,a\}$ and $Y_{23} = Y_{32} = 1$ while
$Y_{ij} = 0$ unless $\{i,j\} = \{2,3\}$.  The result is
$$
\psi(X,Y) = 4\,\eta^{12}\w\eta^{23}\w\eta^{3a} 
          - 4\,\eta^{13}\w\eta^{23}\w\eta^{2a}.
$$
Since~$\eta^{2a} = 0$ on~$E$, the vanishing of~$\psi(X,Y)$ on~$E$
forces~$\eta^{3a}$ to be a linear combination 
of~$\eta^{12}$ and~$\eta^{23}$.
Since~$\eta^{3a} = A^{3a}\,\eta^{31}$, 
this forces~$A^{3a} = 0$, i.e.,
$\eta^{3a}$ vanishes on~$E$.

If~$m=4$, it has now been demonstrated that, 
for any integral $3$-dimensional
integral element of~$\bbS_{(3,1,1,1)}\bigl((\C{m})^*\bigr)$, there is 
a basis~$\vb_1,\dots,\vb_m$ so that~$E$ is defined by the 
equations~$\eta^{ab} = 0$ when~$a>3$.

If~$m>4$ assume~$a>b>3$ and consider~$\psi(X,Y)$ where~$X_{12ab}=1$ 
while~$X_{ijkl}=0$ unless~$\{i,j,k,l\} = \{1,2,a,b\}$ and $Y_{33}=1$ 
while $Y_{ij} = 0$ unless $i = j = 3$.  Then, on~$E$, $\psi(X,Y) 
= 4\,\eta^{13}\w\eta^{23}\w\eta^{ab} = 0$.  Now, permuting~$(1,2,3)$ 
in this construction shows that
$$
\eta^{13}\w\eta^{23}\w\eta^{ab} = \eta^{23}\w\eta^{12}\w\eta^{ab} 
= \eta^{12}\w\eta^{13}\w\eta^{ab} = 0.
$$
This implies that~$\eta^{ab} = 0$, on~$E$, as desired.

Finally, if there is a basis~$\vb_1,\dots,\vb_m$ so that~$E$
is defined by~$\eta^{ab} = 0$ when~$a>3$, 
then it is clear that~$\psi(X,Y)$
vanishes on~$E$ for all~$X$ and~$Y$, so that~$E$ is, indeed, an integral
element of~$\bbS_{(3,1,1,1)}\bigl((\C{m})^*\bigr)$.
\end{proof}

\subsubsection{Integral varieties}
\label{sssec: isotropic deg 3 int var}
 
The next two propositions describe the integral manifolds of 
the exterior differential systems~$\cI_{(2,2,2)}$ and $\cI_{(3,1,1,1)}$.

Before stating the first of these two propositions, I need
to describe a family of projective spaces~$\bbP^{m-1}$
that are embedded in~$N^+_m$.

\begin{example}[The chiral double fibration]
\label{ex: chiral dbl fib}
If~$H\subset\C{2m}$ is any isotropic $(m{-}1)$-plane, 
it lies in two distinct isotropic $m$-planes: 
$H^+$, of positive chirality, and~$H^-$, of negative chirality.
Thus, if~$N_{m-1}(\C{2m})\subset\Gr(m{-}1,2m)$ denotes the
space of isotropic $(m{-}1)$-planes in~$\C{2m}$, this space
is the apex of a double fibration
$$
\begin{array}{ccccc}
&& N_{m-1}(\C{2m})\\
&\swarrow & &\searrow\\
 N^+_m&&&& N^-_m
\end{array}.
$$
Let~$S\subset\C{2m}$ be any isotropic $m$-plane 
of negative chirality.  Then there is a canonical embedding~$\iota_S$ 
of~$\bbP(S^*)\simeq\bbP^{m-1}$ into~$N^+_m$ defined by sending each
hyperplane~$H\in \bbP(S^*) = \Gr(m{-}1,S)$ to its positive chirality
extension~$H^+ = \iota_S(H)$. 
\end{example}

\begin{proposition}\label{prop:222-integrals}
For any~$S\in N^-_m$, the projective space 
$\iota_S\bigl(\bbP(S^*)\bigr)$ is a {\upshape(}maximal\/{\upshape)}
integral manifold of~$\cI_{(2,2,2)}$. 

Conversely, if~$X\subset N^+_m$ is an irreducible variety of 
dimension at least~$3$ that is an integral variety of~$\cI_{(2,2,2)}$,
then there exists a unique~$S\in N^-_m$ so that~$X$ is contained 
in~$\iota_S\bigl(\bbP(S^*)\bigr)$.
\end{proposition}

\begin{proof}
The first task (which will be needed in the next proposition as well),
is to establish the equations of the moving frame for submanifolds
of~$N^+_m$.

Define~$\SO(2m,\bbC)$ be the subgroup of~$\SL(2m,\bbC)$
consisting of the matrices~$\us$ that satisfy
\begin{equation}\label{eq:defining-SO2m}
{}^t\us \begin{pmatrix}0_{m}&\I_{m}\\ \I_{m}&0_{m}\end{pmatrix} \us
= \begin{pmatrix} 0_{m} & \I_{m}\\ \I_{m} & 0_{m}\end{pmatrix}.
\end{equation} 
Also, let~$F\subset\GL(2m,\bbC)$ denote the set of matrices~$\vs$
that satisfy
\begin{equation}\label{eq:defining-F}
{}^t\vs\,\vs
= \begin{pmatrix} 0_{m} & \I_{m}\\ \I_{m} & 0_{m}\end{pmatrix}
\end{equation}
and~$\det(\vs) = \iC^m$.  Evidently, $F$ is an orbit of~$\SO(2m,\bbC)$
acting on~$\GL(2m,\bbC)$ on the right.  
I will regard~$\vs:F\to\GL(2m,\bbC)$
as a matrix-valued function and denote its columns as
$$
\vs = (\vs_1\quad\dots\quad \vs_{m} \quad \vs^1\quad\dots\quad \vs^{m})
$$
where~$\vs_i,\vs^i:F\to\C{2m}$ are regarded as (holomorphic) mappings.

Define
$$
\pi(\vs) = [\vs_1\w\dots\w\vs_m]\,,
$$
so that~$\pi$ is a surjective submersion~$\pi:F\to N^+_m$.
The fibers of~$\pi$ are the orbits of the parabolic subgroup~$P\subset
\SO(2m,\bbC)$ consisting of elements of the form
\begin{equation}\label{eq: elements of P}
\us = \begin{pmatrix} A & AB\\ 0_m & {}^t\!A^{-1}\end{pmatrix}
\qquad\text{for~$A\in \GL(m,\bbC)$ and~$B = -{}^t\!B\in\C{m,m}$.}
\end{equation}
Thus,~$\pi:F\to N^+_m$ is a principal right~$P$-bundle over~$N^+_m$.

In accordance with the usual moving frame conventions, 
write the structure equations as
\begin{equation} \label{eq:F-structure-equations}
\d\vs = \d (\vs_i\ \vs^i) = (\vs_j\ \vs^j)
\begin{pmatrix}\alpha^j_i & \gamma^{ji}\\ 
                \beta_{ji}& -\alpha^i_j \end{pmatrix}
= \vs\,\theta
\end{equation}
where~
\begin{equation}\label{eq:F-structure-relations}
\beta_{ji} = -\beta_{ij}\quad\text{and}\quad\gamma^{ji}=-\gamma^{ij}\,,
\end{equation}
but the components of~$\alpha$, $\beta$, and~$\gamma$ are otherwise
linearly independent. The relations~\eqref{eq:F-structure-relations} 
follow in the usual way from the exterior derivative 
of~\eqref{eq:defining-F}.  The structure equation 
$\d\theta = -\theta\w\theta$ holds since~$\theta = \vs^{-1}\,\d \vs$.
These expand to
\begin{equation}\label{eq:F-structure-abc-relations}
\begin{split}
\d\alpha^i_j  &= -\alpha^i_k\w\alpha^k_j-\gamma^{ik}\w\beta_{kj}\,,\\
\d\beta_{ij}  &= -\beta_{ik}\w\alpha^k_j +\alpha^k_i\w\beta_{kj}\,,\\
\d\gamma^{ij} &= -\alpha^i_k\w\gamma^{kj} +\gamma^{ik}\w\alpha^j_k\,.
\end{split}
\end{equation}

Now suppose that~$X\subset N^+_m$ is an irreducible integral
variety of~$\cI_{(2,2,2)}$ of dimension~$d\ge 3$, 
and let~$X^\circ\subset X$ denote its
smooth locus, which is an embedded submanifold of~$N^+_m$.
For every~$V\in X^\circ$, the 
tangent space~$T_VX$ is an integral element of~$\cI_{(2,2,2)}$ of
dimension~$d\ge3$. By Lemma~\ref{lem:integrals-of-ideals-in-L3N+m}, 
it follows that, for every~$V\in X^\circ$, there exists a~$\vs\in F$
so that
\begin{enumerate}
\item $V$ is spanned by~$\vs_1,\dots,\vs_m$, and
\item $T_VX$ is spanned by 
$\dl\vs^1\dr\w\dl\vs^2\dr,\dl\vs^1\dr\w\dl\vs^3\dr,\dots,
 \dl\vs^1\dr\w\dl\vs^{d+1}\dr$.
\end{enumerate}
Let~$F(X^\circ)\subset F$ denote the set of such~$\vs$ as~$V$
ranges over~$X^\circ$.  Then~$\pi:F(X^\circ)\to X^\circ$ is
a principal~$G$-bundle over~$X^\circ$, where~$G\subset P$ is the 
subgroup consisting of the matrices of the form~\eqref{eq: elements of P}
with~${}^t\!A^{-1}$ in $P_1\cap P_{d+1}\subset\GL(m,\bbC)$.

By construction, the forms~$\beta_{12},\dots,\beta_{1(d{+}1)}$ are 
linearly independent on~$F(X^\circ)$ and span the $\pi$-semibasic 
$1$-forms, while $\beta_{1a}=0$ for~$a>d{+}1$ and~$\beta_{ij}=0$
when both~$i$ and~$j$ are bigger than~$1$.

This paragraph of the argument is necessary only if~$d<m{-}1$,
so suppose this is so for the moment.  Choose a pair~$(i,a)$ 
satisfying~$2\le i\le d{+}1<a\le m$ and differentiate 
the relation~$\beta_{ia}=0$.  By the structure equations, this is
$$
0 = \d \beta_{ia} = -\beta_{i1}\w\alpha^1_a\,.
$$
Since~$d\ge3$, and since~$\beta_{12},\dots,\beta_{1(d{+}1)}$ are
linearly independent,~$\alpha^1_a=0$ for~$a>d{+}1$.

Now choose a pair~$(i,j)$ with~$2\le i \not= j\le d{+}1$ 
and differentiate~$\beta_{ij}=0$.  The structure equations
give that
$$
0 = \d \beta_{ij} = -\beta_{i1}\w\alpha^1_j + \alpha^1_i\w\beta_{1j}\,.
$$
Equivalently, 
\begin{equation}\label{eq:F-ab-symmetry}
\alpha^1_i\w\beta_{1j} = \alpha^1_j\w\beta_{1i}\,.
\end{equation}  
Wedging this relation with~$\beta_{1i}$ gives 
$\alpha^1_i\w\beta_{1i}\w\beta_{1j}=0$ for all~$2\le i\not=j\le d{+}1$.
Again, because~$d\ge3$ and because 
$\beta_{12},\dots,\beta_{1(d{+}1)}$ are
linearly independent, it follows that~$\alpha^1_i\w\beta_{1i}=0$ 
for~$2\le i\le d{+}1$.  In particular, there exist functions~$\lambda_i$
on~$F(X^\circ)$ so that~$\alpha^1_i = \lambda_i\,\beta_{1i}$.
Substituting this back into \eqref{eq:F-ab-symmetry} and again using
the linear independence of~$\beta_{1j}$ and~$\beta_{1i}$, it follows
that~$\lambda_i+\lambda_j=0$ for all $2\le i\not=j\le d{+}1$.  Again,
since~$d\ge3$, this implies that~$\lambda_i=0$ for~$2\le i\le d{+}1$.

In other words, $\alpha^1_i=0$ for~$2\le i\le d{+}1$.  Since the
previous paragraph showed that~$\alpha^1_a=0$ for all~$a>d{+}1$, this
combines to give that~$\alpha^1_i=0$ for all~$i>1$.  This
vanishing together with the fact that~$\beta_{ij}=0$ for all~$i,j\ge2$
yield the congruences
$$
\d\vs_2\equiv\dots\equiv \d\vs_m\equiv \d\vs^1
  \equiv 0 \mod \vs_2\,,\dots,\vs_m\,,\vs^1\,.
$$
In other words, the mapping~$\sigma:F(X^\circ)\to N^-_m$ defined by
$\sigma(\vs) = [\vs_2\w\dots\w\vs_m\w\vs^1]$ is constant.  Let~$S\in N^-_m$
be this constant $m$-plane.    

By construction~$\pi(\vs) = [\vs_1\w\dots\w\vs_m]$ lies 
in~$\iota_S\bigl(\bbP(S^*)\bigr)$, so it follows that~$X^\circ$, and,
hence, $X$ lie in~$\iota_S\bigl(\bbP(S^*)\bigr)$, as desired.

That~$\iota_S\bigl(\bbP(S^*)\bigr)$ really is an integral 
variety of~$\cI_{(2,2,2)}$ follows immediately from the proof of the
first part.
\end{proof}

Now, by Lemma~\ref{lem:integrals-of-ideals-in-L3N+m}, 
there are no integral 
manifolds of~$\cI_{(3,1,1,1)}$ of dimension greater than~$3$.  
The following proposition classifies all of the $3$-dimensional integrals.

First, a definition.  For any isotropic subspace~$A\subset\C{2m}$, let
$N^+_m(A)\subset N^+_m$ denote the set of~$P\in N^+_m$ that contain~$A$.
Note that, if~$a = \dim A < m$, then $N^+_m(A)$ is a smooth subvariety 
of~$N^+_m$ that is isomorphic to~$N^+_{m-a}$.  

\begin{proposition}\label{prop:3111-integrals}
Let~$X\subset N^+_m$ be an irreducible variety of dimension~$3$
that is an integral variety of~$\cI_{(3,1,1,1)}$.  Then there exists
an isotropic $(m{-}3)$-plane~$A\subset\C{2m}$ so that~$X \subset N^+_m(A)$.  
In particular, if~$X$ is closed, then $X\simeq N^+_3\simeq\bbP^3$.
\end{proposition}

\begin{proof}
Recall the moving frame notation and constructions 
from the first part of the proof of Proposition~\ref{prop:222-integrals}.

Suppose now that~$X\subset N^+_m$ 
is an irreducible $3$-dimensional integral variety of of~$\cI_{(3,1,1,1)}$
and let~$X^\circ\subset X$ be its smooth locus, which is connected
since~$X$ is irreducible. 
By Lemma~\ref{lem:integrals-of-ideals-in-L3N+m}, for every~$V\in X^\circ$,
there exists a~$\vs\in F$ so that
\begin{enumerate}
\item $V$ is spanned by~$\vs_1,\dots,\vs_m$, and
\item The tangent space $T_VX^\circ$ is spanned by
       $\dl\vs^2\dr\w\dl\vs^3\dr,\>\dl\vs^3\dr\w\dl\vs^1\dr,
       \>\dl\vs^1\dr\w\dl\vs^2\dr$.
\end{enumerate}
Let~$F(X^\circ)\subset F$ denote the set of such~$\vs$ as~$V$
ranges over~$X^\circ$.  Then~$\pi:F(X^\circ)\to X^\circ$ is a 
principal right $G$-bundle over~$X^\circ$ where $G\subset P$ is the
subgroup consisting of the matrices of the form~\eqref{eq: elements of P}
with~${}^t\!A^{-1}$ in $P_3\subset\GL(m,\bbC)$.  Since~$G$ and~$X^\circ$
are each connected, it follows that~$F(X^\circ)$ is also connected. 

By construction, the $1$-forms~$\beta_{23},\beta_{31},\beta_{12}$
are linearly independent on~$F(X^\circ)$ and span the $\pi$-semibasic
$1$-forms, while $\beta_{ij} = 0$ if either~$i$ or~$j$ is
greater than~$3$.

Let~$i>3$ be fixed and differentiate the identities~$\beta_{i1}=\beta_{i2}
=\beta_{i3}=0$ using the structure equations.  The result is
equations of the form
$$
\begin{pmatrix} \alpha^1_i& \alpha^2_i & \alpha^3_i \end{pmatrix}
\w
\begin{pmatrix} 0 & -\beta_{12} & \beta_{31}\\
               \beta_{12} & 0 & -\beta_{23}\\
               -\beta_{31} & \beta_{23} & 0 \end{pmatrix}
= \begin{pmatrix} 0& 0 & 0 \end{pmatrix} .
$$
By the linear independence of~$\beta_{23},\beta_{31},\beta_{12}$,
it follows that~$\alpha^1_i = \alpha^2_i = \alpha^3_i = 0$.

This vanishing for all~$i>3$ implies 
$$
\d\vs_4\equiv\dots\equiv\d\vs_m \equiv 0 \mod \vs_4,\dots,\vs_m\,,
$$
i.e., the $(m{-}3)$-plane~$[\vs_4\w\dots\vs_m]$ is locally
constant on $F(X^\circ)$.  Since~$F(X^\circ)$ is connected, this
map must be constant.  Thus, let~$A\in \Gr(m{-}3,2m)$ be the isotropic
plane so that~$ [\vs_4\w\dots\vs_m]\equiv A$.
By construction,~$A\subset V$ for all~$V\in X^\circ$, so it follows
that~$X^\circ$ and, hence,~$X$ are subsets of~$N^+_m(A)$, as desired.
\end{proof}

These propositions allow characterizations of the extremal
classes in~$H_6(N^+_m)$ that are analogous 
to that of Schubert cycles in Grassmannians:

\begin{theorem}\label{thm:rigid3-cycles-in-isotropic-Grassmannians}
Suppose~$m\ge 4$.  Let~$C\subset\C{2m}$ be an
isotropic plane of dimension~$m{-}3$.  Fix~$A\in N^-_m$ and 
let~$P\subset \bbP(A^*)$ be a linearly embedded projective~$3$-space.
Define two $3$-dimensional subvarieties of~$N^+_m$ by
\begin{equation}
  X = \iota_A(P)
 \qquad\text{and}\qquad
  Y = N^+_m(C).
\end{equation}
Then~$[X]$ and~$[Y]$ are the generators of~$H^+_6(N^+_m,\bbZ)\simeq\bbZ^2$.  

Any irreducible~$Z\in\cZ^+_3(N^+_m)$ that satisfies~$[Z] = r[X]$ is
of the form~$Z = \iota_S(Z')$ where~$S\in N^-_m$ is fixed 
and~$Z'\subset \bbP(S^*)\simeq\bbP^{m-1}$ is an irreducible variety
of dimension~$3$ and degree~$r$.

Any irreducible~$Z\in\cZ^+_3(N^+_m)$ that satisfies~$[Z] = r[Y]$ 
is of the form~$Z = N^+_m(D)$ for some isotropic~$D\subset\C{2m}$
of dimension~$m{-}3$.
\end{theorem}

\begin{proof}
First of all, it follows by either~\cite{MR21:1586} or 
\eqref{eq: L3 decomp of TN} and the general results of Kostant
mentioned above that $b_6(N^+_m)=2$.  
Let~$\phi_1$ be the $\SO(2m)$-invariant
K\"ahler form on~$N^+_m$ whose cohomology class is a generator
of~$H^2(N^+_m,\bbZ)$.  By \eqref{eq: hss sum of Schur forms},
there is a sum of the form
$$
{\phi_1}^3 = \mu^{(2,2,2)}\,\phi_{(2,2,2)}+\mu^{(3,1,1,1)}\,\phi_{(3,1,1,1)}
$$
where~$\mu^{(2,2,2)}>0$ and $\mu^{(3,1,1,1)}>0$ and 
$\phi_{(2,2,2)}$ and~$\phi_{(3,1,1,1)}$ are positive $\SO(2m)$-invariant  
forms dual to the generalized 
Schubert cycles~$\sigma_{(2,2,2)^*}$ 
and~$\sigma_{(3,1,1,1)^*}$ of complex dimension~$3$ whose 
cohomology classes generate~$H^+_6(N^+_m)$.  

It follows that~$\sigma_{(2,2,2)^*}$ is an irreducible $3$-dimensional 
integral variety of~$\cI_{(3,1,1,1)}$ and, so, 
by Proposition~\ref{prop:3111-integrals}, must be of the form~$N^+_m(C)$ 
for some isotropic $m{-}3$ plane~$C$.  Thus~$\sigma_{(2,2,2)^*}$ is
homologous to~$Y$.  

It also follows that~$\sigma_{(3,1,1,1)^*}$ 
is an irreducible $3$-dimensional
integral variety of~$\cI_{(2,2,2)}$, and so, by
Proposition~\ref{prop:222-integrals}, must lie in 
$\iota_A\bigl(\bbP(A^*)\bigr)$ for some unique~$A\in N^-_m$.  
Since~$\sigma_{(3,1,1,1)^*}$ must be a generator of~$H^+_6(N^+_m)$, it
follows easily that it must be homologous to~$\iota_A(P)$,
where $P\subset\bbP(A^*)$ is a linearly embedded projective~$3$-space.
Thus,~$\sigma_{(3,1,1,1)^*}$ is homologous to~$X$.

Finally, if~$Z\in\cZ^+_3(N^+_m)$ is irreducible and satisfies~$[Z]=r[X]$,
then the integral of~$\phi_{(2,2,2)}$ over~$Z$ must be zero, 
so $\phi_{(2,2,2)}$ must vanish on~$Z$.  Thus~$Z$ is an integral manifold
of~$\cI_{(2,2,2)}$ and Proposition~\ref{prop:222-integrals} applies.

The argument when~$[Z] = r[Y]$ is similar.
\end{proof}

\begin{remark}[Homologies to integrals of~$\cI_{(2,2,2)}$]
When~$A\in N^-_m$, each linear subspace~$P_d\subset\bbP(A^*)$ 
of dimension~$d\ge 3$ determines 
a homology class~$\bigl[\iota_A(P_d)\bigr]\in H^+_{2d}(N^+_m,\bbZ)$
that displays similar quasi-rigidity.  In other words, 
if~$Z\in\cZ^+_{2d}(N^+_m)$ is irreducible and 
satisfies~$[Z] = r\,\bigl[\iota_A(P_d)\bigr]$ for some~$r>0$, then 
$Z = \iota_S(Z')$ where~$S\in N^-_m$ is fixed 
and~$Z'\subset \bbP(S^*)\simeq\bbP^{m-1}$ is an irreducible variety
of dimension~$d$ and degree~$r$.  The argument is left to the reader.
\end{remark}

The extremal cycles of codimension~$3$ in~$N^+_m$ also display rigidity.
I am not going to give all the details of this discussion, since most
of the methods of proof in each of the two cases I am going to consider 
will, by now, be familiar to the reader.  Instead, I will simply highlight
the points at which some interesting or different idea comes into play.

First, I need to recall an elementary fact about intersections
of maximal isotropic planes in~$\C{2m}$.  Namely, if~$P$ and~$Q$ 
lie in~$N_m = N^+_m\cup N^-_m$, then $P$ and~$Q$ lie in the same 
component of~$N_m$ if and only if the dimension of~$P\cap Q$ is
congruent to~$m$ modulo~$2$.  For a proof (which,
in any case, is not difficult), see~\cite[p.~735]{MR80b:14001}.

\begin{theorem}\label{thm:*3111-integrals}
Assume~$m\ge3$ and let $P\subset\C{2m}$ be an isotropic $m$-plane 
that lies in $N^+_m$ if $m$ is odd and in $N^-_m$ if $m$ is even.  
Let~
$$
\Sigma(P)= \left\{\ V\in N^+_m\mid \dim(V\cap P)\ge 3 \ \right\} .
$$
Then~$\Sigma(P)$ is of codimension~$3$ in~$N^+_m$ 
and represents the generalized
Schubert cycle~$\sigma_{(2,2,2)}$.  

Moreover, any irreducible variety
$X\subset N^+_m$ of codimension~$3$ that satisfies~$[X] 
= r\,\bigl[\Sigma(P)\bigr]$ is of the form~$X=\Sigma(Q)$ 
for some isotropic $m$-plane~$Q$ {\upshape(}that lies in $N^+_m$
if $m$ is odd and in $N^-_m$ if $m$ is even{\upshape)}. 
\end{theorem}

\begin{proof}
Using arguments that should, by now, be familiar, one sees that the
homology class of a codimension~$3$ irreducible cycle~$X\subset N^+_m$ 
is some multiple of~$[\sigma_{(2,2,2)}]$ if and only 
if the form~$\phi_{(3,1,1,1)^*}$ vanishes on~$X$, which is the same
as saying that, at any smooth point~$V\in X^\circ$, 
the normal space to~$T_VX$ in~$T_VN^+_m$ is an
integral element of~$\cI_{(3,1,1,1)}$.   

By Lemma~\ref{lem:integrals-of-ideals-in-L3N+m}, it follows that,
for every~$V\in X^\circ$, there is a~$\vs\in F$ so that 
\begin{enumerate}
\item $V$ is spanned by~$\vs_1,\dots,\vs_m$, and
\item $T_VX$ is spanned by the~$\dl v^a\dr\w \dl v^b\dr$,
where $1\le a<b\le m$ and~$b>3$.
\end{enumerate}
The set of all such~$\vs\in F$ as~$V$ ranges over~$X^\circ$ is
a principal right~$G$-bundle~$\pi:F(X^\circ)\to X^\circ$, where
$G\subset P$ is the subgroup consisting of the matrices of the 
form~\eqref{eq: elements of P} with~$A$ in $P_3\subset\GL(m,\bbC)$.  
Since~$G$ and~$X^\circ$ are each connected, 
it follows that~$F(X^\circ)$ is also connected. 

By construction, the $1$-forms~$\beta_{ab}$ with~$1\le a<b\le m$ and
$b>3$ are linearly independent and span the $\pi$-semibasic forms 
on~$F(X^\circ)$ while the forms~$\beta_{12}$, $\beta_{13}$, 
and~$\beta_{23}$ are all identically zero.  

To save writing, adopt the conventions that~$1\le i,j,k\le 3$
while~$4\le a,b,c\le m$.  Differentiating the relations ~$\beta_{ij}=0$
and applying the structure equations then yields the relations
$\alpha^a_i\w\beta_{aj} = \alpha^a_j\w\beta_{ai}$  (summation on~$a$).
Judicious use of Cartan's Lemma, together with the linear independence
of the~$\beta_{ai}$, implies that there exist functions~$T^{ab} = -T^{ba}$
on~$F(X^\circ)$ so that
$$
\alpha^a_i = T^{ab}\,\beta_{bi}\,.
$$

After computing how the functions~$T^{ab}$ vary on the fibers of~$\pi$,
one sees that the equations~$T^{ab}=0$ define a principal right
$G_1$-bundle $F_1\subset F(X^\circ)$ over~$X^\circ$ where~$G_1\subset P$
is the connected subgroup of matrices whose Lie algebra consists of
the matrices of the form
$$
\begin{pmatrix}
x^i_j & x^i_b & y^{ij} & y^{ib}\\
0 & x^a_b & -y^{ja} & 0\\
0 & 0 & -x^j_i & 0\\
0 & 0 & -x^j_a & -x^b_a
\end{pmatrix}.
$$

The relations~$\alpha^a_i=0$ hold on~$F_1$. Differentiating
these relations and applying Cartan's Lemma shows that the relations
$\gamma^{ab} = 0$ must also hold.

These identities combine to show that, on~$F_1$, 
$$
\d\vs_1\equiv\d\vs_2\equiv\d\vs_3\equiv\d\vs^4\equiv\cdots \equiv\d\vs^m 
\equiv 0 \mod \vs_1,\vs_2,\vs_3,\vs^4,\dots,\vs^m.
$$
In other words, the 
map~$[\vs_1\w\vs_2\w\vs_3\w\vs^4\w\dots\w\vs^m]:F_1\to N_m$ is locally
constant.  Since~$F_1$ is connected, this map must be globally constant.
Let~$Q\in N_m$ be its constant value.  Note that~$Q$ lies in~$N^+_m$
if $m$ is odd and in $N^-_m$ if $m$ is even. 

By construction~$\dim(Q\cap V) = 3$ for all~$V\in X^\circ$.  It follows
that~$X^\circ$, and, hence, $X$ lie in~$\Sigma(Q)$, as desired.  

Moreover,
examining the argument given shows that~$\Sigma(Q)$ is indeed
an integral manifold of~$\cI_{(3,1,1,1)^*}$ and has codimension~$3$
in~$N^+_m$.
\end{proof}

\begin{remark}[Singularity of~$\Sigma(P)$]
The variety~$\Sigma(P)$ is singular when~$m\ge 5$, since, in this
case, it will necessarily contain the (non-empty) locus of 
those~$V\in N^+_m$ that satisfy~$\dim(P\cap V)\ge 5$.  However, the proof
above shows that the smooth locus of~$\Sigma(P)$ must consist
of those~$V\in N^+_m$ that satisfy~$\dim(P\cap V)=3$.  
Thus, the singular locus of~$\Sigma(P)$ cannot be empty.
In particular, it follows from~Theorem~\ref{thm:*3111-integrals}
that, when~$m\ge 5$, no multiple 
of the homology class~$[\sigma_{(2,2,2)}]$ 
can be represented by a smooth, effective cycle. 
\end{remark}

Before stating the next proposition, I remind the reader that
$H_{2k}(Q_n,\bbZ)\simeq\bbZ$ when~$0\le 2k < n $.  In particular, 
the notion of \emph{degree} is unambiguous for a $k$-cycle in~$Q_n$
as long as~$2k<n$.  

\begin{theorem}\label{thm:*222-integrals}
Assume~$m\ge 4$.  Let~$Y\subset Q_{2m-2}$ be a subvariety
of dimension~$m{-}4$ and degree~$r$.  Let~$\Psi(Y)\subset N^+_m$ denote 
the set of~$V\in N^+_m$ satisfying~$\bbP(V)\cap Y\not=\emptyset$.

Then~$\Psi(Y)$ has codimension~$3$ in~$N^+_m$ and 
satisfies~$\bigl[\Psi(Y)\bigr]=r\bigl[\sigma_{(3,1,1,1)}\bigr]$.  

Moreover, any codimension~$3$ subvariety~$X\subset N^+_m$ that
satisfies~$[X] = r\bigl[\sigma_{(3,1,1,1)}\bigr]$ is~$\Psi(Y)$ 
for some subvariety~$Y\subset Q_{2m-2}$ of dimension~$m{-}4$ 
and degree~$r$.
\end{theorem}

\begin{proof}
Using arguments that should, by now, be familiar, one sees that the
homology class of a codimension~$3$ irreducible cycle~$X\subset N^+_m$ 
is some multiple of~$[\sigma_{(3,1,1,1)}]$ if and only 
if the form~$\phi_{(2,2,2)^*}$ vanishes on~$X$, which is the same
as saying that, at any smooth point~$V\in X^\circ$, 
the normal space to~$T_VX$ in~$T_VN^+_m$ is an
integral element of~$\cI_{(2,2,2)}$.   

By Lemma~\ref{lem:integrals-of-ideals-in-L3N+m}, it follows that,
for every~$V\in X^\circ$, there is a~$\vs\in F$ so that 
\begin{enumerate}
\item $V$ is spanned by~$\vs_1,\dots,\vs_m$, and
\item $T_VX$ is spanned by the~$\dl v^a\dr\w \dl v^b\dr$,
where $1\le a<b\le m$ and either~$a>1$ or~$b>4$.
\end{enumerate}
The set of all such~$\vs\in F$ as~$V$ ranges over~$X^\circ$ is
a principal right~$G$-bundle~$\pi:F(X^\circ)\to X^\circ$, where
$G\subset P$ is the subgroup consisting of the matrices of the 
form~\eqref{eq: elements of P} with~$A$ 
in~$P_1\cap P_4\subset\GL(m,\bbC)$.  
Since~$G$ and~$X^\circ$ are each connected, 
it follows that~$F(X^\circ)$ is also connected. 

By construction, the $1$-forms~$\beta_{ab}$ with~$1\le a<b\le m$ and
and either~$a>1$ or~$b>4$ are linearly independent 
and span the $\pi$-semibasic forms 
on~$F(X^\circ)$ while the forms~$\beta_{12}$, $\beta_{13}$, 
and~$\beta_{14}$ are all identically zero.  

To save writing, adopt the conventions that~$2\le i,j,k\le 4$
while~$5\le a,b,c\le m$.  Differentiating the relations ~$\beta_{1i}=0$
and applying the structure equations then yields relations of the form
$$
0 =\alpha^j_1\w\beta_{ji}+\alpha^a_1\w\beta_{ai}-\alpha^a_i\w\beta_{a1}\,,
\qquad\text{(summation on~$a$ and $j$).}
$$
for~$i = 2$, $3$, and~$4$.
Judicious use of Cartan's Lemma, 
together with the stated linear independence of the entries of~$\beta$, 
implies that there exist functions~$B^{ja}$, $B^{ab}$, and~$B^{ab}_i$
on~$F(X^\circ)$ so that
\begin{equation*}
\begin{split}
\alpha^i_1 &= B^{ib}\,\beta_{b1}\,,\\
\alpha^a_1 &= B^{ab}\,\beta_{b1}\,,\\
\alpha^a_i &= B^{ja}\,\beta_{ij}-B^{ba}\,\beta_{bi}+B^{ab}_i\,\beta_{b1}\,.
\end{split}
\end{equation*}
In particular, it follows from the structure equations that
$$
d\vs_1 
\equiv\left(\vs^a + B^{ja}\vs_j+B^{ba}\vs_b\right)\,\beta_{a1}\mod\vs_1\,.
$$
Consequently the mapping~$[\vs_1]:F(X^\circ)\to Q_{2m-2}$ 
is constant on the fibers of~$\pi$ and its differential has rank~$m{-}4$
everywhere.  Thus, $[\vs_1] = y\circ\pi$ where~$y:X^\circ\to Q_{2m-2}$ 
is a holomorphic map whose differential has rank~$m{-}4$ everywhere.

When~$X$ is an algebraic variety, it is not hard to argue that~$y$ 
is a rational map and then that the closure of~$y(X^\circ)$ in~$Q_{2m-2}$ 
is an algebraic variety~$Y$ of dimension~$m{-}4$.  At this point, it
is evident that~$X = \Psi(Y)$.  Details will be left to
the reader, along with the verification that the degree~$r$ of~$Y$
satisfies~$[X] = r\bigl[\sigma_{(3,1,1,1)}\bigr]$.
\end{proof}

\begin{remark}
It is not difficult to see that~$\Psi(Y)$ is always singular when~$m\ge 5$.
Thus, the homology classes of the form~$r\bigl[\sigma_{(3,1,1,1)}\bigr]$
cannot be represented by smooth, effective cycles when~$m\ge 5$.
\end{remark}

\subsection{Lagrangian Grassmannians}
\label{ssec:lagrangian-grassmannians}
Fix the standard symplectic form on~$\C{2m}$, namely
\begin{equation}
\Omega := \d z^1\w \d z^{m+1} + \cdots + \d z^m\w \d z^{2m},
\end{equation}
and consider the set~$L_m\subset\Gr(m,2m)$ consisting of the
$\Omega$-Lagrangian $m$-planes in~$\C{2m}$.  This is a compact
manifold of complex dimension~$\frac12m(m{+}1)$ that is
homogeneous under the action of the group~$\Symp(m,\bbC)$.  
The maximal compact subgroup~$\Symp(m)\subset\Symp(m,\bbC)$ also
acts transitively on~$L_m$, with stabilizer isomorphic to~$\Un(m)$.
Thus, 
\begin{equation}
L_m = \frac{\Symp(m)}{\Un(m)}\,,
\end{equation}
which exhibits~$L_m$ as one of the classical Hermitian symmetric
spaces.%
\footnote{In~\cite[Section 16]{MR21:1586}, the notation~$G_m$ is 
used for this symmetric space.}

\subsubsection{Topology}
\label{ssec: lag gr topology}
In~\cite{MR21:1586}, the Poincar\'e polynomial of $L_m$ is found to be
\begin{equation}
p(L_m,t) = (1+t^2)(1+t^4)\dots(1+t^{2m}) 
           = 1 + t^2 + t^4 + 2\,t^6 + \dots,
\end{equation}
so $6$ is the lowest degree in which the rank of a homology group is
greater than~$1$ and this only happens when~$m\ge3$.  For this reason,
I am going to assume that~$m\ge3$ for the rest of this section.

As defined,~$L_m$ is a submanifold of~$\Gr(m,2m)$ and
so inherits bundles~$S$ and~$Q$ by pullback.  For any~$V\in L_m$,
the symplectic structure~$\Omega$ induces an isomorphism~$Q_V\simeq V^*$,
so these bundles satisfy~$S^*=Q$.  

It will be important to understand the tangent space to~$L_m$ at 
a general point~$V\in L_m$.  Now, at~$V\in L_m$, the isomorphism~
$$
T_V\Gr(m,2m)\simeq Q_V\ot V^*\simeq V^*\ot V^* = S^2(V^*)\oplus\L^2(V^*)
$$
is canonical.  Under this isomorphism, the tangent space at~$V$ 
to~$L_m$ corresponds to the subspace~$S^2(V^*)$.  In other words,
$TL_m\simeq S^2(Q)\simeq S^2(S^*)$.

\begin{figure}
\includegraphics[width=\linewidth]{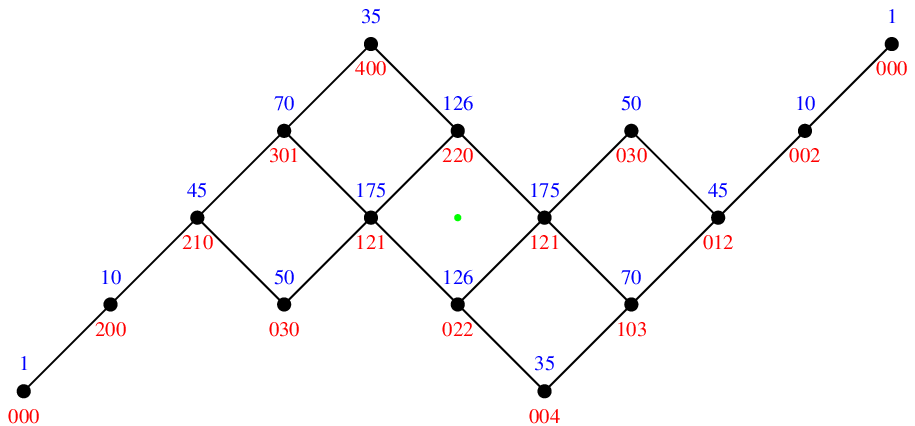}
\caption[The ideal poset for~$L_4 = \Symp(4)/\Un(4)$]
{\label{fig:Sp4poset}  The ideal poset for~$L_4 = \Symp(4)/\Un(4)$.
The upper label on each node is the dimension of the corresponding 
subrepresentation of~$\L^{*,0}\bigl(\eusp(4)/\euu(4)\bigr)$ 
and the lower label is its highest weight 
as a representation of~$\SU(4)$. }
\end{figure}

More detail about the topology and Schubert cell decomposition of~$L_m$
can be found in~\cite{MR55:2941}.  Complete information about 
the irreducible constituents of the exterior powers of the cotangent
bundle of~$L_m$ and the consequent structure of its
ideal poset is collected in a convenient form in~\cite{kKhT99}. 
The corresponding Hasse diagram for the case~$m=4$ is drawn 
in Figure~\ref{fig:Sp4poset}.

\subsubsection{Ideals of degree~$3$}
\label{sssec: lag ideals deg 3}

I have not analyzed the boundary cases in all dimensions
for the Lagrangian Grassmannian, so I will confine myself to studying 
the cases in the first interesting dimension, that of cycles of
dimension or codimension equal to~$3$. 

The first task is to describe the irreducible decomposition
of~$\L^{3,0}(\eum)$ under the action of~$K = \Un(m)$.  Fortunately,
this is relatively easy.  Using the above description of
the tangent bundle of~$L_m$, one can see that, as a
representation of~$\Un(n)$, the space~$\eum$ is isomorphic to the 
representation~$S^2(\C{m}) = \bbS_{(2)}(\C{m})$ 
associated to the standard representation of~$\Un(m)$ on~$\C{m}$.  
Then a little work with  multiplicity formulae from~\cite{MR93a:20069} 
shows that
\begin{equation}\label{eq: L3 decomp of TL}
\L^3\bigl(S^2(\C{m})\bigr) 
\simeq \L^3\bigl(\bbS_{(2)}(\C{m})\bigr)
\simeq \bbS_{(3,3)}(\C{m})\oplus \bbS_{(4,1,1)}(\C{m}).
\end{equation}
These latter two representations are irreducible and their 
dimensions are given by~\cite[Theorem~6.3 or Exercise~6.4]{MR93a:20069} as:
\begin{equation}
\begin{split}
\dim\bbS_{(3,3)}(\C{m})  &=\frac{m^2(m{+}1)^2(m{+}2)(m{-}1)}{144},\\
\dim\bbS_{(4,1,1)}(\C{m})&= \frac{m(m^2{-}1)(m^2{-}4)(m{+}3)}{72}.
\end{split}
\end{equation}

Let~$\cI_{(3,3)}$ and $\cI_{(4,1,1)}$, respectively, 
denote the exterior differential systems on~$L_m$ generated in 
degree~$3$ by the sections of $\bbS_{(3,3)}(S)\subset\L^3(T^*L_m)$ 
and $\bbS_{(4,1,1)}(S)\subset\L^3(T^*L_m)$.

\subsubsection{Integral elements}
\label{ssec: lag ideals deg 3 int elem}

The following linear algebra lemma identifies
the integral elements of dimension three or more for each of the two
$\GL(m,\bbC)$-invariant subspaces of~$\L^3\bigl(S^2((\C{m})^*)\bigr)$. 

\begin{lemma}\label{lem:lagrangian-integral-elements}
Any subspace $E\subset S^2(\C{m})$ of dimension~$3$ or more on which
all of the elements of~$\bbS_{(3,3)}\bigl((\C{m})^*\bigr)$ vanish
is of the form~$E = L{\circ}W$ where~$L\subset\C{m}$ is a line 
and~$W\subset\C{m}$ is a subspace 
whose dimension is the same as that of~$E$.

Any subspace $E\subset S^2(\C{m})$ of dimension~$3$ or more on which
all of the elements of~$\bbS_{(4,1,1)}\bigl((\C{m})^*\bigr)$ vanish
must have dimension~$3$ and be of the form~$E = S^2(W)$ 
where~$W\subset\C{m}$ is a subspace of dimension~$2$.
\end{lemma}

\begin{proof}
This proof is very similar in spirit and structure to the 
proof of Lemma~\ref{lem:integrals-of-ideals-in-L3N+m}, 
so, to save space, I will not go into as much detail 
here as I did there.  Instead, I will limit my discussion 
to the outline, except where some essentially new or different
idea is needed.

Let~$\eta: S^2(\C{m})\to S^2(\C{m})$ be the identity map.
For any basis~$\vb_1,\dots,\vb_m$ of~$\C{m}$, 
write~$\eta = \frac12\eta^{ab}\vb_a{\circ}\vb_b$, where~$\eta^{ab} 
= \eta^{ba}$ are 1-forms on~$S^2(\C{m})$.  
Note that~$\{\eta^{ab}\mid a\le b\}$ is a basis for the dual space.  

Now,~$\bbS_{(3,3)}(\C{m})$ occurs as a constituent
of~$S^3(\C{m})^{\ot 2}$, but~$\bbS_{(4,1,1)}(\C{m})$ 
does not.  Consequently, the $3$-forms of the form
\begin{equation}\label{eq:33-generators}
\psi(X,Y) = -\psi(Y,X) = X_{i_1i_2i_3}Y_{i_4i_5i_6}\,
     \eta^{i_1i_2}\w \eta^{i_3i_4}\w \eta^{i_5i_6}\,,
\end{equation}
when~$X$ and $Y$ are symmetric in their indices,
must lie in the subspace~$\bbS_{(3,3)}\bigl((\C{m})^*\bigr)$
of~$\L^3\bigl(S^2((\C{m})^*)\bigr)$.  Taking,
as a particular example, $X_{111} = Y_{222} = 1$ and 
all other~$X_{ijk}$, $Y_{ijk}$ equal to zero yields
$\psi(X, Y) = \eta^{11}\w\eta^{12}\w\eta^{22}\not=0$, so
the span of the~$\psi(X,Y)$ is nontrivial.  Since
this span is invariant under~$\GL(m,\bbC)$ and since 
$\bbS_{(3,3)}\bigl((\C{m})^*\bigr)$ is irreducible, 
this span must be all of this subspace.  Thus,
$\bbS_{(3,3)}\bigl((\C{m})^*\bigr)$ is the span of the
$3$-forms of the form~\eqref{eq:33-generators}.

Now, to begin, it must be checked that for every line~$L\subset\C{m}$,
all of the forms~$\psi(X,Y)$ vanish on the $m$-dimensional 
subspace~$E = L{\circ}\C{m}\subset S^2(\C{m})$.  By equivariance,
it suffices to check this for the line~$L= \C{}{\cdot}\vb_1$,
which is defined by the equations~$\eta^{ab}=0$ for~$1<a\le b$.
Thus, the claim is equivalent to the claim that the 
forms~$\psi(X,Y)$ all lie in the ideal generated 
by the $1$-forms~$\eta^{ab}$ with~$1<a\le b$.  By obvious reductions
(keeping in mind that~$\psi(X,Y) = -\psi(Y,X)$), it suffices to 
check this for the cases where~$X_{11a}=X_{1a1}=X_{a11}=1$ but
$X_{i_1i_2i_3}=0$ otherwise 
and~$Y_{1bc}=Y_{1cb}=Y_{b1c}=Y_{c1b}=Y_{bc1}
=Y_{cb1}=1$ but~$Y_{i_1i_2i_3}=0$ otherwise.  
In this case, there are $9$ terms in the 
sum~\eqref{eq:33-generators}, and each term either vanishes identically,
cancels in combination with one or two other terms, 
or else is a multiple of some~$\eta^{pq}$ with~$1<p\le q$.  
Thus, the claim is established.

Now suppose that~$E\subset S^2(\C{m})$ is a subspace of dimension~$d\ge3$
on which all of the forms in~$\bbS_{(3,3)}\bigl((\C{m})^*\bigr)$ vanish.
The goal is to show that there is a unique line~$L\subset\C{m}$ so 
that~$E\subset L{\circ}\C{m}$.   Let~$\xi^{ab}$ be the restriction
to~$E$ of~$\eta^{ab}$.

Now,~$E$ will be defined by some set of linear relations among the 
$1$-forms~$\xi^{ab}$.  (By hypothesis, at least three of 
the~$\xi^{ab}$ are linearly independent on~$E$.)
I need to show that one can choose the basis~$\vb_1,\dots,\vb_m$ 
so that these relations include~$\xi^{ab} = 0$ for~$a,b>1$.  

Suppose that the basis~$\vb_1,\dots,\vb_m$ has been chosen so that the 
maximum number, say~$p$, of the forms~$\xi^{11},\dots,\xi^{1m}$ 
are linearly independent.  (This maximum independence will hold
on a Zariski open set of bases~$\vb$.) Clearly,~$1\le p\le \dim E$.
It is not difficult to show that, by making a change of basis in~$\vb$, 
I can assume that~$\xi^{11}\w\dots\w\xi^{1p}\not=0$ but 
that~$\xi^{1a} = 0$ for~$a>p$.  Then, by the same sort of analysis 
as was done in the proof of Lemma~\ref{lem:integrals-of-ideals-in-L3N+m}, 
one can show that maximality implies that
\begin{equation}\label{eq: E rels 1}
\xi^{qr}\equiv 0 \mod \xi^{11},\dots,\xi^{1p}\,,
\qquad\qquad\text{when~$r>p$.}
\end{equation}
Since the~$\xi^{ab}$ must span~$E^*$,
the relations~\eqref{eq: E rels 1} imply that~$p\ge2$.

So far, no use has been made of the assumption that the forms
in~$\bbS_{(3,3)}\bigl((\C{m})^*\bigr)$ vanish on~$E$, i.e., that
\begin{equation}\label{eq: 33 E relations}
0 =  X_{i_1i_2i_3}Y_{i_4i_5i_6}\,
     \xi^{i_1i_2}\w \xi^{i_3i_4}\w \xi^{i_5i_6}\,,
\end{equation}
for all~$X$ and~$Y$ symmetric in their indices.  To make
any further progress, these relations will have to be used.

Since, as has already been seen, the relations~\eqref{eq: 33 E relations}
include~$\xi^{11}\w\xi^{12}\w\xi^{22} = 0$, it follows that~$p\ge 3$.
Now, replacing~$2$ by~$i\le p$ in this relation 
gives~$\xi^{11}\w\xi^{1i}\w\xi^{ii} = 0$, so, in particular,
\begin{equation}\label{eq: E rels 2}
\xi^{ii}\equiv 0 \mod \xi^{11},\dots,\xi^{1p}.
\end{equation} 
On the other hand, polarizing the 
identity $\xi^{11}\w\xi^{1i}\w\xi^{ii} = 0$
gives
\begin{equation}\label{eq: E rels 3}
0 = \xi^{11}\w\xi^{1i}\w\xi^{jk}+\xi^{11}\w\xi^{1j}\w\xi^{ki}
+\xi^{11}\w\xi^{1k}\w\xi^{ij}\,
\end{equation}
for all~$2\le i,j,k\le p$.  Since~$p\ge 3$, taking any pair~$(i,j)$
with~$2\le i<j\le p$ and setting~$k=j$ in the above relation yields
$$
0 = \xi^{11}\w\xi^{1i}\w\xi^{jj}+2\,\xi^{11}\w\xi^{1j}\w\xi^{ij}.
$$
Wedging this relation with~$\xi^{1i}$ 
gives~$\xi^{11}\w\xi^{1i}\w\xi^{1j}\w\xi^{ij}=0$
whenever~$2\le i<j\le p$.  Since~$\xi^{11}\w\xi^{1i}\w\xi^{1j}\not=0$ by
hypothesis, it follows that 
\begin{equation}\label{eq: E rels 4}
\xi^{ij}\equiv 0 \mod \xi^{11},\dots,\xi^{1p}.
\end{equation}

It follows from~\eqref{eq: E rels 1}, \eqref{eq: E rels 2},
and \eqref{eq: E rels 4} 
that $\{\xi^{11},\dots,\xi^{1p}\}$ is a basis for~$E^*$.  
Consequently,~$p = d = \dim E$.  Moreover, 
analysis of~\eqref{eq: E rels 3}, shows that there must exist~$S^i$
so that
$$
\xi^{ij} \equiv S^i\,\xi^{1j} + S^j\,\xi^{1i} \mod\xi^{11},
\qquad\qquad \text{for $2\le i,j\le d$}.
$$
It follows that, by replacing~$\vb_1$ by~$\vb_1 + S^i\,\vb_i$, 
I can arrange that~$\xi^{ij}\equiv 0\mod \xi^{11}$, so assume this.

Now, in the same way that~\eqref{eq: E rels 3} was derived, one can
derive the relations
\begin{equation}\label{eq: E rels 5}
0 = \xi^{11}\w\xi^{1i}\w\xi^{jr}+\xi^{11}\w\xi^{1j}\w\xi^{ir},
\qquad\text{when~$2\le i,j\le d<r$.}
\end{equation}
and
\begin{equation}\label{eq: E rels 6}
0 = \xi^{11}\w\xi^{1i}\w\xi^{qr}, \qquad\text{when~$2\le i\le d < q,r$.}
\end{equation}
The relations~\eqref{eq: E rels 5} imply that there exist~$A^r$
for~$r>d$ so that~$\xi^{jr}\equiv A^r\xi^{1j}\mod \xi^{11}$ 
when~$2\le j\le d<r$, while \eqref{eq: E rels 6} 
implies~$\xi^{jr}\equiv 0\mod \xi^{11}$
when~$d<q,r$.  In particular, it follows that
$$
\xi = \xi^{ab}\,\vb_a\vb_b 
   \equiv 2\,\xi^{1i}\,(\vb_1{+}A^r\vb_r){\circ}\vb_i \mod\xi^{11}.
$$
Thus, the hyperplane~$H=\ker\xi^{11}\subset E$ is of the 
form~$H = L\circ W$ where the line~$L\subset\C{m}$ is spanned by the
vector~$\vb_1' = \vb_1{+}A^r\vb_r$ and~$W$ is spanned 
by~$\vb_2,\dots,\vb_d$.

The analysis so far shows that every~$E\subset S^2(\C{m})$
of dimension~$d\ge3$  on which all the elements 
of~$\bbS_{(3,3)}\bigl((\C{m})^*\bigr)$ vanish contains a hyperplane~$H$
of the form~$L{\circ}W$ where $L$ and~$W$ are subspaces of~$\C{m}$ of
dimensions~$1$ and~$d{-}1$, respectively, that are independent, i.e.,
$L\cap W = 0$.  To finish the characterization of these integral elements,
it suffices to show that, for any basis~$\vb_1,\dots,\vb_m$, every 
integral element of~$\bbS_{(3,3)}\bigl((\C{m})^*\bigr)$ that contains
the $2$-plane spanned by~$\vb_1{\circ}\vb_2$ and~$\vb_1{\circ}\vb_3$
must be itself be a subspace of the 
$m$-plane~$(\C{}{\cdot}\vb_1){\circ}\C{m}$ (which has already been shown
to be an integral element). 

To see this, consider, for every~$a,b>1$, the $3$-forms
\begin{equation*}
\begin{split}
\psi^{ab} 
&= \phantom{{ } + { }}  2\,\eta^{12}\w\eta^{13}\w\eta^{ab} 
      + 2\,\eta^{12}\w\eta^{1a}\w\eta^{b3}
      + 2\,\eta^{12}\w\eta^{1b}\w\eta^{3a}\\
&\phantom{= { } }+\phantom{ 2\,}\,\,\eta^{11}\w\eta^{23}\w\eta^{ab}
      + \phantom{2\,}\eta^{11}\w\eta^{2a}\w\eta^{b3}
      + \phantom{2\,}\eta^{11}\w\eta^{2b}\w\eta^{3a},
\end{split}
\end{equation*}
which manifestly belong to~$\bbS_{(3,3)}\bigl((\C{m})^*\bigr)$.
When~$a,b>3$, 
\begin{equation*}
\begin{split}
\psi^{ab}\bigl(\vb_1{\circ}\vb_2,\vb_1{\circ}\vb_3,\ub\bigr)
  &= 2\,\eta^{ab}(\ub)\,,\\
\psi^{3b}\bigl(\vb_1{\circ}\vb_2,\vb_1{\circ}\vb_3,\ub\bigr)
  &= 4\,\eta^{3b}(\ub)\,,\\
\psi^{2b}\bigl(\vb_1{\circ}\vb_2,\vb_1{\circ}\vb_3,\ub\bigr)
  &= 2\,\eta^{2b}(\ub)\,,
\end{split}
\end{equation*}
while
\begin{equation*}
\begin{split}
\psi^{33}\bigl(\vb_1{\circ}\vb_2,\vb_1{\circ}\vb_3,\ub\bigr)
    &= 6\,\eta^{33}(\ub)\,,\\
\psi^{23}\bigl(\vb_1{\circ}\vb_2,\vb_1{\circ}\vb_3,\ub\bigr)
    &= 4\,\eta^{23}(\ub)\,,\\
\psi^{22}\bigl(\vb_1{\circ}\vb_2,\vb_1{\circ}\vb_3,\ub\bigr)
    &= 2\,\eta^{22}(\ub)\,.\\
\end{split}
\end{equation*}
In particular, if~$\{\vb_1{\circ}\vb_2,\vb_1{\circ}\vb_3,\ub\}$ is
to span an integral element of~$\bbS_{(3,3)}\bigl((\C{m})^*\bigr)$,
then~$\eta^{ab}(\ub)=0$ when~$a,b>1$, i.e., $\ub$ must lie in the
span of~$\{\vb_1{\circ}\vb_1,\dots,\vb_1{\circ}\vb_m\}$, i.e.,
$\ub\in(\C{}{\cdot}\vb_1){\circ}\C{m}$, which is what needed to be shown.

I now turn to the analysis of the integral elements of the 
ideal~$\cI_{(4,1,1)}$.

First,~$\bbS_{(4,1,1)}(\C{m})$ occurs as a constituent
of~$S^3(\C{m})\ot \L^3(\C{m})$, but~$\bbS_{(3,3)}(\C{m})$ 
does not.  Consequently, the $3$-forms of the form
\begin{equation}\label{eq:411-generators}
\psi(X,Y) = X_{i_1i_2i_3}Y_{i_4i_5i_6}\,
     \eta^{i_1i_4}\w \eta^{i_2i_5}\w \eta^{i_3i_6}\,,
\end{equation}
when~$X$ is symmetric in its indices and $Y$ is skewsymmetric 
in its indices,
must lie in the subspace~$\bbS_{(4,1,1)}\bigl((\C{m})^*\bigr)$
of~$\L^3\bigl(S^2((\C{m})^*)\bigr)$.  Taking,
as a particular example, $X_{111} = 1$ with 
all other~$X_{ijk}$ equal to zero, and $Y_{123} = 1$ but~$Y_{ijk}=0$
unless~$\{i,j,k\} = \{1,2,3\}$ yields
$\psi(X, Y) = 6\,\eta^{11}\w\eta^{12}\w\eta^{13}\not=0$.  Thus,
the span of the~$\psi(X,Y)$ is nontrivial.  Since
this span is invariant under~$\GL(m,\bbC)$ and since 
$\bbS_{(4,1,1)}\bigl((\C{m})^*\bigr)$ is irreducible, 
this span must be all of this subspace.  Thus,
$\bbS_{(4,1,1)}\bigl((\C{m})^*\bigr)$ is the span of the
$3$-forms of the form~\eqref{eq:411-generators}.

Now suppose that~$E\subset S^2(\C{m})$ is a subspace of dimension~$d\ge3$
on which all of the forms in~$\bbS_{(4,1,1)}\bigl((\C{m})^*\bigr)$ vanish.
The goal is to show that there is a $2$-plane~$P\subset\C{m}$ so 
that~$E= S^2(P)$.   Let~$\xi^{ab}$ be the restriction
to~$E$ of~$\eta^{ab}$.

Now,~$E$ will be defined by some set of linear relations among the 
$1$-forms~$\xi^{ab}$.  (By hypothesis, at least three of 
the~$\xi^{ab}$ are linearly independent on~$E$.)
I need to show that one can choose the basis~$\vb_1,\dots,\vb_m$ 
so that these relations include~$\xi^{ab} = 0$ when~$b>2$.  

Suppose that the basis~$\vb_1,\dots,\vb_m$ has been chosen so that the 
maximum number, say~$p$, of the forms~$\xi^{11},\dots,\xi^{1m}$ 
are linearly independent.  (This maximum independence will hold
on a Zariski open set of bases~$\vb$.)  Clearly,~$1\le p\le \dim E$.
It is not difficult to show that, by making a change of basis in~$\vb$, 
I can assume that~$\xi^{11}\w\dots\w\xi^{1p}\not=0$ but 
that~$\xi^{1a} = 0$ for~$a>p$.  Then, by the same sort of analysis 
as was done in the proof of Lemma~\ref{lem:integrals-of-ideals-in-L3N+m}, 
one can show that maximality of~$p$ implies that
\begin{equation}\label{eq: 411 E rels 1}
\xi^{qr}\equiv 0 \mod \xi^{11},\dots,\xi^{1p}\,,
\qquad\qquad\text{when~$r>p$.}
\end{equation}
Since the~$\xi^{ab}$ must span~$E^*$,
the relations~\eqref{eq: 411 E rels 1} imply that~$p\ge2$.

So far, no use has been made of the assumption that the forms
in~$\bbS_{(4,1,1)}\bigl((\C{m})^*\bigr)$ vanish on~$E$, i.e., that
\begin{equation}\label{eq: 411 E relations}
0 =  X_{i_1i_2i_3}Y_{i_4i_5i_6}\,
     \xi^{i_1i_4}\w \xi^{i_2i_5}\w \xi^{i_3i_6}\,,
\end{equation}
for all~$X$ symmetric in its indices and~$Y$ skewsymmetric in its indices. 
To make any further progress, these relations will have to be used.

Since, as has already been seen, the relations~\eqref{eq: 411 E relations}
include~$\xi^{11}\w\xi^{12}\w\xi^{13} = 0$, it follows that~$p<3$,
i.e., $p=2$.  Moreover, since~$p=2$, the relations~\eqref{eq: 411 E rels 1}
imply that~$\{\xi^{11},\xi^{12},\xi^{22}\}$ must span~$E^*$.  
Since~$\dim E\ge 3$, by hypothesis, it follows that~$\dim E = 3$ and
that~$\bigl(\xi^{11},\xi^{12},\xi^{22}\bigr)$ must be a basis for~$E^*$.

Now, fix~$a>2$ and let~$Y_{12a}=1$ with~$Y_{ijk}=0$ unless~$\{i,j,k\}
=\{1,2,a\}$.  Letting~$X_{abc}$ be the general symmetric symbol with
$X_{abc}=0$ unless~$1\le a,b,c\le 2$, substituting this into 
the relations~\eqref{eq: 411 E relations}, and using the 
fact that~$\xi^{1a}=0$ yields the relations
$$
\xi^{11}\w\xi^{12}\w\xi^{2a} = \xi^{11}\w\xi^{22}\w\xi^{2a} = 
\xi^{12}\w\xi^{22}\w\xi^{2a} = 0,
$$
which implies~$\xi^{2a} = 0$.  Now, fix~$b>2$.  Letting~$X$
be the symmetric symbol that satisfies~$X_{1,1,b}=1$ but~$X_{ijk}=0$
unless~$(i,j,k)$ is a permutation of~$(1,1,b)$ gives the relation
$\xi^{11}\w\xi^{12}\w\xi^{ab}=0$.   Letting~$X$
be the symmetric symbol that satisfies~$X_{1,2,b}=1$ but~$X_{ijk}=0$
unless~$(i,j,k)$ is a permutation of~$(1,2,b)$ gives the relation
$\xi^{11}\w\xi^{22}\w\xi^{ab}=0$.   Letting~$X$
be the symmetric symbol that satisfies~$X_{2,2,b}=1$ but~$X_{ijk}=0$
unless~$(i,j,k)$ is a permutation of~$(2,2,b)$ gives the relation
$\xi^{12}\w\xi^{22}\w\xi^{ab}=0$.   These three identities imply that
$\xi^{ab}=0$ when~$a,b>2$, which is what remained to be shown.   
\end{proof}

\begin{remark}[Integral element orbit structure]
\label{rem: lag orbit str}
Note an interesting feature of the above description of the integral
elements of~$\bbS_{(3,3)}\bigl((\C{m})^*\bigr)$ when~$m>3$:  
The space of integral elements of a given dimension~$d$ 
in the range~$3\le d<m$ is not homogeneous under the action 
of~$\GL(m,\bbC)$.  After all, such integral elements are of the 
form~$L{\circ}W$ where~$L$ and~$W$ are subspaces of~$\C{m}$ of
dimensions~$1$ and~$d$, respectively.  Generically, one 
will have $L\cap W = 0$ and this describes a $\GL(m,\bbC)$-orbit
that is open and dense in the space of $d$-dimensional integral elements.
However, there are `special' integral elements that 
satisfy~$L\subset W$, and these constitute a closed  $\GL(m,\bbC)$-orbit
of their own.
\end{remark}

\subsubsection{Integral varieties}
\label{ssec: lag ideals deg 3 int var}
 
The next two propositions describe the integral manifolds of 
the exterior differential systems~$\cI_{(3,3)}$ and $\cI_{(4,1,1)}$.  
I remind the reader that~$m\ge 3$ throughout this subsection.

Before stating the first of these two propositions, I need
to describe a family of subvarieties of~$L_m$.

\begin{example}[A Lagrangian Schubert cycle]
\label{ex: lag 411 integral schub}
If~$P\subset\C{2m}$ is any sub-Lagrangian $(m{-}1)$-plane, 
it lies in a $1$-parameter family of Lagrangian $m$-planes:
$$
[P,\C{2m}]_m\cap L_m \simeq \bbP^1.
$$
Let~$S\subset\C{2m}$ be any Lagrangian $m$-plane,
and define a closed subvariety~$\sigma(S)\subset L_m$ by
$$
\sigma(S) = \{\ V\in L_m\ \mid\ \dim(V\cap S)\ge m{-}1\ \}.
$$ 
Of course~$S$ lies in~$\sigma(S)$, and there is a natural submersion
$$
\kappa:\sigma(S)\setminus\{S\}\to 
  \Gr(m{-}1,S)\simeq\bbP(S^*)\simeq\bbP^{m-1}
$$
defined by~$\kappa(V) = S\cap V$ when~$V\in\sigma(S)$ is not~$S$ itself.
By the first statement in this example, the
fibers of~$\kappa$ are biholomorphic to~$\bbC$.  It follows 
that~$\sigma(S)$ is an irreducible subvariety of~$L_m$ of
dimension~$m$ and that~$\sigma(S)$ is smooth away from~$S$.  

It is not difficult to see that the tangent cone to~$\sigma(S)$
at~$S\in L_m$ is the cone consisting of the quadratic forms 
in~$T_SL_m\simeq S^2(S^*)$ that are perfect squares.   
Thus, $S$ is a genuine singular point of~$\sigma(S)$.
\end{example}

\begin{proposition}\label{prop:33-integrals}
For~$S\in L_m$, the variety $\sigma(S)$ is an integral of~$\cI_{(3,3)}$. 

Conversely, if~$X\subset L_m$ is an irreducible variety of 
dimension at least~$3$ that is an integral of~$\cI_{(3,3)}$,
then there is an~$S\in L_m$ so that~$X\subset\sigma(S)$.
\end{proposition}

\begin{proof}
The first task (which will be needed in the next proposition as well),
is to establish the equations of the moving frame for submanifolds
of~$L_m$.

Define~$\Symp(m,\bbC)$ be the subgroup of~$\SL(2m,\bbC)$
consisting of the matrices~$\vs$ that satisfy
\begin{equation}\label{eq:defining-Spm}
{}^t\vs \begin{pmatrix}0_{m}&-\I_{m}\\ \I_{m}&0_{m}\end{pmatrix} \vs
= \begin{pmatrix} 0_{m} & -\I_{m}\\ \I_{m} & 0_{m}\end{pmatrix}.
\end{equation} 
I will regard~$\vs:\Symp(m,\bbC)\to\GL(2m,\bbC)$
as a matrix-valued function and denote its columns as
$$
\vs = (\vs_1\quad\dots\quad \vs_{m} \quad \vs^1\quad\dots\quad \vs^{m})
$$
where~$\vs_i,\vs^i:\Symp(m,\bbC)\to\C{2m}$ 
are regarded as (holomorphic) mappings.

Define
$$
\pi(\vs) = [\vs_1\w\dots\w\vs_m]\,,
$$
so that~$\pi$ is a surjective submersion~$\pi:\Symp(m,\bbC)\to L_m$.
The fibers of~$\pi$ are the orbits of the parabolic subgroup~$P\subset
\Symp(m,\bbC)$ consisting of elements of the form
\begin{equation}\label{eq: lag elements of P}
\vs = \begin{pmatrix} A & AB\\ 0_m & {}^t\!A^{-1}\end{pmatrix}
\qquad\text{for~$A\in \GL(m,\bbC)$ and~$B = {}^t\!B\in\C{m,m}$.}
\end{equation}
Thus,~$\pi:\Symp(m,\bbC)\to L_m$ is a principal 
right~$P$-bundle over~$L_m$.

In accordance with the usual moving frame conventions, 
write the structure equations as
\begin{equation} \label{eq:Spm-structure-equations}
\d\vs = \d (\vs_i\ \vs^i) = (\vs_j\ \vs^j)
\begin{pmatrix}\alpha^j_i & \gamma^{ji}\\ 
                \beta_{ji}& -\alpha^i_j \end{pmatrix}
= \vs\,\theta
\end{equation}
where~
\begin{equation}\label{eq:Spm-structure-relations}
\beta_{ji} = \beta_{ij}\quad\text{and}\quad\gamma^{ji}=\gamma^{ij}\,,
\end{equation}
but the components of~$\alpha$, $\beta$, and~$\gamma$ are otherwise
linearly independent. The relations~\eqref{eq:Spm-structure-relations} 
follow in the usual way from the exterior derivative 
of~\eqref{eq:defining-Spm}.  The structure equation 
$\d\theta = -\theta\w\theta$ holds since~$\theta = \vs^{-1}\,\d \vs$.
These expand to
\begin{equation}\label{eq:Spm-structure-abc-relations}
\begin{split}
\d\alpha^i_j  &= -\alpha^i_k\w\alpha^k_j-\gamma^{ik}\w\beta_{kj}\,,\\
\d\beta_{ij}  &= -\beta_{ik}\w\alpha^k_j +\alpha^k_i\w\beta_{kj}\,,\\
\d\gamma^{ij} &= -\alpha^i_k\w\gamma^{kj} +\gamma^{ik}\w\alpha^j_k\,.
\end{split}
\end{equation}

Now suppose that~$X\subset L_m$ is an irreducible integral
variety of~$\cI_{(3,3)}$ of dimension~$d\ge 3$, 
and let~$X^\circ\subset X$ denote its
smooth locus, which is an embedded submanifold of~$L_m$.
For every~$V\in X^\circ$, the 
tangent space~$T_VX$ is an integral element of~$\cI_{(3,3)}$ of
dimension~$d\ge3$. By Lemma~\ref{lem:lagrangian-integral-elements},
there is a Zariski-open subset~$X^\diamond\subset X^\circ$ (which
may be empty) that consists of the elements~$V\in X^\circ$ such
that~$T_VX = L_V\circ W_V$ where~$L_V$ and~$W_V$ are subspaces of~$V^*$
of dimensions~$1$ and~$d$, respectively, and $L_V\cap W_V=0$.  
(See Remark~\ref{rem: lag orbit str}.)  The proof has to be broken
up into two cases now, depending on whether or not~$X^\diamond$ is
empty.  (Note that~$X^\diamond$ cannot be empty unless~$d<m$.)

The first case is that~$X^\diamond$ is not empty, so assume this.
Note that~$X^\diamond$ is connected since~$X$ is irreducible.
Now, for every~$V\in X^\diamond$, there exists a~$\vs\in \Symp(m,\bbC)$
so that
\begin{enumerate}
\item $V$ is spanned by~$\vs_1,\dots,\vs_m$, and
\item $T_VX$ is spanned by 
$\dl\vs^1\dr{\circ}\dl\vs^2\dr,\dl\vs^1\dr{\circ}\dl\vs^3\dr,\dots,
 \dl\vs^1\dr{\circ}\dl\vs^{d+1}\dr$.
\end{enumerate}
Let~$F(X^\diamond)\subset \Symp(m,\bbC)$ denote the set of such~$\vs$ as~$V$
ranges over~$X^\diamond$.  Then~$\pi:F(X^\diamond)\to X^\diamond$ is
a principal~$G$-bundle over~$X^\diamond$, where~$G\subset P$ is a 
subgroup of the matrices of the form~\eqref{eq: lag elements of P}
with~${}^t\!A^{-1}$ in $P_1\cap P_{d+1}\subset\GL(m,\bbC)$.
The reader can write out the exact conditions defining~$G$ and
verify that it is connected.

By construction, the forms~$\beta_{12},\dots,\beta_{1(d{+}1)}$ are 
linearly independent on~$F(X^\diamond)$ and span the $\pi$-semibasic 
$1$-forms, while $\beta_{11}=\beta^{ab}=0$
when either~$b>d{+}1$ or $a$ and~$b$ are both greater than~$1$.

This paragraph of the argument is necessary only if~$d{+}1<m$,
so suppose this is so for the moment.  Choose a pair~$(i,a)$ 
satisfying~$2\le i\le d{+}1<a\le m$ and differentiate 
the relation~$\beta_{ia}=0$.  By the structure equations, this is
$$
0 = \d \beta_{ia} = -\beta_{i1}\w\alpha^1_a\,.
$$
Since~$d\ge3$, and since~$\beta_{12},\dots,\beta_{1(d{+}1)}$ are
linearly independent,~$\alpha^1_a=0$ for~$a>d{+}1$.

Now choose a pair~$(i,j)$ with~$2\le i,j\le d{+}1$ 
and differentiate~$\beta_{ij}=0$.  The structure equations
give that
$$
0 = \d \beta_{ij} = -\beta_{i1}\w\alpha^1_j + \alpha^1_i\w\beta_{1j}\,.
$$
Equivalently, 
\begin{equation}\label{eq:Spm-ab-symmetry}
\alpha^1_i\w\beta_{1j} + \alpha^1_j\w\beta_{1i} = 0
\end{equation}  
This relation implies that there exists a 
function~$\lambda$ on~$F(X^\diamond)$
so that~$\alpha^1_i = \lambda\,\beta_i$ for~$2\le i\le d{+}1$.  

Computing how the function~$\lambda$ varies on the fibers of~$\pi$
(a standard computation in the technique of the moving frame)
shows that the equation~$\lambda=0$ defines a principal
right $G_1$-bundle~$F_1\subset F(X^\diamond)$ over~$X^\diamond$
where~$G_1\subset G$ is a certain connected Lie subgroup of 
codimension~$1$.

Since the previous paragraph showed that~$\alpha^1_a=0$ 
for all~$a>d{+}1$, it now follows that the 
identities~$\alpha^1_i=0$ for all~$i>1$ hold on~$F_1$.  This
vanishing together with the fact that~$\beta_{ij}=0$ for all~$i,j\ge2$
yield the congruences
$$
\d\vs_2\equiv\dots\equiv \d\vs_m\equiv \d\vs^1
  \equiv 0 \mod \vs_2\,,\dots,\vs_m\,,\vs^1\,.
$$
In other words, the mapping~$\nu:F_1\to L_m$ defined by
$\nu(\vs) = [\vs_2\w\dots\w\vs_m\w\vs^1]$ is locally constant
and hence, by connectedness, globally constant.  Let~$S\in L_m$
be $m$-plane that is the image of~$\nu$.    Of course, it now
follows that~$S\cap \pi(\vs)$ is the $(m{-}1)$-plane
$[\vs_2\w\dots\w\vs_m]$.  Thus, $X^\diamond=\pi(F_1)$ lies in~$\sigma(S)$.  
Of course, this implies that~$X$ itself lies in~$\sigma(S)$ as well,
as desired.

Now, consider the second case, in which~$X^\diamond=\emptyset$.  
Then for every~$V\in X^\circ$, the tangent space~$T_VX$
is of the form $L_V\circ W_V$ where~$L_V$ and~$W_V$ are subspaces of~$V^*$
of dimensions~$1$ and~$d$, respectively, and $L_V\subset W_V$.  
(Again, see Remark~\ref{rem: lag orbit str}.)
Note that~$X^\circ$ is connected since~$X$ is irreducible.
For every~$V\in X^\circ$, there exists a~$\vs\in \Symp(m,\bbC)$
so that
\begin{enumerate}
\item $V$ is spanned by~$\vs_1,\dots,\vs_m$, and
\item $T_VX$ is spanned by 
$\dl\vs^1\dr{\circ}\dl\vs^1\dr,\dl\vs^1\dr{\circ}\dl\vs^2\dr,\dots,
 \dl\vs^1\dr{\circ}\dl\vs^{d}\dr$.
\end{enumerate}
Let~$F(X^\circ)\subset\Symp(m,\bbC)$ denote the set of such~$\vs$ as~$V$
ranges over~$X^\circ$.  Then~$\pi:F(X^\circ)\to X^\circ$ is
a principal~$G$-bundle over~$X^\circ$, where~$G\subset P$ is the 
group consisting of the matrices of the form~\eqref{eq: lag elements of P}
with~${}^t\!A^{-1}$ in $P_1\cap P_{d}\subset\GL(m,\bbC)$.

By construction, the forms~$\beta_{11},\dots,\beta_{1d}$ are 
linearly independent on~$F(X^\circ)$ and span the $\pi$-semibasic 
$1$-forms, while $\beta^{ab}=0$
when either~$b>d$ or $a$ and~$b$ are both greater than~$1$.

This paragraph of the argument is necessary only if~$d<m$,
so suppose this is so for the moment.  Choose a pair~$(i,a)$ 
satisfying~$2\le i\le d< a\le m$ and differentiate 
the relation~$\beta_{ia}=0$.  By the structure equations, this is
$$
0 = \d \beta_{ia} = -\beta_{i1}\w\alpha^1_a\,.
$$
Since~$d\ge3$, and since~$\beta_{12},\dots,\beta_{1d}$ are
linearly independent,~$\alpha^1_a=0$ for~$a>d$.

Now choose a pair~$(i,j)$ with~$2\le i,j\le d$ 
and differentiate~$\beta_{ij}=0$.  The structure equations
give
$$
0 = \d \beta_{ij} = -\beta_{i1}\w\alpha^1_j + \alpha^1_i\w\beta_{1j}\,.
$$
Equivalently, 
\begin{equation}\label{eq:Spm-ab-symmetry 2}
\alpha^1_i\w\beta_{1j} + \alpha^1_j\w\beta_{1i} = 0
\end{equation}  
This relation implies that there exists a 
function~$\lambda$ on~$F(X^\circ)$
so that~$\alpha^1_i = \lambda\,\beta_i$ for~$2\le i\le d$.  

Computing how the function~$\lambda$ varies on the fibers of~$\pi$
(a standard computation in the technique of the moving frame)
shows that the equation~$\lambda=0$ defines a principal
right $G_1$-bundle~$F_1\subset F(X^\circ)$ over~$X^\circ$
where~$G_1\subset G$ is a certain connected Lie subgroup of 
codimension~$1$.

Since the previous paragraph showed that~$\alpha^1_a=0$ 
for all~$a>d$, it now follows that the 
identities~$\alpha^1_i=0$ for all~$i>1$ hold on~$F_1$.  This
vanishing together with the fact that~$\beta_{ij}=0$ for all~$i,j\ge2$
yield the congruences
$$
\d\vs_2\equiv\dots\equiv \d\vs_m\equiv \d\vs^1
  \equiv 0 \mod \vs_2\,,\dots,\vs_m\,,\vs^1\,.
$$
In other words, the mapping~$\nu:F_1\to L_m$ defined by
$\nu(\vs) = [\vs_2\w\dots\w\vs_m\w\vs^1]$ is locally constant
and hence, by connectedness, globally constant.  Let~$S\in L_m$
be $m$-plane that is the image of~$\nu$.    Of course, it now
follows that~$S\cap \pi(\vs)$ is the $(m{-}1)$-plane
$[\vs_2\w\dots\w\vs_m]$.  Thus, $X^\circ=\pi(F_1)$ lies in~$\sigma(S)$.  
Of course, this implies that~$X$ itself lies in~$\sigma(S)$ as well,
as desired.

That~$\sigma(S)$ really is an integral variety of~$\cI_{(3,3)}$ 
follows immediately from the argument in the second case (with $d=m$).
\end{proof}

By Lemma~\ref{lem:lagrangian-integral-elements}, 
there are no integral manifolds of~$\cI_{(4,1,1)}$ 
of dimension greater than~$3$.  
The following proposition classifies the $3$-dimensional integral
varieties of~$\cI_{(4,1,1)}$.

First, a definition.  For any sub-Lagrangian subspace~$A\subset\C{2m}$, 
let
$L_m(A)\subset L_m$ denote the set of~$P\in L_m$ that contain~$A$.
Note that, if~$a = \dim A < m$, then $L_m(A)$ is a smooth subvariety 
of~$L_m$ that is isomorphic to~$L_{m-a}$.  

\begin{proposition}\label{prop:411-integrals}
Let~$X\subset L_m$ be an irreducible variety of dimension~$3$
that is an integral variety of~$\cI_{(4,1,1)}$.  Then there is
a sub-Lagrangian $(m{-}2)$-plane~$A\subset\C{2m}$ 
so that~$X\subset L_m(A)$. 

Thus, if~$X$ is closed, then $X = L_m(A)\simeq L_2\simeq Q_3$.
\end{proposition}

\begin{proof}
Recall the moving frame notation and constructions 
from the first part of the proof of Proposition~\ref{prop:33-integrals}.

Suppose now that~$X\subset L_m$ 
is an irreducible $3$-dimensional integral variety of of~$\cI_{(4,1,1)}$
and let~$X^\circ\subset X$ be its smooth locus, which is connected
since~$X$ is irreducible. 
By Lemma~\ref{lem:lagrangian-integral-elements}, 
for every~$V\in X^\circ$, there exists a~$\vs\in \Symp(m,\bbC)$ so that
\begin{enumerate}
\item $V$ is spanned by~$\vs_1,\dots,\vs_m$, and
\item The tangent space $T_VX^\circ$ is spanned by
       $\dl\vs^1\dr{\circ}\dl\vs^1\dr,\>\dl\vs^1\dr{\circ}\dl\vs^2\dr,
       \>\dl\vs^2\dr{\circ}\dl\vs^2\dr$.
\end{enumerate}
Let~$F(X^\circ)\subset\Symp(m,\bbC)$ denote the set of such~$\vs$ as~$V$
ranges over~$X^\circ$.  Then~$\pi:F(X^\circ)\to X^\circ$ is a 
principal right $G$-bundle over~$X^\circ$ where $G\subset P$ is the
subgroup consisting of the matrices 
of the form~\eqref{eq: lag elements of P}
with~${}^t\!A^{-1}$ in $P_3\subset\GL(m,\bbC)$.  Since~$G$ and~$X^\circ$
are each connected, it follows that~$F(X^\circ)$ is also connected. 

By construction, the $1$-forms~$\beta_{11},\beta_{12},\beta_{22}$
are linearly independent on~$F(X^\circ)$ and span the $\pi$-semibasic
$1$-forms, while $\beta_{ij} = 0$ if either~$i$ or~$j$ is
greater than~$2$.

Let~$i>2$ be fixed and differentiate 
the identities~$\beta_{i1}=\beta_{i2}=0$
using the structure equations.  The result is equations of the form
$$
\begin{pmatrix} \alpha^1_i& \alpha^2_i \end{pmatrix}
\w
\begin{pmatrix} \beta_{11} & \beta_{12}\\
               \beta_{12} &  \beta_{22}
\end{pmatrix}
= \begin{pmatrix} 0 & 0 \end{pmatrix} .
$$
By the linear independence of~$\beta_{11},\beta_{12},\beta_{22}$,
it follows that~$\alpha^1_i = \alpha^2_i = 0 $.

This vanishing for all~$i>2$ implies 
$$
\d\vs_3\equiv\dots\equiv\d\vs_m \equiv 0 \mod \vs_3,\dots,\vs_m\,,
$$
i.e., the $(m{-}2)$-plane~$[\vs_3\w\dots\vs_m]$ is locally
constant on $F(X^\circ)$.  Since~$F(X^\circ)$ is connected, this
map must be constant.  Thus, let~$A\in \Gr(m{-}3,2m)$ be the isotropic
plane so that~$ [\vs_3\w\dots\vs_m]\equiv A$.
By construction,~$A\subset V$ for all~$V\in X^\circ$, so it follows
that~$X^\circ$ and, hence,~$X$ are subsets of~$L_m(A)$, as desired.
\end{proof}

These propositions allow characterizations of the extremal
classes in~$H_6(L_m)$ that are analogous 
to that of Schubert cycles in Grassmannians:

\begin{theorem}\label{thm:rigid3-cycles-in-lag-Grassmannians}
Let~$A\subset\C{2m}$ be a sub-Lagrangian plane of dimension~$m{-}2$.
Fix a Lagrangian $m$-plane~$S\subset\C{2m}$ 
and let~$B\subset S$ be a subspace of dimension~$m{-}3$. 
Define two $3$-dimensional subvarieties of~$L_m$ by
\begin{equation}
  X = \sigma(S)\cap [B,\C{2m}]_m
 \qquad\text{and}\qquad
  Y = L_m(A).
\end{equation}
Then~$[X]$ and~$[Y]$ are the generators of~$H^+_6(L_m,\bbZ)\simeq\bbZ^2$.

An irreducible~$Z\in\cZ^+_3(L_m)$ satisfies~$[Z] = r[X]$ for some~$r>0$
if and only if it lies in~$\sigma(T)$ for some~$T\in L_m$.

An irreducible~$Z\in\cZ^+_3(L_m)$ satisfies~$[Z] = r[Y]$ 
if and only if~$r=1$ and $Z = L_m(D)$ 
for some isotropic~$D\subset\C{2m}$
of dimension~$m{-}2$.
\end{theorem}

\begin{proof}
First of all, it follows by either~\cite{MR21:1586} or 
\eqref{eq: L3 decomp of TL} and the general results of Kostant
mentioned above that $b_6(L_m)=2$.  
Let~$\phi_1$ be the $\SO(2m)$-invariant
K\"ahler form on~$L_m$ whose cohomology class is a generator
of~$H^2(L_m,\bbZ)$.  By \eqref{eq: hss sum of Schur forms},
there is a sum of the form
$$
{\phi_1}^3 = \mu^{(3,3)}\,\phi_{(3,3)}+\mu^{(4,1,1)}\,\phi_{(4,1,1)}
$$
where~$\mu^{(3,3)}>0$ and $\mu^{(4,1,1)}>0$ and 
$\phi_{(3,3)}$ and~$\phi_{(4,1,1)}$ are positive $\SO(2m)$-invariant  
forms that are dual to the generalized Schubert cycles~$\sigma_{(3,3)^*}$ 
and~$\sigma_{(4,1,1)^*}$ of complex dimension~$3$ whose 
cohomology classes generate~$H^+_6(L_m)$.  

It follows that~$\sigma_{(3,3)^*}$ is an irreducible $3$-dimensional 
integral of~$\cI_{(4,1,1)}$ and so, 
by Proposition~\ref{prop:411-integrals}, must be of the form~$L_m(C)$ 
for some isotropic $m{-}3$ plane~$C$.  Thus~$\sigma_{(3,3)^*}$ is
homologous to~$Y$.

It also follows that~$\sigma_{(4,1,1)^*}$ 
is an irreducible $3$-dimensional
integral of~$\cI_{(3,3)}$, and so, by
Proposition~\ref{prop:33-integrals}, must lie in 
$\sigma(T)$ for some~$T\in L_m$.   
When~$m=3$, the dimension of~$\sigma(T)$
is~$3$, so~$\sigma_{(4,1,1)^*}$ must be equal to~$\sigma(T)$,
which is~$X$ in this case.  The definition of~$X$ now
makes it clear that~$X$ represents~$\sigma_{(4,1,1)^*}$ 
when $m>3$ as well.

Finally, if~$Z\in\cZ^+_3(L_m)$ is irreducible and satisfies~$[Z]=r[X]$,
then the integral of~$\phi_{(3,3)}$ over~$Z$ must be zero, 
so $\phi_{(3,3)}$ must vanish on~$Z$.  Thus~$Z$ is an integral manifold
of~$\cI_{(3,3)}$ and Proposition~\ref{prop:33-integrals} applies.

The argument when~$[Z] = r[Y]$ is similar.
\end{proof}

The extremal cycles of codimension~$3$ in~$L_m$ also display rigidity.
I am not going to give all the details of this discussion, since most
of the methods of proof in each of the two cases I am going to consider 
will, by now, be familiar to the reader.  Instead, I will simply mention
the points at which some interesting or different idea comes into play.

\begin{theorem}\label{thm:*33-integrals}
Let~$Y\subset\bbP^{2m-1}$ be a subvariety
of dimension~$m{-}3$ and degree~$r$.  Let~$\Psi(Y)\subset L_m$ denote 
the set of~$V\in L_m$ satisfying~$\bbP(V)\cap Y\not=\emptyset$.

Then~$\Psi(Y)$ has codimension~$3$ in~$L_m$ and 
satisfies~$\bigl[\Psi(Y)\bigr]=r\bigl[\sigma_{(4,1,1)}\bigr]$.  

Moreover, any codimension~$3$ subvariety~$X\subset L_m$ that
satisfies~$[X] = r\bigl[\sigma_{(4,1,1)}\bigr]$ is~$\Psi(Y)$ 
for some subvariety~$Y\subset\bbP^{2m-1}$ of dimension~$m{-}3$ 
and degree~$r$.
\end{theorem}

\begin{proof}
The analysis is very similar to that for the proof of
Theorem~\ref{thm:*222-integrals}, so I will leave the details
to the reader.  Only one aspect of the proof requires comment:
Since there are two types of integral elements of~$\cI_{(3,3)}$,
there are, correspondingly, two types of integral 
elements of~$\cI_{(3,3)^*}$.  This forces a subdivision into two
cases like that of Proposition~\ref{prop:33-integrals}, but this
offers no essential new difficulty. \end{proof}

\begin{remark}
It is not difficult to see that~$\Psi(Y)$ is always singular when~$m\ge4$.
Thus, the homology classes of the form~$r\bigl[\sigma_{(4,1,1)}\bigr]$
cannot be represented by smooth, effective cycles when~$m\ge 4$.
\end{remark}

\begin{theorem}\label{thm:*411-integrals}
Let $P\subset\C{2m}$ be a Lagrangian $m$-plane and let
$$
\Sigma(P)= \left\{\ V\in L_m \mid \dim(V\cap P)\ge 2 \ \right\} .
$$
Then~$\Sigma(P)$ is of codimension~$3$ in~$L_m$ 
and represents the Schubert cycle~$\sigma_{(3,3)}$.  

Any irreducible variety
$X\subset L_m$ of codimension~$3$ that satisfies~$[X] 
= r\,\bigl[\Sigma(P)\bigr]$ is of the form~$X=\Sigma(Q)$ 
for some Lagrangian $m$-plane~$Q$. 
\end{theorem}

\begin{proof}
Follow the pattern of the proof of Theorem~\ref{thm:*3111-integrals}.
\end{proof}

\begin{remark}[Singularity of~$\Sigma(P)$]
The variety~$\Sigma(P)$ is singular when~$m\ge 3$.  
In fact,~$P$ itself is a singular point of~$\Sigma(P)$, as is easily
verified.
In particular, Theorem~\ref{thm:*411-integrals} implies
that, when~$m\ge 3$, no multiple of the homology class~$[\sigma_{(3,3)}]$ 
can be represented by a smooth, effective cycle. 
\end{remark}

\subsection{The space E\,III}
\label{ssec:EIII}
I have not done any deep analysis of the ideals on
the Hermitian symmetric 
space~$\mathsf{E\,III}=\E_6/\bigl(S^1{\cdot}\Spin(10)\bigr)$
(which has (complex) dimension~$16$),
but in this section, I will indicate some of the interesting 
possibilities that are turned up by a preliminary analysis.

\begin{figure}
\includegraphics[width=\linewidth]{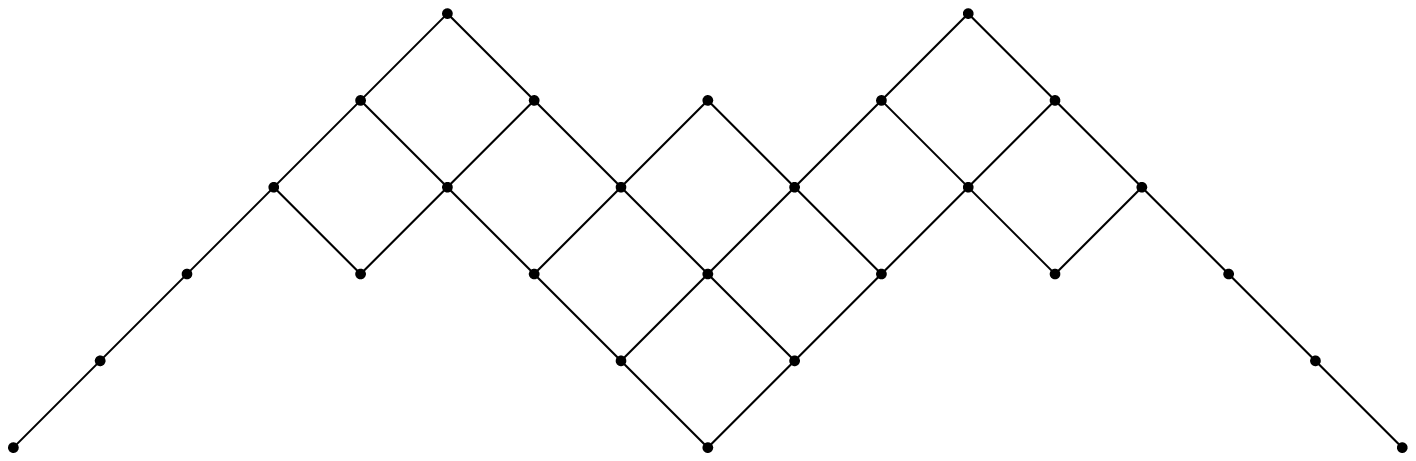}
\caption[The ideal poset 
for~$\mathsf{E\,III}=\E_6/\bigl(S^1{\cdot}\Spin(10)\bigr)$]
{\label{fig: E3idealposet}  
The ideal poset for~$\mathsf{E\,III}=\E_6/\bigl(S^1{\cdot}\Spin(10)\bigr)$.
}
\end{figure}

Using a program such as \texttt{simpLie},
it is not difficult to compute the poset of ideals of~$\mathsf{E\,III}$.  
The Hasse diagram for this poset is drawn in Figure~\ref{fig: E3idealposet}. 

I have not attempted to label the nodes, partly through lack of space.  
However, it will be convenient to be able to refer to some of the nodes,
so I will do so by coordinates.  Thus, the left-most node in the
figure has coordinates~$(0,0)$ and they continue upwards and to
the right as~$(1,1)$, $(2,2)$, $(3,3)$, with the first case of
two nodes at the same level being at level~$4$, namely~$(4,4)$
and~$(4,2)$.  The lowest central node is~$(8,0)$, and so on.
I will refer to the corresponding ideals on~$\mathsf{E\,III}$
by~$\cI_{(p,q)}$.  Thus, for example,~$\cI_{(p,q)}$ is 
generated in degree~$p$ and 
the first interesting ideals are~$\cI_{(4,2)}$
and~$\cI_{(4,4)}$.  In this notation,~$\bigl|(p,q)\bigr| = p$ 
and~$(p,q)^* = (16{-}p,q)$.

These nodes could also be labeled by the the highest weight
of the corresponding irreducible subrepresentation of~$\L^*(\eum)$,
where~$\eum = \eue_6/\bigl(\eut{+}\euso(10)\bigr)$, with respect to
a maximal torus in~$K = S^1{\cdot}\Spin(10)$.  It is traditional
to do this by labeling the nodes of the Dynkin diagram of~$K$ with
the coordinates of the highest weight, so I will do this in a 
compactified form.  Thus, the summand~$\L^1(\eum) = \eum$,
which is, of course, irreducible and of dimension~$16$, will be notated 
as~$1{\cdot}$\raisebox{0.25ex}{$\ss000^{\ss1}_{\ss0}$}.  
It corresponds to the node~$(1,1)$.

The node~$(4,2)$ is the 
representation~$4{\cdot}$\raisebox{0.25ex}{$\ss020^{\ss0}_{\ss0}$}
of dimension~$770$, which is considerably smaller than the
node~$(4,4)$, the 
representation~$4{\cdot}$\raisebox{0.25ex}{$\ss100^{\ss0}_{\ss2}$}
of dimension~$1050$.  Thus, one might expect the integrals 
of the ideal~$\cI_{(4,2)}$ to display more flexibility than those 
of~$\cI_{(4,4)}$.   In fact, as the Hasse diagram makes clear, 
$\cI_{(4,4)}$ contains all forms of degree~$5$ or higher, so its
maximal integrals have dimension~$4$.  It would be interesting to know
whether or not all its irreducible integral varieties are the 
obvious Schubert varieties of dimension~$4$.  In contrast, the 
ideal~$\cI_{(4,2)}$ has integrals of dimension at least~$5$,
since the corresponding form~$\phi_{(4,2)}$ vanishes on the 
$5$-dimensional Schubert cycles that correspond to the node~$(5,5)$.

In fact, the node~$(5,5)$ is very interesting, corresponding
to the 
representation~$5{\cdot}$\raisebox{0.25ex}{$\ss000^{\ss0}_{\ss3}$}
of dimension~$672$, which is less than 20\% 
of the size of its `competitor' at level~$5$, namely the node~$(5,3)$ 
which is 
the representation~$5{\cdot}$\raisebox{0.25ex}{$\ss110^{\ss0}_{\ss1}$} 
of dimension~$3696$.
Again, note that the ideal~$\cI_{(5,3)}$ contains all $6$-forms,
so that its maximal integrals are $5$-dimensional.  Meanwhile, the
ideal~$\cI_{(5,5)}$ fails to contain forms of degree as high as~$8$.
It would be interesting to know whether the maximal
dimension irreducible integrals of this ideal 
are the $8$-dimensional Schubert cycles that correspond
to the node~$(8,0)$.  

As a final note on size disparity among the ideals, the node~$(8,0)$ 
corresponds to the 
representation~$8{\cdot}$\raisebox{0.25ex}{$\ss400^{\ss0}_{\ss0}$}
whose dimension is~$660$, a remarkably small number in comparison with
the dimensions associated to the nearby nodes.  One might expect the
integrals of the ideal~$\cI_{(8,0)}$ to be particularly interesting.
For example, could it be the case that every irreducible integral 
of~$\cI_{(8,0)}$ of dimension at least~$8$ lies in one of the 
$11$-dimensional Schubert varieties associated to the node~$(11,5)$?
(It is evident that these Schubert varieties are integrals of
$\cI_{(8,0)}$ and that this ideal has no integrals of dimension 
greater than~$11$.)

The exploration of these problems is postponed to a later date.

\subsection{The space E\,VII}
\label{ssec:EVII}
I have not done any deep analysis of the ideals on
the Hermitian symmetric 
space~$\mathsf{E\,VII}=\E_7/\bigl(S^1{\cdot}\E_6\bigr)$
(which has (complex) dimension~$27$),
but in this section, I will indicate some of the interesting 
possibilities that are turned up by a preliminary analysis.

\begin{figure}
\includegraphics[width=\linewidth]{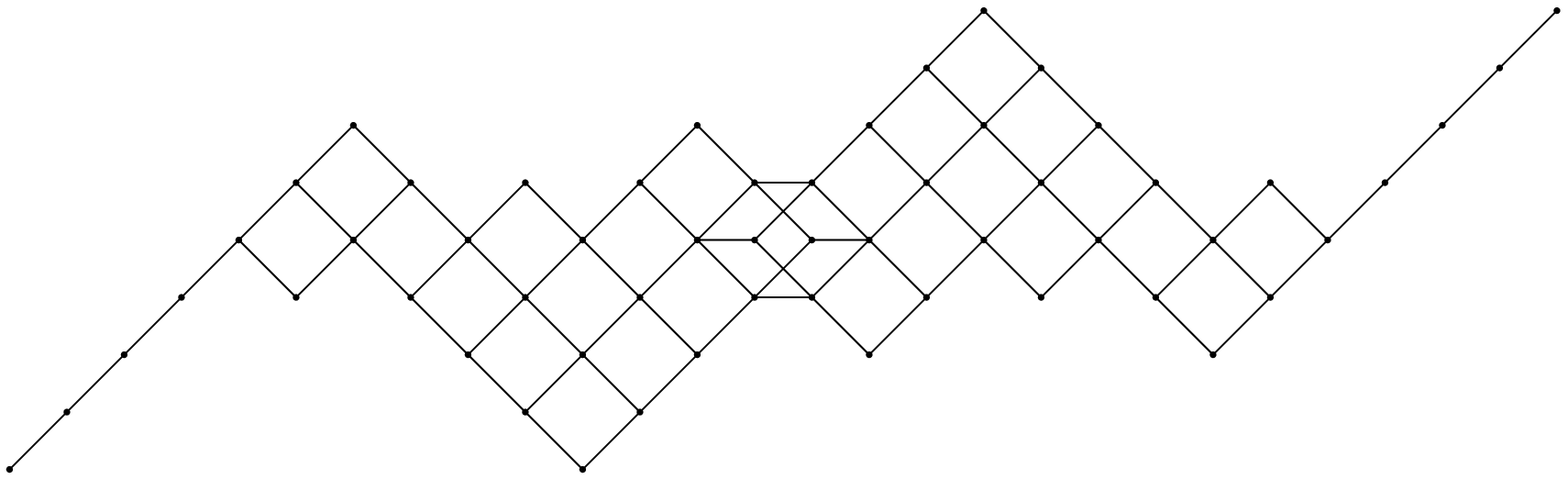}
\caption[The ideal poset 
for~$\mathsf{E\,VII}=\E_7/\bigl(S^1{\cdot}\E_6\bigr)$]
{\label{fig: E7idealposet} 
The ideal poset for~$\mathsf{E\,VII}=\E_7/\bigl(S^1{\cdot}\E_6\bigr)$.
}
\end{figure}

Using a program such as \texttt{simpLie},
it is not difficult to compute 
the poset of ideals of~$\mathsf{E\,VII}$.  
The Hasse diagram for this poset 
is drawn in Figure~\ref{fig: E7idealposet}. 

I have not attempted to label the nodes, partly through lack of space.  
However, it will be convenient to be able to refer to some of the nodes,
so I will do so by coordinates.  Thus, the left-most node in the
figure has coordinates~$(0,0)$ and they continue upwards and to
the right as~$(1,1)$, $(2,2)$, $(3,3)$, with the first case of
two nodes at the same level being at level~$5$, namely~$(5,5)$
and~$(5,3)$.  For each node~$(p,q)$,
I will refer to the corresponding ideal on~$\mathsf{E\,VII}$
by~$\cI_{(p,q)}$.  Thus, for example,~$\cI_{(p,q)}$ is 
generated in degree~$p$ and the first interesting ideals 
are~$\cI_{(5,3)}$
and~$\cI_{(5,5)}$.  In this notation,~$\bigl|(p,q)\bigr| = p$ 
and~$(p,q)^* = (27{-}p,8{-}q)$.

These nodes could also be labeled by the the highest weight
of the corresponding irreducible subrepresentation of~$\L^*(\eum)$,
where~$\eum = \eue_7/\bigl(\eut{+}\eue_6\bigr)$, with respect to
a maximal torus in~$K = S^1{\cdot}\E_6$.  It is traditional
to do this by labeling the nodes of the Dynkin diagram of~$K$ with
the coordinates of the highest weight, so I will do this in a 
compactified form.  Thus, the summand~$\L^1(\eum) = \eum$,
which is, of course, irreducible and of dimension~$27$, will be notated 
as~$1{\cdot}\begin{smallmatrix}&&0\\1&0&0&0&0\end{smallmatrix}$.  
It corresponds to the node~$(1,1)$.

The figure makes clear some of the interesting features of the
ideals, so I will not belabor them here, except to mention two
of the more interesting nodes:

The node~$(6,6)$ is the representation~
$6{\cdot}\begin{smallmatrix}&&3\\0&0&0&0&0\end{smallmatrix}$
of dimension~$43,758$, which is only slightly more than one-eighth of 
the dimension of the representation associated to its `competitor'
node~$(6,4)$.  Note also that~$\cI_{(6,6)}$ has irreducible integrals of 
dimension~$10$, namely the Schubert cycles corresponding to the 
node~$(10,0)$.  It would be interesting to know whether any
irreducible integral of~$\cI_{(6,6)}$ of dimension~$6$ or more is
a subvariety of one of these $10$-dimensional Schubert cycles.

The node~$(10,0)$ is the representation~
$10{\cdot}\begin{smallmatrix}&&0\\0&0&0&0&5\end{smallmatrix}$
of dimension~$100,386$, which is less than one-\emph{eightieth}
of the dimension of~$\L^{10}(\eum)$.  
Note also that~$\cI_{(10,0)}$ has irreducible integrals of 
dimension~$17$, namely the Schubert cycles corresponding to the 
node~$(17,8)$.  It would be interesting to know whether any
irreducible integral of~$\cI_{(10,0)}$ of dimension~$10$ or more is
a subvariety of one of these $17$-dimensional Schubert cycles.

Again, answers to these questions will have to await a detailed study
of the ideals involved.

\bibliographystyle{amsplain}

\end{document}